\documentclass[11pt,a4paper,reqno]{amsart}
\usepackage[paperheight=11in,paperwidth=8.5in]{geometry}
\usepackage[utf8]{inputenc}
\usepackage{amsmath,amsfonts,amsthm,amssymb}
\usepackage{dsfont}
\usepackage{graphicx,enumerate}
\usepackage{mathtools,bm,tikz-cd,esint}
\usepackage{todonotes}
\usepackage{comment}
\usepackage[normalem]{ulem}
\usepackage[dvipsnames]{xcolor}
\usepackage{appendix}

\usepackage{hyperref}
\hypersetup{
    colorlinks=true,
    linkcolor=blue,
    filecolor=magenta,      
    urlcolor=cyan,
    pdftitle={On the resonant Carleson-Radon transform in all dimensions},
    pdfpagemode=FullScreen,
    }

\makeatletter
\def\l@section{\@tocline{1}{12pt plus2pt}{0pt}{}{\bfseries}}
\def\l@subsection{\@tocline{2}{0pt}{2pc}{2pc}{}}
\makeatother
\makeatletter
\def\subsection{\@startsection{subsection}{2}{\z@}%
	{-3.25ex\@plus -1ex \@minus -.2ex}%
	{1.5ex \@plus .2ex}%
	{\normalfont\bfseries\boldmath}}
\def\subsubsection{\@startsection{subsubsection}{3}%
	\z@{.5\linespacing\@plus.7\linespacing}{-.5em}%
	{\normalfont\bfseries\boldmath}}
\renewcommand\paragraph{\@startsection{paragraph}{4}{\z@}%
	{3.25ex \@plus1ex \@minus.2ex}%
	{-1em}%
	{\normalfont\normalsize\bfseries}}
\makeatother

\numberwithin{equation}{section}
\newtheorem{theorem}{Theorem}[section]
\newtheorem*{theorem*}{Theorem}

\newtheorem{observation}[theorem]{Observation}
\newtheorem{lemma}[theorem]{Lemma}
\newtheorem{proposition}[theorem]{Proposition}

\theoremstyle{definition}
\newtheorem{remark}[theorem]{Remark}

\newtheorem{definition}[theorem]{Definition}
\newtheorem{assumption}[theorem]{Assumption}

\newcommand{\R}{\mathbb{R}}
\newcommand{\C}{\mathbb{C}}
\newcommand{\I}{\mathbb{I}}
\newcommand{\J}{\mathbb{J}}

\newcommand{\T}{\mathbb{T}}
\newcommand{\V}{\mathbb{V}}
\newcommand{\Z}{\mathbb{Z}}
\newcommand{\N}{\mathbb{N}}

\newcommand\numberthis{\addtocounter{equation}{1}\tag{\theequation}}

\newcommand{\br}[1]{\left(#1\right)}
\newcommand{\mr}[1]{\left[#1\right]}
\newcommand{\Br}[1]{\left[#1\right)}
\newcommand{\bR}[1]{\left(#1\right]}
\newcommand{\BR}[1]{\left\{#1\right\}}
\newcommand{\abs}[1]{\left\vert#1\right\vert}
\newcommand{\nrm}[1]{\left\Vert#1\right\Vert}
\newcommand{\ang}[1]{\left<#1\right>}

\newcommand{\set}[2]{\left\{ #1 \middle.: #2 \right\}}

\newcommand{\1}{\mathds{1}}

\renewcommand{\P}{\mathbb{P}}

\newcommand{\W}{\mathbb{W}}
\newcommand{\uk}{{\underline{k}}}
\newcommand{\ok}{{\overline{k}}}

\newcommand{\para}{{\operatorname{par}}}
\newcommand{\bmo}{{\operatorname{bmo}}}

\newcommand{\dia}{{\operatorname{diam}}}
\newcommand{\tphi}{{\widetilde{\phi}}}
\newcommand{\tarphi}{{\widetilde{\varphi}}}
\newcommand{\tpsi}{{\widetilde{\psi}}}
\newcommand{\id}{\operatorname{Id}}

\newcommand{\fS}{\mathfrak{S}}
\newcommand{\fs}{\mathfrak{s}}
\newcommand{\fp}{\mathfrak{p}}

\newcommand{\cA}{\mathcal{A}}
\newcommand{\cC}{\mathcal{C}}
\newcommand{\cE}{\mathcal{E}}
\newcommand{\cF}{\mathcal{F}}
\newcommand{\cG}{\mathcal{G}}
\newcommand{\cH}{\mathcal{H}}
\newcommand{\cI}{\mathcal{I}}
\newcommand{\cJ}{\mathcal{J}}
\newcommand{\cK}{\mathcal{K}}
\newcommand{\cL}{\mathcal{L}}
\newcommand{\cM}{\mathcal{M}}

\newcommand{\cP}{\mathcal{P}}
\newcommand{\cQ}{\mathcal{Q}}
\newcommand{\cT}{\mathcal{T}}
\newcommand{\cU}{\mathcal{U}}

\newcommand{\cW}{\mathcal{W}}
\newcommand{\cZ}{\mathcal{Z}}

\newcommand{\va}{{\bm{a}}}
\newcommand{\vb}{{\bm{b}}}
\newcommand{\vc}{{\bm{c}}}
\newcommand{\ve}{{\bm{e}}}

\newcommand{\vu}{{\bm{u}}}
\newcommand{\vv}{{\bm{v}}}
\newcommand{\vw}{{\bm{w}}}
\newcommand{\vx}{{\bm{x}}}
\newcommand{\vy}{{\bm{y}}}
\newcommand{\vz}{{\bm{z}}}
\newcommand{\vt}{{\bm{t}}}
\newcommand{\vC}{{\bm{C}}}
\newcommand{\vI}{{\bm{I}}}
\newcommand{\vJ}{{\bm{J}}}

\newcommand{\vQ}{{\bm{Q}}}
\newcommand{\vR}{{\bm{R}}}
\newcommand{\vS}{{\bm{S}}}
\newcommand{\vX}{{\bm{X}}}

\newcommand{\vZ}{{\bm{Z}}}
\newcommand{\valpha}{{\bm{\alpha}}}
\newcommand{\vbeta}{{\bm{\beta}}}
\newcommand{\vtau}{{\bm{\tau}}}

\newcommand{\vxi}{{\bm{\xi}}}

\newcommand{\vzeta}{{\bm{\zeta}}}

\newcommand{\vomega}{{\bm{\omega}}}
\newcommand{\vnull}{{\bm{0}}}

\newcommand{\vLambda}{{\bm{\Lambda}}}
\newcommand{\vid}{\bm{I\!d}}

\def\beq{\begin{equation}}
	\def\eeq{\end{equation}}

\def\beq{\begin{equation}}
	\def\eeq{\end{equation}}

\def\a{\alpha}
\newcommand{\cR}{\mathcal{R}}

\def\E{\operatornamewithlimits{\mathbb{E}}}
\def\beq{\begin{equation}}
	\def\eeq{\end{equation}}

\def\beq{\begin{equation}}
	\def\eeq{\end{equation}}

\DeclareMathOperator{\supp}{supp}

\DeclareMathOperator{\pv}{p.v.}
\DeclareMathOperator{\Mod}{Mod}
\DeclareMathOperator{\Tr}{Tr}
\DeclareMathOperator{\Dil}{Dil}
\DeclareMathOperator{\dist}{dist.}

\DeclareFontFamily{U}{mathx}{}
\DeclareFontShape{U}{mathx}{m}{n}{<-> mathx10}{}
\DeclareSymbolFont{mathx}{U}{mathx}{m}{n}
\DeclareMathAccent{\widehat}{0}{mathx}{"70}
\DeclareMathAccent{\widecheck}{0}{mathx}{"71}
\allowdisplaybreaks
\begin{document}

\title[Carleson-Radon transform with linear resonance]{\LARGE{On the resonant Carleson-Radon transform in all dimensions}\\\bigskip\Large{The degree one resonant case}}

\author{Martin Hsu}
\author{Victor Lie}

\thanks{\textit{Address}(M. Hsu): Department of Mathematics, Purdue University, 150 N. University St, W. Lafayette, IN 47907, U.S.A}
\thanks{\textit{Address}(V. Lie): Department of Mathematics, Purdue University, 150 N. University St, W. Lafayette, IN 47907, U.S.A. and The ``Simion Stoilow" Institute of Mathematics of the Romanian Academy, Bucharest, RO 70700, P.O. Box 1-764, Romania}

\begin{abstract} 
Since its inception, more than a decade ago, the Carleson-Radon theme revolved around two central \emph{hierarchical structures}---with both displayed below in an antithetical and increasing order of complexity:
\begin{itemize} 
\smallskip
\item the \underline{\emph{higher}} versus \underline{\emph{one-dimensional}} setting;
\smallskip
\item the \underline{\emph{non-resonant}} versus \underline{\emph{resonant}} behavior accounting for the absence/presence of modulation invariance.
\end{itemize}
\smallskip

Consistent with the above---and under suitable hypothesis not detailed here---the historical evolution of the Carleson-Radon theme up to date accounted chronologically for the treatment of the \emph{non-resonant higher-dimensional} regime (\cite{PY19}, \cite{AMPY24}), followed by the corresponding \emph{degree one resonant higher-dimensional} regime (\cite{Bcarlrad}), and, more recently, by the \emph{non-resonant one-dimensional} regime (\cite{BGH24}, \cite{HsuL24}). 

In this paper, we provide the resolution of the \emph{degree one resonant case in all dimensions}. Moreover, confirming the above hierarchies, our methods extend and are able to treat---explicitly or implicitly---all of the above-mentioned regimes and most of the above results.
\bigskip

With the above context settled, our main result reads as follows: for any dimension $D\geq 1$ set $\vX\br{\vt}:=\br{\vt,\abs{\vt}^2},\; \vt\in\R^D$, and let $K(\vt)$ be any suitable translation invariant Calder\'on--Zygmund kernel. If \(\V\leq\R^{D+1}\) is any linear subspace such that $ \exists\:\:\vv_0\in\R^D\times\BR{0}$ nontrivial with $\vv_0\perp\V$ then the following (maximal) Carleson-Radon transform $CR^\ast_\V$ is $L^p(\R^{D+1})-$bounded in the \emph{maximal range} $1<p<\infty$, where
$$CR^\ast_\V f\br{\vx}:=
        \sup_{\substack{
            0<r<R<\infty\\
            \va\in\V
        }}
        \abs{
            \int_{r<\abs{\vt}\leq R}
                f\br{\vx-\vX\br{\vt}}
                e\br{\va\cdot \vX\br{\vt}}
                K\br{\vt}
            d\vt
        }.$$
 The above choice for $\V$ creates a \emph{maximal} linear subspace of $\R^{D+1}$ closed under parabolic scaling for which
\begin{itemize}
\item $CR^\ast_\V$ \emph{\underline{is} degree one resonant}, and
 
\item  $CR^\ast_\V$ \emph{\underline{is not} degree two (or higher) resonant}.
\end{itemize}  
\bigskip
      
The proof of the above result unravels several new manifestations and ideas meant to capture the remarkable features of the resonant Carleson-Radon behavior. Indeed, relying on the deep interplay between the structure of the modulation invariance and the geometry of the underlying manifold, we decompose our operator into three main components, as follows: low frequency (L),  high-frequency non-stationary (HNS), and high-frequency stationary (HS).

In contrast with the canonical situations, a first distinctive manifestation is that all of the above components inherit a linear modulation invariance feature. As a consequence, we develop a new, \emph{formal time-frequency framework}, which, via its \emph{variable resolution multiscale analysis}, is precisely tailored in order to  cover all of the specific attributes of the three components listed above and to exploit the dual nature of our transform:
\begin{itemize} 
\smallskip
\item \underline{\emph{the Carleson behavior}} encapsulated within the \emph{linear modulation invariance feature} relies on multiscale wave-packet analysis, which in turn is based on suitable mass and energy selection algorithms.
\begin{itemize}
\item for (L), (HNS): is expressed via the new concept of a \emph{shifted superposed Carleson operator};

\item for (HS): it involves a \emph{low resolution multiscale analysis} corresponding to wave-packets whose time-frequency information is encoded within super-Heisenberg tiles with expanding, supraunitary volume and variable geometry depending on the dimension and frequency location;
\end{itemize}     
\smallskip
\item \underline{\emph{the Radon behavior}} encapsulated within the \emph{lower dimensional character of the underlying integrant manifold}, it manifests as a cancellation hidden within the curvature of the manifold.
    
\begin{itemize}
\item for (L), (HNS): is exploited via kernel oscillation and non-stationary phase principle;

\item for (HS): is extracted via the LGC-method by appealing to a \emph{high resolution single scale analysis} applied within each super-Heisenberg tile to a Heisenberg tile partition.
\end{itemize}
\end{itemize}

\end{abstract}


\dedicatory{Dedicated to our parents,\\ those recently departed and those among us, \\for their enduring impact on our lives.}

\date{\today}


\maketitle
\newpage
\setcounter{tocdepth}{3}
\tableofcontents
\section{Introduction}

In this paper, we present the resolution of the \emph{degree one resonant} Carleson-Radon transform \emph{in all dimensions}. More precisely, in what follows, we will prove the \emph{maximal $L^p$-boundedness range} for a singular integral maximal operator that exhibits a dual nature:
\begin{itemize} 
\item a \emph{Carleson behavior} encoded in a \emph{linear modulation invariance} property;  

\item a \emph{Radon behavior}  encoded in the \emph{lower dimensional} manifold over which the integration of the operator's kernel is performed.
\end{itemize}

The motivation and historical evolution of this theme will be discussed in \textsc{Section \ref{hyst}}. With these being said, we are now ready to state our

\subsection{Main results}\label{mainres}

\begin{theorem}\label{thm_main1D}[\textsf{Degree one resonance: one-dimensional setting}]
\smallskip

Define the degree one resonant Carleson-Radon transform in one dimension as\footnote{Strictly speaking, a more precise definition would be the ``degree one resonant maximally truncated Carleson-Radon transform along parabola". Also, here and throughout the paper, given $a\in\R$, we use the notation $e(a):=e^{2 \pi i a}$.}
\begin{equation}\label{par}
  CR^{*}_{I}f\br{x_0,x_1}:=
    \sup_{\substack{
            0<r<R<\infty\\
            a\in\R
        }}
        \abs{
            \int_{r<\abs{t}\leq R}
                f\br{x_0-t,x_1-t^2}
                e\br{a t^2 }
            \frac{dt}{t}
        }\,.
\end{equation}
Then $CR^{*}_{I}$ is of strong type \(\br{p,p}\) for all \(p\in\br{1,\infty}\).
\end{theorem}
\medskip

Moreover, we have the following extension addressing the general dimensional case: 
\begin{itemize}
\item let \(D\in\N\) and set the parametric equation for the paraboloid in \(\R^{D+1}\) as
\begin{equation}\label{eq_paraboloid}
    \vX\br{\vt}:=\br{\vt,\abs{\vt}^2},\quad \vt\in\R^D\,.
\end{equation}
\item also consider a singular kernel \(K\) on \(\R^D\setminus\BR{0}\) satisfying
\begin{equation}\label{eq_kernel_cond}
    \int_{r<\abs{\vt}\leq R}
        K\br{\vt}
    d\vt
    =0,\quad
    \abs{\partial^{\valpha}K\br{\vt}}\lesssim \abs{\vt}^{-D-\abs{\valpha}},\quad
    \valpha=\br{\alpha_0,\dots,\alpha_{D-1}},\quad
    0\leq \alpha_i\leq A
\end{equation}
with fixed sufficiently large \(A\in\N\) that only depends on \(D\). 
\end{itemize}

With these settled, we have the following

\begin{theorem}\label{thm_main}[\textsf{Degree one resonance: the general dimensional setting}]
\smallskip

In the above setting, given a linear subspace\footnote{In effect, there is no need to \emph{apriori} enforce on $\V$ a vector space structure as long as one requires $\V$ to be closed under the natural symmetries dictated by \eqref{Modl2} and the structure of \eqref{def_C_ast_V_K}, that is, for $\V$ to be closed under addition and parabolic scaling. For a more detailed discussion on these issues, the reader is invited to consult \textsc{Section \ref{optV}} in the Final Remarks Section.} \(\V\leq\R^{D+1}\), we set
    \begin{equation}\label{def_C_ast_V_K}
        CR^\ast_\V f\br{\vx}:=
        \sup_{\substack{
            0<r<R<\infty\\
            \va\in\V
        }}
        \abs{
            \int_{r<\abs{\vt}\leq R}
                f\br{\vx-\vX\br{\vt}}
                e\br{\va\cdot \vX\br{\vt}}
                K\br{\vt}
            d\vt
        }.
    \end{equation}

Assume now that there exists a nontrivial vector
 \begin{equation}\label{degone}
 \vv_0\in\R^D\times\BR{0}\quad\textrm{with}\quad\vv_0\perp\V\,.
 \end{equation}

Then, for any $f\in\mathcal{S}\br{\R^{D+1}}$ and any $p\in\br{1,\infty}$, the following holds:
    \begin{equation}
        \nrm{CR^\ast_\V  f}_{L^p}
        \underset{p}{\lesssim}
        \nrm{f}_{L^p}\,.
    \end{equation}
\end{theorem}

\begin{remark}[\textsf{Main results: relationship}] Notice that we trivially have that Theorem  \ref{thm_main} is a generalization of Theorem \ref{thm_main1D}. Indeed, taking first $D=1$ followed by $ \vv_0=(1,0)$ and $\V=\{0\}\times\R$ in \eqref{degone} and finally $K(t)=\frac{1}{t}$ we deduce that \eqref{def_C_ast_V_K} identifies with \eqref{par}.
\end{remark}

\begin{remark}[\textsf{Linear resonance}] For $D\geq 1$ fixed, $j,\,l\in \N$ with $0\leq j\leq D$, $\textbf{a}=(a_0,\ldots, a_D)\in \R^{D+1}$, $\vx=(x_0,\ldots, x_D)\in \R^{D+1}$, and $f\in L^1_{loc}(\R^{D+1})$, throughout this paper, we let
\begin{itemize}
\item the \emph{generalized modulation of order $l$ in the $j^{th}$ coordinate}:
 \begin{equation}\label{Mlj}
M^j_{l,a_j} f(\vx)=e^{i\,a_j x_j^l}\,f(\vx)\,;
\end{equation}
\item  the \emph{generalized modulation of order $l$:\footnote{Here we use the convention  $\vx^l:=(x_0^l,\ldots, x_D^l)$.}}
 \begin{equation}\label{Ml}
M_{l,\textbf{a}} f(\vx):=\left(\prod_{j=0}^D M^j_{l,a_j}\right) f(\vx)=e^{i\,\textbf{a}\cdot \vx^l}\,f(\vx)\,.
 \end{equation}
\end{itemize}
With this notation, we notice that both $CR^{*}_{I}$ and $CR^\ast_\V$ are \emph{invariant under (suitable) linear modulations}. Indeed, we have
\begin{equation}\label{Modl}
CR^{*}_{I}(M_{1,a}^1 f)=CR^{*}_{I}(f) \qquad \forall\:a\in\R,\,f\in\mathcal{S}\br{\R^{2}}\,,
\end{equation}
and, more generally 
\begin{equation}\label{Modl2}
CR^\ast_\V(M_{1,\textbf{a}} f)=CR^\ast_\V(f) \qquad \forall\:\textbf{a}\in\V,\,f\in\mathcal{S}\br{\R^{D+1}}\,.
\end{equation}
\end{remark}

\begin{remark}[\textsf{Motivation for the choice of $\V$: linear versus higher order resonance}]\label{MotV} Our Theorem  \ref{thm_main} addresses the boundedness of a Carleson-Radon transform which exhibits modulation invariance sets restricted to a subspace \(\V\leq\R^{D+1}\). Our choice to define \(\V\) in terms of \eqref{degone} is made in order to achieve a \emph{maximal subspace} of $\R^{D+1}$  which \emph{simultaneously}:\footnote{For a justification of this claim, please see \textsc{Section \ref{optV}}.}  
\begin{itemize}
\item  \underline{respects} \textsf{parabolic rescaling}, \emph{i.e.}, $\V$ is \textsf{invariant} under the anisotropic \textsf{paraboloid dilation} symmetry;
\item \underline{provides} \textsf{degree one resonances}, \emph{i.e.}, \textsf{linear modulation invariance} features as  described by \eqref{Modl2};
\smallskip
\item \underline{rules out} \textsf{degree two resonances}, \emph{i.e.}, $CR^\ast_\V$ does not have \textsf{higher order modulation invariance} features.
\end{itemize}
\end{remark}

\begin{remark}[\textsf{Overcoming a dimensional barrier}]\label{apdimCR} 
Until the publication of the current paper, all the treated Carleson-Radon cases involved restrictions to a higher-dimensional setting ($D>1$) or/and non-resonant regimes. In particular, the general point of view relying on a $T T^{*}$ kernel regularization faced a dimensional singularity barrier when attempting to treat the $D=1$ case. One of the main contributions of our present work is to show that, with a proper perspective on the interplay between the zero (Carleson) and non-zero (Radon) curvature behavior, this \emph{dimensional barrier becomes illusive}. In particular, the methods involved in the treatment of the (degree one) resonant one-dimensional case extend in all dimensions, confirming the intuition stated in the Introduction\footnote{See (H3) in \textsc{Section 1.2} therein.} of \cite{HsuL24}. 
\end{remark}

\begin{remark}[\textsf{Extensions of the previous results}]\label{ext}
Finally, reinforcing the two hierarchical structures stated in the abstract, the methods developed in this paper extend/provide a treatment for all of the following\footnote{For the second and third bullet points below, please see the similar comments in \textsc{Section 1.3} of \cite{Bcarlrad}.}: 
\begin{itemize}
\item (\textsf{explicitly}) the \emph{degree one resonant higher-dimensional} Carleson-Radon transform in (\cite{Bcarlrad});

\item (\textsf{implicitly}) the \emph{degree one resonant one-dimensional}  Carleson operator in \cite{c1};

\item (\textsf{implicitly})  the \emph{degree one resonant general dimension}  Carleson operator in \cite{sj2};

\item (\textsf{implicitly}) the \emph{non-resonant one-dimensional} Carleson-Radon transform in \cite{BGH24} and \cite{HsuL24};

 \item (\textsf{implicitly}) in part\footnote{Here some restrictions apply. For more details, please see \textsc{Section \ref{finrem}}.}, the \emph{non-resonant higher-dimensional} Carleson-Radon results in \cite{PY19} and \cite{AMPY24}.
\end{itemize}
For more on all of the above items, we invite the reader to consult the Final Remarks section.
\end{remark}

\subsection{Motivation, historical evolution and future directions}\label{hyst}

In this section, following in part the presentation in \cite{HsuL24}, we focus on a brief antithetical discussion of the evolution of the two central and interconnected themes that generate the main subject of our present paper:
\begin{itemize}  
\item the \emph{Polynomial Carleson operator}:
\begin{equation}\label{CdD}
       C_{d,D}f(x):=\sup_{P\in\mathcal{Q}_{d,D}}\left|\int_{\R^D} f(x-t)\,e^{i P(t)}\,K(t)\,dt\right|,\qquad\,x\in \R^D\,,
\end{equation}

\item the \emph{Polynomial Carleson-Radon transform}:
\begin{equation}\label{CRdD}
       CR_{d,D}f(x,y):=\sup_{P\in\mathcal{Q}_{d,D}}\left|\int_{\R^D} f(x-t,y-|t|^2)\,e^{i P(t)}\,K(t)\,dt\right|,\qquad\,(x,y)\in \R^D\times\R\,,
\end{equation}
\end{itemize}
where in the above $d,D\in\N$ with $\mathcal{Q}_{d,D}$ the set of all real-coefficient polynomials in $D$ variables of degree less than or equal to $d$ and $K$ is a standard translation invariant Calder\'on--Zygmund kernel on $\R^D$.

\subsubsection{Hierarchies}\label{Hier}

Our presentation will be delivered in the increasing order of complexity according to the following hierarchies:
\begin{itemize}
\item \emph{non-resonant} regime, \emph{i.e., no modulation invariance} case;

\item \emph{degree one resonant} regime,  \emph{i.e., linear modulation invariance} case;

\item \emph{degree two (or higher) resonant} regime,  \emph{i.e., higher order modulation invariance} case;
\end{itemize}
All of the above regimes will be analyzed through the lenses offered by a \emph{dimensional dichotomy}\footnote{In hindsight, this dimensional dichotomy appeared as a manifestation of the employed methodologies rather than being of intrinsic nature. For more details, please see \textbf{Remark \ref{apdimCR}} and \textsc{Section \ref{apdimC}}.}
$$ \textrm{\emph{one} dimensonal versus \emph{higher} dimensional setting}.$$

\subsubsection{Polynomial Carleson operator}\label{CEv}
$\newline$

We start our discussion with the following

\begin{remark}[\textsf{Chronology and motivation}]\label{ChronC}
In contrast to the hierarchical order listed above, \emph{historically}, the origin of the Carleson theme goes back to the linear resonant one-dimensional regime. Indeed,  the operator corresponding in \eqref{CdD} to the case $d=D=1$, is referred to as the (classical) Carleson operator and its study was motivated by a century-old conjecture of N. Luzin, \cite{Luz}, stating that the Fourier series of a square integrable function converges almost everywhere. The study of the non-resonant regime as well as the general resonant regime are rooted in Stein's conjecture on the Polynomial Carleson operator, \cite{s2}, \cite{sw}. This is motivated by two distinct directions of investigation: \emph{(i)} the study of oscillatory integrals/operators in Euclidean setting due to E. Stein and S. Wainger, \cite{SW70}, and, \emph{(ii)} the study of singular oscillatory integrals on nilpotent groups due to F. Ricci and E. Stein, see \cite{RS86}, \cite{RS87} and \cite{RS89}.
\end{remark}

With these being said, we resume our historical evolution according to the hierarchies presented above:

\begin{itemize}
\item \emph{non-resonant one-dimensional} regime:
\smallskip
\begin{itemize}
\item This was addressed in \cite{s2} in the particular setting $D=1$, $d=2$ and the supremum in \eqref{CdD} 
 taken over polynomials with no linear term. The approach involved an asymptotic analysis of the kernel $\frac{e^{i y^2}}{y}$.  
 
\item A more general non-resonant setting dealing with arbitrary fewnomials with no linear term was solved via the  LGC-method in \cite{LVunif}.
\end{itemize} 
\smallskip

\item \emph{non-resonant higher/general dimensional} regime: 
\smallskip
\begin{itemize}
\item This corresponds to the situation when the supremum in \eqref{CdD} is taken over polynomials with no linear term, and it was treated in \cite{sw} relying on Van der Corput estimates and $TT^{*}$-method.
\end{itemize}

\smallskip

\item \emph{linear resonant one-dimensional} regime:
\smallskip
\begin{itemize}
\item As mentioned in \textbf{Remark \ref{ChronC}}, this lies at the foundation of theme \eqref{CdD} and corresponds to the classical Carleson operator defined by taking  $d=D=1$. The first proof of the $L^2$ boundedness of this operator was presented in \cite{c1} and set the foundation of the time-frequency analysis area. 
\item The $L^p$, $1<p<\infty$, counterpart was provided by R. Hunt in \cite{hu}. 
\item Two other seminal proofs of Carleson's result\footnote{For a more illuminating discussion on the particularities of each of the three known methods addressing the boundedness of the classical Carleson operator, the reader is invited to consult the longer itemized paragraph above \textbf{Remark \ref{Temp}}.}  are due to
\begin{itemize}
\item C. Fefferman, \cite{f}: the fundamental concepts of \emph{tiles}, \emph{tress}, \emph{mass of a tile} are introduced.

\item M. Lacey and C. Thiele, \cite{lt3}: building on the pioneering contributions in \cite{c1} and \cite{f} as well as on their influential work on the Bilinear Hilbert transform, \cite{LT97}, \cite{LT99}, the authors develop a streamlined approach centered around the \emph{mass and energy selection algorithms}.  
\end{itemize} 
\item More recently, while using Carleson's result as a black-box,  \cite{GALVHybcurves} treats via the LGC-method a more general version in which the supremum in \eqref{CdD} is taken over fewnomials with no quadratic, but possibly linear terms.
\end{itemize}
\smallskip

\item \emph{linear resonant higher/general dimensional} regime: 
\smallskip
\begin{itemize}
\item Corresponding to the $d=1,\,D\geq 1$ case this was treated by P. Sj\"olin in \cite{sj2} via an adaptation to the $D-$dimensional setting  of the ideas in \cite{c1}.
 \end{itemize}   
\smallskip

\item \emph{higher order resonant one-dimensional} regime: 
\smallskip
\begin{itemize}
\item The $d=1,\, D=1$ case was solved in \cite{LVQuadCarl}, while the general one-dimensional version of \eqref{CdD}  was settled in \cite{LVPolynCarl}. The main ideas introduced in this context are: (i) a new perspective on the higher-order wave-packet theory under the name of \emph{relational time-frequency} analysis inspired by \cite{Feucp}, (ii) a specialized \emph{tile selection algorithm} addressing the increased boundary effect of the higher order wave-packet discretization, and (iii) a novel \emph{local analysis} aimed towards the mass selection algorithm implying the removal of the so-called exceptional sets.
\end{itemize} 
\smallskip
\item \emph{higher order resonant higher/general dimensional} regime: 
\smallskip
\begin{itemize}
\item The general $d,\,D\geq 1$  case involves a keen adaptation of the one-dimensional setting according to various scenarios:
\begin{itemize}
\item the standard translation invariant Calder\'on-Zygmund kernel case,  \cite{LVHidimPolyCarl};

\item the general H\"older continuous Calder\'on–Zygmund kernels case,  \cite{ZK};

\item the doubling metric measure spaces case, \cite{BDJST}.

\end{itemize}
\end{itemize}    
\end{itemize}

\subsubsection{Polynomial Carleson-Radon transform}\label{CREv}
$\newline$

Mirroring the previous section, we initiate our discussion with the following

\begin{remark}[\textsf{Chronology and motivation}]\label{ChronCR0}

In part, the chronology of the Carleson-Radon theme follows in reverse that of the Carleson counterpart.\footnote{For more on this see \textsc{Section \ref{apdimC}}.}  Indeed, the former
\begin{itemize}
\item follows the hierarchical order non-resonant/resonant regime presented in Section  \ref{Hier};

\item develops first in the  $D>1$ setting as opposed to the $D=1$ setting.
\end{itemize}

Historically, the topic represented by \eqref{CRdD} was inspirationally formulated by L. Pierce and P.-L. Yung in \cite{PY19}. Its main motivation relies precisely on the new phenomena exhibited by this object, which require, in particular, a deep understanding of the interrelationship between maximal modulation features (Carleson behavior) and the geometry of the underlying lower-dimensional manifold (Radon behavior). Because of all these, the Carleson-Radon theme has a higher complexity than that of the Carleson theme, and in particular, the $L^p$ boundedness of \eqref{CRdD} implies the corresponding bounds on \eqref{CdD} (see also \textsc{Section \ref{finrem}}).

Besides the analogy with the Polynomial Carleson operator theme, further motivation is served by the works of Stein and Wainger, \cite{sw}, and  Ricci and Stein, \cite{RS87}.
\end{remark}

With these done, we proceed with the historical evolution of the Carleson-Radon theme up to the publication of the present result:

\begin{itemize}
\smallskip
\item \emph{non-resonant higher-dimensional} regime: 
\smallskip
\begin{itemize}
\item This was addressed in \cite{PY19} in the particular setting $D>1$, $d>1$ and the supremum in \eqref{CRdD} 
 taken over polynomials with no linear terms, no quadratic resonance terms, and having a uniform homogeneous degree behavior. 
 
\item An extension of this result to more general manifolds arising from non-degenerate quadratic forms but under similar restrictions on dimension and on the supremum appeared in \cite{AMPY24}.
 
 \item Both these works involve some kernel regularizations via $TT^{*}-$method and a refinement of the techniques in \cite{sw}.
\end{itemize}
\smallskip

\item \emph{degree one resonant higher-dimensional} regime: 
\smallskip
\begin{itemize}
\item The earliest studies treating some related toy-models go back to \cite{R19} and \cite{B24}.

\item The main contribution to this regime is made in \cite{Bcarlrad}. There L. Becker proves that for $D>1$ the linear modulation invariant Carleson-Radon transform $ CR^\ast_\V $ defined by \eqref{def_C_ast_V_K}  is $L^p$ bounded as long as $\frac{D^2+4D+2}{(D+1)^2}<p< 2(D + 1)$ and $\V\leq\R^{D+1}$ subspace satisfying a suitable subcase of \eqref{degone}. This is the first important work that implements Carleson operator wave-packet methods into the Carleson-Radon realm.
\end{itemize}
\smallskip
\item \emph{non-resonant one-dimensional} regime: 
\smallskip
\begin{itemize}
\item Some particular toy models confined to one-variable dependence of the linearizing supremum function appear in \cite{GPRY17}. Further resonant one-dimensional toy models---but only using as a black-box Carleson's theorem---are discussed in both \cite{GPRY17} and \cite{R21}.

\item The one-dimensional case in the monomial phase non-resonant setting is finally independently solved in two papers: 
\begin{itemize}
\item The case $P$ monomial of degree three in \eqref{CRdD} is settled in \cite{BGH24} using a local smoothing
estimate for a variable-coefficient Schr\"odinger operator via small cap decoupling theory together with a more elaborate machinery developed in a companion paper \cite{CGGHIW24}. A crucial ingredient for the proof's validity is the so-called Nikodym non-compression hypothesis, which must be satisfied by the phase of the transform's multiplier.
    
\item A more general case corresponding to arbitrary non-resonant phase/manifold monomials is treated in \cite{HsuL24}. The proof represents the stepping stone for our present paper and is based on the LGC method, which has the advantage of being well adapted to both resonant and higher-dimensional regimes. 
\end{itemize}
\end{itemize}
\end{itemize}
We end this section with the following 

\begin{remark}[\textsf{Limitations of the previous methods in approaching the linear resonant one-dimensional regime}]\label{ChronCR}
In what follows, leaving aside the non-trivial technical challenges, we summarize the \emph{conceptual obstructions} of the previous methods in approaching the (degree one) resonant one-dimensional regime. Our discussion follows the order of the contributions given in the above presentation: 
\begin{itemize}
\item The \emph{non-resonant higher-dimensional} regime: the methods employed here face both dimensional and modulation invariance obstructions due to their reliance on\footnote{For more revealing details, please see \textsc{Section \ref{PY}}.} 
\begin{itemize}
\item $T T^{*}$ kernel regularization arguments only available in higher dimensions, and
\item Calder\'on-Zygmund techniques, which favor the special role of the zero frequency and (thus) lack the adaptability to a modulation invariance setting.
\end{itemize}
\smallskip
\item The \emph{degree one resonant higher-dimensional} regime: the methods involved here have dimensional but no modulation invariance obstructions. Indeed, while suited for the linear modulation invariance since they incorporate Carleson operator techniques in the spirit of \cite{f}, these methods are not effective in the one-dimensional setting due to a coarser discretization approach which affects among others the treatment of the oscillatory behavior of the high-frequency anti-chain component.\footnote{For an illuminating discussion on these issues, the interested reader is invited to consult \textsc{Section \ref{LB}}.}
    
 \smallskip
\item The \emph{non-resonant one-dimensional} regime: the methods involved in this regime have modulation invariance obstructions but no dimensional obstructions. Indeed:
\begin{itemize}
\item the work in \cite{BGH24} relies on a local smoothing estimate under the Nikodym non-compression hypothesis \cite{CGGHIW24} as well as on decoupling theory, with these being employed in order to go beyond the universal estimates for oscillatory integrals associated to H\"{o}rmander-type phases given by Stein’s oscillatory integral theorem \cite{MR864375}. While proficient in the context of purely curvature/non-resonant regime, these methods are not adequate for the harder, resonant regimes, that is, they are able to address the Radon-type behavior but not the genuine Carleson-type behavior.\footnote{For a more detailed discussion, please see \textsc{Section \ref{Dec}}.}

\item strictly speaking, the work in \cite{HsuL24} also exploits the special role played by the zero-frequency, which, as already mentioned, is lost in the modulation invariant regime. However, given the special structure of the degree one resonance in the Carleson-Radon context, we reveal the following:
\end{itemize}
\end{itemize}
\end{remark}

\begin{remark}[\textsf{Incorporating the non-resonant LGC approach into the linear resonant one-dimensional regime}]\label{IncorpCR}
In contrast with the above approaches, the versatility of the LGC-method allows us to ``localize and shift" our entire argument in \cite{HsuL24}  in order to obtain a so-called single tile decay estimate. Thus, loosely speaking, our whole work in the non-resonant one-dimensional regime can be recovered and incorporated into the resonant case within the process of controlling the behavior of each elementary building block embodied by what we will refer to as a super-Heisenberg tile. From this point on, the crux of the linear-resonant approach requires the propagation of the single super-Heisenberg tile decay estimate within the full collection of super-Heisenberg tiles involved in the resonant Carleson-Radon transform discretization, all these while preserving the fabric of the mass and energy selection algorithms which address the linear modulation invariance, hence the Carleson-type behavior.   
\end{remark}

\subsubsection{Revisiting the dimensional dichotomy}\label{apdimC}

In view of both the chronology and methodology development discussed in \textsc{Section \ref{CEv}} and \textsc{Section \ref{CREv}}, it is interesting to notice that \emph{even though the two main themes started from the opposite end-points of the dimensional spectrum, after decades of dedicated study they evolved towards a common point of equilibrium.} More precisely, we have

\begin{itemize}

\item \textsf{The Carleson setting}: For both non-resonant and resonant regimes, this theme started with the $D=1$ setting and then moved towards the $D>1$ setting. This is explained in part by the fact that, for a while---especially in the one and then higher degree resonant regime---the general perceived point of view was that the combinatorics of the wave-packets/tiles is more challenging in higher dimensions due to the increased number of degrees of freedom.  
\smallskip
\item \textsf{The Carleson-Radon setting}: In contrast with the above,  for both non-resonant and resonant regimes, this second theme started with the $D>1$ setting and then evolved towards the $D=1$ setting. This discrepancy is due to the singular and geometric properties of the kernel, which in  \eqref{CRdD} is integrated over a manifold of co-dimension one. Indeed, the previously known higher-dimensional approaches used the fact that for $D>1$ one can exploit a $TT^{*}$ regularization effect of the kernel, a feature which breaks down for $ D=1$.
\end{itemize}

However, in view of the methods developed for approaching
\begin{enumerate}
\item the degree one resonant Carleson regime in one/general dimension---see \cite{c1} and \cite{sj2},
\item the general degree resonant Carleson regime in one/general dimension---see \cite{LVPolynCarl}, \cite{LVHidimPolyCarl} and \cite{ZK},
\item the degree one resonant Carleson-Radon regime in all dimensions---see the present paper,
\end{enumerate}
the following is revealed: once we are able to capture the quintessence of the one-dimensional reasonings, the higher-dimensional setting follows from a natural, though often non-trivial, adaptation. 
\smallskip

\emph{Thus, we conclude that at the very foundation, the key underlying ideas within the time-frequency analysis of the above two themes have a general applicability character irrespective of the dimension.} 

\subsubsection{The major remaining conceptual challenge: the higher degree resonance regime}\label{GoalCR} $\smallskip$

Comparing the advancements corresponding to the two themes---see \textsc{Section \ref{CEv}} and \textsc{Section \ref{CREv}}, we notice that the last major step in the resolution of the Carleson-Radon theme in \eqref{CRdD} is given by the higher degree resonant regime which, as we already know by now, corresponds to the situation when our main object of study enjoys quadratic or higher order modulation invariance features. 
 
 As explained in \cite{GALVHybcurves}, we recall here that there is a conceptual, fundamental difference between the nature of the methods addressing  problems which exhibit
\begin{itemize}
\item \emph{linear modulation invariance}: in this setting, one employs \emph{classical} time-frequency analysis, in the spirit of  \cite{c1}, \cite{f}, and \cite{lt3} or, in the case in which also curvature features are involved---see \emph{e.g.} \cite{LVBilHilbCarl}, \cite{GALVHybcurves} or the present paper, one adds into the picture the LGC method introduced in \cite{LVunif}. Notice that both these approaches rely on \emph{linear wave-packet analysis}.
\smallskip 
\item \emph{quadratic or higher modulation invariance}: in this setting, none of the above methods work, and thus, one is naturally led to consider the \emph{relational} time-frequency analysis introduced in \cite{LVQuadCarl} and \cite{LVPolynCarl} and relying on \emph{higher order wave-packet analysis}.  
\end{itemize}

Moreover, up to the present day, the Polynomial Carleson operator represented by \eqref{CdD} is \emph{the only operator in the literature whose boundedness behavior is understood in the higher degree resonance regime}.
\smallskip

In view of the above considerations, we have the following:
\begin{itemize}
\item relying on \textbf{Remark \ref{MotV}}, we deduce that our \textbf{Theorem \ref{thm_main}} is an \emph{optimal result within the linear resonant regime covered by the classical time-frequency and LGC methods}.

\item as also confirmed by the discussion of the regime hierarchies in \textsc{Section \ref{Hier}}, the quadratic (or higher) resonant Carleson-Radon transform
\begin{itemize} 
\item is significantly more challenging than its linear resonant counterpart;

\item it requires a new conceptual approach that has to include and further develop the relational time-frequency method.  
\end{itemize}
\end{itemize}
We end this section with a brief discussion---in increasing order of complexity---of the two key problems modeling the essence of the quadratic resonant Carleson-Radon realm. In view of \textsc{Section \ref{apdimC}}, we only focus on the one-dimensional setting:
\medskip

\noindent\underline{\textsf{Problem 1}}: \emph{Investigate the $L^p$ boundedness of the (type 1)-degree two resonant Carleson-Radon transform}
\begin{equation}\label{q1}
    CR_{II,1}f\br{x_0,x_1}:=
    \sup_{a\in\R}
        \abs{
            \textrm{p.v.}\,\int_{\R}
                f\br{x_0-t,x_1-t^2}
                e\br{a t }
            \frac{dt}{t}
        }\,.
\end{equation}
\indent Notice that in this setting we have that $CR_{II,1}$ obeys:\footnote{Recall notation \eqref{Mlj}.}
\begin{itemize}
\item a linear modulation invariance on the first component:
\begin{equation}\label{qi11}
CR_{II,1}(M_{1,a}^0)=CR_{II,1}\quad \textrm{for any}\:a\in\R\,.
\end{equation}

\item a \emph{dependent linear-quadratic resonance}:
\begin{equation}\label{qi12}
CR_{II,1}(M_{2,b}^0 M_{1,b}^1)=CR_{II,1}\quad \textrm{for any}\:b\in\R\,.
\end{equation}
\end{itemize}
\medskip

\noindent\underline{\textsf{Problem 2}}: \emph{Investigate the $L^p$ boundedness of the (type 2)-degree two resonant Carleson-Radon transform}
\begin{equation}\label{q2}
    CR_{II,2}f\br{x_0,x_1}:=
    \sup_{a,b\in\R}
        \abs{
            \textrm{p.v.}\,\int_{\R}
                f\br{x_0-t,x_1-t^2}
                e\br{a t+bt^2}
            \frac{dt}{t}
        }\,.
\end{equation}
In this increased complexity setting $CR_{II,2}$ obeys the \emph{independent linear-quadratic resonance}
\begin{equation}\label{qi21}
CR_{II,2}(M_{1,a}^0 M_{2,b}^0 M_{1,c}^1)=CR_{II,2}\quad \textrm{for any}\:a,\,b,\,c\in\R\,.
\end{equation}

\begin{remark}[\textsf{Quadratic resonance hierarchy}]\label{QCR}
Based on the intuition provided by the Polynomial Carleson theme, it is natural to conjecture that both $CR_{II,1}$ and $CR_{II,2}$ are $L^p$ bounded operators whenever $1<p<\infty$. Also, it is worth mentioning the following:
\begin{itemize}
\item Comparing \eqref{Modl} with \eqref{qi11} and \eqref{qi12} with \eqref{qi21} we notice that $CR_{II,1}$ lies ``halfway" between $CR^{*}_{I}$ and $CR_{II,2}$ representing a natural stepping stone between the degree one resonant regime and the full degree two resonant regime;

\item In analogy with the Polynomial Carleson theme in \eqref{CdD}, we have that $CR_{II,2}$ captures essentially the full complexity of the one-dimensional Polynomial Carleson-Radon theme in \eqref{CRdD}.
\end{itemize} 
\end{remark}

\subsection{Structure of the paper and main ideas}\label{mainid}

Following the presentation of our main results in \textsc{Section \ref{mainres}}, our paper develops as follows:
\begin{itemize}
\item The one-dimensional linear resonant regime, which is the subject of \textbf{Theorem \ref{thm_main1D}}, is covered from \textsc{Section \ref{prelred}} to \textsc{Section \ref{sec_pf_key_lemmas}}.

\item The extension to the general dimensional linear resonant regime, which is the content of \textbf{Theorem \ref{thm_main}}, is presented in \textsc{Section  \ref{sec_H_dim}} to \textsc{Section \ref{sec_sub_lev_est}}.
    
\item The relationship between the general character of our methods and the various approaches employed for proving previous related results---see \textbf{Remark \ref{ext}}, is explored in \textsc{Section \ref{finrem}}.

\end{itemize}

The very essence of our approach resides, as expected, in the two critically defining features of the linear resonant Carleson-Radon transform:
\begin{itemize}
\item (C) the \emph{linear resonant Carleson behavior} encapsulated within the key modulation invariance feature described by \eqref{Modl}, or more generally, by \eqref{Modl2};

\item (R) the \emph{Radon behavior} encoded in the geometric and analytic properties of the underlying integrant manifold (in our case the paraboloid) paired with the singular properties of the kernel $K$ whose synthesis is captured in the preliminaries of \textsc{Section \ref{prelred}} by considering auxiliary Borel measures supported on the paraboloid, see \eqref{eq_mu_k_def}.
\end{itemize} 
For expository reasons, in what follows, we focus first on the one-dimensional case. After a standard linearization procedure,  
the interaction of the above two features now paired with the essential hypothesis \eqref{degone}, gives via \eqref{eq_C_VKa_symb_rep} the following crucial correspondence for the multiplier's behavior: 
\begin{itemize}
\item (C1) the Carleson feature (C) turns into a \emph{translation invariance property in the second variable $\xi_1$};

\item (R1) the Radon feature (R) manifests through the curvature in 
\begin{itemize}
\item (R1a) an \emph{anisotropic scaling behavior}, and 

\item (R1b) the \emph{preferential role of the zero frequency in the first variable $\xi_0$} (no translation invariance). 
\end{itemize}

\end{itemize}
These, in turn, naturally create the following \emph{phase plane analysis} strategy:
\begin{itemize}
\item (R1a) requires a  \emph{multiscale analysis} based on a 
\begin{equation}\label{parabres}
\textrm{family of parabolic scalings of the form}\quad\left\{\left(2^{k} \xi_0,\, 2^{2k} \xi_1\right)\right\}_{k\in\Z}
\end{equation}
which is a direct consequence of the dilation symmetries obeyed by \eqref{par};
    
\item (R1b) invites for a quantification in \eqref{parabres} of the frequency distance from the first coordinate to the origin---this becomes the ``curvature" parameter $|\lambda|$ with $\lambda\in \R$;

\item (C1) induces a partition of the frequency universe covered by the second coordinate in \eqref{parabres} into intervals of size roughly $|\lambda|$.
\end{itemize}
In view of the above argumentation, in \textsc{Section \ref{subsec_pf_main_planar}} we implement the following decomposition: for each fixed $x\in\R^2$ and $\lambda\in\R$, taking $k\in\Z$ and $|2^k\xi_0|\approx |\lambda|$, we split our analysis into three regimes depending on the relative coordinate interaction within the pair
\begin{equation}\label{parabres1}
\left(2^{k} \xi_0,\, 2^{2k} (\xi_1-a(x))\right)
\end{equation}
representing the point of reference for our multiplier after linearization---see \eqref{eq_C_VKa_symb_rep}. More precisely, our operator \eqref{par} is subdivided into three components, as follows:
 \begin{itemize}
\item (L) a low-frequency component: $|\lambda|\lesssim 1$

\item (HNS) a high-frequency non-stationary component: $|\lambda|\gtrsim 1$ and $|2^k\xi_0|\approx |\lambda|\not\approx | 2^{2k} (\xi_1-a(x))|$,
   
\item (HS) a high-frequency stationary component: $|\lambda|\gtrsim 1$ and $|2^k\xi_0|\approx |\lambda|\approx | 2^{2k} (\xi_1-a(x))|$.
\end{itemize} 

With these, we reduce the initial statement in \textbf{Theorem \ref{thm_main1D}} to providing an $L^{p}$-control on each of the above components. The latter is the subject of \textbf{Theorem \ref{prop_4cases}} and \textbf{Theorem \ref{prop_4casesdecay}} stated in \textsc{Section \ref{subsec_pf_main_planar}}.

Once here, we reach a main novel point: in contrast with the standard non-zero or hybrid curvature settings\footnote{See \emph{e.g.} the comparison with the treatment for each of (L), (HNS) and (HS) terms in the non-resonant Carleson-Radon setting proved in \cite{HsuL24} or in other general non-resonant settings for several main harmonic analysis themes as for example in \cite{LVunif}, or even in the hybrid setting corresponding to the Bilinear Hilbert transform treated in \cite{GALVHybcurves}.}, \emph{all} of the three resulting components involve zero-curvature features, and thus there is no easy or direct treatment for either of them. More precisely, confirming once again the novel features brought to light by the resonant Carleson-Radon theme, one has to develop a time-frequency modulation invariant analysis adapted to \emph{each} and \emph{all} of the three cases above!

In order to achieve the latter, we employ a \emph{new, unitary approach} consisting of four main stages:
\medskip

The \emph{first stage}, treated in \textsc{Section \ref{subsec_class}}, consists of the reduction of \textbf{Theorem \ref{prop_4cases}} to \textbf{Theorem \ref{thm_phy_shift_est_unit}} via the newly introduced concept of \emph{shifted superposed Carleson operator} according to the following procedure:
\begin{itemize} 
\item Firstly, each of the three components (L), (HNS) and (HS) may be realized as both a frequency or space parameter spatially-shifted superposed Carleson operators, that is, as a linear combination of single scale spatially-shifted Carleson-type operators with the shift parameter encoding a frequency or a space location---see \eqref{eq_cC_4ier_dec} and \eqref{eq_cC_int_dec}, respectively;\footnote{In what follows, for notational simplicity, we abuse the language and refer to a frequency or space parameter spatially-shifted superposed Carleson operators as simply frequency or space shifted superposed Carleson operators.}

\item Secondly, exploiting the essence of their definition, one shows that for the (L) and (HNS) components the coefficients arising in the linear combination decomposition of the attached frequency shifted superposed Carleson operator verify a decay estimate in the curvature parameter $|\lambda|$---this is the content of \textbf{Lemma \ref{lem_mikhlin_symb}}.   

\item Thirdly, employing a rescaling argument---see  \eqref{resc} and making use of the above two items, one reduces the $L^p$ control of the (L) and (HNS) components to the $L^p$-control of some normalized \emph{frequency shifted} superposed Carleson operators involving coefficients of essentially unitary size.
    
\item  While a similar frequency shift argument fails for the (HS) component, one notices that via the dual manifestation evoked at the first item above, the attached \emph{space shifted} superposed Carleson operator also involves coefficients of essentially unitary size.
    
\item  Using all of the above, one deduces that \textbf{Theorem \ref{prop_4cases}} follows if one can obtain the $L^p$ control of a generic shifted superposed Carleson operator in terms of the $l^{\infty}(\Z)$-norm of the coefficients involved in its superposition decomposition together with a log-loss in the shifting parameter. This latter statement is precisely the content of \textbf{Theorem \ref{thm_phy_shift_est_unit}}, completing our first stage. 
\end{itemize}

Once at this point, we remark on two things: 

\noindent \emph{i)}  on the one hand, exploiting the innate extra cancellation within the (L) and (HNS) components, we are able via \textbf{Theorem \ref{thm_phy_shift_est_unit}} to obtain in \textbf{Theorem \ref{prop_4cases}}---see \eqref{eq_LF} and \eqref{eq_HF_nst}---a $\lambda$ power decay control on their corresponding $L^p$ norms;

\noindent \emph{ii)}  on the other hand, confirming the higher complexity level of the stationary component (HS), the mechanism within \textbf{Theorem \ref{thm_phy_shift_est_unit}} can only provide a tame, $\lambda$ logarithmic blow-up of its corresponding $L^p$ norm. Consequently, one has to design a different approach---generating the proof of \textbf{Theorem \ref{prop_4casesdecay}}---in order to obtain the desired $\lambda$ power decay control.

Thus, summarizing our work up to this point, our \textbf{Theorem \ref{thm_main1D}} has been reduced to two main statements which both involve modulation invariance behavior:
\begin{itemize}
\item \textbf{Theorem \ref{thm_phy_shift_est_unit}} addressing $L^p$ bounds for a normalized\footnote{The curvature parameter $\lambda$ is no longer present} shifted superposed Carleson operator; 

\item \textbf{Theorem \ref{prop_4casesdecay}} targeting a $\lambda$ power decay in the $L^2$-bound of the stationary component (HS).
\end{itemize}

We are thus naturally brought to the \emph{second stage}: the creation of a \emph{unitary time-frequency framework} able to incorporate both settings corresponding to the two theorems above. This is performed in two steps according to the following description:
\begin{itemize}
\item at the first step, addressed in \textsc{Section \ref{subsec_Freq_disc}} and \textsc{Section \ref{subsec_Freq_disc1}}, we perform a preliminary frequency discretization of the two operators we have to control, and, via some averaging arguments determined by their classes of symmetries, we reduce our analysis to suitable discretized model operators. In the newly created context, \textbf{Theorem \ref{prop_4casesdecay}} is rephrased into \textbf{Theorem \ref{thm_HF_st_disc}} while \textbf{Theorem \ref{thm_phy_shift_est_unit}} is reduced to \textbf{Theorem  \ref{thm_phy_shift_disc}}.

\item at the second step, addressed in \textsc{Section \ref{sec_formal_tf_ana}}, we create a formal time-frequency framework adapted to \emph{variable resolution}\footnote{See also \textbf{Remark \ref{varres}}.} families of tiles with distinct shapes and volumes. This is required in order to accommodate the different hypotheses as well as aims corresponding to \textbf{Theorem \ref{thm_HF_st_disc}} and \textbf{Theorem  \ref{thm_phy_shift_disc}}. In order to achieve this level of versatility we formulate a \emph{general template}---see \textbf{Assumption \ref{ass_phy_loc_meas}} and \textbf{Assumption \ref{ass_tf_proj_emb}}, design the concept of a \emph{shifted mass} of a tile---see \textbf{Definition \ref{def_shift_mass}}, and, crucially, introduce a \emph{``Heisenberg" norm scale} associated with a generic tile which captures the trade-off between spatial localization and smoothing (with the latter being an expression of the frequency localization)---see \textbf{Definition \ref{def_sing_tile_est}}. Moreover, in the same \textsc{Section \ref{sec_formal_tf_ana}}, a general plan of dealing with \textbf{Theorem \ref{thm_HF_st_disc}} and \textbf{Theorem  \ref{thm_phy_shift_disc}} is laid out. This plan consists of four key lemmas whose proofs are postponed for \textsc{Section \ref{sec_pf_key_lemmas}} and four main propositions which are the subject of the next two stages. 
\end{itemize}

The \emph{third stage} addresses \textbf{Theorem \ref{thm_HF_st_disc}}, which, we recall, focuses on obtaining a $\lambda$ power decay in the $L^2$-bound of a model operator associated to the stationary component (HS). The resolution of this theorem relies on, yet again, two steps:
\begin{itemize}
\item a \emph{low resolution multiscale analysis} addressing the Carleson-type behavior of the (HS) component and involving a combinatorial analysis of families of wave-packets whose time-frequency localization are adapted to super-Heisenberg tiles defined in \eqref{fattile} (notice that the volume of such a tile is of order $\lambda^2$). The fact that this low-resolution multiscale analysis falls within the abstract time-frequency framework conceived at the second stage above is one of the key aspects treated in \textsc{Section \ref{sec_pf_HF_st_disc}} and has as a pinnacle the proof of \textbf{Proposition \ref{thm_simp_w22}}.
    
\item a \emph{high resolution single scale analysis} addressing the Radon-type behavior of the (HS) component and involving a \emph{smoothing estimate} which captures the behavior of each single super-Heisenberg tile. This is the subject of \textbf{Proposition \ref{singtileestimdec}} whose proof is presented in \textsc{Section \ref{sec_est_of_sing_tile}} and relies crucially on the LGC-method. It employs \emph{i)} the partition of the ambient universe represented by a super-Heisenberg tile into a family of Heisenberg tiles (volume one) with linearizing effect on the curvature of the underlying manifold appearing in the definition of the single super-Heisenberg tile operator, \emph{ii)} a sparse-uniform dichotomy, and  \emph{iii)} a level set analysis.
\end{itemize}
Thus, in a nutshell, the proof of \textbf{Theorem \ref{thm_HF_st_disc}} is achieved by decoupling the modulation invariance, that is the Carleson behavior/multiscale analysis, from the curvature, that is the Radon behavior/single tile analysis, and feeding the resulting outputs into the abstract machinery designed at the second stage.

Finally, the \emph{fourth stage} addresses \textbf{Theorem  \ref{thm_phy_shift_disc}} concerned with the  $L^p$ bounds of a normalized shifted superposed Carleson model operator. The proof of this result relies on  \textbf{Propositions \ref{thm_endpoint_infty_1}} and \textbf{\ref{thm_endpoint_1_infty}} and is delivered in \textsc{Section \ref{sec_pf_physift_disc}}. In an antithetical note to the reasonings presented at the third stage, we have:
\begin{itemize}
\item the proof of \textbf{Theorem  \ref{thm_phy_shift_disc}} involves a single\footnote{The low and high resolution from the third stage become indistinguishable in the $\lambda=1$ setting present in the fourth stage.}, \emph{standard resolution multiscale analysis} addressing the Carleson-type behavior of the operator under analysis. This analysis involves Heisenberg tiles (of volume one) and is performed by feeding its particular features into the  abstract time-frequency framework conceived at the second stage;

\item however, the single resolution simplification is counterbalanced by the presence of a shifting parameter\footnote{This is absent in the reasonings presented at the third stage.} whose careful consideration requires the development of some new techniques.
\end{itemize}
\medskip

The treatment of the general dimensional case certifies the versatility and robustness of our methods: while facing a key new challenge to be discussed below, the required modifications---though far from trivial---are all fitting within the general four-stage approach presented in the one-dimensional setting.   

Indeed, after some standard reductions---see \eqref{red}, and departing from the intuition presented in \eqref{parabres}, the natural correspondent of \eqref{parabres1} becomes
\begin{equation}\label{parabres10}
\left(2^{k} \xi_0,\, 2^{k} (\xi_1-a_1(x)),\,\ldots,\,2^{k} (\xi_{D-1}-a_{D-1}(x)),\,2^{2k} (\xi_D-a_D(x))\right)\,,
\end{equation}
while the analog three-component subdivision of \eqref{def_C_ast_V_K} is given by 
\begin{itemize}
\item (L) a low-frequency component: $|\lambda|\lesssim 1$

\item (HNS) a high-frequency non-stationary component: $|\lambda|\gtrsim 1$ and $|\lambda|\lesssim \left| 2^{2k} (\xi_D-a_D(x))\right|$,
   
\item (HS) a high-frequency stationary component: $|\lambda|\gtrsim 1$ and $|\lambda|\gtrsim \left| 2^{2k} (\xi_D-a_D(x))\right|$.
\end{itemize} 
With these settled, in perfect symmetry with the one-dimensional procedure of reducing  \textbf{Theorem \ref{thm_main1D}} to \textbf{Theorem \ref{prop_4cases}} and \textbf{Theorem \ref{prop_4casesdecay}}, in the general dimensional case---see \textsc{Section \ref{Gendimred}}---we maintain the same strategy by aiming to reduce \textbf{Theorem \ref{thm_main}} to \textbf{Theorem \ref{prop_Hdim_4cases}} and \textbf{Theorem \ref{prop_Hdim_4cases1}}.

Indeed, following the same spirit within the one-dimensional four-stage approach described earlier, at the \emph{first stage}---completed in \textsc{Section \ref{ssco}}---we reduce \textbf{Theorem \ref{prop_Hdim_4cases}} to \textbf{Theorem \ref{thm_log_shift_est_D_dim}} via the same concept of a shifted superposed Carleson operator. The proof of  \textbf{Theorem \ref{thm_log_shift_est_D_dim}} follows the same ideas as those presented for proving its one-dimensional version in the form of \textbf{Theorem \ref{thm_phy_shift_est_unit}}. Thus, the first and fourth stages from the $D=1$ case transfer with only minor changes to the general dimensional setting. Moreover, the \emph{second stage} corresponding to the  \emph{formal time-frequency framework} remains valid in all dimensions. Consequently, the main remaining task is the completion of the \emph{third stage} corresponding to the treatment of the (HS) component, which is the subject of  \textbf{Theorem \ref{prop_Hdim_4cases1}}. 

This is the stage where one encounters a genuine novel challenge manifesting only in dimension $D>1$, which, beyond the expected higher-dimensional technical adaptations, will require some new ideas. However, in order to be able to provide a more conceptual and intuitive description of this new enemy, we have to make a short detour and provide a very brief philosophical outlook on the three known proofs of the classical Carleson operator:
\begin{itemize}
\item \emph{the original proof of Carleson}, \cite{c1}, used in an implicit manner an energy selection algorithm, namely, a quantification of the size of some suitable weighted local Fourier coefficients corresponding to the input function; the above can be re-interpreted in the form of a hierarchy of BMO-type norms attached to shifted, rescaled and modulated versions of the above-mentioned input function. 

\item \emph{the next proof due to Fefferman}, \cite{f}, introduced the modern language of trees and tiles and, in contrast with Carleson's approach, focused on the study of the linearizing function (derived from linearizing the supremum) as opposed to the input function. Specifically, the author designs a mass selection algorithm with the concept of mass of a tile capturing the ``amount" of the graph of the linearizing function contained in the corresponding tile.

\item \emph{finally, the third proof due to Lacey and Thiele},  \cite{lt3}, combines and streamlines the above two approaches via a dualization argument: indeed, their approach employs both a mass and an energy selection argument. During this process, a standard Heisenberg tile is split into an (+) upper half subtile and an (-) lower half subtile, with the former capturing the time-frequency localization of the input function---the key ingredient within the energy concept, and the latter encoding the time-frequency localization of the linearizing function---the key ingredient within the mass concept. Consequently, the necessary analysis of the tree structures is split accordingly into $(+)-$trees and $(-)-$trees, respectively. 
\end{itemize}

We end our brief historical detour with the following general comment

\begin{remark}[\textsf{Template utilized in the present paper}]\label{Temp}
While many of our initial advancements in the realization of this paper were conceived within the Fefferman framework, our final steps, together with the current format of our presentation, are delivered within the Lacey-Thiele framework for two main reasons: 
\begin{itemize}
\item a \emph{specific} one: as noticed in \textbf{Remark \ref{IncorpCR}}, our approach incorporates the non-resonant treatment provided in  \cite{HsuL24} in the form of the single (super-Heisenberg) tile decay estimate completed in \textsc{Section \ref{sec_est_of_sing_tile}}; this treatment is better adapted to the  Lacey-Thiele framework due to the format of \eqref{decste} which alludes to the mass-energy duality;
\item a \emph{general} one: the symbiosis between the Carleson and Radon behaviors captured within the formal time-frequency analysis developed in \textsc{Section \ref{sec_formal_tf_ana}} encodes, yet again, an effective mass-energy duality.     
\end{itemize}
\end{remark}

With all these settled, we are now ready to reveal the new higher-dimensional challenge, which, for expository reasons, is analyzed antithetically:
\begin{itemize}
\item in the $D=1$ case,  the behavior of the (HS) component captured in \eqref{eq_HF_st} focuses on the $j=0$ setting in \eqref{eq_cC_j_lambda_def}, which in turn, via the structure of \eqref{parabres1}, involves in the time-frequency analysis of (HS) only  $(-,0)-$trees---for more on this see \emph{e.g.} \textbf{Observation \ref{obs_total_L2_energy_bd}};

\item in the $D>1$ case, the behavior of the (HS) component captured by \eqref{eq_prop_Hdim_HF_st} focuses on \emph{all} $0\leq j<D$ situations in \eqref{eq_Hdim_C_V_decomp}, which, via the structure of \eqref{parabres10}, involves in the time-frequency analysis of (HS) both $(-,j)-$trees and $(+,j)-$trees whenever $j\not\in\{0,D\}$---for more on this see \emph{e.g.} the introductory paragraph in \textsc{Section \ref{subsec_Hdim_issue_j=1}}.\footnote{Due to the symmetries of our operator, the study of the (HS) component needs only to address the $j=0$ and $j=1$ cases, and hence, the main new challenge may be reduced to the $j=1$ case.} 
\end{itemize}

Far from being a mere technicality, the task of dealing with $(+,1)-$trees requires a reconstruction with old and new elements of the entire edifice on which our one-dimensional approach was built: 
\begin{itemize}
\item the ``original sin" is captured within the essential fact that the main statements provided by \textbf{Lemma \ref{lem_gen_tree_c_half}} and \textbf{Proposition \ref{thm_simp_w22}} fail in the new context, in the sense that, the upper bounds in those estimates (which involve a suitable quantity $\nrm{\Lambda}_{\frac{1}{2}}$) blow up as $|\lambda|$ grows to infinity---see \eqref{blowup}.

\item in order to solve this difficulty one has to introduce a dual high-resolution analysis within the low-resolution analysis of the (HS) component---see \textbf{Remark \ref{dualhres}} and also \textbf{Remark \ref{varres}}.

\item indeed, on top of the high-resolution analysis employed via the LGC method for the single scale decay estimate---see the proof of \textbf{Proposition \ref{thm_Hdim_sing_tile_ests}} and \textbf{Lemma \ref{lem_sub_lev_2_2}} in \textsc{Section \ref{singltilesestimgen}} and \textsc{Section \ref{sec_sub_lev_est}}, respectively---we design a second high-resolution analysis relying on the new concept of \emph{frequency sub-systems}---see \textbf{Definition \ref{def_Hdim_freq_sub}}---which allows us to divide each super-Heisenberg tile in \eqref{superHeis} into an outer and central component---see \eqref{eq_Hdim_freq_proj_for_subsys}---\eqref{eq_cent_n_out_nrm}.
    
\item via this second high-resolution analysis, we are able to solve the issue presented at the first item by formulating the proper extension of \textbf{Lemma \ref{lem_gen_tree_c_half}} as given by \textbf{Lemma \ref{lem_Hdim_gen_tree}}, with the proof of the latter displayed in \textsc{Section \ref{sgtreejone}}

\item the final step for treating the $D>1$ case of the (HS) component in \textbf{Theorem \ref{thm_Hdim_HF_st_disc}} requires a new form of the energy selection algorithm---see \textbf{Lemma \ref{lem_Hdim_energy_sel}}---which has to reflect the introduction of a second high-resolution analysis; the treatment of the latter is offered in \textsc{Section \ref{sec_pf_gen_energy_est}}.
\end{itemize}

We end our introductory section with the following panoramic remark concerning a key concept\footnote{The first occurrence of a low-high resolution as part of a multiscale analysis problem is in  \cite{LVBilHilbCarl} followed by \cite{GALVHybcurves}.} in the present paper: \emph{the variable (adaptive) resolution} and its evolution within the dimensional paradigm:

\begin{remark}[\textsf{Variable resolutions: one versus higher-dimensional setting}]\label{varres}
Our variable resolution multiscale analysis, developed in accordance with the formal time-frequency framework, encompasses, in an antithetical description, the following resolutions:
\begin{itemize}    
\item \emph{in the $D=1$ case, we have three resolutions} with the first one focusing\footnote{Strictly speaking, the first resolution also provides tame $L^p$-bounds for the (HS) component.} on the (L) and (HNS) components while the second and third ones addressing the (HS) component:\smallskip  
\begin{enumerate}
\item the \emph{standard resolution}\footnote{For the scale $k=0$ the generic tile $P=\vI\times\vomega$ is given by $\vI=I_1\times I_2$ and $\vomega=\omega_1\times\omega_2$ with $I_j=[\a_j,\a_j+1]$ and $\omega_1=[-\frac{1}{2},\frac{1}{2}]$,  $\omega_2=[w,w+1]$ with $\a_1,\,\a_2,\,w\in\Z$.} involving Heisenberg tiles of volume one---defined by taking in \eqref{fattile} the parameter $\lambda=1$---is employed for proving \textbf{Theorem \ref{prop_4cases}}, or, more precisely, its reduced form \textbf{Theorem  \ref{thm_phy_shift_disc}}---see \textsc{Section \ref{sec_pf_physift_disc}};
\smallskip     
\item the \emph{$\lambda-$low resolution}\footnote{For the scale $k=0$ the generic tile $P=\vI\times\vomega$ is given by $\vI=I_1\times I_2$ and $\vomega=\omega_1\times\omega_2$ with $I_j=[\a_j,\a_j+1]$ and $\omega_1=[-\frac{\lambda}{2},\frac{\lambda}{2}]$,  $\omega_2=[w\lambda,(w+1)\lambda]$ with $\a_1,\,\a_2,\,w\in\Z$.} involving super-Heisenberg tiles of volume $\lambda^2$---defined directly by \eqref{fattile}---is used for treating the Carleson behavior/multiscale component in \textbf{Theorem \ref{prop_4casesdecay}} which is performed in \textsc{Section \ref{sec_pf_HF_st_disc}};\smallskip      
\item the \emph{$\sqrt{\lambda}-$high resolution}\footnote{For the scale $k=0$ the generic tile $P=\vI\times\vomega$ is given by $\vI=I_1\times I_2$ and $\vomega=\omega_1\times\omega_2$ with $I_j=[\frac{\a_j}{\sqrt{\lambda}},\frac{\a_j+1}{\sqrt{\lambda}}]$ and $\omega_1=[v \sqrt{\lambda},(v+1)\sqrt{\lambda}]$,  $\omega_2=[w\sqrt{\lambda},(w+1)\sqrt{\lambda}]$ with $\a_1,\,\a_2,\,v,\,w\in\Z$ and $-\frac{\sqrt{\lambda}}{2}\leq v< \frac{\sqrt{\lambda}}{2}$.} (or simply the \emph{LGC-resolution}) involving Heisenberg tiles of volume one---defined by the equipartition of each fixed super-Heisenberg tile into $\lambda^2$ tiles whose shape match the time-frequency localization dictated by \eqref{eq_def_varphi_vomega_xyuv}-\eqref{eq_planar_tf_proj}---is employed for addressing the Radon behavior/single scale component in \textbf{Theorem \ref{prop_4casesdecay}}, namely the single (super-Heisenberg) tile decay estimate in \textbf{Proposition \ref{singtileestimdec}} treated in \textsc{Section \ref{sec_est_of_sing_tile}}. 
\smallskip    
\end{enumerate}

\item \emph{in the $D>1$ case, we have four resolutions} with the first one applied as before for the (L) and (HNS) components and the other three all addressing the (HS) component. While the first three resolutions are simply the general dimensional versions of those presented in the one-dimensional setting above, the fourth one addresses precisely the higher-dimensional challenge concerning  $(+,1)-$tree structures. Indeed, we have:
\begin{enumerate}    
\item the \emph{standard resolution} involving Heisenberg tiles of volume one---defined by taking in \eqref{superHeis} the parameter $\lambda=1$---is employed for proving \textbf{Theorem \ref{prop_Hdim_4cases}}, or, more precisely, its reduced form \textbf{Theorem  \ref{thm_log_shift_est_D_dim}}---see \textsc{Section \ref{ssco}};
    
\item  the \emph{$\lambda-$low resolution} involving super-Heisenberg tiles of volume $\lambda^{D+1}$---defined by \eqref{eq_Hdim_W_n_Wk}, \eqref{spinterv} and \eqref{superHeis}---is used for treating the Carleson behavior/multiscale component in \textbf{Theorem \ref{prop_Hdim_4cases1}}, or, equivalently after some reductions, \textbf{Theorem \ref{thm_Hdim_HF_st_disc}}, see  \textsc{Section \ref{redmain}} and \textsc{Section \ref{sec_pf_Hdim_HF_st_mod}};\smallskip

\item the \emph{$\sqrt{\lambda}-$high resolution} (the \emph{LGC-resolution}) involving Heisenberg tiles of volume one---defined by the equipartition of each fixed super-Heisenberg tile into $\lambda^{D+1}$ tiles with a shape matching the time-frequency localization of the wave-packet \eqref{wvpa}---is employed for addressing the Radon behavior/single scale component in \textbf{Theorem \ref{thm_Hdim_HF_st_disc}}, namely the single (super-Heisenberg) tile decay estimate in \textbf{Proposition \ref{thm_Hdim_sing_tile_ests}} treated in \textsc{Section \ref{singltilesestimgen}}.

\item the \emph{$(1,+)-$high resolution} involving super-Heisenberg tiles of volume $\lambda^{D-1}$---defined by the equipartition of each fixed original super-Heisenberg tile (see the second item above) into $\lambda^{2}$ tiles dictated by \eqref{eq_sticks}---is employed for addressing the presence of $(+)-$tree structures within the analysis of the higher-dimensional (HS) component and is applied in the proof of the single tree estimate---see \textbf{Lemma \ref{lem_Hdim_gen_tree}} in \textsc{Section \ref{sgtreejone}} as well as in the energy selection algorithm---see \textbf{Lemma \ref{lem_Hdim_energy_sel}} and its proof in \textsc{Section \ref{sec_pf_gen_energy_est}}. 
\end{enumerate}

\end{itemize}
    
\end{remark}

\subsection*{Acknowledgments}
We thank Christoph Thiele for a nice, motivational discussion on the one-dimensional resonant Carleson-Radon transform which took place more than a decade ago when the second author visited University of Bonn. Also we thank Lars Becker for helpful comments and suggestions on an earlier version of this paper.

Both authors were supported by the NSF grant DMS-2400238 and by the Simons Travel grant MPS-TSM-00008072.

\section{The one-dimensional setting: Preliminary decompositions and reductions}\label{prelred}
\subsection{Some introductory notations and definitions}

\begin{definition}[Symmetry operations]\label{def_sym_ops}
    Given \(\vx,\vxi\in \R^{D+1}\), \(\vC\in GL_{D+1}\br{\R}\), and \(p\in\bR{0,\infty}\), we define for a function \(\Psi:\R^{D+1}\to\C\) the following:
    \begin{equation}
        \Dil^p_\vC
        \Psi\br{\vz}:=
        \abs{\det \vC}^{-1/p}\Psi\br{\vC^{-1}\vz},
    \end{equation}
    \begin{equation}
        \Tr_\vx \Psi\br{\vz}:=
        \Psi\br{\vz-\vx},
        \quad
        \Mod_\vxi \Psi\br{\vz}:=
        e\br{\vxi\cdot\vz}\Psi\br{\vz}.
    \end{equation}
    We also introduce for \(\vR\in\br{0,\infty}^{D+1}\) the matrix:
    \begin{equation}
        \vLambda_\vR:=
        \begin{pmatrix}
            R_0 & \vnull & 0\\
            \vnull & \ddots & \vnull\\
            0 & \vnull & R_D
        \end{pmatrix}
    \end{equation}
    and the shorthand notation:
    \begin{equation}
        \Dil^p_\vR \Psi\br{\vz}:=
        \Dil^p_{
            \vLambda_\vR
        }\Psi\br{\vz}
        =
        \br{\prod^D_{j=0} R^{-1/p}_j}\Psi\br{z_0/R_0,\dots,z_D/R_D}.
    \end{equation}
\end{definition}
\begin{definition}[\textsf{Functions adapted to an interval}]\label{def_Psi_I}
    Let \(\Psi\) be a function on \(\R^{D+1}\), \(\psi\) be a function on \(\R\), and \(\vI\subset \R^{D+1}\) be an interval centered at \(\vc:=\br{c_0,\dots,c_D}\) with dimension \(l_0\times\cdots\times l_D\). Unless specified otherwise, the two expressions \(\Psi_\vI\) and \(\psi_\vI\) mean:
    \begin{equation*}
        \Psi_\vI\br{\vz}:=\Psi\br{\frac{z_0-c_0}{l_0},\dots,\frac{z_D-c_D}{l_D}},\quad
        \psi_\vI:=\br{\psi^{\otimes (D+1)}}_\vI.
    \end{equation*}
\end{definition}
\begin{definition}[\textsf{Japanese bracket and various weights}]\label{def_jap_n_weight}
    Let \(\ang{\cdot}:=\sqrt{1+\abs{\cdot}^2}\) denote the Japanese bracket on \(\R\). For \(\vz\in\R^{D+1}\), we write \(\ang{\vz}_\otimes:=\prod_{j=0}^D\ang{z_j}\). With these, we have
\begin{itemize}   
\item for $z\in\R$  we let \(\chi\br{z}:=\ang{z}^{-1}\) be a one-dimensional positive weight  on \(\R\);
\item for \(\vz\in\R^{D+1}\) we let $\chi\br{\vz}:=(\ang{\vz}_\otimes)^{-1}$;
\item using Definition \ref{def_Psi_I}, for \(\vI\subset \R^{D+1}\) and \(\vz\in\R^{D+1}\) we make sense of the $\vI-$adapted weight $\chi_{\vI}\br{\vz}$.
\end{itemize} 
\end{definition}
\begin{definition}[\textsf{Half intervals}]\label{def_half_I}
    Let \(I=\bR{\inf I,\sup I}\).
    We denote the center of interval \(c_I:=\frac{\inf I+\sup I}{2}\) and:
    \begin{equation*}
        I^+:=\bR{c_I,\sup I}
        ,\quad
        I^-:=\bR{\inf I,c_I}.
    \end{equation*}
    In general, for interval \(\vI:=\prod_{j=0}^DI_j\subset \R^{D+1}\) and \(k\in\{0,\dots,D\}\), we set:
    \begin{equation*}
    c_\vI:=\br{c_{I_0},\dots,c_{I_D}},\quad
        \vI^{\pm,k}:=\prod_{j<k}I_j\times I^\pm_k\times\prod_{j>k}I_j.
    \end{equation*}
\end{definition}

\begin{definition}[\textsf{The diameter of a set}]\label{def_diam}
     If \(D\in\N\) and \(E\subset\R^{D}\) is a measurable set we denote by $\dia(E)$ the diameter of the set $E$, that is
\begin{equation*}
        \dia (E):=\sup_{x,y\in E} |x-y|\,.
    \end{equation*}
\end{definition}

\subsection{Kernel discreization: sequence of measures on paraboloid}\label{subsec_kernel_decomp}
In this section, for expository reasons, we carry out our discussion in the general $D\geq 1$ dimensional setting.

We start by considering a partition of unity on \(\br{0,\infty}\) generated by \(\rho\in C^\infty_c\br{\br{0,\infty}}\) and satisfying
\begin{equation}\label{eq_parti_uni}
    \1_{\br{0,\infty}}=\sum_{k\in\Z}\rho_k,\quad
    \rho_k\br{t}:=\rho\br{t/2^k},\quad
    \supp \rho\subset \br{1/2,2}.
\end{equation}

Recalling now the definition of the maximal operator along the paraboloid
    \begin{equation}\label{mparab}
    M_\para f\br{\vx}
    :=
        \sup_{R>0}
        \fint_{\abs{\vt}\leq R}
            \abs{f}\br{\vx-\vX\br{\vt}}
        d\vt
    \end{equation}
and using \eqref{def_C_ast_V_K} and \eqref{mparab} a standard argument shows that
$$CR^\ast_\V  f\lesssim 
        \cC \cR^\ast_\V  f
        +M_\para f\,$$
where \(\cC \cR^\ast_\V \) is a scale-discretized variant
    \begin{equation*}
        \cC \cR^\ast_\V  f\br{\vx}:=
        \sup_{\substack{\uk<\ok\\ \va\in\V}}
            \abs{
                \sum_{k=\uk}^{\ok}
                    \int_{\R^D}
                        f\br{\vx-\vX\br{\vt}}
                        e\br{\va\cdot \vX\br{\vt}}
                        \rho_k\br{\abs{\vt}/\sqrt{D}}
                        K\br{\vt}
                    d\vt
            }.
    \end{equation*}

Moreover, by considering the following auxiliary Borel measure supported on the paraboloid:
\begin{equation}\label{eq_mu_k_def}
    \mu_k\br{E}:=\int_{\R^D} \1_E\br{\vX\br{\vt}}\rho_k\br{\abs{\vt}/\sqrt{D}}K\br{\vt}
    d\vt,\quad E\subset \R^{D+1},
\end{equation}
we can further reformulate \(\cC\cR^\ast_{\V}\) in a manner resembling the classical setting:
\begin{equation*}
    \cC \cR^\ast_\V f\br{\vx}=
    \sup_{\substack{\uk<\ok\\ \va\in\V}}
    \abs{
        \int
            f\br{\vx-\vtau}
            e\br{\va\cdot\vtau}
        \sum_{k=\uk}^{\ok}
            d\mu_k\br{\vtau}
    }
    .
\end{equation*}

Standard reasonings reduce now the study of $\cC R^\ast_\V$ to that of

\begin{equation}\label{eq_1st_linearized_sum}
    \cC\cR^{\br{\va}}_\V  f\br{\bm{x}}
    :=
    \int
        f\br{\vx-\vtau}
        e\br{\va\br{\vx}\cdot\vtau}
    \sum_{k=\uk}^{\ok\br{\bm{x}}}
        d\mu_k\br{\vtau},
\end{equation}
where here both $\ok\br{\cdot}: \R^{D+1}\to\R$ and  $\va(\cdot):\R^{D+1}\to\V$ are arbitrary measurable functions while \(\uk\in\Z\) is a fixed parameter.

We thus conclude now that both Theorems \ref{thm_main1D} and \ref{thm_main} follow if we are able to show that for all singular kernels \(K\) satisfying \eqref{eq_kernel_cond}, all linear subspaces \(\V\leq\R^{D+1}\) satisfying \eqref{degone} and any $\va(\cdot):\R^{D+1}\to\V$ measurable function one has
\begin{equation}\label{ref}
\nrm{\cC\cR^{\br{\va}}_\V f}_{L^p}\lesssim\nrm{f}_{L^p}
\end{equation}
with the implicit constant in \eqref{ref} independent of the choice of the functions \(\ok\br{\cdot}, \va\br{\cdot}\) and of \(\uk\in\Z\). 

\subsection{Symbol discretization}\label{subsec_symb_decomp}

In this section, we focus on the planar setting, leaving for later--- see \textsc{Sections \ref{sec_H_dim}--\ref{sec_pf_gen_energy_est}}---the necessary adaptations to the general dimensional case \(D\geq 1\). Consequently, from now on, we fix \(D=1\) and notice that in view of \eqref{degone} the only interesting case is when \(\V=\BR{0}\times\R\). 

Consider thus the multiplier/symbol formulation of \(\cC\cR^{\br{\va}}_{\BR{0}\times\R} \):
\begin{equation}\label{eq_C_VKa_symb_rep}
    \cC\cR^{\br{\va}}_{\BR{0}\times\R}  f\br{\bm{x}}=
    \int\widehat{f}\br{\vxi}
    \left(\sum_{k=\uk}^{\ok\br{\vx}}
        \widehat{\mu_k}\br{\vxi-\va\br{\vx}}\right)
    e\br{\vxi \cdot\vx}
    d\vxi
    ,\quad \va\br{\vx}=\br{0,a_1\br{\vx}}=\br{0,a\br{x_0,x_1}}.
\end{equation}
Following the natural intuition, we start by performing a decomposition on the frequency information according to the stationary/non-stationary phase contribution of the oscillatory integral expression
\begin{equation}
    \widehat{\mu_k}\br{\vxi}=
    \int   
        \overline{e\br{\xi_0 t +\xi_1 t^2}}\,
    \rho_k\br{\abs{t}} K\br{t} dt.
\end{equation}
This is achieved by implementing a cone decomposition:
\begin{equation}\label{conedec}
\end{equation}
\begin{itemize}
\item fix an even bump function \(\phi\in C^\infty_c\br{\R}\) satisfying \(\1_{\mr{-0.005,0.005}}\leq\phi \leq \1_{\br{-0.01,0.01}}\);

\item define the function \(\varphi:=\nrm{\phi}_{L^1}^{-1}\cdot\1_{\mr{-\frac{1}{2},\frac{1}{2}}}\ast \phi\in C^\infty_c\br{\R}\) and notice that
\(\1_{\br{-0.49,0.49}}\leq \varphi \leq \1_{\br{-0.51,0.51}}\);

\item finally, consider some auxiliary smooth functions \(\tarphi,\tpsi\) obeying
\begin{equation}\label{eq_tar_tpsi_supp}
    \supp \tarphi\subset \br{-0.26,0.26},\quad
    \supp \tpsi\subset \br{0.49,0.51}\,.
\end{equation}
\end{itemize}
With these, we have
\begin{equation}\label{eq_cone_dec}
    \1_{\R^2\setminus\BR{0}}\br{\vxi}:=
    \int_{\R\setminus\BR{0}}
   \left( c_0\,\phi\br{\frac{2\xi_0}{\lambda}-\frac{1}{2}}\varphi\br{\frac{\xi_1}{\lambda}}
    + c_1\,\tarphi\br{\frac{\xi_0}{\lambda}}\tpsi\br{\frac{\xi_1}{\lambda}} \right)\,\frac{d\lambda}{\abs{\lambda}},
\end{equation}
where here \(c_0,c_1\eqsim 1\) are some suitable normalizing constants. 

Since the exact choice of \(\phi,\varphi,\tarphi,\tpsi\) is not relevant for the present argument, for expository reasons, whenever \(\phi_0,\phi_1\) are supported on a small neighborhood of \(0\), we set
\begin{equation*}
    \chi_j\br{\vxi}:=\left(\phi_0\otimes\phi_1\right)\br{\vxi-C_j\ve_j},\quad
    \ve_0:=\br{1,0},\quad
    \ve_1:=\br{0,1},\quad C_j\in\br{0,1}\,.
\end{equation*}

With these we rewrite \eqref{eq_cone_dec} in the simplified form
\begin{equation}\label{eq_cone_abst}
     \1_{\R^2\setminus\BR{0}}\br{\vxi}:=
    \sum_{j\in\{0,1\}}
    c_j
    \int_{\R\setminus\BR{0}}
    \chi_j\br{\frac{\vxi}{\lambda}}
    \frac{d\lambda}{\abs{\lambda}}.
\end{equation}
Accordingly, involving a suitable parabolic scaling factor, we decompose the symbol
\begin{equation}\label{eq_symb_lamb_dec}
    \widehat{\mu}_k=
    \sum_{j\in\{0,1\}}
    c_j
    \int_{\R\setminus\BR{0}}
    \widehat{\mu_{k,\lambda}^{\br{j}}}
    \frac{d\lambda}{\abs{\lambda}},
    \quad
    \widehat{\mu_{k,\lambda}^{\br{j}}}\br{\vxi}
    :=
    \widehat{\mu}_k\br{\vxi}
    \chi_j\br{\frac{2^k\xi_0}{\lambda},\frac{4^k\xi_1}{\lambda}}
\end{equation}
with the natural operator correspondent
\begin{equation}\label{eq_C_V_decomp}
    \cC\cR^{\br{\va}}_{\BR{0}\times\R} f\br{\vx}=
    \sum_{j\in\{0,1\}}
    c_j
    \int_{\R\setminus\BR{0}}
        \cC\cR^{\br{\va}}_{\lambda,j} f\br{\vx}
    \frac{d\lambda}{\abs{\lambda}}
\end{equation}
where
\begin{equation}\label{eq_cC_j_lambda_def}
    \cC\cR^{\br{\va}}_{\lambda,j} f\br{\vx}:=
    \int \widehat{f}\br{\vxi}\left(\sum_{k=\uk}^{\ok\br{\vx}}\widehat{\mu^{\br{j}}_{k,\lambda}}\br{\vxi-\va\br{\vx}}\right)
    e\br{\vxi\cdot\vx}
    d\vxi.
\end{equation}

\subsection{Further discretization of \texorpdfstring{\(\cC\cR^{\br{\va}}_{\BR{0}\times\R}\)}{}: Reduction of Theorem \ref{thm_main1D} to Theorems  \ref{prop_4cases} and \ref{prop_4casesdecay}.}\label{subsec_pf_main_planar}

In this section we decompose $\cC\cR^{\br{\va}}_{\BR{0}\times\R}$ into three sub-operators:
\begin{itemize}
\item a low-frequency component,

\item a high-frequency non-stationary component,

\item a high-frequency stationary component,
\end{itemize}
and reduce our Theorem \ref{thm_main1D} to the two results stated below, encapsulating the key estimates for each of the above components:

\begin{theorem}\label{prop_4cases}
    There exists a universal constants \(\epsilon>0\) such that for any measurable functions \(\va:=\br{0,a}\), \(\ok\), integers \(\uk\in\Z\) and index \(1<p<\infty\), the following estimates hold:
    \begin{itemize}
        \item \textsf{Low-frequency} component---polynomial decay: given \(\abs{\lambda}\lesssim 1\) one has
        \begin{equation}\label{eq_LF}
            \nrm{
            \cC\cR^{\br{\va}}_{\lambda,j}f
            }_{L^p}
            \underset{p}{\lesssim} \abs{\lambda}^\epsilon  
            \nrm{f}_{L^p},\qquad\forall\:\:j\in\{0,1\}.
        \end{equation}
        \item \textsf{High-frequency non-stationary} component---polynomial decay: given \(\abs{\lambda}\gtrsim 1\) one has
        \begin{equation}\label{eq_HF_nst}
                \nrm{
                \cC\cR^{\br{\va}}_{\lambda,1}f
                }_{L^p}
                \underset{p}{\lesssim}  \abs{\lambda}^{-\epsilon} \nrm{f}_{L^p}.
        \end{equation}
        \item \textsf{High-frequency stationary} component---tame log growing estimate: given \(\abs{\lambda}\gtrsim 1\) one has
      
        \begin{equation}\label{eq_HF_log}
                \nrm{
                \cC\cR^{\br{\va}}_{\lambda,0}f
                }_{L^p}
                \underset{p}{\lesssim} \log^2\br{e+\abs{\lambda}} \nrm{f}_{L^p}.
        \end{equation} 
        \end{itemize}
\end{theorem}

\begin{theorem}\label{prop_4casesdecay}
   Under the same notations and conditions as above, the \textsf{high-frequency stationary} component obeys the following $L^2$-polynomial decay estimate: given \(\abs{\lambda}\gtrsim 1\) one has
   \begin{equation}\label{eq_HF_st}
                \nrm{
                \cC\cR^{\br{\va}}_{\lambda,0}f
                }_{L^{2,\infty}} 
                \lesssim \abs{\lambda}^{-\epsilon} \nrm{f}_{L^2}.
   \end{equation}
\end{theorem}

Assuming for a moment that both \textbf{Theorems \ref{prop_4cases}} and \textbf{\ref{prop_4casesdecay}} hold we notice that indeed \textbf{Theorem \ref{thm_main1D}} follows immediately:
\begin{proof}[\textbf{Proof of} \textbf{Theorem \ref{thm_main1D}}] Recalling \eqref{eq_C_V_decomp}, triangle inequality yields:
    \begin{equation*}
        \nrm{
        \cC\cR^{\br{\va}}_{\BR{0}\times\R}f
        }_{L^p}
        \lesssim
        \int_{0<\abs{\lambda}\lesssim 1}
            \sum_{j\in\{0,1\}}
            \nrm{
                \cC\cR^{\br{\va}}_{\lambda,j}f
            }_{L^p}
        \frac{d\lambda}{\abs{\lambda}}
        +
        \sum_{j\in\{0,1\}}
        \int_{\abs{\lambda}\gtrsim 1}
            \nrm{
                \cC\cR^{\br{\va}}_{\lambda,j}f
            }_{L^p}
        \frac{d\lambda}{\abs{\lambda}}=I_1+I_2.
    \end{equation*}
    To treat \(I_1\), we apply \eqref{eq_LF} and obtain:
    \begin{equation*}
        I_1\underset{p}{\lesssim}
        \int_{0<\abs{\lambda}\lesssim 1}
            \abs{\lambda}^\epsilon \nrm{f}_{L^p}
        \frac{d\lambda}{\abs{\lambda}}
        \underset{\epsilon}{\eqsim}
        \nrm{f}_{L^p}.
    \end{equation*}
    As for the treatment of \(I_2\), from \eqref{eq_HF_nst}--\eqref{eq_HF_st} and real interpolation (for some \(\epsilon_p>0\)), we have that 
    \begin{equation*}
        \nrm{
            \cC\cR^{\br{\va}}_{\lambda,j}f
        }_{L^p}
        \underset{p}{\lesssim}
        \abs{\lambda}^{-\epsilon_p}
        \nrm{f}_{L^p},\qquad \forall\:\:\abs{\lambda}\gtrsim 1\:\:\textrm{and}\:\:j\in\{0,1\}
    \end{equation*}
and thus     
\begin{equation*}
        I_2\underset{p}{\lesssim}
        \int_{\abs{\lambda}\gtrsim 1}
            \abs{\lambda}^{-\epsilon_p}
            \nrm{
                f
            }_{L^p}
        \frac{d\lambda}{\abs{\lambda}}
        \underset{\epsilon_p}{\eqsim}
        \nrm{f}_{L^p},
    \end{equation*}
    which completes the proof.
\end{proof}
Thus, in what follows, we focus on proving \textbf{Theorems \ref{prop_4cases}} and \textbf{\ref{prop_4casesdecay}}.

\subsection{Reduction of Theorem \ref{prop_4cases} to Theorem \ref{thm_phy_shift_est_unit}: shifted superposed Carleson operators}\label{subsec_class}
In this section we utilize the smoothness of the symbol \(\widehat{\mu^{\br{j}}_{k,\lambda}}\) in order to provide a global treatment of the three statements \eqref{eq_LF}, \eqref{eq_HF_log} and \eqref{eq_HF_nst} within \textbf{Theorem \ref{prop_4cases}} via the boundedness of some general shifted superposed Carleson operator model with the latter being the content of \textbf{Theorem \ref{thm_phy_shift_est_unit}} (or equivalently \textbf{Theorem \ref{thm_log_shift_est}}).

Now, in order to achieve the above-mentioned goal, we apply two decomposition types:
\begin{itemize}
\item \textsf{a frequency decomposition} via a Fourier series argument adapted to the support\footnote{Recall \eqref{eq_tar_tpsi_supp}. One may choose \(\chi_j\) to be different from those in \eqref{eq_cone_dec} and \eqref{eq_cone_abst}. Also the equality in \eqref{eq_sym_4ier_sum} must be understood in terms of a superposition of similar looking objects---these are standard Fourier analysis adaptations that are left to the interested reader.} of the symbol: 
\begin{equation}\label{eq_sym_4ier_sum}
    \widehat{\mu^{\br{j}}_{k,\lambda}}
    \br{\vxi}
    =
    \sum_{\vu\in \Z^2}
    m^{\br{j}}_{k,\lambda,\vu}\,
    \overline{e\br{\frac{2^k\xi_0 u_0}{\lambda}+\frac{4^k\xi_1 u_1}{\lambda}}}\,
    \chi_j\br{\frac{2^k \xi_0}{\lambda},\frac{4^k\xi_1}{\lambda}}.
\end{equation}
\item \textsf{a spatial decomposition} employing the integral definition \eqref{eq_mu_k_def}:
\begin{equation}\label{eq_sym_int_dec}
    \widehat{\mu^{\br{j}}_{k,\lambda}}\br{\vxi}
    =\int
        2^k K\br{2^k t}\,
        \overline{e\br{2^k\xi_0 t+4^k\xi_1 t^2}}\,
        \chi_j\br{\frac{2^k \xi_0}{\lambda},\frac{4^k\xi_1}{\lambda}}
    \rho\br{\abs{t}}
    dt.
\end{equation}
\end{itemize}

Combining the above two equivalent decompositions---as a direct manifestation of the time-frequency dual analysis---and invoking a parabolic scaling, we obtain the operator decomposition
\begin{equation}\label{eq_cC_2ways_2_decomp}
    \cC\cR^{\br{\va}}_{\lambda,j}= \sum_{\vu\in \Z^2}
    \cC\cR^{\br{\va}}_{\lambda,j,\vu}= 
    \int
    \cC\cR^{\br{\va}}_{\lambda,j,t}\,
    \rho\br{\abs{t}}
    dt,
\end{equation}
where the elementary building blocks \(\cC\cR^{\br{\va}}_{\lambda,j,\vu}\) and \(\cC\cR^{\br{\va}}_{\lambda,j,t}\) are the natural correspondent to the above dual time-frequency representation:
\begin{itemize}
\item \textsf{the frequency component} \(\cC\cR^{\br{\va}}_{\lambda,j,\vu}\) given by
\begin{align}\label{eq_cC_4ier_dec}
    \cC\cR^{\br{\va}}_{\lambda,j,\vu}f\br{\vx}:= &
    \sum_{k=\uk}^{\ok\br{\vx}} 
    m^{\br{j}}_{k,\lambda,\vu}
    \cdot e\br{\frac{4^k a\br{\vx} u_1}{\lambda}} \nonumber\\
    \cdot
    \int 
        \chi_j & \br{\frac{2^k \xi_0}{\lambda},\frac{4^k\br{\xi_1-a\br{\vx}}}{\lambda}}
        \widehat{f}\br{\vxi}
    e\br{\xi_0\br{x_0-\frac{2^k u_0}{\lambda}}+\xi_1\br{x_1-\frac{4^k u_1}{\lambda}}}d\vxi,
\end{align}
\item \textsf{the spatial  component}  \(\cC\cR^{\br{\va}}_{\lambda,j,t}\) given by
\begin{align}\label{eq_cC_int_dec}
    \cC\cR^{\br{\va}}_{\lambda,j,t}f\br{\vx}:= &
    \sum_{k=\uk}^{\ok\br{\vx}} 
    2^k K\br{2^k t}
    \cdot e\br{4^k a\br{\vx} t^2} \nonumber\\
    \cdot
    \int 
        \chi_j & \br{\frac{2^k \xi_0}{\lambda},\frac{4^k\br{\xi_1-a\br{\vx}}}{\lambda}}\,
        \widehat{f}\br{\vxi}\,
    e\br{\xi_0\br{x_0-2^k t}+\xi_1\br{x_1-4^k t^2}}\,d\vxi.
\end{align}
\end{itemize}

Combining the above two forms, we are naturally led to the following:

\begin{definition}\label{def_op_w_phy_shift} [\textsf{Shifted superposed Carleson operator}]
    Given a sequence of coefficients \(\br{\varepsilon_k}_{k\in\Z} \in \C^\Z\), a vector \(\vu\in\R^2\), and a tuple \(\vR:=\br{R_0,R_1}\) with positive entries, we define its associated \textsf{shifted superposed Carleson operator} as:
    \begin{align}\label{eq_op_w_phy_shift}
        \cC^{\br{\va,\varepsilon,\uk,\ok}}_{\vR,j,\vu}f\br{\vx}:= &
        \sum_{k=\uk}^{\ok\br{\vx}} 
        \varepsilon_k
        \cdot e\br{\frac{4^k a\br{\vx} u_1}{R_1}} \nonumber\\
        \cdot
        \int 
            \chi_j & \br{\frac{2^k \xi_0}{R_0},\frac{4^k\br{\xi_1-a\br{\vx}}}{R_1}}
            \widehat{f}\br{\vxi}\,
        e\br{\xi_0\br{x_0-\frac{2^k u_0}{R_0}}+\xi_1\br{x_1-\frac{4^k u_1}{R_1}}}d\vxi.
    \end{align}
\end{definition}

With these done, we claim that the following holds:

\begin{theorem}[\textsf{Log-type dependence on the physical-shift}]\label{thm_log_shift_est}
For any $p\in\br{1,\infty}$ one has
    \begin{equation}
        \nrm{\cC^{\br{\va,\varepsilon,\uk,\ok}}_{\vR,j,\vu}f
        }_{L^p}
        \underset{p,\chi_j}{\lesssim}
        \log^2\br{e+\abs{\vu}}
        \nrm{\varepsilon}_{\ell^\infty}
        \nrm{f}_{L^p}\,.
    \end{equation}
\end{theorem}

However, via a rescaling argument, it turns out that in order to prove \textbf{Theorem \ref{thm_log_shift_est}}, it is enough to analyze the \(\vR=\br{1,1}\) case. Indeed, recalling \eqref{eq_op_w_phy_shift}, we observe that
\begin{align*}
    \Dil^p_\vR\cC^{\br{\va,\varepsilon,\uk,\ok}}_{\vR,j,\vu}\br{\Dil^p_\vR}^{-1}f\br{\vx}= &
    \sum_{k=\uk}^{\Dil^\infty_\vR \ok\br{\vx}} 
    \varepsilon_k
    \cdot e\br{4^k \cdot \frac{\Dil^\infty_\vR a\br{\vx}}{R_1} \cdot u_1}\\
    \cdot
    \int 
        \chi_j & \br{2^k \xi_0,4^k\xi_1}
        \widehat{f}\br{\vxi}
    e\br{\xi_0\br{x_0-2^k u_0}+\xi_1\br{x_1-4^k u_1}}d\vxi.
\end{align*}
Setting now \(\va_\vR:=\Dil^1_\vR\va\) and \(\ok_\vR:=\Dil^\infty_\vR\ok\), we obtain the simple identity
\begin{equation}\label{resc}
    \Dil^p_\vR\cC^{\br{\va,\varepsilon,\uk,\ok}}_{\vR,j,\vu}\br{\Dil^p_\vR}^{-1}
    =
    \cC^{\br{\va_\vR,\varepsilon,\uk,\ok_\vR}}_{\br{1,1},j,\vu}.
\end{equation}
Thus, as claimed, \textbf{Theorem \ref{thm_log_shift_est}} is a direct consequence of the following

\begin{theorem}[\textsf{Log-type dependence on physical-shift, the normalized version}]\label{thm_phy_shift_est_unit}
    Set \(\cC^{\br{\va,\varepsilon}}_{j,\vu}:=\cC^{\br{\va,\varepsilon,\uk,\ok}}_{\br{1,1},j,\vu}\). Then, for any $p\in\br{1,\infty}$, the following estimate holds:
    \begin{equation}\label{shiftct}
        \nrm{
            \cC^{\br{\va,\varepsilon}}_{j,\vu}f
        }_{L^p}
        \underset{p,\chi_j}{\lesssim}
        \log^2\br{e+\abs{\vu}}\nrm{\varepsilon}_{\ell^\infty}\nrm{f}_{L^p}\,.
    \end{equation}
\end{theorem}

\begin{remark}
    One may refine the result on anisotropic variants of the Carleson operator presented in the Ph.D. dissertation by J. Roos. The \(\log\)-type bound implies the sharp Mikhlin regularity condition. In our case, having the \(\log\) type dependence is crucial for proving \eqref{eq_HF_log} as it will be seen in a moment.
\end{remark}

With these being said, we claim that \textbf{Theorem \ref{prop_4cases}} is a direct consequence of \textbf{Theorem \ref{thm_phy_shift_est_unit}}. To see this, assuming in what follows that \eqref{shiftct} holds, we first prove the following auxiliary result: 

\begin{lemma}[\textsf{Mikhlin-type symbol condition}]\label{lem_mikhlin_symb}
Let \(\chi_j\) in \eqref{eq_symb_lamb_dec} be chosen so that \eqref{eq_cone_abst} matches \eqref{eq_cone_dec}. Let also \(k\in\Z\) and \(j\in\{0,1\}\). For multi-indices\footnote{Here $A$ is the parameter appearing in relation \eqref{eq_kernel_cond}.} \(\valpha\in \BR{0,\dots, \lfloor\frac{A-1}{2}\rfloor }^2\), the following holds.
\begin{itemize}
    \item Low-oscillatory case: given \(\abs{\lambda}\lesssim 1\) one has
    \begin{equation}\label{eq_LF_par}
        \abs{\partial^\valpha_\vxi
        \widehat{\mu^{\br{j}}_{k,\lambda}}\br{\frac{\xi_0}{2^k},\frac{\xi_1}{4^k}}}
        \lesssim
        \abs{\lambda}
        \cdot
        \abs{\vxi}^{-\abs{\valpha}}.
    \end{equation}
    \item High oscillatory, non-stationary case: given \(\abs{\lambda}\gtrsim 1\) one has
    \begin{equation}\label{eq_HF_nst_par}
        \abs{\partial^\valpha_\vxi
        \widehat{\mu^{\br{1}}_{k,\lambda}}\br{\frac{\xi_0}{2^k},\frac{\xi_1}{4^k}}}
        \lesssim
        \abs{\lambda}^{-1}
        \cdot
        \abs{\vxi}^{-\abs{\valpha}}.
    \end{equation}
\end{itemize}
\end{lemma}

\begin{proof}
Recall \eqref{eq_symb_lamb_dec} and \eqref{eq_sym_int_dec}. The normalization gives:
\begin{equation*}
    \widehat{\mu^{\br{j}}_{k,\lambda}}
    \br{\frac{\xi_0}{2^k},\frac{\xi_1}{4^k}}
    =
    \widehat{\mu_k}\br{\frac{\xi_0}{2^k},\frac{\xi_1}{4^k}}
    \cdot
    \chi_j\br{\frac{\xi_0}{\lambda},\frac{\xi_1}{\lambda}}
    =
    \left(\int
        \overline{e\br{\xi_0 t+\xi_1 t^2}}\,
    2^k\, K\br{2^k t}
    \rho\br{\abs{t}}
    dt\right)
    \cdot
    \chi_j\br{\frac{\xi_0}{\lambda},\frac{\xi_1}{\lambda}}
    .
\end{equation*}
By Leibniz's rule, it suffices to verify the Mikhlin-type symbol condition for the oscillatory integral expression \(\widehat{\mu_k}\br{\frac{\xi_0}{2^k},\frac{\xi_1}{4^k}}\) when \(\vxi\in \supp \Dil^\infty_\lambda \chi_j\). Starting with the treatment of \eqref{eq_LF_par}, we observe that the truncated kernal \(2^k K\br{2^k t}
    \rho\br{t}\) satisfies the mean-zero condition:
\begin{equation*}
    \int 2^k K\br{2^k t} \rho\br{\abs{t}} dt
    = 2
    \int_{2^{-1}}^2
    \rho'\br{t}
        \int_{2^{k-1}<\abs{s}\leq 2^k t}
                K\br{s}
        ds
    dt
    =0.
\end{equation*}
As a direct consequence, we have:
\begin{equation}\label{eq_symb_cancel}
    \abs{\widehat{\mu_k}\br{\frac{\xi_0}{2^k},\frac{\xi_1}{4^k}}}
    =
    \abs{
    \int
        \br{
        \overline{e\br{\xi_0 t+\xi_1 t^2}}
        -1}
    2^k K\br{2^k t}
    \rho\br{\abs{t}}
    dt
    }
    \lesssim
    \abs{\xi_0 t+\xi_1 t^2}
    \lesssim
    \abs{\lambda}
\end{equation}
whenever \(\vxi\in \supp \Dil^\infty_\lambda \chi_j\). On the other hand, a trivial estimate gives:
\begin{equation}\label{eq_symb_triv}
    \abs{\partial^\valpha_\vxi \widehat{\mu_k}\br{\frac{\xi_0}{2^k},\frac{\xi_1}{4^k}}}
    \underset{\valpha}{\eqsim}
    \abs{
    \int
        \overline{e\br{\xi_0 t+\xi_1 t^2}}
    2^k K\br{2^k t}
    t^{\alpha_0+2\alpha_1}\rho\br{\abs{t}}
    dt
    }
    \underset{\valpha}{\lesssim} 1.
\end{equation}
If \(\abs{\lambda}\lesssim 1\), \(\vxi\in \supp \Dil^\infty_\lambda \chi_j\), and \(\va\neq \vnull\), we may further dominate the above with \(\underset{\va}{\lesssim} \abs{\lambda}^{1-\abs{\valpha}}\underset{\valpha}{\eqsim}
    \abs{\lambda}\cdot\abs{\vxi}^{-\abs{\valpha}}\).
In combination with \eqref{eq_symb_cancel}, we verify \eqref{eq_LF_par}. To demonstrate \eqref{eq_HF_nst_par}, we first let:
\begin{equation}\label{eq_symb_xi_range}
    \vxi=\br{\xi_0,\xi_1}\in\supp \Dil^\infty_\lambda\chi_1
    = \supp\Dil^\infty_\lambda \tarphi
    \times \supp\Dil^\infty_\lambda \tpsi
    \subset \br{-0.26\lambda,0.26\lambda}\times \br{0.49\lambda,0.51\lambda}.
\end{equation}
We now integrate by parts the middle expression of \eqref{eq_symb_triv} to obtain:
\begin{equation}\label{eq_symb_by_part}
    \abs{\partial^\valpha_\vxi \widehat{\mu_k}\br{\frac{\xi_0}{2^k},\frac{\xi_1}{4^k}}}
    \underset{\valpha}{\eqsim}
    \abs{
    \int
        \overline{e\br{\xi_0 t+\xi_1 t^2}}
        \br{
            \frac{\partial}{\partial t}
            \cdot
            \frac{1}{\xi_0+2\xi_1 t}\cdot
        }^{1+\abs{\valpha}}
    2^k K\br{2^k t}
    t^{\alpha_0+2\alpha_1}\rho\br{\abs{t}}
    dt
    }.
\end{equation}
By design, the conditions \eqref{eq_parti_uni} and \eqref{eq_symb_xi_range} guarantee that \(\abs{\xi_0+2\xi_1 t}\geq 2\abs{\xi_1}\cdot\abs{t}-\abs{\xi_0}\geq \br{0.49-0.26}\abs{\lambda}\). As a result, we may dominate the above with \(\underset{\valpha}{\lesssim}\abs{\lambda}^{-1-\abs{\valpha}}\underset{\valpha}{\eqsim} \abs{\lambda}^{-1}\cdot\abs{\vxi}^{-\abs{\valpha}}\).
This concludes \eqref{eq_HF_nst_par}.
\end{proof}

\begin{remark}
The treatment of \eqref{eq_HF_st} cannot follow the same strategy since the estimate \(\abs{m^{\br{0}}_{k,\lambda,\vu}}\lesssim \abs{\lambda}^{-c}\br{1+\abs{\vu}^2}^{-1-c}\)
fails for all \(c>0\) when \(\abs{\lambda}\gg 1\). This is the direct consequence of having stationary phase contribution in the oscillatory integral expression \(\widehat{\mu_k}\br{\vxi}\). One must utilize the oscillatory nature of the symbol to control the expression.
\end{remark}

With these done, we are now ready to show 

\begin{proof}[\textsf{Proof of \textbf{Theorem \ref{prop_4cases}} as a consequence of \textbf{Theorem \ref{thm_phy_shift_est_unit}} (and hence \textbf{Theorem \ref{thm_log_shift_est}}).}] 
$\newline$

\noindent 1) \underline{\textsf{The proof of \eqref{eq_HF_log}}}. 
\medskip

We start by observing that
\begin{equation*}
    \cC\cR^{\br{\va}}_{\lambda,j,t}=\cC^{\br{\va,\varepsilon\br{t},\uk,\ok}}_{\br{\lambda,\lambda},j,\lambda\br{t,t^2}},\quad
    \varepsilon_k\br{t}:=2^kK\br{2^k t}.
\end{equation*}
By \textbf{Theorem \ref{thm_log_shift_est}}, we deduce for \(\abs{t}\eqsim 1\) the estimate:
\begin{equation*}
    \nrm{\cC\cR^{\br{\va}}_{\lambda,j,t} f}_{L^p}\underset{p}{\lesssim} 
    \log^2\br{e+\lambda\abs{\br{t,t^2}}}
    \nrm{
        2^k K\br{2^k t}
    }_{\ell^\infty\br{k}}
    \nrm{f}_{L^p}
    \eqsim
    \log^2\br{e+\lambda}
    \nrm{f}_{L^p}.
\end{equation*}
Combined with \eqref{eq_cC_2ways_2_decomp}, we conclude:
\begin{equation*}
    \nrm{\cC\cR^{\br{\va}}_{\lambda,j} f}_{L^p}\leq
    \int
        \nrm{\cC\cR^{\br{\va}}_{\lambda,j,t} f}_{L^p}
    \rho\br{\abs{t}}
    dt
    \underset{p}{\lesssim}
    \log^2\br{e+\lambda}
    \nrm{f}_{L^p}.
\end{equation*}

\noindent 2) \underline{\textsf{The proofs of \eqref{eq_LF} and \eqref{eq_HF_nst}}}. 
\medskip

In this situation, we have
\begin{equation}\label{eq_LF_HF_nst_2_shift_mod}
    \cC\cR^{\br{\va}}_{\lambda,j,\vu}=\cC^{\br{\va,\varepsilon^{\br{j}},\uk,\ok}}_{\br{\lambda,\lambda},j,\vu},\quad \varepsilon^{\br{j}}_k:=m^{\br{j}}_{k,\lambda,\vu},
\end{equation}
and hence from \textbf{Theorem \ref{thm_log_shift_est}} and \eqref{eq_cC_2ways_2_decomp} we deduce 
\begin{equation}\label{eq_eq_LF_HF_nst_2_shift_est}
    \nrm{\cC\cR^{\br{\va}}_{\lambda,j}}_{L^p}\underset{p}{\lesssim}
    \sum_{\vu\in\Z^2}
    \log^2\br{e+\abs{\vu}}
    \nrm{m^{\br{j}}_{k,\lambda,\vu}}_{\ell^\infty\br{k}}
    \nrm{f}_{L^p}.
\end{equation}

Thus, it remains to obtain a suitable control on the coefficients \(m^{\br{j}}_{k,\lambda,\vu}\), which is a consequence of a standard argument and \textbf{Lemma \ref{lem_mikhlin_symb}}:
\begin{itemize}
        \item low-oscillatory case: given \(\abs{\lambda}\lesssim 1\) one has
        \begin{equation*}
            \abs{m^{\br{j}}_{k,\lambda,\vu}}
            \lesssim
            \abs{\lambda}
            \br{1+\abs{\vu}^2}^{-2};
        \end{equation*}
        \item high-oscillatory non-stationary case: given \(\abs{\lambda}\gg 1\) one has
        \begin{equation*}
            \abs{m^{\br{1}}_{k,\lambda,\vu}}
            \lesssim
            \abs{\lambda}^{-1}
            \br{1+\abs{\vu}^2}^{-2}.
        \end{equation*}
\end{itemize}
Consequently, putting everything together, we obtain
\begin{itemize}
\item for \(\abs{\lambda}\lesssim 1\) the following holds:
    \begin{equation*}
        \sum_{\vu\in\Z^2}
        \log^2\br{e+\abs{\vu}}
        \nrm{m^{\br{j}}_{k,\lambda,\vu}}_{\ell^\infty\br{k}}
        \lesssim
        \abs{\lambda}
        \sum_{\vu\in\Z^2}
        \log^2\br{e+\abs{\vu}}
        \br{1+\abs{\vu}^2}^{-2}
        \eqsim\abs{\lambda}
    \end{equation*}
\item for \(\abs{\lambda}\gtrsim 1\) the following holds:
    \begin{equation*}
        \sum_{\vu\in\Z^2}
        \log^2\br{e+\abs{\vu}}
        \nrm{m^{\br{1}}_{k,\lambda,\vu}}_{\ell^\infty\br{k}}
        \lesssim
        \abs{\lambda}^{-1}
        \sum_{\vu\in\Z^2}
        \log^2\br{e+\abs{\vu}}
        \br{1+\abs{\vu}^2}^{-2}
        \eqsim\abs{\lambda}^{-1}.
    \end{equation*}
\end{itemize}    
This concludes our proof.
\end{proof}

Thus with these we have reduced our main \textbf{Theorem \ref{thm_main1D}} to \textbf{Theorems \ref{prop_4casesdecay}} and \textbf{\ref{thm_phy_shift_est_unit}}.

\subsection{Preliminary frequency discretization: reduction of Theorems \ref{prop_4casesdecay} and \ref{thm_phy_shift_est_unit} to Theorems \ref{thm_HF_st_disc} and \ref{thm_phy_shift_disc}, respectively}\label{subsec_Freq_disc}

In this subsection, employing the time-frequency analytic framework developed by Lacey and Thiele, we focus on obtaining suitable discretized models for \(\cC\cR^{\br{\va}}_{\lambda,0}\) and \(\cC\cR^{\br{\va}}_{j,\vu}\).

\subsubsection{Model sums} Our first aim here is to decouple the frequency variable \(\xi_1\) from the measurable function expression \(a\br{\vx}\) in both \eqref{eq_cC_j_lambda_def} and \eqref{eq_op_w_phy_shift}.

For this, we start by recalling \eqref{conedec}, and, assuming without loss of generality \(\lambda>0\), we let
\begin{equation}\label{eq_chi_01_nice}
    \chi_0\br{\vxi}:=
    \phi\br{2\xi_0-\frac{1}{2}}
    \varphi\br{\xi_1},\quad
    \chi_1\br{\vxi}:=
    \phi\br{\xi_0}\varphi\br{2\xi_1-1}.
\end{equation}
Our declared aim above invites us to consider the following representation:
\begin{align*}
    \varphi\br{\xi-a}= &
    \int
        \phi\br{\xi-\zeta}
        \1_{\bR{-\frac{1}{2},\frac{1}{2}}}\br{a-\zeta}
    \frac{d\zeta}{\nrm{\phi}_{L^1}}=\int
        \phi\br{\xi-\zeta-\frac{1}{2}}
        \1_{\bR{-\frac{1}{2},\frac{1}{2}}}\br{a-\zeta-\frac{1}{2}}
    \frac{d\zeta}{\nrm{\phi}_{L^1}}\\
    = &
    \fint_0^M
        \sum_{z\in\Z}
        \phi\br{\xi-\zeta-\frac{1}{2}-z}
        \1_{\bR{-\frac{1}{2},\frac{1}{2}}}\br{a-\zeta-\frac{1}{2}-z}
    \frac{d\zeta}{\nrm{\phi}_{L^1}},\quad M\in\N.
\end{align*}
Using now Definition \ref{def_Psi_I}, we may rewrite the above identity as
\begin{equation}\label{eq_varphi_conv_expan}
    \varphi\br{\xi-a}
    = 
    \fint_0^M
        \sum_{z\in\Z}
        \phi_{\bR{z,z+1}}\br{\xi-\zeta}
        \1_{\bR{z,z+1}}\br{a-\zeta}
    \frac{d\zeta}{\nrm{\phi}_{L^1}},\quad M\in\N.
\end{equation}
On the other hand, an almost identical argument produces the following identity:
\begin{align*}
    \varphi\br{2\br{\xi-a}-1}
    = &
    \int \phi\br{2\xi-\zeta-\frac{3}{2}}
    \1_{\bR{-\frac{1}{2},\frac{1}{2}}}\br{2a-\zeta-\frac{1}{2}}\frac{d\zeta}{\nrm{\phi}_{L^1}}\\
    = &
    2
    \int \phi\br{2\br{\xi-\zeta-\frac{3}{4}}}
    \1_{\bR{-\frac{1}{2},\frac{1}{2}}}\br{2\br{a-\zeta-\frac{1}{4}}}\frac{d\zeta}{\nrm{\phi}_{L^1}}\\
    = &
    2 \fint_0^M \sum_{z\in\Z} \phi_{\bR{z+\frac{1}{2},z+1}}\br{\xi-\zeta}\1_{\bR{z,z+\frac{1}{2}}}\br{a-\zeta}\frac{d\zeta}{\nrm{\phi}_{L^1}},\quad
    M\in\N.
\end{align*}

Incorporating now the notation in Definition \ref{def_half_I}, we rewrite the above identity as
\begin{equation}\label{eq_varphi_w_shift_conv_expan}
    \varphi\br{2\br{\xi-a}-1}
    =
    2 \fint_0^M \sum_{z\in\Z} \phi_{\bR{z,z+1}^+}\br{\xi-\zeta}\1_{\bR{z,z+1}^-}\br{a-\zeta}\frac{d\zeta}{\nrm{\phi}_{L^1}},\quad
    M\in\N.
\end{equation}
Using relation \(\phi\br{2\xi-\frac{1}{2}}=\phi_{\bR{-\frac{1}{2},\frac{1}{2}}^+}\br{\xi}\) and \eqref{eq_varphi_conv_expan}, we take \(M=4^{k-\uk}\in\N\) and deduce
\begin{align*}
    &
    \chi_0\br{\frac{2^{k+1}\xi_0}{\lambda},\frac{4^k\br{\xi_1-a\br{\vx}}}{\lambda}}\\
    = &
    \fint_0^{4^{k-\uk}} 
    \sum_{z\in\Z}
    \phi_{\bR{-\frac{1}{2},\frac{1}{2}}^+}\br{\frac{2^k\xi_0}{\lambda}}
    \phi_{\bR{z,z+1}}\br{\frac{4^k\xi_1}{\lambda}-\zeta}
    \1_{\bR{z,z+1}}\br{\frac{4^ka\br{\vx}}{\lambda}-\zeta}
    \frac{d\zeta}{\nrm{\phi}_{L^1}}\\
    = &
    \fint_0^{\lambda/4^\uk} 
    \sum_{z\in\Z}
    \phi_{\bR{-\frac{\lambda}{2^{k+1}},\frac{\lambda}{2^{k+1}}}^+}\br{\xi_0}
    \phi_{\bR{\frac{\lambda z}{4^k},\frac{\lambda\br{z+1}}{4^k}}}\br{\xi_1-\zeta}
    \1_{\bR{\frac{\lambda z}{4^k},\frac{\lambda\br{z+1}}{4^k}}}\br{a\br{\vx}-\zeta}
    \frac{d\zeta}{\nrm{\phi}_{L^1}}.
\end{align*}
The above identity suggests the introduction of the following multi-scale collection of frequency rectangles:
\begin{definition}[\textsf{Frequency $2-$intervals (rectangles)}]\label{def_W_n_Wk}
    \begin{equation}\label{eq_def_W_n_Wk}
    \W\br{\lambda}:=\bigsqcup_{k\in\Z}\W_k\br{\lambda},\quad
    \W_k\br{\lambda}:=
    \BR{
        \bR{
            -\lambda/2^{k+1},
            \lambda/2^{k+1}
        }
        \times
        \bR{
            \lambda w/4^k,
            \lambda\br{w+1}/4^k
        }
        \subset \R^2
        \::\:
        w\in\Z
    }.
\end{equation}
When the choice of \(\lambda\) is clear from the context, for simplicity, we write \(\W:=\W\br{\lambda}\) and \(\W_k:=\W_k\br{\lambda}\).
\end{definition}
We thus have the succinct expression:
\begin{equation}\label{eq_chi_0_exp_succinct}
    \chi_0\br{\frac{2^{k+1}\xi_0}{\lambda},\frac{4^k\br{\xi_1-a\br{\vx}}}{\lambda}}
    =\fint^{\lambda/4^\uk}_0
    \sum_{\vomega\in \W_k\br{\lambda}}
    \phi_{\vomega^{\br{+,0}}}\br{\vxi-\br{0,\zeta}}
    \1_{\vomega^{\br{-,0}}}\br{\va\br{\vx}-\br{0,\zeta}}
    \frac{d\zeta}{\nrm{\phi}_{L^1}}.
\end{equation}
In view of \eqref{eq_varphi_w_shift_conv_expan}, we also obtain:
\begin{equation}\label{eq_chi_1_exp_succinct}
    \chi_1\br{\frac{2^{k+1}\xi_0}{\lambda},\frac{4^k\br{\xi_1-a\br{\vx}}}{\lambda}}
    =2\fint^{\lambda/4^\uk}_0
    \sum_{\vomega\in \W_k\br{\lambda}}
    \phi_{\vomega^{\br{+,1}}}\br{\vxi-\br{0,\zeta}}
    \1_{\vomega^{\br{-,1}}}\br{\va\br{\vx}-\br{0,\zeta}}
    \frac{d\zeta}{\nrm{\phi}_{L^1}}.
\end{equation}
The above two identities together with the expressions \eqref{eq_cC_j_lambda_def} and \eqref{eq_op_w_phy_shift} suggest introducing the following 

\begin{definition}[\textsf{Model sum associated to} \(\cC\cR^{\br{\va}}_{\lambda,0}\)]\label{def_cC_Delta_lambda_0}
    \begin{equation*}
        \cC^{\br{\va}}_{\Delta,\lambda,0}f\br{\vx}:=
        \sum_{k=\uk}^{\ok\br{\vx}}
        \sum_{\vomega\in\W_k\br{\lambda}} 
        \1_{\vomega^{\br{-,0}}}\br{\va\br{\vx}}
        \int
            \widehat{\mu_k}\br{\vxi-\va\br{\vx}}
        \phi_{\vomega^{\br{+,0}}}\br{\vxi}
        \widehat{f}\br{\vxi}
        e\br{\vx\cdot\vxi}
        d\vxi.
    \end{equation*}
\end{definition}
\begin{definition}[\textsf{Model sum associated to} \(\cC^{\br{\va,\varepsilon}}_{j,\vu}\)]\label{def_cC_Delta_j_vu}
    \begin{align*}
        \cC^{\br{\va,\varepsilon}}_{\Delta,j,\vu}f\br{\vx}:=
        \sum_{k=\uk}^{\ok\br{\vx}}
        \sum_{\vomega\in\W_k\br{1}} 
        &
        \varepsilon_k
        \cdot
        \1_{\vomega^{\br{-,j}}}\br{\va\br{\vx}}
        e\br{4^ka\br{\vx}u_1}\\
        &\cdot \int
        \phi_{\vomega^{\br{+,j}}}\br{\vxi}
        \widehat{f}\br{\vxi}
        e\br{\xi_0\br{x_0-2^k u_0}+\xi_1\br{x_1-4^k u_1}}
        d\vxi.
    \end{align*}
\end{definition}

\subsubsection{Reformulation of our goal and reductions}

For the next three sections, our main goal will be to prove the following two key results:

\begin{theorem}[\textsf{Model sum estimate for the high-frequency stationary-phase contribution}]\label{thm_HF_st_disc}
There is a universal constant \(\epsilon>0\) such that for any \(\lambda\gtrsim 1\) and \(\uk\in\Z\), function \(f\in L^2\br{\R^2}\), and measurable functions \(\va:\R^2\to\BR{0}\times\R\), \(\ok:\R^2\to\Z\), the following estimate holds:
    \begin{equation*}
        \nrm{
            \cC^{\br{\va}}_{\Delta,\lambda,0}
            f
        }_{L^{2,\infty}}
        \lesssim \lambda^{-\epsilon}
        \nrm{f}_{L^2}.
    \end{equation*}
\end{theorem}
\begin{theorem}[\textsf{Model sum estimate with Log-type dependence on physical-shift}]\label{thm_phy_shift_disc}
Let \(p\in\br{1,\infty}\) and \(j\in\BR{0,1}\). For any sequence \(\br{\varepsilon_k}_{k\in\Z}\in\C^\Z\), integer \(\uk\in\Z\), function \(f\in L^2\br{\R^2}\), and measurable functions \(\va:\R^2\to\BR{0}\times\R\), \(\ok:\R^2\to\Z\), the following estimate holds:
    \begin{equation*}
        \nrm{
        \cC^{\br{\va,\varepsilon}}
        _{\Delta,j,\vu}f
        }_{L^p}
        \underset{p}{\lesssim}
        \log^2\br{e+\abs{\vu}}\nrm{\varepsilon}_{L^\infty}
        \nrm{f}_{L^p}.
    \end{equation*}
\end{theorem}
\begin{proof}[\textsf{Proof of \textbf{Theorem \ref{prop_4casesdecay}} and \textbf{Theorem \ref{thm_phy_shift_est_unit}} assuming that both \textbf{Theorems
\ref{thm_HF_st_disc}} and \textbf{\ref{thm_phy_shift_disc}} hold.}]
$\newline$

\noindent By \eqref{eq_chi_0_exp_succinct} and \eqref{eq_chi_1_exp_succinct}, we obtain the following two averaging identities:
\begin{equation*}
    \cC\cR^{\br{\va}}_{\lambda,0}f\br{\vx}=
    c_0\fint^{\lambda/4^\uk}_0
    \Mod_{\br{0,\zeta}}\cC^{\br{\va-\br{0,\zeta}}}_{\Delta,\lambda,0}\Mod_{-\br{0,\zeta}}f\br{\vx}
    d\zeta;
\end{equation*}
\begin{equation*}
    \cC^{\br{\va,\varepsilon}}_{j,\vu}f\br{\vx}
    =c_j\fint^{4^{-\uk}}_0
    \Mod_{\br{0,\zeta}}\cC^{\br{\va-\br{0,\zeta},\varepsilon}}_{\Delta,j,\vu}\Mod_{-\br{0,\zeta}}
    f\br{\vx}
    d\zeta,
\end{equation*}
with \(c_j:=2^j/\nrm{\phi}_{L^1}\) and \(\Mod_\vxi f\br{\vx}:=e\br{\vxi \cdot\vx}f\br{\vx}\). We now apply \textbf{Theorems \ref{thm_HF_st_disc}} and \textbf{\ref{thm_phy_shift_disc}} and triangle inequality to finish the proof. 
\end{proof}

\subsection{A final reformulation of our main object of study}\label{subsec_Freq_disc1}

In this final preliminary subsection, we will introduce some further notions in order to help recast Definitions \ref{def_cC_Delta_lambda_0} and \ref{def_cC_Delta_j_vu} in a unified language.

We start with the following convenient notation:
\begin{definition}[\textsf{Shifted intervals}]\label{def_shift_I}
    Let \(I=\bR{\inf I,\sup I}\subset \R\). Given \(v\in\R\), we define:
    \begin{equation*}
        v\boxplus I:=v\abs{I}+I=\bR{v\abs{I}+\inf I,v\abs{I}+\sup I}.
    \end{equation*}
    For interval \(\vI\subset \R^{D+1}\) and vector \(\vv\in\R^{D+1}\), we define the higher-dimensional analog accordingly
    \begin{equation*}
        \vv\boxplus \vI:=\prod_{j=0}^D v_j\boxplus I_j,\quad
        \vI=\prod_{j=0}^D I_j,\quad
        \vv=\br{v_j}_{j=0}^D.
    \end{equation*}
\end{definition}
This provides an alternate way to express \textbf{Definition \ref{def_W_n_Wk}}, that is
\begin{equation*}
    \W_k\br{\lambda}=\BR{
        \br{0,w}\boxplus 
        \lambda
        \begin{pmatrix}
            2 & 0\\
            0 & 4
        \end{pmatrix}^{-k}
        \hspace{-2ex}
        \cdot
        \bR{-1/2,1/2}\times\bR{0,1}\subset \R^2
    \::\:
        w\in\Z
    }
\end{equation*}
Next, the presence of \(\phi_{\vomega^{\br{+,j}}}\widehat{f}\) and \(\sum_{k=\uk}^\ok\1_{\vomega^{\br{-,j}}}\br{\va}\dots\) in both \textbf{Definition \ref{def_cC_Delta_lambda_0} and \ref{def_cC_Delta_j_vu}} motivates the following two definitions:
\begin{definition}[\textsf{Frequency projection}]\label{def_Freq_proj}
    Given an interval \(\vomega\subset \R^{D+1}\) we let \(\widehat{\pi_\vomega f}:=\phi_\vomega \widehat{f}\).
 \end{definition}
\begin{definition}[\textsf{Level set decomposition}]\label{def_A_omega}
Given measurable functions \(\va: \R^{D+1}\to \V\) and \(\ok: \R^{D+1} \to \Z\cap \Br{\uk,\infty}\), we define for each \(\vomega\in\W_k\br{\lambda}\) the following:
    \begin{equation}\label{eq_def_A_omega}
        A_\vomega:=
            \va^{-1}\vomega\cap \ok^{-1}\bR{k,\infty}
        ,\quad
        A^{\br{-,j}}_\vomega:=
            \va^{-1}\vomega^{\br{-,j}}\cap \ok^{-1}\bR{k,\infty}
        .
    \end{equation}
    We also write \(A^-_\vomega:=A^{\br{-,j}}_\vomega\) when the choice of \(j\) is clear from the context.
\end{definition}
We can thus rewrite the main expression in \textbf{Definition \ref{def_cC_Delta_lambda_0}} as
\begin{align}\label{eq_def_cC_Delta_lambda_0}
    \cC^{\br{\va}}_{\Delta,\lambda,0}f\br{\vx}= &
    \sum_{k\in\Z}
    \sum_{\vomega\in\W_k\br{\lambda}}
        \1_{A^{\br{-,0}}_\vomega}\br{\vx}
        \int
            \widehat{\mu_k}\br{\vxi-\va\br{\vx}}
            \widehat{\pi^{\br{+,0}}_\vomega f}
            \br{\vxi}
            e\br{\vx\cdot\vxi}
        d\vxi\nonumber\\
    = & \sum_{k\in\Z}
    \sum_{\vomega\in\W_k\br{\lambda}}
        \1_{A^{\br{-,0}}_\vomega}\br{\vx}
        \int
            \pi^{\br{+,0}}_\vomega f\br{\vx-\vtau}
            e\br{\va\br{\vx}\cdot\vtau}
        d\mu_k\br{\vtau}
\end{align}
and, by setting \( \alpha_k\br{\vtau}:=\varepsilon_k \delta\br{\vtau-\br{2^k u_0,4^k u_1}}\), also the main object in \textbf{Definition \ref{def_cC_Delta_j_vu}}:
\begin{align}\label{eq_def_cC_Delta_j_vu}
    \cC^{\br{\va,\varepsilon}}_{\Delta,j,\vu}f\br{\vx}= &
    \sum_{k\in\Z}
    \sum_{\vomega\in\W_k\br{1}}
        \varepsilon_k
        \1_{A^{\br{-,j}}_\vomega}\br{\vx}
        e\br{4^k a\br{\vx} u_1}
        \pi^{\br{+,j}}_\vomega f
        \br{x_0-2^ku_0 , x_1-4^k u_1}
        \nonumber\\
    = & \sum_{k\in\Z}
    \sum_{\vomega\in\W_k\br{1}}
        \1_{A^{\br{-,j}}_\vomega}\br{\vx}
        \int
            \pi^{\br{+,j}}_\vomega f\br{\vx-\vtau}
            e\br{\va\br{\vx}\cdot\vtau}
        d\alpha_k\br{\vtau}.
\end{align}
The similarity between \eqref{eq_def_cC_Delta_lambda_0} and \eqref{eq_def_cC_Delta_j_vu} naturally suggests a general framework:
\begin{definition}\label{def_cC_nu}[\textsf{A unified object of study}]
Given a sequence of complex Borel measures \(\br{\nu_k}_{k\in\Z}\) we define:
    \begin{equation}\label{eq_def_cC_nu}
        \cC_{\nu,\lambda,j}^{\br{\va}}
        f:=
        \sum_{k\in\Z}
        \sum_{\vomega\in\W_k\br{\lambda}}
            \1_{A^{\br{-,j}}_\vomega}
            \cC^{\br{\va}}_{\nu_k}
            \pi^{\br{+,j}}_\vomega f
        ,\quad
        \cC_{\nu_k}^{\br{\va}}f\br{\vx}:=
        \int
            f\br{\vx-\vtau}
            e\br{\va\br{\vx}\cdot\vtau}
        d\nu_k\br{\vtau}
        .
    \end{equation}
    We write \(\cC_{\lambda,j}:=\cC_{\nu,\lambda,j}^{\br{\va}}\) and \(\cC_k:=\cC_{\nu_k}^{\br{\va}}\) whenever \(\br{\nu_k}_{k\in\Z}\), and \(\va:\R^{D+1}\to\V\) are clear from the context.
\end{definition}
Remark from the above that we have the following two equalities 
\begin{equation}\label{eq_2_identity_2_abs}
    \cC^{\br{\va}}_{\Delta,\lambda,0}
    =\cC^{\br{\va}}_{\mu,\lambda,0},\quad
    \cC^{\br{\va,\varepsilon}}_{\Delta,j,\vu}
    =\cC^{\br{\va}}_{\alpha,1,j}.
\end{equation}
Consequently, in what follows, we develop a general theory for \(\cC^{\br{\va}}_{\nu,\lambda,j}\) given suitable assumptions on \(\br{\nu_k}_{k\in\Z}\).

\section{Formal time-frequency analysis}\label{sec_formal_tf_ana}
The goal of this section is to develop a formal time-frequency analysis framework for the operator \(\cC^{\br{\va}}_{\nu,\lambda,j}\) introduced in \textbf{Definition \ref{def_cC_nu}} in order to address the estimate 
$$\nrm{\cC^{\br{\va}}_{\nu,\lambda,j} f}_{L^p} \leq C_{\nu,\lambda,p} \nrm{f}_{L^p}\,.$$

Specifically, in this section, we present the following:
\begin{itemize}
    \item An abstract environment for the notion of tiles, trees, energy, and mass---see Sections \ref{assum} and \ref{relevterm};
    \item The statements corresponding to the single tile estimates---see \textbf{Propositions \ref{singtileestimdec}} and \textbf{\ref{singtileestimbd}};
    \item The statements corresponding to the single tree estimates---see \textbf{Lemmas \ref{lem_-_tree}} and \textbf{\ref{lem_gen_tree_c_half}}; 
    \item  The statements of the mass selection algorithm---see \textbf{Lemma \ref{lem_u_shift_mass_sel}} and of the energy selection algorithm---see \textbf{Lemma \ref{lem_cE_2_bd}};
  \item The statements of the three key \textbf{Propositions \ref{thm_simp_w22}, \ref{thm_endpoint_infty_1}} and \textbf{\ref{thm_endpoint_1_infty}}.
\end{itemize}
Throughout the section, we fix and suppress the dependence on \(\br{\nu_k}_{k\in\Z}\), and \(\va:\R^{D+1}\to\V\) whenever they are clear from the context. For simplicity, we will write \(\cC_{\lambda,j}:=\cC_{\nu,\lambda,j}^{\br{\va}}\) and \(\cC_k:=\cC_{\nu_k}^{\br{\va}}\). Since \(\lambda\) will also be fixed, we will suppress its dependency whenever it does not affect the understanding of statements and thus, in particular, we occasionally write \(\W:=\W\br{\lambda}\) and \(\W_k:=\W_k\br{\lambda}\).


\subsection{Assumptions}\label{assum}
One must impose assumptions on \(\br{\nu_k}_{k\in\Z}\) to guarantee boundedness of the model operator \(\cC_{\lambda,j}\). In our case, the support of the measure plays an important role. Indeed, recall \eqref{eq_def_cC_Delta_lambda_0} and \eqref{eq_def_cC_Delta_j_vu}.
We observe that there is some large enough constant \(C\lesssim 1\) such that for all \(k\in\Z\) the following holds:
\begin{equation*}
    \begin{pmatrix}
        2 & 0\\
        0 & 4
    \end{pmatrix}^{-k}
    \hspace{-2ex}
    \cdot
    \supp \mu_k
    \subset B\br{\vnull,C}
    ,\quad
    \begin{pmatrix}
        2 & 0\\
        0 & 4
    \end{pmatrix}^{-k}
    \hspace{-2ex}
    \cdot
    \supp \alpha_k
    \subset B\br{\vu,C}.
\end{equation*}
We will thus make the following assumptions for our formal settings:
\begin{assumption}[\textsf{Support assumptions on measures}]\label{ass_phy_loc_meas}
    There is a vector \(\vu\) and a positive constant \(C\lesssim 1\) such that for all \(k\in\Z\) we have \(\begin{pmatrix}
        2 & 0\\
        0 & 4
    \end{pmatrix}^{-k}
    \hspace{-2ex}
    \cdot
    \supp \nu_k
    \subset B\br{\vu,C}\).
\end{assumption}
Guided by Heisenberg's uncertainty principle, the above suggests that we further decompose \(\pi^{\br{+,j}}_\vomega f\) in \eqref{eq_def_cC_nu} on the physical side to respect the scaling and the shape of the support of the measure \(\nu_k\). To make this precise, we introduce the system of physical (time) rectangles
\begin{equation}
    \I:=\bigsqcup_{k\in\Z}\I_k,\quad
    \I_k:=\BR{
        \valpha \boxplus
        \begin{pmatrix}
            2 & 0\\
            0 & 4
        \end{pmatrix}^k
        \hspace{-1ex}\cdot
            \bR{0,1}^2
        \subset \R^2\::\:
        \valpha\in\Z^2
    }.
\end{equation}
and the system of time-frequency tiles:
\begin{equation}\label{fattile}
    \P\br{\lambda}:=\bigsqcup_{k\in\Z} \P_k\br{\lambda},\quad
    \P_k\br{\lambda}:=\BR{
        \vI\times\vomega\::\:
        \br{\vI,\vomega}\in\I_k\times\W_k\br{\lambda}
    }.
\end{equation}
Similar to \(\W,\W_k\), we will occasionally suppress the \(\lambda\) dependency and write \(\P:=\P\br{\lambda}\) and \(\P_k:=\P_k\br{\lambda}\).
\medskip

With these settled, we next employ a \textsf{Gabor-type representation} of \(f\) to produce an identity of the form
\begin{equation}\label{eq_time_decomp}
    \pi^{\br{+,j}}_\vomega  f=
    \sum_{\vI\in\I_k}f^{\br{+,j}}_{\vI\times\vomega},\quad
    \vomega\in\W_k.
\end{equation}
Since the exact representation will depend on the specific context, we will not specify the precise meaning of \eqref{eq_time_decomp} at this moment. Instead, we focus on the functional properties we expect for such construction; more precisely, we view the map \(\br{f,P}\mapsto f^{\br{+,j}}_P\) as a functional \(F:L^2\br{\R^{D+1}}\times\P\to L^2\br{\R^{D+1}}\) and make the following:
\begin{assumption}[\textsf{Properties of the tile operator}]\label{ass_tf_proj_emb}
For \(\vI\times\vomega\in\P\) and \(f\in L^2\br{\R^{D+1}}\) the following hold:
\begin{itemize}
\item \textsf{Frequency localization}: Representation \eqref{eq_time_decomp} preserves the frequency localization, that is\footnote{Recall \eqref{conedec}.}
\begin{equation}\label{eq_ass_freq_supp}
    \supp \widehat{f^{\br{+,j}}_{\vI\times\vomega}}
    \subset \supp \br{\phi\ast\phi}_{\vomega^{\br{+,j}}}.
\end{equation}

\item \textsf{Tile operator norm and induced spatial localization}: There exists a quantity \(\nrm{f^{\br{+,j}}}_P\geq 0\) such that 
\begin{itemize}
    \item Given \(\vomega\in\W_k\), the following estimate holds:
    \begin{equation}\label{eq_ass_wp_t_l2}
        \nrm{\nrm{f^{\br{+,j}}}_{\vI\times\vomega}}_{\ell^2\br{\vI\in\I_k}}\lesssim
        \nrm{\pi^{\br{+,j}}_\vomega  f}_{L^2}.
    \end{equation}
    \item Given \(f,g\in L^2\br{\R^{D+1}}\) and \(\vI_1\times\vomega_1, \vI_2\times\vomega_2\in\P:=\P\br{\lambda}\), if \(\abs{\vI_2}\leq\abs{\vI_1}\), the following estimate holds:\footnote{Recall Definition \ref{def_jap_n_weight}.}
    \begin{equation}\label{eq_ass_wp_t_in_prod_V}
        \abs{
            \ang{
                f^{\br{+,j}}_{\vI_1\times\vomega_1},
                g^{\br{+,j}}_{\vI_2\times\vomega_2}
            }
        }
        \underset{N}{\lesssim} 
        \nrm{f^{\br{+,j}}}_{\vI_1\times\vomega_1}
        \nrm{g^{\br{+,j}}}_{\vI_2\times\vomega_2}
        \min\br{1,\sqrt{\frac{\lambda^{D+1}\abs{\vI_2}}{\abs{\vI_1}}}}
        \nrm{\chi_{\vI_1}^N}_{L^\infty\br{\vI_2}}.
    \end{equation}
\end{itemize}
\end{itemize}
\end{assumption}
\begin{remark}
    As a direct consequence of \eqref{eq_ass_wp_t_l2} and \eqref{eq_ass_wp_t_in_prod_V}, for \(\cI\subset \I_k\) and \(\vomega\in\W_k\), we have the following:
    \begin{equation}\label{eq_ass_sum_I}
        \bigg\Vert
        \sum_{\vI\in\cI}
            f^{\br{+,j}}_{\vI\times\vomega}
        \bigg\Vert_{L^2}
        \lesssim 
        \nrm{
            \nrm{f^{\br{+,j}}}_{\vI\times\vomega}
        }_{\ell^2\br{\vI\in\cI}}\lesssim \nrm{\pi^{\br{+,j}}_\vomega  f}_{L^2}.
    \end{equation}
\end{remark}
Substituting now the formal identity \eqref{eq_time_decomp} into \textbf{Definition \ref{def_cC_nu}} and dualizing the resulting expression we deduce
\begin{equation}\label{dec}
    \ang{\cC_{\lambda,j} f,g}=
        \sum_{k\in\Z}
        \sum_{\vI\times\vomega\in\P_k\br{\lambda}}
            \ang{
                \cC_k
                    f^{\br{+,j}}_{\vI\times\vomega},
                \1_{A^{\br{-,j}}_\vomega} g
            }.
\end{equation}
In what follows, it will be convenient to isolate the individual term in \eqref{dec}:
\begin{definition}[\textsf{Single tile bilinear form}]\label{def_sing_tile_bili} With the previous notations we set
    \begin{equation}\label{sgtile}
        \Lambda^{\br{j}}_{\vI\times\vomega}\br{f,g}:=
        \ang{
            \cC_k
                f^{\br{+,j}}_{\vI\times\vomega},
            \1_{A^{\br{-,j}}_\vomega} g
        },\quad
        \vI\times\vomega\in\P_k.
    \end{equation}
\end{definition}
Finally, we notice that as a trivial consequence of the triangle inequality, we have 
$$\left|\ang{\cC_{\lambda,j} f,g}\right|\leq\nrm{\Lambda^{\br{j}}_P\br{f,g}}_{\ell^1\br{P\in\P\br{\lambda}}}\,,$$
and thus the \(L^p\) boundedness of \(\cC_{\lambda,j}\) reduces to the following estimate: for any $E'\subset E$ with $\abs{E'}\geq \abs{E}/2$ one has 
\begin{equation}\label{eq_nu_wpp}
    \nrm{\Lambda^{\br{j}}_P\br{\1_F,\1_{E'}}
    }_{\ell^1\br{P\in\P\br{\lambda}}}\leq C_{\nu,\lambda,p} \abs{F}^{\frac{1}{p}}\abs{E}^{\frac{p-1}{p}}\,.
\end{equation}

\subsection{Relevant terminologies: order relations among tiles, trees, mass and energy}\label{relevterm}

In what follows, with the proper adaptations to our context, we recall the foundational concepts employed for proving the celebrated theorem of Carleson on the pointwise convergence of Fourier Series (\cite{c1}, \cite{f}, \cite{lt3}):

\begin{definition}[\textsf{Order relation}]\label{def_tile_leqs}
    Given \(P_i=\vI_i\times\vomega_i\in\P\), we define the following relations:
    \begin{align*}
        P_0\leq P_1\iff & \vI_0 \subseteq \vI_1\text{ and }\vomega_0\supseteq\vomega_1\\
        P_0\leq_{\pm,j} P_1\iff & \vI_0 \subseteq \vI_1\text{ and }\vomega^{\br{\pm,j}}_0\supseteq\vomega^{\br{\pm,j}}_1
    \end{align*}
    Note that \(\leq,\leq_{+,j},\leq_{-,j}\) are all partial orderings on \(\P\).
\end{definition}

\begin{definition}[\textsf{Trees}]\label{def_tree}
    A tree \(\cT\subset \P\) is a
    collection that comes with a tree top \(P_\cT:=\vI_\cT\times\vomega_\cT\in\P\) such that \(P\leq P_\cT\) for all \(P\in\cT\). If additionally, \(P\leq_{\pm,j} P_\cT\) for all \(P\in\cT\), we call \(\cT\) a \(\br{\pm,j}\)-tree.
\end{definition}
\begin{observation}[\textsf{Frequency-scale correlation}]\label{obs_tree_freq_scale_rel}
    Let \(\vomega'\in\W_{k'}\). For any \(k\leq k'\), there is a unique \(\vomega\in\W_k\) such that \(\vomega'\subset \vomega\).
    As a direct consequence, 
    for a tree \(\cT\subset\P\) and \(k\in\Z\) with \(\cT_k:=\cT\cap \P_k \neq \varnothing\), there is a unique \(\vomega_{\cT_k}\in\W_k\) such that for all \(\vI\times\vomega\in\cT_k\) we have \(\vomega=\vomega_{\cT_k}\). 
\end{observation}

\begin{definition}[\(L^2\) \textsf{and} \(\mathrm{BMO}\) \textsf{Energy}]\label{def_L2_BMO_energy}
Given a sub-collection of tiles \(\cP\subset \P\) and a function  \(f\in L^2\) we set
\begin{equation}
    \nrm{f^{\br{+,j}}}_{L^2\br{\cP}}:=
    \nrm{
        \nrm{f^{\br{+,j}}}_P
    }_{\ell^2\br{P\in\cP}}
    ,\quad
    \nrm{f^{\br{+,j}}}_{\bmo_2\br{\cP}}:=
    \sup_{\substack{
        \cT\subset \cP\\
        \cT\text{ be }-,j\text{-tree}
    }}
    \nrm{f^{\br{+,j}}}_{L^2\br{\cT}}
    /\abs{\vI_\cT}^{\frac{1}{2}}.
\end{equation}
\end{definition}

\begin{definition}[\textsf{Weighted measure adapted to an interval}]
    Let \(N\gg 1\) and \(\vI\subset \R^{D+1}\) be an interval. We define the weighted measure \(d\mu^{\br{N}}_\vI\br{\vx}:=\abs{\vI}^{-1}\chi_\vI^N\br{\vx}d\vx\) on \(\R^{D+1}\).
\end{definition}

\begin{remark}\label{rmk_N_freedom}
    The above choice of \(N\gg 1\) is irrelevant for our purposes and is allowed to change from line to line. We thus, for notational simplicity, suppress the sup-script \(\br{\cdot}^{\br{N}}\) and write \(\mu_\vI=\mu^{\br{N}}_\vI\) henceforth.
\end{remark}

\begin{definition}[\textsf{Shifted Mass}]\label{def_shift_mass} Given a sub-collection of tiles \(\cP\subset \P\) and a measurable set \(E\subset \R^2\) we let
\begin{equation}
        \cM_\cP\br{E}:=\sup_{\substack{
            P\in\cP
        }}\sup_{\substack{
            \vI\times\vomega\in\P\\
            P\leq\vI\times\vomega
        }}
        \max_{\substack{
            \vu_\ast\in\cU
        }}
        \mu_{\vu_\ast\boxplus\vI}\br{E\cap A_\vomega},
    \end{equation}
    where \(\cU\) is defined as follows:
    \begin{equation}\label{eq_cU_cV_def}
        \cU:=\BR{
            \br{u_0/2^j,u_1/4^j}\in\R^2\setminus\mr{-1,1}^2\::\:
            j\in\N\sqcup\BR{0}
        }\sqcup\BR{-1,0,1}^2.
    \end{equation}
\end{definition}

\subsection{Time-frequency localization and single tile estimates}
As is suggested by the notation we have chosen, the function \(f^{\br{+,j}}_{\vI\times\vomega}\) shall carry information about the physical location \(\vI\), and \textbf{Assumption \ref{ass_tf_proj_emb}} does so in a \emph{relational time-frequency} sense---see \cite{LVQuadCarl}, \cite{LVPolynCarl}, and more recently, \cite{GALVHybcurves}. Indeed, the following statement implicitly measures the concentration of the time-frequency localization captured within expression \eqref{sgtile}:
\begin{definition}[\textsf{A Heisenberg norm scale: $\{\nrm{\Lambda}_c\}_{c\in\mr{0,1/2}}$}]\label{def_sing_tile_est}
$\smallskip$

Let \(\vu\) be as in \textbf{Assumption \ref{ass_phy_loc_meas}}. For \(j\in\{0,1\}\), \(c\in\mr{0,1/2}\) and \(\lambda\geq 1\), define \(\nrm{\Lambda}_c:=\nrm{\Lambda}_{c,\lambda,j} \geq 0\) to be the smallest quantity such that for all \(f,g\in L^2\) and all \(\vI\times\vomega\in\P\br{\lambda}\) the following inequality holds:
    \begin{equation}\label{eq_sing_tile_est}
        \abs{\Lambda^{\br{j}}_{\vI\times\vomega}\br{f,g}}
        \leq \nrm{\Lambda}_c \mu_{\vu\boxplus \vI}^c\br{\supp g\cap A^{\br{-,j}}_\vomega}
        \nrm{f^{\br{+,j}}}_{\vI\times\vomega}
        \nrm{\1_{A^{\br{-,j}}_\vomega}g}_{L^2\br{d\mu_{\vu\boxplus\vI}}}\abs{\vI}^{\frac{1}{2}}.
    \end{equation}
\end{definition}
\begin{remark}
    The shift in \(\mu_{\vu\boxplus\vI}\) can be seen by setting \(c=0\). To be more precise, recall \textbf{Definitions \ref{def_cC_nu} and \ref{def_sing_tile_bili}}. For \(\vI\times\vomega\in\P_k\), we expect the operator \(\cC_k\) to shift the information physically by increments of the corresponding scale in \(\vu\) direction.
\end{remark}

A first key result that will be proved later is the following:

\begin{proposition}[\textsf{High resolution analysis: single tile decay estimate addressing $\cC\cR^{\br{\va}}_{\Delta,\lambda,0}$}]\label{singtileestimdec}
$\smallskip$

There exists \(\sigma>0\) such that for $j=0$,  \(\vu=\vnull\), \(\br{\nu_k}_{k\in\Z}=\br{\mu_k}_{k\in\Z}\) and any $\lambda\geq 1$ one has 
\begin{equation}\label{sgtiledecest}
 \nrm{\Lambda}_0\lesssim \lambda^{-\sigma},\quad\textrm{that is},\quad \abs{\Lambda^{\br{j}}_{\vI\times\vomega}\br{f,g}}
    \lesssim \lambda^{-\sigma}\nrm{f}_{\vI\times\vomega}\nrm{\1_{A^{\br{-,0}}_\vomega}g}_{L^2\br{d\mu_\vI}}.
\end{equation}
\end{proposition}

We will also need the easier counterpart

\begin{proposition}[\textsf{Low resolution analysis: single tile $O(1)$-estimate addressing both $\cC\cR^{\br{\va}}_{\Delta,\lambda,0}$ and $ \cC^{\br{\va,\varepsilon}}_{\Delta,j,\vu}$}]\label{singtileestimbd}
$\smallskip$

\noindent i) For $j=0$,  \(\vu=\vnull\), \(\br{\nu_k}_{k\in\Z}=\br{\mu_k}_{k\in\Z}\) and any $\lambda\geq 1$ one has in \eqref{eq_sing_tile_est}
\begin{equation}\label{sgtileest1}
  \nrm{\Lambda}_{1/4}\lesssim 1\,.
\end{equation}

\noindent ii) For $j\in\{0,1\}$,  \(\vu\in\Z^2\),  \(\br{\nu_k}_{k\in\Z}=\br{\alpha_k}_{k\in\Z}\)  and  $\lambda=1$ one has in \eqref{eq_sing_tile_est}
\begin{equation}\label{sgtileest2}
  \nrm{\Lambda}_{1/2}\lesssim 1\,.
\end{equation}

\end{proposition}

\subsection{Single tree estimates}\label{sgtree}

For integer \(k\in\Z\), we consider two forms of Littlewood-Paley frequency projections:
\begin{equation}\label{eq_LP_freq_proj}
    \widehat{\pi^k f}\br{\vxi}:= 
    \varphi^{\otimes 2}\br{
        \begin{pmatrix}
            2 & 0\\
            0 & 4
        \end{pmatrix}^k
        \cdot
        \frac{\vxi}{2\lambda}
    }
    \widehat{f}\br{\vxi},\quad
    \pi_k:=\pi^k-\pi^{k+5}
\end{equation}
We define two related operators that are relevant to our analysis as below:
\begin{definition}[\textsf{Singular integral type and maximal type operators;  the $\cL$ correction factor}]\label{def_sio_max}
    \begin{equation}
        T^{\cZ}_\nu f\br{\vx}:=
        \sum_{k\in\cZ} \pi_k f \ast d\nu_k\br{\vx},\quad
        M_\nu f\br{\vx}:=
        \sup_{k\in\Z}
            \abs{\pi_k f}\ast \abs{d\nu_k}\br{\vx},\quad
            \cZ\subset\Z.
    \end{equation}
\end{definition}
Let \(\nrm{T^{\cZ}_\nu}_2,\nrm{M_\nu}_2\geq 0\) denote the strong type \(\br{2,2}\) operator norms and consider the following quantities:
\begin{equation}\label{eq_T_nu_sup_CZO}
    \nrm{T_\nu}_2:=\sup_{\cZ\subset \Z}
    \nrm{T^\cZ_\nu}_2
\end{equation}
and the \(\log\) size correction factor
\begin{equation}\label{defL}
    \cL:=\log\br{e+\nrm{M_\nu}_2/\nrm{T_\nu}_2+\lambda\br{1+\abs{\vu}}}.
\end{equation}

\medskip

With these settled, we now state two key single tree estimates---their proofs are given in \textsc{Section \ref{sec_pf_key_lemmas}}:
\begin{lemma}[\textsf{The improved \(\br{-,j}\)-tree estimate: gain over the ($l^2$-$l^2$) dual pair case}]\label{lem_-_tree}
For a \(\br{-,j}\)-tree \(\cT\subset \P\), we have:
    \begin{equation}
        \nrm{\Lambda^{\br{j}}_P\br{f,\1_E}}_{\ell^1\br{P\in\cT}}
        \lesssim \br{\nrm{\Lambda}_c+\nrm{T_\nu}_2}\,\cL\,\nrm{f^{\br{+,j}}}_{L^2\br{\cT}}M_\cT^{\frac{1}{2}+c}\br{E}\abs{\vI_\cT}^{\frac{1}{2}}.
    \end{equation}
\end{lemma}
\begin{lemma}[\textsf{The general tree estimate: the ($l^\infty$-$l^1$) dual pair case}]\label{lem_gen_tree_c_half}
For a tree \(\cT\subset \P\), we have:
    \begin{equation}
        \nrm{\Lambda^{\br{j}}_P\br{f,\1_E}}_{\ell^1\br{P\in\cT}}
        \lesssim \br{\nrm{\Lambda}_{\frac{1}{2}}+\nrm{T_\nu}_2}\,\cL\,\nrm{f^{\br{+,j}}}_{\bmo_2\br{\cT}}M_\cT\br{E}\abs{\vI_\cT}.
    \end{equation}
\end{lemma}

\subsection{Mass selection algorithm}\label{subsec_mass_sel}
The following mass selection lemma is a slight modification of the standard one implemented in Lacey and Thiele's approach to the Carleson operator. Its proof will be postponed for \textsc{Section \ref{sec_pf_key_lemmas}}.
\begin{lemma}[\textsf{Mass selection}]\label{lem_u_shift_mass_sel}
    For all \(\varsigma>0\), \(\cP\subset \P\) finite, and \(E\subset \R^2\) measurable, there is finite disjoint collection of trees \(\BR{\cT_i}_i\) with \(\cT_i\subset \cP\) such that the following two inequalities hold:
    \begin{equation}
        \cM_{\cP\setminus\bigsqcup_i\cT_i}\br{E}\leq \varsigma
        \qquad\textrm{and}\qquad
        \sum_{i}\abs{\vI_{\cT_i}}
        \lesssim \frac{\log\br{e+\abs{\vu}}\abs{E}}{\varsigma}.
    \end{equation}
\end{lemma}

\subsection{Energy selection algorithm}\label{energy_select}

In this subsection, we present the energy selection algorithm, the proof of which will be displayed in \textsc{Section \ref{sec_pf_key_lemmas}}.

\begin{definition}[\textsf{Strong disjointness}]\label{def_strong_disj}
    A sequence of trees \(\BR{\cT_i}_i\) are strongly disjoint if for any two tiles \(\vI_m\times\vomega_m\in\cT_m\) and \(\vI_n\times\vomega_n\in\cT_n\) with \(\vomega^{\br{+,1}}_m\cap\vomega^{\br{+,1}}_n\neq\varnothing\), the following holds:
    \begin{equation*}
        \br{m<n
        \text{ or }
        \vomega_m
        \subsetneq 
        \vomega_n }
        \implies 
        3\vI_{\cT_m}\cap \vI_n=\varnothing.
    \end{equation*}
The definition presented is a slight modification of the \textbf{Definition 6.8} in \cite{muscalu2013classicalv2}.
\end{definition}
With the above notion of strong disjointness, we state the following key estimate:
\begin{lemma}[\textsf{Bessel-type inequalities}]\label{lem_bessel}
    Let \(\BR{\cT_i}_i\) be a sequence of strongly disjoint \(\br{-,1}\)-trees. If for all \(i\) and \(\vI\times\vomega\in\cT_i\), the two estimates \(\nrm{f^{\br{+,1}}}_{\vI\times\vomega}\lesssim\abs{\vI}^{\frac{1}{2}}\) and \(\nrm{f^{\br{+,1}}}_{L^2\br{\cT_i}}\eqsim\abs{\vI_{\cT_i}}^{\frac{1}{2}}\) hold, we have \(\sum_i\abs{\vI_{\cT_i}}\lesssim \nrm{f}_{L^2}^2\).
\end{lemma}

\begin{remark}
The above lemma appears implicitly in most of the literature that utilizes modern time-frequency analysis techniques. The interested reader can refer to \textsc{Section 6.5} in \cite{muscalu2013classicalv2} for further reference.\footnote{Our present formulation corresponds to the \(n=0\) case therein.} However, we will in effect prove a generalization of the above--- see \textbf{Lemma \ref{lem_Hdim_bessel}} in \textsc{Section \ref{sec_pf_gen_energy_est}}.
\end{remark}

With these, we have the following:

\begin{lemma}[\textsf{Energy selection}]\label{lem_cE_2_bd}
    For all \(\cP\subset \P\) finite and \(f\in L^2\), there is a finite collection \(\T\) of mutually disjoint trees in \(\cP\) satisfying the following two estimates:
    \begin{equation}
        \nrm{f^{\br{+,j}}}_{\bmo_2\br{\cP\setminus\bigsqcup\T}}\leq \frac{1}{2}\nrm{f^{\br{+,j}}}_{\bmo_2\br{\cP}}\quad\qquad\textrm{and}\qquad\quad
        \sum_{\cT\in\T}\abs{\vI_\cT}
        \lesssim
        \br{\frac{\nrm{f}_{L^2}}{\nrm{f^{\br{+,j}}}_{\bmo_2\br{\cP}}}}^2.
    \end{equation}
\end{lemma}

\subsection{The final three main propositions}\label{subsec_mass_sel_ez_22_bd}

We end this chapter by stating (the last) three key propositions that will prove quintessential in addressing our main result in the one-dimensional setting:
\begin{itemize}
\item \textbf{Proposition \ref{singtileestimdec}} from before and \textbf{Proposition \ref{thm_simp_w22}} below are the cornerstones in the resolution of \textbf{Theorem \ref{thm_HF_st_disc}}. The proof of the former proposition will be presented in \textsc{Section \ref{sec_est_of_sing_tile}} while the proof of the latter will be delivered in \textsc{Section \ref{sec_pf_HF_st_disc}}. 

\item \textbf{Propositions \ref{thm_endpoint_infty_1}} and \textbf{\ref{thm_endpoint_1_infty}} are essential in the resolution of \textbf{Theorem \ref{thm_phy_shift_disc}} and both these results will be proved in \textsc{Section \ref{sec_pf_physift_disc}}.
\end{itemize}

With these clarifications, we are now ready to state our key propositions: 

\begin{proposition}[\textsf{Weak} \(\br{2,2}\) \textsf{estimate}]\label{thm_simp_w22}
    For all \(c\in\bR{0,\frac{1}{2}}\), \(\lambda\geq 1\), \(f\in L^2\) and \(E\subset \R^2\) measurable, the following estimate holds:
    \begin{equation*}
        \nrm{\Lambda^{\br{0}}_P\br{f,\1_E}}_{\ell^1\br{P\in\P\br{\lambda}}}
        \lesssim
        \frac{
            \nrm{\Lambda}_c+\nrm{T_\nu}_2
        }{\sqrt{c}}\cdot\cL^{\frac{3}{2}}
        \nrm{f}_{L^2}\abs{E}^{\frac{1}{2}}.
    \end{equation*}
\end{proposition}

The general restricted weak type \(\br{p,p}\) estimates for general \(j\) require a more refined argument, for which we first need to introduce the following two definitions:
\begin{definition}[\textsf{Outer} \(L^{2,\infty}\) \textsf{embedding estimate}]\label{def_out_L2_infty_energy_emb}
    Let \(\cE_2:=\cE_{2,\lambda,j}\geq 0\) be the smallest constant such that for all \(\varsigma>0\), \(\cP\subset\P:=\P\br{\lambda}\) finite, and \(f\in L^2\), there is finite disjoint collection of trees \(\BR{\cT_i}_i\) with \(\cT_i\subset \cP\) satisfying the following two estimates:
    \begin{equation}\label{eq_def_out_L2_infty_energy_emb}
        \nrm{f^{\br{+,j}}}_{\bmo_2\br{\cP\setminus\bigsqcup_i\cT_i}}\leq \varsigma\qquad\textrm{and}\qquad
        \sum_i\abs{\vI_{\cT_i}}\leq \frac{\cE_2^2\nrm{f}_{L^2}^2}{\varsigma^2}.
    \end{equation}
\end{definition}

\begin{definition}[\textsf{Outer} \(L^\infty\) \textsf{embedding estimate}]\label{def_out_Linfty_energy_emb}
    Let \(\cE_\ast:=\cE_{\ast,\lambda,j}\geq 0\) be the smallest constant such that for all \(\cP\subset\P:=\P\br{\lambda}\) finite and \(f\in L^1\cap L^\infty\),
    \begin{equation}\label{eq_Linfty_emb_sup}
        \nrm{f^{\br{+,j}}}_{\bmo_2\br{\cP}}\leq \cE_\ast
        \sup_{\substack{\varnothing\neq\cT\subset\cP\\
        \cT\text{ be tree}\\
        \vI_\cT\subset \vI_\cP
        }}\nrm{f}_{L^1\br{d\mu_{\vI_\cT}}},
    \end{equation}
    where \(\vI_\cP\in\I\) denotes the minimal interval that contains \(\bigcup_{\vI\times\vomega}\vI\). 
\end{definition}

With these settled, we conclude our section with the following two endpoint estimates:
\begin{proposition}[\textsf{Estimate near} \(L^\infty\times L^1\)]\label{thm_endpoint_infty_1}
    For any $j\in\{0,1\}$ and \(E,F\subset \R^2\) with finite measure, we have:
    \begin{equation}
        \nrm{\Lambda^{\br{j}}_P\br{\1_F,\1_E}}_{\ell^1\br{P\in\P}}\lesssim 
        \br{\nrm{\Lambda}_{\frac{1}{2}}+\nrm{T_\nu}_2}
        \cE_\ast\,\cL^2
        \log\br{e+\frac{\cE_2^2\abs{F}}{\cE_\ast^2\abs{E}}}\abs{E}.
    \end{equation}
\end{proposition}

\begin{proposition}[\textsf{Estimate near} \(L^1\times L^\infty\)]\label{thm_endpoint_1_infty}
    For any $j\in\{0,1\}$ and \(E,F\subset \R^2\) with finite measure, there is a majorant subset \(E'\subset E\) satisfying \(\abs{E'}\geq \abs{E}/2\) such that:
    \begin{equation}
        \nrm{\Lambda^{\br{j}}_P\br{\1_F,\1_{E'}}}_{\ell^1\br{P\in\P}}\lesssim 
        \br{\nrm{\Lambda}_{\frac{1}{2}}+\nrm{T_\nu}_2}
        \cE_\ast\,\cL^2
        \log\br{e+\frac{\cE_2^2\abs{E}}{\cE_\ast^2\abs{F}}}\abs{F}.
    \end{equation}
\end{proposition}

\section{Reduction of \textbf{Theorem \ref{thm_HF_st_disc}} to \textbf{Proposition \ref{singtileestimdec}}}\label{sec_pf_HF_st_disc}

In this section, our goal is to reduce the proof of  \textbf{Theorem \ref{thm_HF_st_disc}} to the single tile decay estimate evoked in \textbf{Proposition \ref{singtileestimdec}} subject to the key lemmas stated in \textsc{Section \ref{sec_formal_tf_ana}}. In view of our goal, we will appeal to the formal setting introduced in the latter section and, recalling  \textbf{Definition \ref{def_cC_nu}}, throughout this section 
$$\textrm{we set}\:\:\: j=0,\:\:\vu=\vnull,\:\:\textrm{and}\:\:\br{\nu_k}_{k\in\Z}=\br{\mu_k}_{k\in\Z}\,.$$

In what follows, we present \textsf{the main steps of our program} in achieving the above-stated goal:
\begin{itemize}
    \item Construct a well-behaved Gabor decomposition verifying \eqref{eq_time_decomp};
    \item Verify \textbf{Assumption \ref{ass_tf_proj_emb}};
    \item Prove \(\nrm{\Lambda}_c:=\nrm{\Lambda}_{c,\lambda,0}\lesssim \lambda^{-\sigma}\) for some \(c,\sigma\in\br{0,\frac{1}{2}}\) subject to \textbf{Proposition \ref{singtileestimdec}};
    \item Prove \(\nrm{T_\nu}_2\lesssim \lambda^{-\frac{1}{2}}\) and \(\nrm{M_\nu}_2\lesssim 1\);
    \item Prove \textbf{Proposition \ref{thm_simp_w22}} subject to \textbf{Lemma \ref{lem_-_tree}} and \textbf{Lemma \ref{lem_u_shift_mass_sel}}.
\end{itemize}

With these settled, we can immediately derive \textbf{Theorem \ref{thm_HF_st_disc}} as a consequence of the above program.

Indeed, assuming \textbf{Proposition \ref{thm_simp_w22}}, we deduce
\begin{equation*}
    \abs{
        \ang{
            \cC^{\br{\va}}_{\Delta,\lambda,0}f,
            \1_E
        }
    }
    =
    \abs{
        \ang{
            \cC^{\br{\va}}_{\nu,\lambda,0}f,
            \1_E
        }
    }
    \leq
    \nrm{\Lambda^{\br{0}}_P\br{f,\1_E}}_{\ell^1\br{P\in\P\br{\lambda}}}
    \lesssim
    \frac{
        \nrm{\Lambda}_c+\nrm{T_\nu}_2
    }{\sqrt{c}}\cdot\cL^{\frac{3}{2}}
    \nrm{f}_{L^2}\abs{E}^{\frac{1}{2}}
\end{equation*}
for all \(E\subset \R^2\) measurable. Recalling now \eqref{defL}, by direct calculation, we further dominate the above by
\begin{equation*}
    \lesssim \frac{\lambda^{-\sigma}+\lambda^{-\frac{1}{2}}}{\sqrt{c}}
    \log^{\frac{3}{2}}\br{e+1/\lambda^{-\frac{1}{2}}+\lambda} \nrm{f}_{L^2}\abs{E}^{\frac{1}{2}}
    \underset{c}{\lesssim}\lambda^{-\sigma/2}\nrm{f}_{L^2}\abs{E}^{\frac{1}{2}}.
\end{equation*}
We thus conclude \textbf{Theorem \ref{thm_HF_st_disc}} with \(\epsilon=\sigma/2\). It remains to complete the above-mentioned program.

\subsection{Gabor decompositions \texorpdfstring{\eqref{eq_time_decomp}}{(time decomp)} adapted to \texorpdfstring{$\cC\cR^{\br{\va}}_{\Delta,\lambda,0}$}{CR} : integral representations}
Directed by the LGC-method treatment of the non-resonant variant of the problem, see \cite{HsuL24}, we decompose the frequency projection operator \(\pi^{\br{+,0}}_\vomega f\) at the scale that \emph{linearizes} the phase information of the (original) Carleson-Radon transform. To achieve this, we observe the following:
\begin{proposition}\label{prop_Gabor_dec_conti}
    Given \(f\in L^2\br{\R^{D+1}}\) and \(\Phi\in \mathcal{S}\br{\R^{D+1}}\) with \(\nrm{\Phi}_{L^2}=1\),
    the following formula holds:
    \begin{align*}
        f= &
        \int_{\R^{D+1}\times\R^{D+1}}
            \ang{f,\Mod_\vxi\Tr_\vx\Phi}
            \Mod_\vxi\Tr_\vx\Phi
        d\vx d\vxi\\
        =&
        \int_{\R^{D+1}\times\R^{D+1}}
            \ang{f,\Tr_\vx\Mod_\vxi\Phi}
            \Tr_\vx\Mod_\vxi\Phi
        d\vx d\vxi.
    \end{align*}
    In fact, for any \(\Psi\in L^1_{\operatorname{loc}}\br{\R^{D+1}}\), one has the following localized analog:
    \begin{align*}
        \br{\abs{\Phi}^2\ast\Psi} f= &
        \int_{\R^{D+1}\times\R^{D+1}}
            \ang{f,\Mod_\vxi\Tr_\vx\Phi}
            \Mod_\vxi\Tr_\vx\Phi
            \cdot
            \Psi\br{\vx}
        d\vx d\vxi\\
        =&
        \int_{\R^{D+1}\times\R^{D+1}}
            \ang{f,\Tr_\vx\Mod_\vxi\Phi}
            \Tr_\vx\Mod_\vxi\Phi
            \cdot
            \Psi\br{\vx}
        d\vx d\vxi.
    \end{align*}
\end{proposition}
Take now \(k\in\Z\) and \(\vI\in \I_k\). We set:
\begin{equation}
    \Phi=\varphi_k:=\br{\Dil^2_{2^k /\sqrt{\lambda}}\widecheck{\phi}/\nrm{\phi}_{L^2}}\otimes
    \br{\Dil^2_{4^k /\sqrt{\lambda}}\widecheck{\phi}/\nrm{\phi}_{L^2}},\quad
    \Psi=\1_\vI
\end{equation}
and define
\begin{equation}\label{eq_def_varphi_vomega_xyuv}
    \varphi_{k,\vx,\vxi}:=
    \Tr_\vx\Mod_\vxi\varphi_k,\quad
    \chi^\Phi_\vI:=\abs{\Phi}^2\ast\Psi=\abs{\varphi_k}^2\ast\1_\vI.
\end{equation}
\textbf{Proposition \ref{prop_Gabor_dec_conti}} suggests that we set:
\begin{equation}\label{eq_planar_tf_proj}
    f^{\br{+,0}}_{\vI\times\vomega}:=\chi^\Phi_\vI \pi^{\br{+,0}}_\vomega  f
    =
    \iint_{\vI}
        \iint_{\R^2}
            \ang{
                \pi^{\br{+,0}}_\vomega  f,
                \varphi_{k,\vx,\vxi}
            }
            \varphi_{k,\vx,\vxi}
        d\vxi
    d\vx,\quad
    \vI\times\vomega\in\P_k:=\P_k\br{\lambda}.
\end{equation}
Note that by \textbf{Proposition \ref{prop_Gabor_dec_conti}}, we trivially have \eqref{eq_time_decomp}. On the same note, for \(\vI\times\vomega\in\P_k\), we set
\begin{equation}\label{eq_wp_nrm}
    \nrm{f^{\br{+,0}}}_{\vI\times\vomega}
    :=
    \nrm{\pi^{\br{+,0}}_\vomega  f\br{\vy}}_{L^2\br{\chi^\Phi_\vI\br{\vy}d\vy}}
    =
    \nrm{
        \ang{
            \pi^{\br{+,0}}_\vomega  f,
            \varphi_{k,\vx,\vxi}
        }
    }_{L^2\br{d\vx d\vxi,\vI\times\R}}.
\end{equation}

\subsection{Verification of \texorpdfstring{\eqref{eq_ass_freq_supp}}{}, \texorpdfstring{\eqref{eq_ass_wp_t_l2}}{}, and \texorpdfstring{\eqref{eq_ass_wp_t_in_prod_V}}{} in Assumption \ref{ass_tf_proj_emb}}\label{subsec_ver_tf_proj_emb}
Starting with \eqref{eq_ass_freq_supp}, we notice the frequency localization of \(f^{\br{+,0}}_{\vI\times\vomega}\) is captured by the set of parameters
\begin{equation*}
    \Xi:=
    \bigcup_{\vx\in\vI}
    \BR{
        \vxi\in\R^2\::\:
        \ang{\pi^{\br{+,0}}_\vomega  f,
            \varphi_{k,\vx,\vxi}}\neq 0
    }
    ,\quad
    \vI\times\vomega\in\P_k\br{\lambda}
\end{equation*}
which have non-trivial contribution to the integral in \eqref{eq_planar_tf_proj}.
Direct calculation shows the following relations:
\begin{align*}
    \Xi\subset &
    \bigcup_{\vx\in\vI}
    \BR{
        \vxi\in\R^2\::\:
        \supp\widehat{\pi^{\br{+,0}}_\vomega  f}\cap
        \supp\widehat{\varphi_{k,\vx,\vxi}}
        \neq \varnothing
    }
    \\
    \subset &
    \BR{
        \vxi\in\R^2\::\:
            \supp\phi_{
                \vomega^{\br{+,0}}}\cap
            \supp
            \br{
                \Tr_{\xi_0}\Dil^\infty_{\frac{\sqrt{\lambda}}{2^k}}\phi \otimes 
                \Tr_{\xi_1}\Dil^\infty_{\frac{\sqrt{\lambda}}{4^k}}\phi
            }
        \neq \varnothing
    }
    \subset  \supp \br{\phi\ast\Dil^1_{1/\sqrt{\lambda}}\phi}_{\vomega^{\br{+,0}}}.
\end{align*}
As a direct consequence, we have:
\begin{equation*}
    \supp \widehat{f^{\br{+,0}}_{\vI\times\vomega}}\subset 
    \supp \br{\phi\ast\Dil^1_{1/\sqrt{\lambda}}\phi\ast\Dil^1_{1/\sqrt{\lambda}}\phi}_{\vomega^{\br{+,0}}}
    \subset \supp\br{\phi\ast\phi}_{\vomega^{\br{+,0}}}
\end{equation*}
whenever \(\lambda \geq 4\). This confirms the frequency support assumption \eqref{eq_ass_freq_supp}.
To show \eqref{eq_ass_wp_t_l2}, we apply the first identity of \eqref{eq_wp_nrm} and expand the left hand side of \eqref{eq_ass_wp_t_l2}:
\begin{align*}
    \nrm{\nrm{f^{\br{+,0}}}_{\vI\times\vomega}}_{\ell^2\br{\vI\in\I_k}}
    = &
    \nrm{
        \nrm{
            \pi^{\br{+,0}}_\vomega f\br{\vx}
        }_{L^2\br{\chi^\Phi_\vI\br{\vx}d\vx}}
    }_{\ell^2\br{\vI\in\I_k}}\\
    = &
    \nrm{
        \pi^{\br{+,0}}_\vomega f\br{\vx}
    }_{L^2\br{\sum_{\vI\in\I_k}\chi^\Phi_\vI\br{\vx}d\vx}}
    =\nrm{\pi^{\br{+,0}}_\vomega f}_{L^2}.
\end{align*}
Finally, to prove \eqref{eq_ass_wp_t_in_prod_V}, we let \(\vI_1\times\vomega_1\in\P_{k_1}\) and \(\vI_2\times\vomega_2\in\P_{k_2}\) with \(k_2\leq k_1\) and expand the inner product expression using the integral identity in \eqref{eq_planar_tf_proj}:
\begin{align}\label{eq_gaba_wp_inner_prod}
    \abs{\ang{
        f^{\br{+,0}}_{\vI_1\times\vomega_1},
        g^{\br{+,0}}_{\vI_2\times\vomega_2}
    }}
    \leq 
    \int &
    \abs{\ang{\pi^{\br{+,0}}_{\vomega_1}f,\varphi_{k_1,\vx_1,\vxi_1}}}
    \1_{\vI_1}\br{\vx_1}\1_{\vomega^{\br{+,0}}_1}\br{\vxi_1}\nonumber\\
    \cdot &
    \abs{\ang{\pi^{\br{+,0}}_{\vomega_2}g,
    \varphi_{k_2,\vx_2,\vxi_2}}}
    \1_{\vI_2}\br{\vx_2}\1_{\vomega^{\br{+,0}}_2}\br{\vxi_2}\nonumber\\
     \cdot &
     \abs{\ang{
         \varphi_{k_1,\vx_1,\vxi_1},
         \varphi_{k_2,\vx_2,\vxi_2}
     }}
     d\vx_1 d\vxi_1 d\vx_2 d\vxi_2.
\end{align}
By Cauchy-Schwarz inequality and the trivial bound \(\1_{\vI_k}\underset{N}{\lesssim} \chi_{\vI_k}^N\), we dominate the above with
\begin{equation}\label{eq_ass_wp_pf_f_g_phi}
\!
    \underset{N}{\lesssim}
    \nrm{f^{\br{+,0}}}_{\vI_1\times\vomega_1}
    \nrm{g^{\br{+,0}}}_{\vI_2\times\vomega_2}
    \nrm{
        \ang{
            \varphi_{k_1,\vx_1,\vxi_1},
            \varphi_{k_2,\vx_2,\vxi_2}
        }
        \prod_{k=0,1}
        \chi_{\vI_k}^N\br{\vx_k}
        \1_{\vomega^{\br{+,0}}_k}\br{\vxi_k}
    }_{L^2\br{d\vx_1 d\vxi_1
    d\vx_2 d\vxi_2}}
    .
\end{equation}
Now, on the one hand, due to the frequency localization of the two wave-packets, we deduce that the inner product contributes to the third factor in the above expression only when the following relations hold simultaneously:
\begin{equation*}
    \abs{\xi_{01}-\xi_{02}}\lesssim 2^{-k_2}\sqrt{\lambda},\quad
    \abs{\xi_{11}-\xi_{12}}\lesssim 4^{-k_2}\sqrt{\lambda},\quad
    \vxi_1:=\br{\xi_{01},\xi_{11}}
    \in\vomega^{\br{+,0}}_1.
\end{equation*}
On the other hand, by the standard fact
\begin{equation}\label{eq_1dim_std_fact}
    \int_\R
        \ang{\frac{z-c_1}{l_1}}^{-N}
        \ang{\frac{z-c_2}{l_2}}^{-N}
    dz
    \underset{N}{\lesssim}
    \min\br{l_1,l_2}\cdot
    \ang{\frac{c_1-c_2}{\max\br{l_1,l_2}}}^{-N}, \quad c_1,c_2\in\R,\quad l_1,l_2>0,
\end{equation}
arguing purely on the spatial side yields the following estimate:
\begin{align*}
    \abs{\ang{
         \varphi_{k_1,\vx_1,\vxi_1},
         \varphi_{k_2,\vx_2,\vxi_2}
    }}
    \underset{N}{\lesssim} &
    \sqrt{
        \frac{\br{2^{k_2}/\sqrt{\lambda}}\cdot \br{4^{k_2}/\sqrt{\lambda}}}{\br{2^{k_1}/\sqrt{\lambda}}\cdot\br{4^{k_1}/\sqrt{\lambda}}}
    }
    \ang{
        \frac{x_{01}-x_{02}}{2^{k_1}/\sqrt{\lambda}}
    }^{-N}
    \ang{
        \frac{x_{11}-x_{12}}{4^{k_1}/\sqrt{\lambda}}
    }^{-N}\\
    = &
    \sqrt{
        \frac{\abs{\vI_2}}{\abs{\vI_1}}
    }
    \ang{
        \frac{x_{01}-x_{02}}{2^{k_1}/\sqrt{\lambda}}
    }^{-N}
    \ang{
        \frac{x_{11}-x_{12}}{4^{k_1}/\sqrt{\lambda}}
    }^{-N}
    .
\end{align*}
In combination, we derive the following estimate:
\begin{align*}
    &
    \nrm{
        \ang{
             \varphi_{k_1,\vx_1,\vxi_1},
             \varphi_{k_2,\vx_2,\vxi_2}
         }
         \prod_{k=0,1}
         \chi^N_{\vI_k}\br{\vx_k}
         \1_{\vomega^{\br{+,0}}_k}\br{\vxi_k}
    }_{L^2\br{d\vx_1 d\vxi_1
    d\vx_2 d\vxi_2}}\\
    \underset{N}{\lesssim} &
    \abs{\vomega_1}^{1/2}
    \sqrt{ 
        2^{-k_2}\sqrt{\lambda}\cdot
        4^{-k_2}\sqrt{\lambda}
    }
    \cdot
    \sqrt{
        \frac{\abs{\vI_2}}{\abs{\vI_1}}
    }
    \cdot
    \nrm{
        \ang{
            \frac{x_{01}-x_{02}}{2^{k_1}/\sqrt{\lambda}}
        }^{-N}
        \ang{
            \frac{x_{11}-x_{12}}{4^{k_1}/\sqrt{\lambda}}
        }^{-N}
        \prod_{k=0,1}
        \chi_{\vI_k}^N\br{\vx_k}
    }_{L^2\br{
        d\vx_1d\vx_2
    }}.
\end{align*}
For the last factor in the above expression, we apply \eqref{eq_1dim_std_fact} twice to obtain:
\begin{align*}
    &
    \nrm{
        \ang{
            \frac{x_{01}-x_{02}}{2^{k_1}/\sqrt{\lambda}}
        }^{-N}
        \ang{
            \frac{x_{11}-x_{12}}{4^{k_1}/\sqrt{\lambda}}
        }^{-N}
        \prod_{k=0,1}
        \chi^N_{\vI_k}\br{\vx_k}
    }_{L^2\br{
        d\vx_1 d\vx_2
    }}\\
    \underset{N}{\lesssim} &
    \sqrt{\br{2^{k_1}/\sqrt{\lambda}}\cdot\br{4^{k_1}/\sqrt{\lambda}}}
    \nrm{
        \chi^N_{\vI_1}\br{\vx_2}
        \chi^N_{\vI_2}\br{\vx_2}
    }_{L^2\br{d\vx_2}}\\
    \underset{N}{\lesssim} &
    \sqrt{\br{2^{k_1}/\sqrt{\lambda}}\cdot\br{4^{k_1}/\sqrt{\lambda}}}
    \abs{\vI_2}^{1/2}
    \nrm{
        \chi^N_{\vI_1}
    }_{L^\infty\br{\vI_2}}
    .
\end{align*}
We thus deduce the following estimate:
\begin{align*}
    &
    \nrm{
        \ang{
             \varphi_{k_1,\vx_1,\vxi_1},
             \varphi_{k_2,\vx_2,\vxi_2}
         }
         \prod_{k=0,1}
         \chi^N_{\vI_k}\br{\vx_k}
         \1_{\vomega^{\br{+,0}}_k}\br{\vxi_k}
    }_{L^2\br{d\vx_1 d\vxi_1
    d\vx_2 d\vxi_2}}\\
    \underset{N}{\lesssim} &
    \abs{\vomega_1}^{1/2}
    \cdot
    \sqrt{ 
        2^{-k_2}\sqrt{\lambda}\cdot
        4^{-k_2}\sqrt{\lambda}
    }
    \cdot
    \sqrt{
        \frac{\abs{\vI_2}}{\abs{\vI_1}}
    }
    \cdot
    \sqrt{\br{2^{k_1}/\sqrt{\lambda}}\cdot\br{4^{k_1}/\sqrt{\lambda}}}
    \cdot
    \abs{\vI_2}^{1/2}
    \nrm{
        \chi^N_{\vI_1}
    }_{L^\infty\br{\vI_2}}\\
    = &
    \abs{\vomega_1}^{1/2}
    \cdot
    \abs{\vI_2}^{1/2}
    \nrm{
        \chi^N_{\vI_1}
    }_{L^\infty\br{\vI_2}}
    = 
    \sqrt{
        \frac{\lambda^2 \abs{\vI_2}}{\abs{\vI_1}}
    }
    \nrm{
        \chi^N_{\vI_1}
    }_{L^\infty\br{\vI_2}}.
\end{align*}
In other words, we have
\begin{equation*}
    \abs{\ang{
        f^{\br{+,0}}_{\vI_1\times\vomega_1},
        g^{\br{+,0}}_{\vI_2\times\vomega_2}
    }}
    \underset{N}{\lesssim}
    \nrm{f^{\br{+,0}}}_{\vI_1\times\vomega_1}
    \nrm{g^{\br{+,0}}}_{\vI_2\times\vomega_2}
    \sqrt{
        \frac{\lambda^{D+1}\abs{\vI_2}}{\abs{\vI_1}}
    }
    \nrm{
        \chi^N_{\vI_1}
    }_{L^\infty\br{\vI_2}}.
\end{equation*}
As for the other estimate, we apply H\"older, \eqref{eq_planar_tf_proj}, \eqref{eq_wp_nrm} and the fact that \(0\leq \chi^\Phi_\vI\underset{N}{\lesssim} \chi^N_\vI\) to conclude
\begin{align*}
    \abs{\ang{
        f^{\br{+,0}}_{\vI_1\times\vomega_1},
        g^{\br{+,0}}_{\vI_2\times\vomega_2}
    }}
    \leq &
    \nrm{
        \sqrt{\chi^\Phi_{\vI_1}}
        \pi^{\br{+,0}}_{\vomega_1}f
    }_{L^2}
    \nrm{
        \sqrt{\chi^\Phi_{\vI_2}}
        \pi^{\br{+,0}}_{\vomega_2}g
    }_{L^2}
    \nrm{
        \sqrt{\chi^\Phi_{\vI_1}}
        \sqrt{\chi^\Phi_{\vI_2}}
    }_{L^\infty}\\
    \underset{N}{\lesssim} &
    \nrm{f^{\br{+,0}}}_{\vI_1\times\vomega_1}
    \nrm{g^{\br{+,0}}}_{\vI_2\times\vomega_2}
    \nrm{
        \chi^N_{\vI_1}
        \chi^N_{\vI_2}
    }_{L^\infty}
    \leq
    \nrm{f^{\br{+,0}}}_{\vI_1\times\vomega_1}
    \nrm{g^{\br{+,0}}}_{\vI_2\times\vomega_2}
    \nrm{\chi^N_{\vI_1}}_{L^\infty\br{\vI_2}}.
\end{align*}

\subsection{Estimate of the quantity \texorpdfstring{\(\nrm{\Lambda}_c\)}{} subject to Proposition \ref{singtileestimdec}}\label{subsec_planar_sing_tile_est_bd_Lambda_c}
We start by recalling the concepts introduced in \textbf{Definition \ref{def_sing_tile_bili} and \ref{def_sing_tile_est}}. Our goal here is to prove the validity of \textbf{Proposition \ref{singtileestimbd}, \emph{i)}}. Indeed, once we show this, since we assume for now also the validity of \textbf{Proposition \ref{singtileestimdec}}, we would immediately deduce that
\begin{equation*}
    \nrm{\Lambda}_{1/4}\lesssim 1,\quad \nrm{\Lambda}_0\lesssim \lambda^{-\sigma}
\end{equation*}
for some \(\sigma>0\), and thus, by taking the geometric mean, we would be able to conclude that 
$$\nrm{\Lambda}_{1/8}\lesssim \lambda^{-\sigma/2}\,.$$

With these said, we first perform a reduction procedure to eliminate the dependency of a few parameters. By the symmetry of the formulation, we may assume without loss of generality that \(\vI\times\vomega\in\P_k\) with
\begin{equation*}
    k=0
    \quad\text{and}\quad
    \vI\times\vomega=
    \bR{0,1}^2\times
    \bR{-\lambda/2,\lambda/2}\times\bR{0,\lambda}\in\P_0\br{\lambda}.
\end{equation*}
Starting now the proof of \(\nrm{\Lambda}_{1/4}\lesssim 1\), we recall that
\begin{equation}\label{eq_planar_sing_tile_id}
    \Lambda^{\br{0}}_{\vI\times\vomega}\br{f,g}=
    \ang{\cC_0 f^{\br{+,0}}_{\vI\times\vomega},\1_{A^{\br{-,0}}_\vomega}g}=
    \ang{\cC_0 \br{\chi^\Phi_\vI \pi^{\br{+,0}}_\vomega f},\1_{A^{\br{-,0}}_\vomega}g}.
\end{equation}
By ignoring the oscillatory term within the expression, we obtain the following trivial bound:
\begin{equation}\label{eq_planar_sing_tile_triv}
    \abs{\Lambda^{\br{0}}_{\vI\times\vomega}\br{f,g}}
    \lesssim
    \nrm{
        \nrm{
            \chi^\Phi_\vI\br{x_0-t,x_1-t^2} \pi^{\br{+,0}}_\vomega f\br{x_0-t,x_1-t^2}
        }_{L^1\br{dt,\abs{t}\lesssim 1}}
        \1_{A^{\br{-,0}}_\vomega}\br{\vx}g\br{\vx}
    }_{L^1\br{d\vx}}.
\end{equation}
We make two key observations:
\begin{equation}\label{eq_jap_simp_approx}
    0\leq \chi^\Phi_\vI\underset{N}{\lesssim} \chi^N_\vI,\quad
    \abs{a}\lesssim 1 \eqsim \abs{b} \implies \ang{x+a}\eqsim \ang{x} \eqsim \ang{bx},\quad a,b,x\in\R.
\end{equation}
In combination, we have:
\begin{equation*}
    \chi^\Phi_\vI\br{x_0-t,x_1-t^2}\underset{N}{\lesssim}\chi_\vI^N\br{x_0-t,x_1-t^2}\underset{N}{\eqsim}\chi_\vI^N\br{\vx},\quad \abs{t}\lesssim 1.
\end{equation*}
As a result, we obtain:
\begin{equation}\label{eq_planar_sing_tile_triv_for_Lp_improv}
    \eqref{eq_planar_sing_tile_triv}\underset{N}{\lesssim}
    \nrm{
        \nrm{
            \sqrt{\chi^\Phi_\vI\br{x_0-t,x_1-t^2}}
            \pi^{\br{+,0}}_\vomega f\br{x_0-t,x_1-t^2}
        }_{L^1\br{dt,\abs{t}\lesssim 1}}
        \chi_\vI^N\br{\vx}\1_{A^{\br{-,0}}_\vomega}\br{\vx}g\br{\vx}
    }_{L^1\br{d\vx}}.
\end{equation}
We notice that the \(t\) variable evolves along the parabola and thus, we are naturally led to consider the following: 

\begin{theorem}\label{thm_Lp_improv}
Let \(D\geq 1\). For \(f\in L^{\frac{D+2}{D+1}}\br{\R^{D+1}}\), the following estimate holds:
    \begin{equation}
        \nrm{
            \int_{\R^D}
                f\br{\vx-\vX\br{\vt}}
            d\vt
        }_{L^{D+2}\br{d\vx}}
        \lesssim \nrm{f}_{L^{\frac{D+2}{D+1}}}.
    \end{equation}
    Additionally, via interpolation with trivial estimates, it is easy to see that:
    \begin{equation}\label{eq_thm_Lp_improv_p_mult}
        \nrm{
            \nrm{
                f\br{\vx-\vX\br{\vt}}
            }_{L^p\br{d\vt,\abs{\vt}\lesssim 1}}
        }_{L^q\br{d\vx}}
        \lesssim \nrm{f}_{L^{\frac{q}{D+1}}},\quad 0<\br{D+2}p\leq q \leq \infty.
    \end{equation}
\end{theorem}
With the above theorem in mind, we perform \(L^4\br{d\vx}\)--\(L^{4/3}\br{d\vx}\) H\"{o}lder's inequality:
\begin{equation*}
    \eqref{eq_planar_sing_tile_triv_for_Lp_improv}
    \lesssim
    \nrm{
        \nrm{
            \sqrt{\chi^\Phi_\vI\br{x_0-t,x_1-t^2}}
            \pi^{\br{+,0}}_\vomega f\br{x_0-t,x_1-t^2}
        }_{L^1\br{dt,\abs{t}\lesssim 1}}
    }_{L^4\br{d\vx}}
    \nrm{
        \chi_\vI^N\1_{A^{\br{-,0}}_\vomega}g
    }_{L^{4/3}}
\end{equation*}
and apply \textbf{Theorem \ref{thm_Lp_improv}} to the first factor to dominate the above with:
\begin{equation*}
    \lesssim
    \nrm{
        \sqrt{\chi^\Phi_\vI}
        \pi^{\br{+,0}}_\vomega f
    }_{L^2}
    \nrm{
        \chi_\vI^N\1_{A^{\br{-,0}}_\vomega}g
    }_{L^{4/3}}
    .
\end{equation*}
Finally, recall \eqref{eq_wp_nrm} and \textbf{Definition \ref{def_Psi_I} and \ref{def_jap_n_weight}} and apply H\"{o}lder's inequality on \(\1_{A^{\br{-,0}}_\vomega}g=\1_{A^{\br{-,0}}_\vomega\cap\supp g}\cdot \1_{A^{\br{-,0}}_\vomega}g\) to further dominate the above with
\begin{equation*}
    \lesssim
    \mu_\vI^{1/4}\br{A^{\br{-,0}}_\vomega\cap \supp g}
    \nrm{f}_{\vI\times\vomega}
    \nrm{\1_{A^{\br{-,0}}_\vomega}g}_{L^2\br{d\mu_\vI}}
    .
\end{equation*}
This confirms that \(\nrm{\Lambda}_{1/4}\lesssim 1\) thus concluding our proof.

\subsection{Control over \texorpdfstring{\(\nrm{T_\nu}_2\)}{} and \texorpdfstring{\(\nrm{M_\nu}_2\)}{}}\label{subsec_est_T_M_4_st_pt}
We recall the specific setting \(\br{\nu_k}_{k\in\Z}=\br{\mu_k}_{k\in\Z}\).
We aim to show the following two estimates:
\begin{equation*}
    \nrm{T_\nu}_2\lesssim \lambda^{-\frac{1}{2}},\quad
    \nrm{M_\nu}_2\lesssim 1.
\end{equation*}
Recalling \textbf{Definition \ref{def_sio_max}} we fix \(\cZ\subset \Z\), and, starting  with the estimate for \(T^{\cZ}_\nu\), we observe the following multiplier representation:
\begin{equation*}
    \widehat{T^{\cZ}_\nu f}\br{\vxi}=
    \sum_{k\in\cZ}
    \Psi\br{
        \begin{pmatrix}
            2 & 0\\
            0 & 4
        \end{pmatrix}^k
        \vxi/\lambda
    }
    \widehat{\mu_k}\br{\vxi}
    \widehat{f}\br{\vxi},
\end{equation*}
where the smooth bump function \(\Psi\) is defined as
\begin{equation}\label{eq_def_Psi_for_pi_k}
    \Psi\br{\vxi}:=
    \varphi^{\otimes 2}\br{\vxi/2}
    -
    \varphi^{\otimes 2}
    \br{\begin{pmatrix}
        2 & 0 \\
        0 & 4
    \end{pmatrix}^5
    \vxi/2
    }
    .
\end{equation}
Most importantly, \(\vxi\in\supp \Psi\) implies \(\abs{\vxi}\eqsim 1\). To estimate \(T^{\cZ}_\nu f\), we apply the \(L^2\) multiplier theory:
\begin{equation*}
    \nrm{T^{\cZ}_\nu}_2=
    \nrm{
        \sum_{k\in \cZ}
            \Psi\br{
                \begin{pmatrix}
                    2 & 0\\
                    0 & 4
                \end{pmatrix}^k
                \vxi/\lambda
            }
        \widehat{\mu_k}\br{\vxi}
    }_{L^\infty\br{d\vxi}}.
\end{equation*}
However, the almost disjointness of the support for the expression in the summand allows us to dominate the above with the following:
\begin{equation*}
    \lesssim
    \sup_{k\in\Z}
        \nrm{
            \Psi\br{
                \begin{pmatrix}
                    2 & 0\\
                    0 & 4
                \end{pmatrix}^k
                \vxi/\lambda
            }
            \widehat{\mu_k}\br{\vxi}
        }_{L^\infty\br{d\vxi}}
    \leq
    \sup_{k\in\Z}
        \nrm{
            \widehat{\mu_k}
            \br{
                \frac{\xi_0}{2^k}
                \frac{\xi_1}{4^k}
            }
        }_{L^\infty\br{\abs{\vxi}\eqsim \lambda}}
\end{equation*}
It remains obtain an estimate for \(\widehat{\mu_k}\br{
                \frac{\xi_0}{2^k}
                \frac{\xi_1}{4^k}
            }\). This can be achieved by expanding the formula:
\begin{equation*}
    \widehat{\mu_k}
    \br{
        \frac{\xi_0}{2^k}
        \frac{\xi_1}{4^k}
    }
    =
    \int 
    \overline{
        e\br{\xi_0 t+\xi_1 t^2}
    }
    2^kK\br{2^kt}
    \rho\br{\abs{t}}
    dt
\end{equation*}
and applying the Van der Corput lemma on oscillatory integral to obtain for \(\abs{\vxi}\eqsim \lambda\) the following:
\begin{equation*}
    \abs{
        \widehat{\mu_k}
        \br{
            \frac{\xi_0}{2^k}
            \frac{\xi_1}{4^k}
        }
    }
    \lesssim 
    \lambda^{-\frac{1}{2}}
    \br{
        \nrm{2^kK\br{2^kt}
        \rho\br{\abs{t}}}_{L^\infty\br{dt}}
        +
        \nrm{
            \frac{d}{dt}
            2^kK\br{2^kt}
            \rho\br{\abs{t}}
        }_{L^1\br{dt}}
    }
    \lesssim \lambda^{-\frac{1}{2}}.
\end{equation*}
We thus conclude \(\nrm{T^{\cZ}_\nu}_2\lesssim\lambda^{-\frac{1}{2}}\) and thus \(\nrm{T_\nu}_2\lesssim \lambda^{-\frac{1}{2}}\).
As for \(M_\nu\), we estimate directly via \eqref{eq_kernel_cond} to deduce
\begin{equation}\label{eq_M_nu_st_pw_est}
    M_\nu f\br{\vx}:=
    \sup_{k\in\Z}
    \int
        \abs{\pi_k f}\br{x_0-t,x_1-t^2}
        \rho_k\br{\abs{t}}\abs{K\br{t}}
    dt
    \lesssim
    M_\para M^{\otimes 2} f\br{\vx}.
\end{equation}
As a direct consequence,  we obtain the estimate \(\nrm{M_\nu}_2\lesssim 1\).

\subsection{The proof of Proposition \ref{thm_simp_w22}}

We start our subsection with the following

\begin{observation}[\textsf{Total \(L^2\) energy bound in the $j=0$ case}]\label{obs_total_L2_energy_bd}
Due to the absence of resonances on the first frequency variable \(\xi_0\), the energy-related reasonings are significantly simpler in the $j=0$ case relative to the \(j=1\) case. Indeed, we notice that in the former situation, the collection of intervals \(\BR{\vomega^{\br{+,0}}\::\: \vomega\in\W}\) is disjoint. This implies the orthogonality among all time-frequency projections on tiles. To be more precise, the following becomes a direct consequence of \eqref{eq_ass_sum_I}: 
\begin{equation}\label{eq_obs_total_L2_energy_bd}
        \nrm{f^{\br{+,0}}}_{L^2\br{\P}}=
        \nrm{\nrm{f^{\br{+,0}}}_{P}}_{\ell^2\br{P\in\P}}\lesssim 
        \nrm{
            \nrm{\pi^{\br{+,0}}_\vomega f}_{L^2}
        }_{\ell^2\br{\vomega\in\W}}
        \lesssim \nrm{f}_{L^2}.
    \end{equation}
Moreover, we note that in the \textbf{$j=0$ case} given any  \(P,P'\in\P\) one has \(P\leq P'\) if and only if  \(P\leq_{-,0} P'\). In other words, \textbf{all trees are \((-,0)-\)}trees.
\end{observation}
With these settled, we utilize \textbf{Lemma \ref{lem_-_tree}} and  \textbf{Lemma \ref{lem_u_shift_mass_sel}} in order to prove our weak \(\br{2,2}\) estimate:

\begin{proof}[\textbf{Proof of} \textbf{Proposition \ref{thm_simp_w22}}]
Due to the positive nature of the expression on the left-hand side, it suffices to prove the statement with \(\P\) replaced by a finite sub-collection \(\cP\subset \P:=\P\br{\lambda}\) and with the implicit constant independent of the choice of \(\cP\). We initiate the argument by decomposing \(\cP\) into a collection of trees \(\BR{\cT_{n,i}}_{n,i}\):
\begin{enumerate}
    \item Set \(\cP_0:=\cP\).
    \item Given \(\cP_{n-1}\), we apply \textbf{Lemma \ref{lem_u_shift_mass_sel}} to obtain a disjoint collection of trees \(\cT_{n,i}\subset \cP_{n-1}\) such that
    \begin{equation}
        \cM_{\cP_{n-1}\setminus\bigsqcup_i\cT_{n,i}}\br{E}\lesssim 2^{-n},\quad
        \sum_i\abs{\vI_{\cT_{n,i}}}\lesssim 2^n \abs{E}.
    \end{equation}
    \item Set \(\cP_n:=\cP_{n-1}\setminus\bigsqcup_i\cT_{n,i}\).
\end{enumerate}
We iterate until the algorithm exhausts the collection of tiles. By construction, we have:
\begin{equation}
    \cP=\bigsqcup_{n,i}\cT_{n,i},\quad
    \cM_{\cT_{n,i}}\br{E}\leq\cM_{\cP_{n-1}}\br{E}\lesssim 2^{-n}.
\end{equation}
We now apply the decomposition to rewrite the model sum:
\begin{equation}
    \nrm{\Lambda^{\br{0}}_P\br{f,\1_E}}_{\ell^1\br{P\in\cP}}=\sum_n \sum_i \nrm{\Lambda^{\br{1}}_P\br{f,\1_E}}_{\ell^1\br{P\in\cT_{n,i}}}.
\end{equation}
Recall that all trees are \(\br{-,0}\)-trees. We may apply \textbf{Lemma \ref{lem_-_tree}} to dominate the above with:
\begin{align*}
    \lesssim &
    \br{\nrm{\Lambda}_c+\nrm{T_\nu}_2}
    \cL
    \sum_n\sum_i
        \nrm{f^{\br{+,0}}}_{L^2\br{\cT_{n,i}}}
        \cM^{\frac{1}{2}+c}_{\cT_{n,i}}\br{E}\abs{I_{\cT_{n,i}}}^{\frac{1}{2}}\\
    \lesssim &
    \br{\nrm{\Lambda}_c+\nrm{T_\nu}_2}\cL
    \sum_n
        2^{-\frac{n}{2}-cn}
        \sum_i
            \nrm{f^{\br{+,0}}}_{L^2\br{\cT_{n,i}}}
            \abs{I_{\cT_{n,i}}}^{\frac{1}{2}}.
\end{align*}
Applying Cauchy-Schwarz on \(\sum_i\) and \textbf{Lemma \ref{lem_u_shift_mass_sel}}, one dominates the above by
\begin{align*}
    \lesssim &
    \br{\nrm{\Lambda}_c+\nrm{T_\nu}_2}\cL
    \sum_n
        2^{-\frac{n}{2}-cn}
        \nrm{
            \nrm{f^{\br{+,0}}}_{L^2\br{\cT_{n,i}}}
        }_{\ell^2\br{i}}
        \br{
        \sum_i
            \abs{I_{\cT_{n,i}}}
            }^{\frac{1}{2}}\\
    \lesssim &
    \br{\nrm{\Lambda}_c+\nrm{T_\nu}_2}\cL
    \sum_n
        2^{-\frac{n}{2}-cn}
        \nrm{f^{\br{+,0}}}_{L^2\br{\bigsqcup_i\cT_{n,i}}}
        2^{\frac{n}{2}}
        \abs{E}^{\frac{1}{2}}\\
    \lesssim &
    \br{\nrm{\Lambda}_c+\nrm{T_\nu}_2}\cL^{\frac{3}{2}}
    \sum_n
        2^{-cn}
        \nrm{f^{\br{+,0}}}_{L^2\br{\cP_n}}
        \abs{E}^{\frac{1}{2}}.
\end{align*}
Another Cauchy-Schwarz on \(\sum_n\) allows us to dominate the above with:
\begin{equation*}
    \leq
    \br{\nrm{\Lambda}_c+\nrm{T_\nu}_2}
    \cL^{\frac{3}{2}}
    \br{
    \sum_n
        2^{-2cn}
    }^{\frac{1}{2}}
    \nrm{f^{\br{+,0}}}_{L^2\br{\P}}
        \abs{E}^{\frac{1}{2}}
\end{equation*}
Finally, using Observation \ref{obs_total_L2_energy_bd} and in particular \eqref{eq_obs_total_L2_energy_bd}, we conclude our proof as follows:
\begin{equation*}
    \nrm{\Lambda^{\br{0}}_P\br{f,\1_E}}_{\ell^1\br{P\in\cP}}
    \lesssim
    \frac{
        \nrm{\Lambda}_c+\nrm{T_\nu}_2
    }{\sqrt{c}}\cdot\cL^{\frac{3}{2}}
        \nrm{f}_{L^2}
        \abs{E}^{\frac{1}{2}}.
\end{equation*}
\end{proof}


\section{High resolution analysis---single tile decay estimate: Proof of Proposition \ref{singtileestimdec} subject to Lemma \ref{lem_sub_lev_2_2}}\label{sec_est_of_sing_tile}

In this section our goal is to show that \(\exists\:\sigma>0\) s.t. \(\nrm{\Lambda}_0\lesssim \lambda^{-\sigma}\), or equivalently, via \textbf{Definition \ref{def_sing_tile_est}}
\begin{equation}\label{decste}
    \abs{\Lambda^{\br{j}}_{\vI\times\vomega}\br{f,g}}
    \lesssim \lambda^{-\sigma}\nrm{f}_{\vI\times\vomega}\nrm{\1_{A^{\br{-,0}}_\vomega}g}_{L^2\br{d\mu_\vI}}.
\end{equation}

\begin{remark}[\textsf{Inclusion of the one-dimensional non-resonant regime}]\label{nonrezinclus}
The proof of the above statement will be done via the LGC method and essentially encapsulates all the key ideas implemented in the treatment of the non-resonant Carleson-Radon transform treated by the same authors in \cite{HsuL24}.
\end{remark}

\subsection{Preliminaries}

We first remark here that due to scaling and translation invariance, it is enough to prove our single tile decay estimate in the case $k=0$ and $I=[0,1]$.
Next, we recall \eqref{eq_planar_tf_proj}---\eqref{eq_planar_sing_tile_id}, and, for the sake of matching dimension when performing interpolation arguments, we renormalize the wave-packet via the change of variable \(\vxi=\sqrt{\lambda}\vzeta\), that is
\begin{equation}\label{eq_planar_def_phi_xyuv_norm}
    \varphi_{\vx,\vzeta}:=
    \Tr_\vx\Dil^1_{1/\sqrt{\lambda}}
    \Mod_\vzeta\widecheck{\phi}^{\otimes 2}
    =
    \nrm{\phi}_{L^2}^2
    \sqrt{\lambda}
    \varphi_{0,\vx,\sqrt{\lambda}\vzeta}.
\end{equation}
For simplicity, we also use the approximation 
\begin{equation}\label{eq_planar_def_dmu_I_norm}
    d\mu\br{\vx}:=
    \ang{\vx}^{-N}_\otimes d\vx
    =
    \ang{x_0}^{-N}
    \ang{x_1}^{-N} d\vx
    \eqsim
    d\mu_\vI\br{\vx}.
\end{equation}
Additionally, we shall ignore the internal structure of the coefficient in \eqref{eq_planar_tf_proj} and consider:
\begin{equation}\label{eq_planar_def_Fxyuv}
    F\br{\vx,\vzeta}:=
    \ang{\pi^{\br{+,0}}_\vomega f,
    \varphi_{\vx,\vzeta}}\1_\vI\br{\vx}.
\end{equation}
As a direct consequence, \eqref{eq_planar_tf_proj} and \eqref{eq_wp_nrm} become
\begin{equation*}
    f^{\br{+,0}}_{\vI\times\vomega}=
    \nrm{\phi}_{L^2}^{-4}
    \int
        F\br{\vx,\vzeta}
        \varphi_{\vx,\vzeta}
    d\vx d\vzeta,\quad
    \nrm{f}_{\vI\times\vomega}\eqsim \nrm{F}_{L^2}.
\end{equation*}
On the other hand, we can make the harmless assumption that 
\begin{equation*}
    \supp g\subset 
    \va^{-1}\vomega^{\br{-,0}}\cap \ok^{-1}\bR{k,\infty}=
    a^{-1}\bR{0,\lambda}\cap \ok^{-1}\bR{0,\infty}.
\end{equation*}
At this stage of the reduction, the statement \(\nrm{\Lambda}_0\lesssim \lambda^{-\sigma}\) becomes equivalent to:
\begin{equation*}
    \abs{
        \int
                F\br{\vx,\vzeta}
                \ang{\cC_0 \varphi_{\vx,\vzeta},
                    g
                }
        d\vx d\vzeta
    }
    \lesssim \lambda^{-\sigma}
    \nrm{F}_{L^2}
    \nrm{g}_{L^2\br{d\mu}}.
\end{equation*}
As it turns out, in what follows, we will prove the following slightly stronger version of \textbf{Proposition \ref{singtileestimdec}}:
\begin{proposition}\label{thm_tile_est_for_C_lambda}
    There is a universal constant \(\sigma>0\) such that for \(\lambda\gg 1\), \(F\in L^2\br{\R^4}\), and \(g\in L^2\br{\R^2}\) with\footnote{Recall \bf{Definition \ref{def_diam}}.} \(\dia\br{a\br{\supp g}}\lesssim \lambda\), the following estimate holds:
    \begin{equation}\label{eq_the_tile}
        \nrm{
            F\br{\vx,\vzeta}
           \ang{\cC_0\varphi_{\vx,\vzeta},g}
        }_{L^1\br{\substack{
            \abs{\vx}\lesssim 1,\\
            \abs{\zeta_1}\lesssim \zeta_0\eqsim \sqrt{\lambda}
        }}}
        \lesssim
         \lambda^{-\sigma}
        \nrm{F}_{L^2}
        \nrm{g}_{L^2\br{d\mu}}.
    \end{equation}
\end{proposition}

Throughout the section, we make the following harmless assumption:
\begin{equation*}
    \supp F \subset 
    \BR{
        \br{\vx,\vzeta}\in\R^4\::\:
        \abs{\vx}\lesssim 1,\quad
        \abs{\zeta_1}\lesssim \zeta_0\eqsim \sqrt{\lambda}
    }.
\end{equation*}
\subsection{LGC method---Step 1: Time-frequency representation at linearizing scale}\label{subsec_tf_rep_lin}
The goal of this subsection is to derive a phase-linearized model for
\begin{align*}
     \cC_0\varphi_{\vx,\vzeta}\br{\vz}:= &
    \int
        \varphi_{\vx,\vzeta}
        \br{
            \vz-\br{t,t^2}
        }
        e\br{a\br{\vz}t^2}
        \rho\br{\abs{t}}
        K\br{t}
    dt\\
    =&
    \int
        \phi^{\otimes 2}\br{\frac{\vxi}{\sqrt{\lambda}}-\vzeta}
        \int
            \overline{e\br{\xi_0 t+\br{\xi_1-a\br{\vz}}t^2}}
            \rho\br{\abs{t}}
            K\br{t}
        dt
        \cdot
        e\br{\br{\vz-\vx}\cdot\vxi}
    d\vxi
\end{align*}
whenever \(\abs{a\br{\vz}}\lesssim \lambda\). By assumption, for \(\abs{\zeta_1}\lesssim \zeta_0\eqsim \sqrt{\lambda}\)
\begin{align*}
    \supp \phi\br{\frac{\cdot}{\sqrt{\lambda}}-\zeta_0}\subset \BR{\xi_0\in\R\::\:\xi_0\eqsim \lambda},&\quad \abs{\supp \phi\br{\frac{\cdot}{\sqrt{\lambda}}-\zeta_0}}\eqsim \sqrt{\lambda}\\
    \supp \phi\br{\frac{\cdot}{\sqrt{\lambda}}-\zeta_1}\subset \BR{\xi_1\in\R\::\:\abs{\xi_1}\lesssim \lambda},&\quad \abs{\supp \phi\br{\frac{\cdot}{\sqrt{\lambda}}-\zeta_1}}\eqsim \sqrt{\lambda}.
\end{align*}
This is exactly the correct scale for \(\vxi\) to perform phase linearization. As for \(a\br{\vz}\), we shall introduce a partition of unity: \(1=\int\phi\br{\frac{a}{\sqrt{\lambda}}-w}\frac{dw}{\nrm{\phi}_{L^1}}\). The formula of \(\cC_0\varphi_{\vx,\vzeta}\br{\vz}\) now reads:
\begin{align*}
    & \cC_0\varphi_{\vx,\vzeta}\br{\vz}
    =
    \int
            \phi\br{\frac{a\br{\vz}}{\sqrt{\lambda}}-w}
    \int
        \phi^{\otimes 2}\br{\frac{\vxi}{\sqrt{\lambda}}-\vzeta}
        \\
        & \phantom{\int}\cdot
        \int
            \overline{e\br{\xi_0 t+\br{\xi_1-a\br{\vz}}t^2}}
            \rho\br{\abs{t}}
            K\br{t}
        dt
        \cdot
        e\br{\br{\vz-\vx}\cdot\vxi}
    \frac{d\vxi
    dw}{\nrm{\phi}_{L^1}}.
\end{align*}
Fix now \(\vxi,a=a\br{\vz}\) and \(\vzeta,w\) as in the above expression. We now address the oscillatory integral within the above expression. Observe that:
\begin{align*}
    &
    \int
        \overline{e\br{\xi_0 t+\br{\xi_1-a}t^2}}
        \rho\br{\abs{t}}
        K\br{t}
    dt\\
    =&
    \int
        \overline{e\br{\xi_0 \br{t+\frac{\tau}{\sqrt{\lambda}}}+\br{\xi_1-a}\br{t+\frac{\tau}{\sqrt{\lambda}}}^2}}
        \rho\br{\abs{t+\frac{\tau}{\sqrt{\lambda}}}}
        K\br{t+\frac{\tau}{\sqrt{\lambda}}}
        \phi\br{\tau}
    \frac{d\tau
    dt}{\nrm{\phi}_{L^1}}\\
    =&
    \int
    \overline{e\br{
            \xi_0 t+\br{\xi_1-a}t^2
        }}
            \int
            \overline{
                e\br{
                    \br{
                       \zeta_0+ 2t\br{\zeta_1-w}
                    }
                    \tau
                }
            }
            \mathcal{E}\br{
                \tau,\vxi,a,t,\vzeta,w
            }
            \rho\br{\abs{t+\frac{\tau}{\sqrt{\lambda}}}}
            K\br{t+\frac{\tau}{\sqrt{\lambda}}}
            \phi\br{\tau}
        \frac{d\tau
        dt}{\nrm{\phi}_{L^1}},
\end{align*}
where \(\mathcal{E}(\cdot)\), obtained via a Taylor decomposition, stands for the expression:
\begin{align*}
    \mathcal{E}
    \br{\cdot}
    = &
    e\br{
        -
        \br{
            \frac{\xi_0}{\sqrt{\lambda}}-\zeta_0
        }\tau
    }
    e\br{
        -2\br{
            \frac{\xi_1}{\sqrt{\lambda}}-\zeta_1
        }t\cdot\tau
    }
    e\br{
        2\br{
            \frac{a}{\sqrt{\lambda}}-w
        }t\cdot\tau
    }\nonumber\\
    \cdot &
    e\br{
        -\lambda^{-\frac{1}{2}}
            \br{
                \frac{\xi_1}{\sqrt{\lambda}}-\zeta_1
            }
        \cdot\tau^2
    }
    e\br{
        -\frac{
            \zeta_1
        }{
            \sqrt{\lambda}
        }
        \cdot\tau^2
    }
    e\br{
        \lambda^{-\frac{1}{2}}\br{\frac{a}{\sqrt{\lambda}}-w}\cdot\tau^2
    }
    e\br{
        \frac{w}{\sqrt{\lambda}}\cdot\tau^2
    }.
\end{align*}
On the one hand, the support of the function \(\phi\) imposes the following condition:
    \begin{equation*}
        \abs{\tau},\;
        \abs{
            \frac{\xi_0}{\sqrt{\lambda}}-\zeta_0
        },\;
        \abs{
            \frac{\xi_1}{\sqrt{\lambda}}-\zeta_1
        },\;
        \abs{
            \frac{a}{\sqrt{\lambda}}-w
        }\leq 0.01.\nonumber
    \end{equation*}
    On the other hand, since \(\abs{\xi_1},\abs{a}\lesssim \lambda \eqsim \xi_0\) and \(\abs{t+\frac{\tau}{\sqrt{\lambda}}}\eqsim 1\), we further deduce:
    \begin{equation*}
        \abs{t},\,
        \frac{\zeta_0}{\sqrt{\lambda}},\,
        \abs{
            \frac{\zeta_1}{\sqrt{\lambda}}
        },\,
        \abs{
            \frac{w}{\sqrt{\lambda}}
        }
        \lesssim 1 \lesssim
        \abs{t},\,\frac{\zeta_0}{\sqrt{\lambda}}.
    \end{equation*}
    Using Taylor expansions for all seven exponential terms appearing in $\mathcal{E}
        \br{\cdot}$ and making the notations $\vec{k}=(k_1,\ldots, k_7)$ and \(\vec{k}!:=\prod_{j=1}^7\br{k_l!}\) we deduce
    \begin{align*}
        \mathcal{E}
        \br{\cdot}
        &=
        \sum_{\vec{k}}
            \frac{
                \br{-1}^{k_1+k_2+k_4+k_5}
                2^{k_2+k_3}
                \lambda^{-\frac{k_4+k_6}{2}}
                \br{2\pi i}^{\sum_{l=1}^7 k_l}
            }{\vec{k}!}
            t^{
                k_2+k_3
            }
            \br{
                \frac{\zeta_1}{\sqrt{\lambda}}
            }^{
                k_5
            } \nonumber\\
        \cdot &
            \br{
                \frac{w}{\sqrt{\lambda}}
            }^{
                k_7
            }
            \br{
            \frac{\xi_0}{\sqrt{\lambda}}-\zeta_0}^{k_1}
        \br{
            \frac{\xi_1}{\sqrt{\lambda}}-\zeta_1
        }^{k_2+k_4}
        \br{
            \frac{a}{\sqrt{\lambda}}-w
        }^{k_3+k_6}
        \tau^{k_1+k_2+k_3+2k_4+2k_5+2k_6+2k_7}.
    \end{align*}
  With these, setting
    \begin{align*}
        C_{\vec{k}}:= &
        \br{-1}^{k_1+k_2+k_4+k_5}
        2^{k_2+k_3}
        \lambda^{-\frac{k_4+k_6}{2}}
        \br{2\pi i}^{\sum_{l=1}^7 k_l}
        \br{\vec{k}!}^{-1/5}
        \\
        D_{\vec{k}}\br{t,\zeta_1,w}:= &
        t^{
            k_2+k_3
        }
        \br{
            \frac{\zeta_1}{\sqrt{\lambda}}
        }^{
            k_5
        }
        \br{
            \frac{w}{\sqrt{\lambda}}
        }^{
            k_7}
        \br{\vec{k}!}^{-1/5}\\
        \phi^{\br{1,\vec{k}}}\br{z}:=&\br{\vec{k}!}^{-1/5}z^{k_1}\phi\br{z},\quad
        \phi^{\br{2,\vec{k}}}\br{z}:=\br{\vec{k}!}^{-1/5}z^{k_2+k_4}\phi\br{z}\\
        \phi^{\br{3,\vec{k}}}\br{z}:=&z^{k_3+k_6}\phi\br{z},\quad
        \phi^{\br{4,\vec{k}}}\br{z}:=\br{\vec{k}!}^{-1/5}z^{k_1+k_2+k_3+2k_4+2k_5+2k_6+2k_7}\phi\br{z}
    \end{align*}
    we deduce
    \begin{align*}
        \cC_0\varphi_{\vx,\vzeta}\br{\vz}
        = &
        \sum_{\vec{k}}
        C_{\vec{k}}
        \nrm{\phi}_{L^1}^{-2}
        \iint
        D_{\vec{k}}\br{t,\zeta_1,w}
        e\br{
            a\br{\vz}t^2
        }
        \phi^{\br{3,\vec{k}}}\br{\frac{a\br{\vz}}{\sqrt{\lambda}}-w}
        \\
        & \cdot
        \int
            \overline{e\br{
                \br{
                    \zeta_0+2t\br{\zeta_1-w}
                }
                \tau
            }
            }
            \rho\br{\abs{t+\frac{\tau}{\sqrt{\lambda}}}}
            K\br{t+\frac{\tau}{\sqrt{\lambda}}}
            \phi^{\br{4,\vec{k}}}\br{\tau}
        d\tau\nonumber\\
        \cdot &
        \underbrace{
        \iint
            \phi^{\br{1,\vec{k}}}
            \otimes
            \phi^{\br{2,\vec{k}}}
            \br{\frac{\vxi}{\sqrt{\lambda}}-\vzeta}
            e\br{
                \br{\vz-\vx-\br{t,t^2}}\cdot\vxi
            }
        d\vxi
        }_{
            =:\varphi^{\br{\vec{k}}}_{\vx+\br{t,t^2},\vzeta}\br{\vz}
        }
        dtdw.
    \end{align*}
As a direct consequence, we obtain the following approximation:
\begin{equation*}
    \abs{\ang{\cC_0\varphi_{\vx,\vzeta},g}}\lesssim 
    \E_{\vec{k}}
    \nrm{
    \frac{
        \ang{
            e\br{
                at^2
            }
            \phi^{\br{3,\vec{k}}}\br{\frac{a}{\sqrt{\lambda}}-w}
            g,
            \varphi^{\br{\vec{k}}}_{\vx+\br{t,t^2},\vzeta}
        }
    }{
    \ang{\zeta_0+2t\br{\zeta_1-w}}^A
    }}_{L^1\br{\abs{t}\eqsim 1, \abs{w}\lesssim \sqrt{\lambda}}}.
\end{equation*}
Note the important fact that \(\abs{\phi^{\br{i,\vec{k}}}}\lesssim \1_{\mr{-0.01,0.01}}\) for \(i=1,2,3,4\) and that \(\abs{\br{\phi^{\br{i,\vec{k}}}}^{\vee}}\underset{N}{\lesssim}\ang{\cdot}^{-10N}\) for \(i=1,2,4\) with implicit constants only depending on \(N\). By setting 
\begin{equation}
    G^{\br{\vec{k}}}\br{\vx,t,\vzeta,w}:=\abs{\ang{
        e\br{
            at^2
        }
        g_w^{\br{\vec{k}}},
        \varphi^{\br{\vec{k}}}_{\vx,\vzeta}
    }},\quad
    g_w^{\br{\vec{k}}}:=
    \phi^{\br{3,\vec{k}}}
    \br{\frac{a}{\sqrt{\lambda}}-w}
    g,
\end{equation}
we may estimate \eqref{eq_the_tile} in the following manner:
\begin{align*}
    \nrm{
        F\br{\vx,\vzeta}
        \ang{\cC_0\varphi_{\vx,\vzeta},g}
    }&_{L^1\br{\substack{ \abs{\vx}\lesssim 1,\\ \abs{\zeta_1}\lesssim \zeta_0\eqsim \sqrt{\lambda}}}}
    \lesssim 
    \E_{\vec{k}} 
    \nrm{
    \frac{
        F\br{\vx-\br{t,t^2},\vzeta}
        G^{\br{\vec{k}}}\br{\vx,t,\vzeta,w}
    }{
    \ang{\zeta_0+2t\br{\zeta_1-w}}^A
    }
    }_{L^1\br{\substack{\abs{\vx}\lesssim \abs{t}\eqsim 1,\\ \abs{\zeta_1},\abs{w}\lesssim \zeta_0\eqsim \sqrt{\lambda}}}}\\
    \lesssim 
    \sup_{\vec{k}} &
    \nrm{
    \frac{
        F\br{\vx-\br{t,t^2},\vzeta}
        G^{\br{\vec{k}}}\br{\vx,t,\vzeta,w}
    }{
    \ang{\zeta_0+2t\br{\zeta_1-w}}^A
    }
    }_{L^1\br{\substack{\abs{\vx}\lesssim \abs{t}\eqsim 1,\\ \abs{\zeta_1},\abs{w}\lesssim \zeta_0\eqsim \sqrt{\lambda}}}}.
\end{align*}
We can thus take \(\vec{k}\) that extremizes the expression, drop the superscript \(\br{\cdot}^{\br{\vec{k}}}\), and assume without loss of generality that \(t>0\). It remains to show:
\begin{equation}\label{eq_the_tile_mod}
    \overbrace{
        \nrm{
        \frac{
            F\br{\vx-\br{t,t^2},\vzeta}
            G\br{\vx,t,\vzeta,w}
        }{
        \ang{\zeta_0+2t\br{\zeta_1-w}}^A
        }
        }_{L^1\br{\substack{\abs{\vx}\lesssim t\eqsim 1,\\ \abs{\zeta_1},\abs{w}\lesssim \zeta_0\eqsim \sqrt{\lambda}}}}
    }^{=:\Lambda\br{F,G}}
    \hspace{-1ex}
    \lesssim
    \lambda^{-\sigma}
    \nrm{F}_{L^2}
    \nrm{g}_{L^2\br{d\mu}}.
\end{equation}


\subsection{Proof of Proposition \ref{thm_tile_est_for_C_lambda}: a blueprint}\label{subsec_3_endpoints}

\subsubsection{LGC method---Step 2: A sparse-uniform dichotomy analysis}\label{dichotomy}
$\newline$

We start with some very brief remarks: 
\begin{itemize}
\item in order to properly characterize the spatial localization of the wave packets, it is helpful to introduce the auxiliary measure
\begin{equation}\label{eq_planar_aux_meas}
    \mu_\vx\br{S}:=
    \int
        \1_S\br{\vx+\frac{\vz}{\sqrt{\lambda}}}\,
    d\mu\br{\vz}\,,
\end{equation}
which is a finite measure with the natural bounds
\begin{equation}
    \mu_\vx \br{S}\leq \mu_\vx \br{\R^2}\eqsim 1,\quad
    \nrm{\mu_\vx \br{S}}_{L^1\br{\abs{\vx}\lesssim 1}}\lesssim \mu\br{S}.
\end{equation}
\item notice that
\begin{equation}
    \ang{f,\varphi_{\vx,\vzeta}}=
    \int
        f\br{\vx-\frac{\vz}{\sqrt{\lambda}}}
        \widehat{\phi}^{\otimes 2}\br{\vz}
        \overline{e\br{\vz\cdot \vzeta}}\,
    d\vz\,,
\end{equation}
and thus, by the Hausdorff-Young's inequality, we deduce that for any $ p\in \mr{1,2}$ the following hold:
\begin{equation}\label{cor_G_bds}
        \nrm{\ang{f,\varphi_{\vx,\vzeta}}}_{L^{p'}\br{d\vzeta }}
        \lesssim\nrm{f}_{L^p\br{d\mu_\vx }}\qquad\textrm{and}\qquad    \nrm{
            G\br{\vx,t,\vzeta,w}
        }_{L^{p'}\br{d\vzeta }}
        \lesssim
        \nrm{g_w}_{L^p\br{d\mu_\vx }}\,.
\end{equation}
\end{itemize}

We are now ready to comment on our strategy. In order to show \eqref{eq_the_tile_mod}, we will apply an interpolation argument among the following regimes:
\begin{itemize}
\item the \textsf{sparse-uniform} regime described by
 \begin{equation}\label{eq_S_U}
  \Lambda\br{F,G}\lesssim 
    \nrm{
        \nrm{
            F\br{\vx,\vzeta}
        }_{L^1\br{d\vzeta }}
    }_{L^2\br{d\vx }}
    \nrm{
        \nrm{\nrm{g_w}_{L^1\br{d\mu_\vx }}}_{L^\infty\br{dw}}
    }_{L^2\br{\abs{\vx}\lesssim 1}}\,.
\end{equation}

\item the \textsf{sparse-sparse} regime described by
\begin{align}\label{eq_S_S}
\Lambda\br{F,G}\lesssim 
        &\nrm{
            \ang{
                \zeta_0\br{\vx}+
                2t\br{\zeta_1\br{\vx}-w\br{\vx+\br{t,t^2}}}
            }^{-A}
        }_{L^4\br{\abs{\vx}\lesssim t\eqsim 1}}\nonumber\\
        & \cdot 
        \nrm{
            \nrm{
                F\br{\vx,\vzeta}
            }_{L^1\br{d\vzeta }}
        }_{L^2\br{d\vx }}
        \nrm{
            \nrm{\nrm{g_w}_{L^1\br{d\mu_\vx }}}_{L^1\br{dw}}
        }_{L^2\br{\abs{\vx}\lesssim 1}}
\end{align}
with \(\zeta_0\br{\cdot},\zeta_1\br{\cdot},w\br{\cdot}\) some measurable functions satisfying \(\abs{\zeta_1\br{\cdot}},\abs{w\br{\cdot}}\lesssim \zeta_0\br{\cdot}\eqsim \sqrt{\lambda}\).

\item the \textsf{uniform-uniform} regime described by

 \begin{equation}\label{eq_U_S}
\Lambda\br{F,G} \lesssim  \log^{\frac{1}{2}}\br{e+\lambda}
    \nrm{
        \nrm{
            F\br{\vx,\vzeta}
        }_{L^4\br{d\vzeta }}
    }_{L^2\br{d\vx }}
    \nrm{
        \nrm{
            \mu_\vx ^{\frac{1}{4}}\br{\supp g_w}
            \nrm{g_w}_{L^2\br{d\mu_\vx }}
        }_{L^1\br{dw}}
    }_{L^2\br{\abs{\vx}\lesssim 1}}\,.
 \end{equation}   
\end{itemize}

The proofs of all the above estimates will be presented in Subsection \ref{proofspunifdich}.

\subsubsection{LGC method---Step 3: Control over the time-frequency correlation set}

The final ingredient is given by the following key decay estimate necessary for obtaining a good control over \eqref{eq_S_S}:

\begin{lemma}\label{lem_sub_lev_2_2}[\textsf{Sub-level set estimate}]

There is a universal constant \(\delta\in\br{0,1}\) such that for any \(\lambda \gtrsim 1\), \(N\geq 1\), and any three measurable functions \(u,v,w\) with \( \abs{u}\gtrsim \lambda\), the following estimate holds:
    \begin{equation}\label{levset}
        \nrm{
            \ang{u\br{x,y}+t\br{v\br{x,y}+w\br{x+t,y+t^2}}}^{-N}
        }_{L^1\br{\abs{x},\abs{y},\abs{t}\lesssim 1}}
        \lesssim \lambda^{-\delta}.
    \end{equation}
\end{lemma}
The proof of this lemma is left to \textsc{Section \ref{sec_sub_lev_est}}---there we will prove in fact a more general form of \eqref{levset}. 

\subsubsection{Putting all the pieces together: an interpolation argument among \texorpdfstring{\eqref{eq_S_U}}{}, \texorpdfstring{\eqref{eq_S_S}}{}, and \texorpdfstring{\eqref{eq_U_S}}{}} \label{subsec_interpol}

In direct correspondence with the sparse-uniform dichotomy in Subsection \ref{dichotomy}, we introduce the following sets:
\begin{align*}
    \cF_j:= &
    \BR{
        \br{\vx,\vzeta}\in \R^4\::\:
        2^j
        \leq 
        \frac{\abs{F\br{\vx,\vzeta}}}{\nrm{F\br{\vx,\vzeta'}}_{L^2\br{d\vzeta'}}}<
        2^{j+1}
    }\\
    \cG_k:= &
    \BR{
        \br{\vx,w}\in \R^3\::\:
        2^k \leq
        \frac{
            \mu_\vx ^{\frac{1}{4}}
            \br{\supp g_w}\nrm{g_w}_{L^2\br{d\mu_\vx }}
        }{
            \mu_\vx ^{\frac{1}{4}}
            \br{\supp g}\nrm{g}_{L^2\br{d\mu_\vx }}
        }
        <2^{k+1}
    }\\
    \cG^l:= &
    \BR{
        \vx\in \R^2\::\:
        2^l \leq
        \frac{
            \nrm{g}_{L^2\br{d\mu_\vx }}
        }{
            \nrm{g}_{L^2\br{d\mu}}
        }
        <2^{l+1}
        ,\quad
        \abs{\vx}\lesssim 1
    }\\
    \cF_{j,\vx}:= &
    \BR{
        \vzeta\in\R^2\::\:
        \br{\vx,\vzeta}\in\cF_j
    },\quad
    \cG_{k,\vx}:= 
    \BR{
        w\in\R\::\:
        \br{\vx,w}\in\cG_k
    }
\end{align*}
for \(j,k,l\in \Z\) and decompose the functions \(F,G\) accordingly:
\begin{equation*}
    F_j\br{\vx,\vzeta}:=
    F\br{\vx,\vzeta}\cdot\1_{\cF_j}\br{\vx,\vzeta},\quad
    G^l_k\br{\vx,t,\vzeta,w}:=
    G\br{\vx,t,\vzeta,w}\cdot\1_{\cG_k\cap \br{\cG^l\times\R}}\br{\vx,w}.
\end{equation*}
The triangle inequality produces the following:
\begin{equation*}
    \Lambda\br{F,G}
    \leq
    \sum_{j,k,l\in \Z}
    \Lambda\br{F_j,G^l_k},
\end{equation*}
and thus it suffices to derive suitable estimate on \(\Lambda\br{F_j,G^l_k}\).

 We start with the observation that:
\begin{equation*}
    \nrm{f}_{L^1\br{d\mu_\vx }}\lesssim 
    \mu_\vx ^{\frac{1}{2}}\br{\supp f}
    \nrm{f}_{L^2\br{d\mu_\vx }}
    \lesssim
    \mu_\vx ^{\frac{1}{4}}\br{\supp f}
    \nrm{f}_{L^2\br{d\mu_\vx }}
    \lesssim 
    \nrm{f}_{L^2\br{d\mu_\vx }},\quad f\in L^1_{\mathrm{loc}}\br{\R^2}.
\end{equation*}
Using the above together with the three estimate \eqref{eq_S_U}, \eqref{eq_U_S}, \eqref{eq_S_S} and \textbf{Lemma \ref{lem_sub_lev_2_2}} we derive:
\begin{equation*}
    \Lambda\br{F_j,G^l_k}
    \lesssim 
    2^{j+k+l}
    \nrm{
        F
    }_{L^2}
    \nrm{g}_{L^2\br{d\mu}}
    \cdot
    \left\{
    \begin{aligned}
        &
        \sup_{\vx\in\R^2}\abs{\cF_{j,\vx}}
        \cdot
        \abs{\cG^l}^{\frac{1}{2}}, & \br{i}\\
        &
        \lambda^{-\frac{\delta}{8}}
        \sup_{\vx\in\R^2}\abs{\cF_{j,\vx}}
        \sup_{\vx\in\cG^l}\abs{\cG_{k,\vx}}
        \cdot
        \abs{\cG^l}^{\frac{1}{2}}, &\br{i\!i}\\
        &
        \log^{\frac{1}{2}}\br{e+\lambda}
        \sup_{\vx\in\R^2}
        \abs{\cF_{j,\vx}}^{\frac{1}{4}}
        \sup_{\vx\in\cG^l}
        \abs{\cG_{k,\vx}}
        \cdot
        \abs{\cG^l}^{\frac{1}{2}}. &\br{i\!i\!i}
    \end{aligned}
    \right.
\end{equation*}
We then estimated via Chebyshev's inequality the sizes of the measurable sets involved:
\begin{equation*}
    \abs{\cF_{j,\vx}}\lesssim 2^{-2j}\land \lambda,\quad
    \abs{\cG_{k,\vx}}\lesssim 
    2^{-\frac{4k}{3}}\land \sqrt{\lambda},\quad
    \abs{\cG^l}\lesssim 
    2^{-2l}\land 1
\end{equation*}
for almost all \(\vx\in \R^2\) with \(\abs{\vx}\lesssim 1\). For all \(\epsilon>0\), we claim there exists \(0<\alpha,\beta,\gamma<1\) satisfying
\begin{equation}\label{eq_interpol_index_rel}
    \begin{pmatrix}
        1 &  1 & \frac{1}{4} \\
        0 &  1 & 1 \\
        1 & 1 & 1
    \end{pmatrix}
    \cdot
    \begin{pmatrix}
        \alpha  \\
        \beta \\
        \gamma
    \end{pmatrix}
    =
    \begin{pmatrix}
        \frac{1}{2}+\epsilon  \\
        \frac{3}{4} \\
        1
    \end{pmatrix}
\end{equation}
such that the following holds:
\begin{align*}
    2^{j+k+l}
    \br{i}^\alpha
    \br{i\!i}^\beta
    \br{i\!i\!i}^\gamma
    = &
    \lambda^{-\frac{\delta\beta}{8}}\log^{\frac{\gamma}{2}}\br{e+\lambda}
    \cdot
    2^j
    \sup_{\vx\in\R^2}\abs{\cF_{j,\vx}}^{\frac{1}{2}+\epsilon}
    \cdot
    2^k
    \sup_{\vx\in\cG^l}\abs{\cG_{k,\vx}}^{\frac{1}{2}}
    \cdot
    2^l
    \abs{\cG^l}^{\frac{1}{2}}
    \\
    \lesssim &
    \lambda^{\epsilon-\frac{\delta\beta}{8}}\log^{\frac{\gamma}{2}}\br{e+\lambda}
    \br{
        \br{2^j\sqrt{\lambda}}^{-2\epsilon}\land
        2^j\sqrt{\lambda}
    }
    \br{
        1\land 2^k\lambda^{\frac{3}{8}}
    }
    \br{
        1\land 2^l
    }.
\end{align*}
Pairing with the direct consequence:
\begin{equation*}
    \abs{\cG_{k,\vx}}\neq 0\implies 2^k\lesssim 1,\quad
    \abs{\cG^l}\neq 0\implies 2^l\lesssim \lambda^{\frac{1}{2}}
\end{equation*}
of the following fact:
\begin{align*}
    \mu_\vx ^{\frac{1}{4}}\br{\supp g_w}\nrm{g_w}_{L^2\br{d\mu_\vx }}
    \lesssim 
    \mu_\vx ^{\frac{1}{4}}\br{\supp g}\nrm{g}_{L^2\br{d\mu_\vx }}
    \lesssim 
    \nrm{g}_{L^2\br{d\mu_\vx }}
    \lesssim
    \lambda^{\frac{1}{2}}
    \nrm{g}_{L^2\br{d\mu}},\quad
    \abs{\vx}\lesssim 1,
\end{align*}
we conclude the following:
\begin{align*}
    \Lambda\br{F,G}\leq
    \sum_{j,k,l\in\Z}\Lambda\br{F_j,G^l_k}
    \lesssim &
    \lambda^{\epsilon-\frac{\delta\beta}{8}}\log^{\frac{\gamma}{2}}\br{e+\lambda}
    \nrm{F}_{L^2}\nrm{g}_{L^2\br{d\mu}}\\
    \cdot &
    \sum_{j\in\Z}
    \br{
        \br{2^j\sqrt{\lambda}}^{-2\epsilon}\land
        2^j\sqrt{\lambda}
    }
    \sum_{\substack{
        k\in\Z\\
        2^k\lesssim 1
    }}
    \br{
        1\land 2^k\lambda^{\frac{3}{8}}
    }
    \sum_{\substack{
        l\in\Z\\
        2^l\lesssim \lambda^{\frac{1}{2}}
    }}
    \br{
        1\land 2^l
    }\\
    \underset{\epsilon}{\lesssim} &
    \lambda^{\epsilon-\frac{\delta\beta}{8}}\log^{2+\frac{\gamma}{2}}\br{e+\lambda}
    \nrm{F}_{L^2}\nrm{g}_{L^2\br{d\mu}}.
\end{align*}
This proves \eqref{eq_the_tile_mod} and thus, \textbf{Proposition \ref{thm_tile_est_for_C_lambda}} as long as \(\epsilon\in\br{0,\frac{\delta\beta}{8}}\).
To complete the interpolation argument, it remains to justify \eqref{eq_interpol_index_rel} with a suitable choice of \(\epsilon>0\).
Taking \(\epsilon=\frac{\delta\beta}{16}\), \eqref{eq_interpol_index_rel} becomes:
\begin{equation*}
    \begin{pmatrix}
        1 &  1-\frac{\delta}{16} & \frac{1}{4} \\
        0 &  1 & 1 \\
        1 & 1 & 1
    \end{pmatrix}
    \cdot
    \begin{pmatrix}
        \alpha  \\
        \beta \\
        \gamma
    \end{pmatrix}
    =
    \begin{pmatrix}
        1/2 \\
        3/4 \\
        1
    \end{pmatrix}
    \implies
    \begin{pmatrix}
        \alpha  \\
        \beta \\
        \gamma
    \end{pmatrix}
    =
    \begin{pmatrix}
        1/4 \\
        \frac{1}{12-\delta} \\
        \frac{3}{4}-\frac{1}{12-\delta}
    \end{pmatrix}
    .
\end{equation*}
This demonstrates that \(0<\alpha,\beta,\gamma<1\) whenever \(\delta\in\br{0,1}\).

\subsection{The sparse-uniform dichotomy}\label{proofspunifdich}

In this section, we present the proof of the estimates stated in Subsection \ref{dichotomy}.

\subsubsection{Some notation and basic preparatives}

Following \eqref{eq_jap_simp_approx}, we introduce a lemma that allows us to perform similar arguments to the Weyl differentiation in the absence of exponentiation:
\begin{lemma}\label{lem_jap_mul_diff}  For any $a,\,b,\,x,\,y\in\R$ we have
    \begin{equation}
        \ang{x}^{-1}\ang{y}^{-1}\leq \ang{\frac{ax+by}{\abs{a}+\abs{b}}}^{-1}\br{\ang{x}^{-1}+\ang{y}^{-1}}.
    \end{equation}
\end{lemma}
As a consequence, we have the following estimate. 
\begin{lemma}\label{lem_mu_tempering}
    Given \(w\) measurable and \(\omega\) a weight, we have the following estimate:
    \begin{equation}
        \nrm{
            \nrm{
                \ang{
                    w\br{t}+\zeta_1 t+\zeta_0
                }^{-2}
            }_{L^1\br{dt,\abs{t}\lesssim 1}}
        }_{
            L^2\br{
                \omega\br{\zeta_1}d\vzeta 
            }
        }
        \lesssim 
        \sqrt{
            \nrm{\omega}_{L^\infty}
            \log\br{e+\frac{\nrm{\omega}_{L^1}}{\nrm{\omega}_{L^\infty}}}
        }.
    \end{equation}
\end{lemma}

\begin{proof}
    \begin{align*}
        &
        \nrm{
            \nrm{
                \ang{
                    w\br{t}+\zeta_1 t+\zeta_0
                }^{-2}
            }_{L^1\br{dt,\abs{t}\lesssim 1}}
        }_{
            L^2\br{
                \omega\br{\zeta_1}d\vzeta 
            }
        }^2
        \\
        =&
        \iint
            \iint_{\abs{t_0},\abs{t_1}\lesssim 1}
                \abs{\omega\br{\zeta_1}}
                \prod_{j=0,1}
                \ang{
                    w\br{t_j}+ \zeta_1 t_j+\zeta_0 
                }^{-2}
            dt_0 dt_1
        du d\zeta_1\\
        \lesssim &
        \iint_{\abs{t_0},\abs{t_1}\lesssim 1}
            \int
                \int
                    \ang{
                        w\br{t_0} +\zeta_1 t_0+\zeta_0
                    }^{-2}
                du
        \cdot 
        \abs{\omega\br{\zeta_1}}
        \ang{
             \left.w\right\vert_{t_0}^{t_1}+\zeta_1\br{t_1-t_0}
        }^{-2}
        d\zeta_1
        dt_0 dt_1\\
        \lesssim & 
        \iint_{\abs{t_0},\abs{t_1}\lesssim 1}
            \nrm{\omega}_{L^1} \wedge \frac{\nrm{\omega}_{L^\infty}}{\abs{t_1-t_0}}
        dt_0 dt_1
        \lesssim 
        \int_{\abs{\Delta t}\lesssim 1} 
            \nrm{\omega}_{L^1} \wedge \frac{\nrm{\omega}_{L^\infty}}{\abs{t_\Delta}}
        dt_\Delta\\
        \eqsim & 
        \int_0^{\frac{\nrm{\omega}_{L^\infty}}{\nrm{\omega}_{L^1}}\wedge 1}
            \nrm{\omega}_{L^1} 
        dt_\Delta+
        \int_{\frac{\nrm{\omega}_{L^\infty}}{\nrm{\omega}_{L^1}}\wedge 1}^1
            \frac{\nrm{\omega}_{L^\infty}}{\abs{t_\Delta}}
        dt_\Delta
        \lesssim
        \nrm{\omega}_{L^\infty}
        \log\br{e+\frac{\nrm{\omega}_{L^1}}{\nrm{\omega}_{L^\infty}}}.
    \end{align*}
\end{proof}

\subsubsection{Treatment of the sparse-uniform regime: proof of \texorpdfstring{\eqref{eq_S_U}}{}}\label{subsec_S_U}
We start with an \(L^1\)--\(L^\infty\) H\"{o}lder's inequality on \(d\vzeta dw\).
\begin{align*}
    \Lambda\br{F,G}\leq & 
    \nrm{
        \nrm{
            \frac{F\br{\vx-\br{t,t^2},\vzeta}}{\ang{\zeta_0+2t\br{\zeta_1-w}}^A}
        }_{L^1\br{d\vzeta dw}}\hspace{-6ex}\cdot
        \nrm{
            G\br{\vx,t,\vzeta,w}
        }_{L^\infty\br{d\vzeta dw}}
    }_{L^1\br{\abs{\vx}\lesssim t \eqsim 1}}\\
    \lesssim &
    \nrm{
        \nrm{
            F\br{\vx-\br{t,t^2},\vzeta}
        }_{L^1\br{d\vzeta }}
        \nrm{\nrm{g_w}_{L^1\br{d\mu_\vx }}}_{L^\infty\br{dw}}
    }_{L^1\br{\abs{\vx}\lesssim t \eqsim 1}}.
\end{align*}
We then apply Cauchy-Schwarz on \(d\vx \) to finish the proof of \eqref{eq_S_U}:
\begin{align*}
    \leq &
    \nrm{
        \nrm{
            \nrm{
                F\br{\vx-\br{t,t^2},\vzeta}
            }_{L^1 \br{ \substack{ d\vzeta }}}
        }_{L^2\br{d\vx }}
    }_{L^1\br{t\eqsim 1}}
    \nrm{
        \nrm{
            \nrm{g_w}_{L^1\br{d\mu_\vx }}}_{L^\infty\br{dw}
        }
    }_{L^2\br{\abs{\vx}\eqsim 1}}\\
    \lesssim &
    \nrm{
        \nrm{
            F\br{\vx,\vzeta}
        }_{L^1 \br{ \substack{ d\vzeta }}}
    }_{L^2\br{d\vx }}
    \nrm{
        \nrm{
            \nrm{g_w}_{L^1\br{d\mu_\vx }}}_{L^\infty\br{dw}
        }
    }_{L^2\br{\abs{\vx}\eqsim 1}}.
\end{align*}

\subsubsection{Treatment of the sparse-sparse regime: proof of \texorpdfstring{\eqref{eq_S_S}}{}}\label{subsec_S_S}
In this case, we start with an \(L^1\)--\(L^\infty\) H\"{o}lder's inequality on \(d\vzeta \) and juggle the \(\nrm{\cdot}_{L^1\br{dt}}\) to the inner layer.
\begin{align*}
    \Lambda\br{F,G}\leq & 
    \nrm{
        \nrm{
            \frac{F\br{\vx-\br{t,t^2},\vzeta}}{\ang{\zeta_0+2t\br{\zeta_1-w}}^A}
        }_{L^1\br{\abs{\zeta_1}\lesssim \zeta_0\eqsim \sqrt{\lambda}}}
        \nrm{
            G\br{\vx,t,\vzeta,w}
        }_{L^\infty\br{d\vzeta }}
    }_{L^1\br{\substack{\abs{w}\lesssim \sqrt{\lambda}\\ \abs{\vx}\lesssim t\eqsim 1}}}\\
    \lesssim &
    \nrm{
        \nrm{
            \frac{F\br{\vx-\br{t,t^2},\vzeta}}{\ang{\zeta_0+2t\br{\zeta_1-w}}^A}
        }_{L^1\br{\substack{\abs{\zeta_1}\lesssim \zeta_0\eqsim \sqrt{\lambda}\\ t\eqsim 1}}}
        \nrm{g_w}_{L^1\br{d\mu_\vx }}
    }_{L^1\br{\substack{\abs{w}\lesssim \sqrt{\lambda}\\ \abs{\vx}\lesssim 1}}}
\end{align*}
Taking \(L^\infty\)--\(L^1\) H\"{o}lder's inequality on \(dw\) and linearizing the \(\nrm{\cdot}_{L^\infty\br{dw}}\) allow us to dominate the above with
\begin{align*}
    \lesssim &
    \nrm{
        \nrm{
            \frac{F\br{\vx-\br{t,t^2},\vzeta}}{\ang{\zeta_0+2t\br{\zeta_1-w\br{\vx}}}^A}
        }_{L^1\br{\substack{\abs{\zeta_1}\lesssim \zeta_0\eqsim \sqrt{\lambda}\\ t\eqsim 1}}}
        \nrm{\nrm{g_w}_{L^1\br{d\mu_\vx }}}_{L^1\br{dw}}
    }_{L^1\br{\abs{\vx}\lesssim 1}}\\
    = &
    \nrm{
        \frac{
            F\br{\vx-\br{t,t^2},\vzeta}
            \nrm{\nrm{g_w}_{L^1\br{d\mu_\vx }}}_{L^1\br{dw}}
        }{\ang{\zeta_0+2t\br{\zeta_1-w\br{\vx}}}^A}
    }_{L^1\br{\substack{
        \abs{\zeta_1}\lesssim \zeta_0 \eqsim \sqrt{\lambda}\\
        \abs{\vx}\lesssim t \eqsim 1
    }}}\\
    = &
    \nrm{
        \nrm{
            F\br{\vx,\vzeta}
            \nrm{
                \frac{
                    \nrm{\nrm{g_w}_{L^1\br{d\mu_{\vx+\br{t,t^2}}}}}_{L^1\br{dw}}
                }{\ang{\zeta_0+2t\br{\zeta_1-w\br{\vx+\br{t,t^2}}}}^A}
            }_{L^1\br{t\eqsim 1}}
        }_{L^1\br{\abs{\zeta_1}\lesssim \zeta_0\eqsim \sqrt{\lambda}}}
    }_{L^1\br{\abs{\vx}\lesssim 1}}.
\end{align*}
We again take \(L^\infty\)--\(L^1\) H\"{o}lder's inequality on \(d\vzeta \) and linearize the \(\nrm{\cdot}_{L^\infty\br{d\vzeta }}\) to dominate the above with
\begin{align*}
    \leq &
    \nrm{
        \frac{
            \nrm{
                F\br{\vx,\vzeta}
            }_{L^1\br{d\vzeta }}
            \nrm{\nrm{g_w}_{L^1\br{d\mu_{\vx+\br{t,t^2}}}}}_{L^1\br{dw}}
        }{\ang{\zeta_0\br{\vx}+2t\br{\zeta_1\br{\vx}-w\br{\vx+\br{t,t^2}}}}^A}
    }_{L^1\br{\abs{\vx}\lesssim t\eqsim 1}}.
\end{align*}
We apply \(L^4\)--\(L^{\frac{4}{3}}\) H\"{o}lder's inequality on \(d\vx dt\) to dominate the above with
\begin{align*}
    \leq & 
    \nrm{
        \ang{
            \zeta_0\br{\vx}+
            2t\br{\zeta_1\br{\vx}-w\br{\vx+\br{t,t^2}}}
        }^{-A}
    }_{L^4\br{\abs{\vx}\lesssim t\eqsim 1}}
    \numberthis\label{eq_SS_cor}\\
    \cdot &
    \nrm{
        \nrm{
            F\br{\vx,\vzeta}
        }_{L^1\br{d\vzeta }}
        \nrm{\nrm{g_w}_{L^1\br{d\mu_{\vx+\br{t,t^2}}}}}_{L^1\br{dw}}
    }_{L^{\frac{4}{3}}\br{\abs{\vx}\lesssim t\eqsim 1}}.\numberthis\label{eq_SS_func_imp}
\end{align*}
This produces the correct factor \eqref{eq_SS_cor} in \eqref{eq_S_S}. For \eqref{eq_SS_func_imp}, we juggle the \(\nrm{\cdot}_{L^{\frac{4}{3}}\br{dt}}\) to the inner layer, apply \(L^2\)--\(L^4\) H\"{o}lder's inequality on \(d\vx \), and apply \(L^p\) improving \textbf{Theorem \ref{thm_Lp_improv}} to finish the proof:
\begin{align}\label{eq_SS_last}
    \eqref{eq_SS_func_imp}= &
    \nrm{
        \nrm{
            F\br{\vx,\vzeta}
        }_{L^1\br{d\vzeta }}
        \nrm{
            \nrm{
                \nrm{g_w}_{L^1\br{d\mu_{\vx+\br{t,t^2}}}}
            }_{L^1\br{dw}}
        }_{L^{\frac{4}{3}}\br{t\eqsim 1}}
    }_{L^{\frac{4}{3}}\br{\abs{\vx}\lesssim 1}}\nonumber
    \\
    \leq &
    \nrm{
        \nrm{
            F\br{\vx,\vzeta}
        }_{L^1\br{d\vzeta }}
    }_{L^2\br{d\vx }}
    \nrm{
        \nrm{
            \nrm{
                \nrm{g_w}_{L^1\br{d\mu_{\vx+\br{t,t^2}}}}
            }_{L^1\br{dw}}
        }_{L^{\frac{4}{3}}\br{t\eqsim 1}}
    }_{L^4\br{\abs{\vx}\lesssim 1}}
    \\
    \lesssim &
    \nrm{
        \nrm{
            F\br{\vx,\vzeta}
        }_{L^1\br{d\vzeta }}
    }_{L^2\br{d\vx }}
    \nrm{
        \nrm{
            \nrm{g_w}_{L^1\br{d\mu_\vx }}
        }_{L^1\br{dw}}
    }_{L^2\br{\abs{\vx}\lesssim 1}}.\nonumber
\end{align}

\subsubsection{Treatment of the uniform-uniform regime: proof of \texorpdfstring{\eqref{eq_U_S}}{}}\label{subsec_U_S}
We distribute the time-frequency correlation factor and apply the Cauchy-Schwarz inequality on \(d\vzeta \).
\begin{equation*}
    \Lambda\br{F,G}\leq 
    \nrm{
        \nrm{
            \frac{F\br{\vx-\br{t,t^2},\vzeta}}{
                \ang{\zeta_0+2t\br{\zeta_1-w}}^{A/2}
            }
        }_{L^2\br{\abs{\zeta_1}\lesssim \zeta_0\eqsim \sqrt{\lambda}}}
        \nrm{
            \frac{
                G\br{\vx,t,\vzeta,w}
            }{
                \ang{\zeta_0+2t\br{\zeta_1-w}}^{A/2}
            }
        }_{L^2\br{\abs{\zeta_1}\lesssim \zeta_0\eqsim \sqrt{\lambda}}}
    }_{L^1\br{\substack{dw\\ \abs{\vx}\lesssim t\eqsim 1}}}.
\end{equation*}
Focus on the second factor in \(\nrm{\cdot}_{L^1\br{d\vx dtdw}}\) for the moment. We note that due to the frequency localization, we may dominate the below with the smooth frequency truncation variant:
\begin{equation*}
    \nrm{
        \frac{
            G\br{\vx,t,\vzeta,w}
        }{
            \ang{\zeta_0+2t\br{\zeta_1-w}}^{A/2}
        }
    }_{L^2\br{\abs{\zeta_1}\lesssim \zeta_0\eqsim \sqrt{\lambda}}}
    \lesssim 
    \nrm{
        \frac{
            G\br{\vx,t,\vzeta,w}
        }{
            \ang{\zeta_0+2t\br{\zeta_1-w}}^{A/2}
        }
    }_{L^2\br{\phi\br{\lambda^{-\frac{1}{2}}\zeta_1}d\vzeta }}.
\end{equation*}
We now compute the doubling:
\begin{align*}
    &
    \nrm{
        \frac{
            G\br{\vx,t,\vzeta,w}
        }{
            \ang{\zeta_0+2t\br{\zeta_1-w}}^{A/2}
        }
    }_{L^2\br{\phi\br{\lambda^{-\frac{1}{2}}\zeta_1}d\vzeta }}^2\\
    = &
    \iint
        \iint
                g_w\br{\vx-\lambda^{-\frac{1}{2}}\vx_0}
                \overline{g_w\br{\vx-\lambda^{-\frac{1}{2}}\vx_1}}
                \widehat{\phi}^{\otimes 2}\br{\vx_0}\overline{\widehat{\phi}^{\otimes 2}\br{\vx_1}}
            \\
    &\phantom{\int\iint}
        \cdot e\br{a\br{\vx-\lambda^{-\frac{1}{2}}\vx_0}t^2-a\br{\vx-\lambda^{-\frac{1}{2}}\vx_1}t^2}\\
    &\phantom{\int\iint}
        \cdot
        \br{
            \iint
                \frac{e\br{\vzeta\cdot\br{\vx_1-\vx_0}}}{\ang{\zeta_0+2t\br{\zeta_1-w}}^A}
                \cdot \phi\br{\lambda^{-\frac{1}{2}}\zeta_1}
            d\vzeta 
        }
    d\vx_0 d\vx_1\\
    \leq &
    \iint
        \iint
            \abs{
                \prod_{j=0,1}
                    \frac{
                        g_w\br{\vx-\lambda^{-\frac{1}{2}}\vx_j}
                    }{
                        \ang{\vx_j}^{2N}_\otimes
                    }
            }
            \cdot
            \abs{
                \widehat{\ang{\cdot}^{-A}}\br{x_{00}-x_{01}}
            }
            \\
    &\phantom{\int\iint}
        \cdot
        \sqrt{\lambda}
        \abs{
            \widehat{\phi}
            \br{
                \sqrt{\lambda}
                \br{-2t,1}\cdot\br{\vx_0-\vx_1}
            }
        }
    d\vx_0 d\vx_1\\
    \lesssim &
    \nrm{
        \frac{
            \sqrt{\lambda}
            \prod_{j=0,1}
                g_w\br{\vx-\lambda^{-\frac{1}{2}}\vx_j}
        }
        {
            \ang{
                \sqrt{\lambda}
                \br{-2t,1}\cdot\br{\vx_0-\vx_1}
            }^N
        }
    }_{L^1
    \br{
        \frac{d\vx_0 d\vx_1}{
            \ang{\vx_0}^N_\otimes
            \ang{\vx_1}^N_\otimes
        }
    }}
    .\numberthis\label{eq_UU_doub}
\end{align*}
Substituting back to the original inequality, we have:
\begin{align*}
    \Lambda\br{F,G}
    \lesssim &
    \nrm{
        \nrm{
            \frac{F\br{\vx-\br{t,t^2},\vzeta}}{
                \ang{\zeta_0+2t\br{\zeta_1-w}}^{A/2}
            }
        }_{L^2\br{\abs{\zeta_1}\lesssim \zeta_0\eqsim \sqrt{\lambda}}}
        \cdot
        \sqrt{\eqref{eq_UU_doub}}
    }_{L^1\br{\substack{dw\\ \abs{\vx}\lesssim t\eqsim 1}}}\\
    \leq &
    \nrm{
        \nrm{
            \frac{F\br{\vx-\br{t,t^2},\vzeta}}{
                \ang{\zeta_0+2t\br{\zeta_1-w}}^{A/2}
            }
        }_{L^2\br{\substack{\abs{\zeta_1}\lesssim \zeta_0 \eqsim \sqrt{\lambda}\\ t\eqsim 1}}}
        \nrm{\sqrt{\eqref{eq_UU_doub}}}_{L^2\br{t\eqsim 1}}
    }_{L^1\br{\substack{dw\\ \abs{\vx}\lesssim 1}}}.
\end{align*}
Again, we focus on the doubling of the second factor and apply Cauchy-Schwarz inequality on the weighted measure \(\frac{d\vx_0 d\vx_1}{\ang{\vx_0}^N_\otimes\ang{\vx_1}^N_\otimes}\). This produces the following chain of inequalities:
\begin{align*}\label{eq_UU_doub_sqrt}
    &
    \nrm{\sqrt{\eqref{eq_UU_doub}}}_{L^2\br{t\eqsim 1}}
    =\nrm{\eqref{eq_UU_doub}}_{L^1\br{t\eqsim 1}}^{1/2},\quad \\
    =&
    \nrm{
        \prod_{j=0,1}
            g_w\br{\vx-\lambda^{-\frac{1}{2}}\vx_j}
        \cdot
        \int_{\abs{t}\lesssim 1}
            \frac{
                \sqrt{\lambda}\1_{\supp g_w}\br{\vx-\lambda^{-\frac{1}{2}}\vx_1}dt
            }{
                \ang{
                    \sqrt{\lambda}
                    \br{-2t,1}\cdot
                    \br{\vx_0-\vx_1}
                }^N
            }
    }_{L^1
    \br{
        \frac{d\vx_0 d\vx_1}{\ang{\vx_0}^N_\otimes\ang{\vx_1}^N_\otimes}
    }}^{1/2}\\
    \leq &
    \nrm{g_w}_{L^2\br{d\mu_\vx }}
    \nrm{\int_{\abs{t}\lesssim 1}
            \frac{
                \sqrt{\lambda}\1_{\supp g_w}\br{\vx-\lambda^{-\frac{1}{2}}\vx_1}dt
            }{
                \ang{
                    \sqrt{\lambda}
                    \br{-2t,1}\cdot \Delta \vx
                }^N
            }
    }_{L^2
    \br{
        \frac{d\Delta \vx d\vx_1}{\ang{\Delta x_0}^N\ang{\vx_1}^N_\otimes}
    }}^{1/2}\\
    = &
    \mu_\vx ^{\frac{1}{4}}\br{\supp g_w}\nrm{g_w}_{L^2\br{d\mu_\vx }}
    \nrm{
        \int_{\abs{t}\lesssim 1}
            \frac{
                \sqrt{\lambda}dt
            }{
                \ang{
                    \sqrt{\lambda}
                    \br{-2t,1}\cdot \Delta \vx
                }^N
            }
    }_{L^2
    \br{
        \frac{d\Delta \vx}{\ang{\Delta x_0}^N}
    }}^{1/2}.\numberthis
\end{align*}
Notice that by \textbf{Lemma \ref{lem_mu_tempering}}, we obtain the following estimate:
\begin{align*}
    &
    \nrm{
        \int_{\abs{t}\lesssim 1}
            \frac{
                \sqrt{\lambda}dt
            }{
                \ang{
                    \sqrt{\lambda}
                    \br{-2t,1}\cdot \Delta \vx
                }^N
            }
    }_{L^2
    \br{
        \frac{d\Delta \vx}{\ang{\Delta x_0}^N}
    }}^{1/2}
    \lesssim
    \nrm{
        \nrm{
            \ang{
                t\zeta_1+\zeta_0
            }^{-N}
        }_{L^1\br{\abs{t}\lesssim 1}}
    }_{L^2
    \br{
        \ang{\lambda^{-\frac{1}{2}}\zeta_1}^{-N}d\vzeta 
    }}^{1/2}\\
    \lesssim &
    \nrm{\ang{\lambda^{-\frac{1}{2}}\zeta_1}^{-N}}_{L^\infty\br{d\zeta_1}}^{1/4}
    \log^{1/4}\br{e+\frac{\nrm{\ang{\lambda^{-\frac{1}{2}}\zeta_1}^{-N}}_{L^1\br{d\zeta_1}}}{\nrm{\ang{\lambda^{-\frac{1}{2}}\zeta_1}^{-N}}_{L^\infty\br{d\zeta_1}}}} \lesssim \log^{\frac{1}{4}}\br{e+\lambda}.
\end{align*}
As a brief summary, we have obtained the following estimate:
\begin{equation*}
    \Lambda\br{F,G}
    \lesssim \log^{\frac{1}{4}}\br{e+\lambda}
    \nrm{
        \nrm{
            \frac{F\br{\vx-\br{t,t^2},\vzeta}}{
                \ang{\zeta_0+2t\br{\zeta_1-w}}^{A/2}
            }
        }_{L^2\br{\substack{\abs{\zeta_1}\lesssim \zeta_0 \eqsim \sqrt{\lambda}\\ t\eqsim 1}}}\cdot
        \mu_\vx ^{\frac{1}{4}}\br{\supp g_w}\nrm{g_w}_{L^2\br{d\mu_\vx }}
    }_{L^1\br{\substack{dw\\ \abs{\vx}\lesssim 1}}}.
\end{equation*}
We apply \(L^\infty\)--\(L^1\) H\"{o}lder's inequality on \(dw\), linearize the \(\nrm{\cdot}_{L^\infty\br{dw}}\), and thus dominate the above with
\begin{equation*}
    \lesssim
    \log^{\frac{1}{4}}\br{e+\lambda}
    \nrm{
        \nrm{
            \frac{F\br{\vx-\br{t,t^2},\vzeta}}{
                \ang{\zeta_0+2t\br{\zeta_1-w\br{\vx}}}^{A/2}
            }
        }_{L^2\br{\substack{\abs{\zeta_1}\lesssim \zeta_0 \eqsim \sqrt{\lambda}\\ t\eqsim 1}}}
        \cdot
        \nrm{
            \mu_\vx ^{\frac{1}{4}}\br{\supp g_w}
            \nrm{g_w}_{L^2\br{d\mu_\vx }}
        }_{L^1\br{dw}}
    }_{L^1\br{\abs{\vx}\lesssim 1}}.
\end{equation*}
To extract the time-frequency correlation, we apply the Cauchy-Schwarz inequality on \(d\vx \)
\begin{align*}\label{eq_UU_1st_fact}
    \leq & \log^{\frac{1}{4}}\br{e+\lambda}
    \nrm{
        \frac{F\br{\vx-\br{t,t^2},\vzeta}}{
            \ang{\zeta_0+2t\br{\zeta_1-w\br{\vx}}}^{A/2}
        }
    }_{L^2\br{\substack{\abs{\zeta_1}\lesssim \zeta_0 \eqsim \sqrt{\lambda}\\ \abs{\vx}\lesssim t\eqsim 1}}}
    \cdot
    \nrm{
        \nrm{
            \mu_\vx ^{\frac{1}{4}}\br{\supp g_w}
            \nrm{g_w}_{L^2\br{d\mu_\vx }}
        }_{L^1\br{dw}}
    }_{L^2\br{\abs{\vx}\lesssim 1}}\\
    = &
    \nrm{
        F\br{\vx,\vzeta}
        \nrm{\ang{\zeta_0+2t\br{\zeta_1-w\br{\vx+\br{t,t^2}}}}^{-A/2}}_{L^2\br{t\eqsim 1}}
    }_{L^2\br{\substack{\abs{\zeta_1}\lesssim \zeta_0 \eqsim \sqrt{\lambda}\\ d\vx }}}\numberthis\\
    &\cdot  \log^{\frac{1}{4}}\br{e+\lambda}
    \nrm{
        \nrm{
            \mu_\vx ^{\frac{1}{4}}\br{\supp g_w}
            \nrm{g_w}_{L^2\br{d\mu_\vx }}
        }_{L^1\br{dw}}
    }_{L^2\br{\abs{\vx}\lesssim 1}}.
\end{align*}
It remains to estimate \eqref{eq_UU_1st_fact}. We apply \(L^4\)--\(L^4\) H\"{o}lder's inequality on \(d\vzeta \) and then \(L^2\)--\(L^\infty\) H\"{o}lder's inequality on \(d\vx \).
\begin{align*}\label{eq_UU_2nd_cor}
   \eqref{eq_UU_1st_fact} \lesssim &
    \nrm{
        \nrm{
            \nrm{\ang{\zeta_0+2t\zeta_1-2tw\br{\vx+\br{t,t^2}}}^{-A/2}}_{L^2\br{\abs{t}\lesssim 1}}
        }_{L^4\br{\phi\br{\lambda^{-\frac{1}{2}}\zeta_1}d\vzeta }}
    }_{L^\infty\br{d\vx }}\numberthis\\
        &\cdot
    \nrm{
        \nrm{
            F\br{\vx,\vzeta}
        }_{L^4\br{d\vzeta }}
    }_{L^2
        \br{
            d\vx 
        }
    }.
\end{align*}
We can thus take \(\vx\) to be the almost maximizer, set \(w\br{t}:=-2tw\br{\vx+\br{t,t^2}}\), and derive
\begin{equation}\label{eq_UU_2nd_cor_simp}
    \eqref{eq_UU_2nd_cor}
    \eqsim
    \nrm{
        \nrm{\ang{\zeta_0+2t\zeta_1+w\br{t}}^{-A/2}}_{L^2\br{\abs{t}\lesssim 1}}
    }_{L^4\br{\phi\br{\lambda^{-\frac{1}{2}}\zeta_1}d\vzeta }}.
\end{equation}
Finally, to estimate \eqref{eq_UU_2nd_cor_simp}, we apply \textbf{Lemma \ref{lem_mu_tempering}} to obtain:
\begin{align*}
    \eqref{eq_UU_2nd_cor_simp}
    = &
    \nrm{
        \nrm{\ang{\zeta_0+2t\zeta_1+w\br{t}}^{-A}}_{L^1\br{\abs{t}\lesssim 1}}
    }_{L^2\br{\phi\br{\lambda^{-\frac{1}{2}}\zeta_1}d\vzeta }}^{1/2}\\
    \lesssim &
    \nrm{\phi\br{\lambda^{-\frac{1}{2}}\zeta_1}}^{1/4}_{L^\infty\br{d\zeta_1}}\log^{1/4}\br{e+\frac{\nrm{\phi\br{\lambda^{-\frac{1}{2}}\zeta_1}}_{L^1\br{d\zeta_1}}}{\nrm{\phi\br{\lambda^{-\frac{1}{2}}\zeta_1}}_{L^\infty\br{d\zeta_1}}}}
    \lesssim \log^{\frac{1}{4}}\br{e+\lambda}.
\end{align*}
With all things combined, we obtain \eqref{eq_U_S}.

\section{Proof of \textbf{Theorem \ref{thm_phy_shift_disc}}}\label{sec_pf_physift_disc}

In this section, our goal is to present the proof of  \textbf{Theorem \ref{thm_phy_shift_disc}} subject to the key lemmas stated in \textsc{Section \ref{sec_formal_tf_ana}}. Recalling \eqref{eq_2_identity_2_abs}, \eqref{eq_time_decomp} and \textbf{Definition \ref{def_cC_nu}} we proceed with the formal setting introduced in \textsc{Section \ref{sec_formal_tf_ana}} and throughout this section 
$$\textrm{we set}\:\:\: j\in\BR{0,1},\:\:\vu\in \Z^2,\:\:\lambda=1,\:\:\textrm{and}\:\:\br{\nu_k}_{k\in\Z}=\br{\alpha_k}_{k\in\Z}\,.$$

In what follows, we present \textsf{the main steps of our program} in achieving the above-stated goal:
\begin{itemize}
    \item Construct a well-behaved Gabor decomposition verifying \eqref{eq_time_decomp};
    \item Verify \textbf{Assumption \ref{ass_tf_proj_emb}};
    \item Prove \(\nrm{\Lambda}_{\frac{1}{2}}:=\nrm{\Lambda}_{\frac{1}{2},1,j}, \nrm{T_\nu}_2,\nrm{M_\nu}_2\lesssim 1\);
    \item Prove \(\cE_2\lesssim 1\) and \(\cE_\ast\lesssim 1\);
    \item Prove \textbf{Proposition \ref{thm_endpoint_infty_1}} and \textbf{Proposition \ref{thm_endpoint_1_infty}} subject to \textbf{Lemma \ref{lem_gen_tree_c_half}} and \textbf{Lemma \ref{lem_u_shift_mass_sel}}.
\end{itemize}

Assuming for now the completion of the above program, we can immediately derive \textbf{Theorem \ref{thm_HF_st_disc}}: 
Indeed, applying \textbf{Proposition \ref{thm_endpoint_infty_1}} and \textbf{Proposition \ref{thm_endpoint_1_infty}} we obtain that for all \(F,E\subset \R^2\) with finite measure there exists a majorant subset \(E'\subset E\) with \(\abs{E'}\geq \abs{E}/2\) such that
\begin{align*}
    \abs{
        \ang{
            \cC\cR^{\br{\va}}_{\Delta,j,\vu}\1_F,
            \1_{E'}
        }
    }
    = &
    \abs{
        \ang{
            \cC\cR^{\br{\va}}_{\nu,1,j}\1_F,
            \1_{E'}
        }
    }
    \leq
    \nrm{\Lambda^{\br{j}}_P\br{\1_F,\1_{E'}}}_{\ell^1\br{P\in\P}}\\
    \lesssim &
    \br{
        \nrm{\Lambda}_{\frac{1}{2}}+\nrm{T_\nu}_2
    }
    \cE_\ast
    \cL^2
    \log\br{
        e+\frac{
            \cE^2_2\max\br{\abs{F},\abs{E}}
        }{
            \cE^2_\ast\min\br{\abs{F},\abs{E}}
        }
    }
    \min\br{\abs{F},\abs{E}}.
\end{align*}
Via direct calculations, cleaning up the various constants, we further dominate the above by
\begin{equation*}
    \lesssim 
    \log^2\br{e+\abs{\vu}}
    \log\br{
        e+\frac{
            \max\br{\abs{F},\abs{E}}
        }{
            \min\br{\abs{F},\abs{E}}
        }
    }
    \min\br{\abs{F},\abs{E}}
    \underset{p}{\lesssim}
    \log^2\br{e+\abs{\vu}}
    \abs{F}^{\frac{1}{p}}\abs{E}^{\frac{p-1}{p}}.
\end{equation*}
The latter, together with a real interpolation argument, completes the proof of \textbf{Theorem \ref{thm_phy_shift_disc}}.

It thus remains to tackle the five steps in our above-mentioned program.

\subsection{Gabor decompositions \texorpdfstring{\eqref{eq_time_decomp}}{} adapted to \texorpdfstring{\(\cC\cR^{\br{\va}}_{\Delta,j,\vu}\)}{}: short-time Fourier series} 
Given \(\vomega=\br{0,w}\boxplus 
        \begin{pmatrix}
            2 & 0\\
            0 & 4
        \end{pmatrix}^{-k}
        \hspace{-2ex}
        \cdot
        \bR{-1/2,1/2}\times\bR{0,1}\in\W_k:=\W_k\br{1}\),
we recall the definition of \(\pi^{\br{+,j}}_\vomega f\) and perform a Fourier series decomposition on the frequency side. Note that we may choose \(\phi:=\tphi^2\) with \(0\leq \tphi\in C^\infty_c\br{\R}\) at the very start.
We thus have the following identity:
\begin{equation}\label{eq_Gabor_frame_starter}
    \widehat{\pi^{\br{+,j}}_\vomega f}
    =
    \tphi_{\vomega^{\br{+,j}}}
    \widehat{f}
    \cdot
    \tphi_{\vomega^{\br{+,j}}}.
\end{equation}
we then equate \(\tphi_{\vomega^{\br{+,j}}}
    \widehat{f}\) with its Fourier series expansion on \(\vomega\):
\begin{equation*}
    \sum_{\valpha\in\Z^2}
        \int
            \tphi_{\vomega^{\br{+,j}}}
            \br{\vxi'}
            \widehat{\pi^{\br{+,j}}_\vomega f}\br{\vxi'}
            e\br{
                2^k \alpha_0\xi'_0+
                4^k \alpha_1\xi'_1
            }
        2^{3k}d\vxi'_0
        \cdot
        \overline{
            e\br{
                2^k \alpha_0\xi_0+
                4^k \alpha_1\xi_1
            }
        }.
\end{equation*}
This produces a wave-packet representation of \(\pi^{\br{+,j}}_\vomega f\) and
suggests that we set for \(\vI=\valpha \boxplus
    \begin{pmatrix}
        2 & 0\\
        0 & 4
    \end{pmatrix}^k
    \hspace{-1ex}\cdot
        \bR{0,1}^2
        \in\I_k
\)
the following wave-packet function:
\begin{equation*}
    \widehat{
    \varphi^{\br{+,j}}_{\vI\times\vomega}
    }
    :=
    2^{3k/2}
    \Mod_{
        -\br{2^k\alpha_0,4^k\alpha_1}
    }
        \tphi_{\vomega^{\br{+,j}}}
\end{equation*}
to simplify the representing formula:
\begin{equation*}
    \pi^{\br{+,j}}_\vomega f=
    \sum_{\vI\in\I_k}
    f^{\br{+,j}}_{\vI\times\vomega},\quad
    f^{\br{+,j}}_{\vI\times\vomega}
    :=
    \ang{
        f,
        \varphi^{\br{+,j}}_{\vI\times\vomega}
    }
    \varphi^{\br{+,j}}_{\vI\times\vomega},\quad
    \vomega\in\W_k.
\end{equation*}
By setting \(\nrm{f^{\br{+,j}}}_{\vI\times\vomega}:=
\abs{\ang{
        f,
        \varphi^{\br{+,j}}_{\vI\times\vomega}
    }}\), it's routine computation to verify \textbf{Assumption \ref{ass_tf_proj_emb}}.

\subsection{Estimating \texorpdfstring{\(\nrm{\Lambda}_{\frac{1}{2}}\)}{}, \texorpdfstring{\(\nrm{T_\nu}_2\)}{}, and \texorpdfstring{\(\nrm{M_\nu}_2\)}{}}
Recall that we are currently working under the assumption that \(\nu_k=\alpha_k\) with the defining formula:
\begin{equation*}
    d\alpha_k\br{\vx}:=
    \epsilon_k\delta\br{x_0 -2^k u_1, x_1-4^k u_2} d\vx.
\end{equation*}
Starting with the estimate on \(\nrm{\Lambda}_{\frac{1}{2}}\),
we take \(\vI\times\vomega\in\P\), expand \textbf{Definition \ref{def_sing_tile_bili}}:
\begin{equation*}
    \abs{\Lambda^{\br{j}}_{\vI\times\vomega}\br{f,g}}
    =
    \nrm{f^{\br{+,j}}}_{\vI\times\vomega}
    \cdot
    \abs{
        \ang{
            \1_{A^{\br{-,j}}_\vomega}g,
            \varphi^{\br{+,j}}_{
                \br{\vu\boxplus\vI}\times
                \vomega
            }
        }
    },
\end{equation*}
and dominate the above with:
\begin{equation*}
    \lesssim 
    \nrm{f^{\br{+,j}}}_{\vI\times\vomega}
    \nrm{\1_{A^{\br{-,j}}_\vomega}g}_{L^1\br{d\mu_{\vu\boxplus\vI}}}
    \abs{\vI}^{\frac{1}{2}}
    \lesssim
    \mu_{\vu\boxplus\vI}^{\frac{1}{2}}\br{A^{\br{-,j}}_\vomega\cap \supp g}
    \nrm{f^{\br{+,j}}}_{\vI\times\vomega}
    \nrm{\1_{A^{\br{-,j}}_\vomega}g}_{L^2\br{d\mu_{\vu\boxplus\vI}}}
    \abs{\vI}^{\frac{1}{2}}.
\end{equation*}
This proves \(\nrm{\Lambda}_{\frac{1}{2}}\lesssim 1\).
To estimate \(\nrm{T_\nu}_2\), we recall \textbf{Definition \ref{def_sio_max}}, fix \(\cZ\subset \Z\), and expand the definition of \(T^{\cZ}_\nu\):
\begin{equation*}
    T^{\cZ}_\nu f\br{\vx}=
    \sum_{k\in \cZ} 
        \epsilon_k 
        \pi_k f\br{x_0-2^k u_1,x_1-4^k u_2}
    .
\end{equation*}
Further appealing to definition \eqref{eq_LP_freq_proj} and \eqref{eq_def_Psi_for_pi_k} we deduce the following multiplier formulation:
\begin{equation*}
    \widehat{T^{\cZ}_\nu f}\br{\vxi}
    =
    \sum_{k\in\cZ}
        \epsilon_k
        \overline{
            e\br{
                \vu^\top
                \begin{pmatrix}
                    2 & 0\\
                    0 & 4
                \end{pmatrix}^k
                \vxi
            }
        }
        \Psi\br{
            \begin{pmatrix}
                2 & 0\\
                0 & 4
            \end{pmatrix}^k
            \vxi
        }
    \widehat{f}\br{\vxi}.
\end{equation*}
Due to the almost disjointness of the Littlewood-Paley pieces on the frequency side, we trivially have
\begin{equation*}
    \nrm{
    \sum_{k\in \cZ}
        \epsilon_k
        \overline{
            e\br{
                \vu^\top
                \begin{pmatrix}
                    2 & 0\\
                    0 & 4
                \end{pmatrix}^k
                \vxi
            }
        }
        \Psi\br{
            \begin{pmatrix}
                2 & 0\\
                0 & 4
            \end{pmatrix}^k
            \vxi
        }
    }_{L^\infty\br{d\vxi}}
    \lesssim 1.
\end{equation*}
As a direct consequence, we deduce \(\nrm{T^{\cZ}_\nu f}_{L^2}\lesssim \nrm{f}_{L^2}\) and thus \(\nrm{T_\nu}_2\lesssim 1\).
To estimate \(\nrm{M_\nu}_2\), we start with a pointwise domination via a suitably defined square function. Namely, we consider the following:
\begin{equation*}
    M_\nu f\br{\vx}
    \lesssim
    \nrm{
        \pi_k f\br{x_0-2^k u_1,x_1-4^k u_2}
    }_{\ell^\infty\br{k\in\Z}}
    \leq 
    \nrm{
        \pi_k f\br{x_0-2^k u_1,x_1-4^k u_2}
    }_{\ell^2\br{k\in\Z}}
    =:
    S_\nu f\br{\vx}.
\end{equation*}
By orthogonality among the Littlewood-Paley projection \(\pi_k f\), we conclude:
\begin{equation*}
    \nrm{M_\nu f}_{L^2}\lesssim
    \nrm{S_\nu f}_{L^2}
    = 
    \nrm{
        \nrm{
            \pi_k f\br{x_0-2^k u_1,x_1-4^k u_2}
        }_{L^2\br{d\vx}}
    }_{\ell^2\br{k\in\Z}}
    = 
    \nrm{
        \nrm{
            \pi_k f
        }_{L^2}
    }_{\ell^2\br{k\in\Z}}
    \lesssim \nrm{f}_{L^2}.
\end{equation*}
This completes the justification for \(\nrm{M_\nu}_2\lesssim 1\).

\subsection{Energy and size estimates}\label{subsec_energy_est}

\subsubsection{Control over \texorpdfstring{\(\cE_2\)}{}}\label{enn}

For this, we simply observe that \(\cE_2\lesssim 1\) follows from iterating the energy selection lemma proved in \textsc{Section \ref{energy_select}}---see \textbf{Lemma \ref{lem_cE_2_bd}} therein.

\subsubsection{Control over \texorpdfstring{\(\cE_\ast\)}{}}\label{sizee}

To prove \(\cE_\ast\lesssim 1\), the reader may refer to \textbf{Lemma 6.13} in \cite{muscalu2013classicalv2} for a statement almost identical to \eqref{eq_Linfty_emb_sup}. We include a sketch of proof of \eqref{eq_Linfty_emb_sup} for readers' convenience. We start with:
\begin{theorem}[Adaptation of \textbf{Theorem 2.7} and \textbf{Corollary 2.8} in \cite{muscalu2013classicalv2}]
     Let \(r\in\br{0,\infty}\). For any finite family of intervals \(\cJ\subset\I\) and any sequence of complex numbers \(\BR{a_\vI}_{\vI\in\cJ}\), we have:
    \begin{equation*}
        \nrm{\br{a_\vI}_{\vI\in\cJ}}_{\bmo_r}:=\hspace{-1ex}
        \sup_{\vJ\in\cJ}\frac{1}{\abs{\vJ}^{\frac{1}{r}}}
        \nrm{
            \br{
                \sum_{\vI\subset \vJ}\frac{\abs{a_\vI}^2}{\abs{\vI}}\1_\vI
            }^{\frac{1}{2}}
        }_{L^r}\underset{r}{\eqsim}
        \sup_{\vJ\in\cJ}\frac{1}{\abs{\vJ}}
        \nrm{
            \br{
                \sum_{\vI\subset \vJ}\frac{\abs{a_\vI}^2}{\abs{\vI}}\1_\vI
            }^{\frac{1}{2}}
        }_{L^{1,\infty}}\hspace{-3ex}
        =:
        \nrm{\br{a_\vI}_{\vI\in\cJ}}_{\bmo_{1,\infty}}.
    \end{equation*}
\end{theorem}
This implies that for any finite \(\cP\subset \P\), the following estimate holds:
\begin{equation*}
    \nrm{f^{\br{+,j}}}_{\bmo_2\br{\cP}}\eqsim 
    \sup_{\substack{
        \cT\subset \cP\\
        \cT\text{ be }-,j\text{-tree}
    }}
    \frac{1}{\abs{\vI_\cT}}
    \nrm{
        \br{
        \sum_{\vI\times\vomega\in\cT}
            \frac{\abs{\ang{
                f,
                \varphi^{\br{+,j}}_{\vI\times\vomega}
            }}^2}{\abs{\vI}}
            \1_{\vI}
        }^{\frac{1}{2}}
    }_{L^{1,\infty}}
    =:\nrm{f}_{\bmo_{1,\infty}\br{\cP}}.
\end{equation*}
Note that due to the maximality, we can restrict the class of trees taken in the supremum and equate the right-hand side of the above with:
\begin{equation*}
    =
    \sup_{\substack{
        \cT\subset \cP\\
        \cT\text{ be }-,j\text{-tree}\\
        \vI_\cT\subset \vI_\cP
    }}
    \frac{1}{\abs{\vI_\cT}}
    \nrm{
        \br{
        \sum_{\vI\times\vomega\in\cT}
            \frac{\abs{\ang{
                f,
                \varphi^{\br{+,j}}_{\vI\times\vomega}
            }}^2}{\abs{\vI}}
            \1_{\vI}
        }^{\frac{1}{2}}
    }_{L^{1,\infty}},
\end{equation*}
where \(\vI_\cP\in\I\) denotes the minimal interval that contains \(\bigcup_{\vI\times\vomega}\vI\).
To estimate the right-hand side, we utilize the following result:
\begin{lemma}[Reinterpretation of equation (2.74) of \textbf{Lemma 2.13} in \cite{muscalu2013classicalv2}]
    For any \(\br{-,j}\)-tree \(\cT\subset \P\), we have:
    \begin{equation*}
        \frac{1}{\abs{\vI_\cT}}
        \nrm{
            \br{
            \sum_{\vI\times\vomega\in\cT}
                \frac{\abs{\ang{
                    f,
                    \varphi^{\br{+,j}}_{\vI\times\vomega}
                }}^2}{\abs{\vI}}
                \1_{\vI}
            }^{\frac{1}{2}}
        }_{L^{1,\infty}}
        \lesssim
        \nrm{f}_{L^1\br{d\mu_{\vI_\cT}}}.
    \end{equation*}
\end{lemma}
In combination, we conclude:
\begin{equation*}
    \nrm{f^{\br{+,j}}}_{\bmo_2\br{\cP}}\eqsim 
    \nrm{f}_{\bmo_{1,\infty}\br{\cP}}
    \lesssim
    \sup_{\substack{
        \varnothing\neq\cT\subset\cP\\
        \cT\text{ be tree}\\
        \vI_\cT\subset \vI_\cP
    }}
    \nrm{f}_{L^1\br{d\mu_{\vI_\cT}}}
\end{equation*}
and thus, \(\cE_\ast \lesssim 1\).

\subsection{Proof of Propositions \ref{thm_endpoint_infty_1} and \ref{thm_endpoint_1_infty}}

The arguments in this section are closely related to those in the approach of the classical Carleson theorem. 

\begin{proof}[\textbf{Proof of} \textbf{Proposition \ref{thm_endpoint_infty_1}}]
Fix \(E,F\subset \R^2\).
Without loss of generality, we assume \(0<\abs{E},\abs{F}<\infty\).
Again, it suffices to show the estimate with \(\P\) replaced by a finite collection of tiles \(\cP\subset \P\).
We initiate the argument by decomposing \(\cP\) into collection of trees \(\BR{\cT_{n,i}}_{n,i}\) in the following manner:
\begin{enumerate}
    \item Set \(\cP_{n_0}:=\cP\) for \(n_0\in\Z\) being the largest number such that \(2^{-n_0}\log\br{e+\abs{\vu}}\abs{E}\geq 1\).
    \item Given \(\cP_{n-1}\), we apply \textbf{Lemma \ref{lem_u_shift_mass_sel}} to obtain a disjoint collection of trees \(\BR{\cT_{n,i}}_{i=1}^m\) in \(\cP_{n-1}\) such that the following two estimates hold:
    \begin{equation*}
        \cM_{\cP_{n-1}\setminus\bigsqcup_i\cT_{n,i}}\br{E}\lesssim 2^{-n}\log\br{e+\abs{\vu}}\abs{E},\quad
        \sum_{i=1}^m\abs{\vI_{\cT_{n,i}}}\lesssim 2^n.
    \end{equation*}
    \item By \textbf{Definition \ref{def_out_L2_infty_energy_emb} and \ref{def_out_Linfty_energy_emb}}, there is a disjoint collection of trees \(\BR{\cT_{n,i}}_{i=m+1}^M\) in \(\cP_{n-1}\setminus\bigsqcup_{i=1}^m\cT_{n,i}\) such that the following two estimates hold:
    \begin{equation*}
        \nrm{\1^{\br{+,j}}_F}_{\bmo_2\br{
            \cP_{n-1}\setminus\bigsqcup_{i=1}^M\cT_{n,i}
        }}
        \lesssim \min\br{\cE_\ast,\cE_2 2^{-n/2}\abs{F}^{\frac{1}{2}}},\quad
        \sum_{i=m+1}^M\abs{\vI_{\cT_{n,i}}}\lesssim 2^n.
    \end{equation*}
    \item Set \(\cP_n:=\cP_{n-1}\setminus\bigsqcup_{i=1}^M\cT_{n,i}\).
\end{enumerate}
We iterate until the algorithm exhausts the collection of tiles. By construction, we have \(\cP=\bigsqcup_{n,i}\cT_{n,i}\) and
\begin{equation}\label{eq_energy_mass_bd_in_proof}
    \nrm{\1^{\br{+,j}}_F}_{\bmo_2\br{\cT_{n,i}}}
    \lesssim \min\br{\cE_\ast,\cE_2 2^{-n/2}\abs{F}^{\frac{1}{2}}},\quad
    \cM_{\cT_{n,i}}\br{E}\lesssim
    2^{-n}\log\br{e+\abs{\vu}}\abs{E}
    .
\end{equation}
We now apply the decomposition to rewrite the model sum:
\begin{equation*}
    \nrm{\Lambda^{\br{j}}_P\br{\1_F,\1_E}}_{\ell^1\br{P\in\cP}}=\sum_n \sum_i \nrm{\Lambda^{\br{j}}_P\br{\1_F,\1_E}}_{\ell^1\br{P\in\cT_{n,i}}}.
\end{equation*}
By \textbf{Lemma \ref{lem_gen_tree_c_half}}, we dominate the above with:
\begin{equation*}
    \lesssim \br{\nrm{\Lambda}_{\frac{1}{2}}+\nrm{T_\nu}_2}\cL
    \sum_{n,i}
        \nrm{\1^{\br{+,j}}_F}_{\bmo_2\br{\cT_{n,i}}}
        \cM_{\cT_{n,i}}\br{E}
        \abs{\vI_{\cT_{n,i}}}.
\end{equation*}
Applying \eqref{eq_energy_mass_bd_in_proof}, we further estimate the above with:
\begin{align*}
    \lesssim &\br{\nrm{\Lambda}_{\frac{1}{2}}+\nrm{T_\nu}_2}\cL
    \sum_{n>n_0}
        \min\br{\cE_\ast,\cE_2 2^{-n/2}\abs{F}^{\frac{1}{2}}}
        2^{-n}\log\br{e+\abs{\vu}}\abs{E}
        \sum_i
            \abs{\vI_{n,i}}\\
    \lesssim &\br{\nrm{\Lambda}_{\frac{1}{2}}+\nrm{T_\nu}_2}\cL^2
    \sum_{n>n_0}
        \min\br{\cE_\ast,\cE_2 2^{-n/2}\abs{F}^{\frac{1}{2}}}
        \abs{E}.
\end{align*}
Let \(n_1\in\Z\) be the largest number such that \(\cE_2 2^{-n/2}\abs{F}^{\frac{1}{2}}\geq \cE_\ast\). A direct calculation completes the proof:
\begin{align*}
    \lesssim &\br{\nrm{\Lambda}_{\frac{1}{2}}+\nrm{T_\nu}_2}\cL^2
    \br{
    \sum_{n_0<n\leq n_1}
        \cE_\ast
        \abs{E}
    +\sum_{n>n_1}
        \cE_2 2^{-n/2}\abs{F}^{\frac{1}{2}}
        \abs{E}
    }\\
    \lesssim &\br{\nrm{\Lambda}_{\frac{1}{2}}+\nrm{T_\nu}_2}
    \cE_\ast\cL^2
    \log\br{e+\frac{\cE_2^2\abs{F}}{\cE_\ast^2\abs{E}}}\abs{E}.
\end{align*}
\end{proof}

\begin{proof}[\textbf{Proof of} \textbf{Proposition \ref{thm_endpoint_1_infty}}]
Suggested by \textbf{Definition \ref{def_out_Linfty_energy_emb}}, the majorant subset \(E'\subset E\) shall be chosen to respect the maximal operator adapted to the system of interval \(\I\):
\begin{equation}
    \cM f\br{\vx}:=
    \sup_{\vx\in \vI\in\I}
    \nrm{f}_{L^1\br{d\mu_\vI}}.
\end{equation}
To be precise, we define:
\begin{equation}
    \Omega:=\br{\cM\1_F}^{-1}\bR{C\abs{F}/\abs{E},\infty}.
\end{equation}
By the weak \(\br{1,1}\) boundedness of \(\cM\), there is a universal constant \(C\eqsim 1\) such that \(\abs{\Omega}<\abs{E}/2\) for all \(E,F\) with finite measure. We thus set \(E':=E\setminus \Omega\) and \(\cP=\bigsqcup_{n=0}^\infty \cP_n\) with respect to the physical distance to \(\Omega^c\):
\begin{equation}
    \cP_0:=\BR{
        \vI\times\vomega\in\cP\::\:
        \vI\not\subset\Omega
    },\quad
    \cP_n:=\BR{
        I_0\times I_1 \times\vomega\in\cP\::\:
        2^{n-1} I_0\times 4^{n-1}I_1\subset\Omega,\,
        2^n I_0\times 4^n I_1\not\subset\Omega
    }.
\end{equation}
\subsubsection{Treatment of \texorpdfstring{\(\cP_0\)}{}}
By design, \textbf{Definition \ref{def_out_Linfty_energy_emb}} gives:
\begin{equation}
    \nrm{\1^{\br{+,j}}_F}_{\bmo_2\br{\cP_0}}\leq 
    \cE_\ast 
    \sup_{\substack{
        \varnothing\neq\cT\subset\cP_0\\
        \cT\text{ be tree}
    }} 
        \mu_{\vI_\cT}\br{F} 
    \leq \cE_\ast \nrm{\cM\1_F}_{L^\infty\br{\Omega^c}} \lesssim \cE_\ast\abs{F}/\abs{E}.
\end{equation}
An almost identical argument as in the proof of \textbf{Theorem \ref{thm_endpoint_infty_1}} with all presence of \(\cE_\ast\) replaced by \(\cE_\ast\abs{F}/\abs{E}\) proves the following estimate:
\begin{equation}
    \nrm{\Lambda^{\br{j}}_P\br{\1_F,\1_{E'}}}_{\ell^1\br{P\in\cP_0}}\lesssim 
        \br{\nrm{\Lambda}_{\frac{1}{2}}+\nrm{T_\nu}_2}
        \cE_\ast\cL^2
        \log\br{e+\frac{\cE_2^2\abs{E}}{\cE_\ast^2\abs{F}}}\abs{F}.
\end{equation}
\subsubsection{Treatment of \texorpdfstring{\(\cP_n\)}{}}
Take \(\vI\times\vomega\in\cP_n\). Expanding \textbf{Definition \ref{def_sing_tile_est} and \ref{def_out_Linfty_energy_emb}} gives:
\begin{align*}
    \abs{\Lambda^{\br{j}}_{\vI\times\vomega}\br{\1_F,\1_{E'}}}
    \lesssim &
    \nrm{\Lambda}_{\frac{1}{2}}
    \nrm{\1^{\br{+,j}}_F}_{\vI\times\vomega}
    \mu_{\vu\boxplus \vI}\br{E'\cap A_\vomega}\abs{\vI}^{\frac{1}{2}}\\
    \lesssim &
    \nrm{\Lambda}_{\frac{1}{2}}
    \nrm{\1^{\br{+,j}}_F}_{\bmo_2\br{\BR{\vI\times\vomega}}}
    \mu_{\vu\boxplus \vI}\br{E'\cap A_\vomega}\abs{\vI}\\
    \lesssim &
    \nrm{\Lambda}_{\frac{1}{2}}\cE_\ast
    \mu_\vI\br{F}
    \mu_{\vu\boxplus \vI}\br{E'\cap A_\vomega}\abs{\vI}.
\end{align*}

Fix \(\vI=I_0\times I_1\) for a brief moment. Summing over \(\vomega\) gives:
\begin{equation}\label{eq_Pn_universal_bd}
    \sum_{\substack{\vomega\in\W\\
    : \vI\times\vomega\in\cP_n
    }}
    \abs{\Lambda^{\br{j}}_{\vI\times\vomega}\br{\1_F,\1_{E'}}}
    \lesssim 
    \nrm{\Lambda}_{\frac{1}{2}}\cE_\ast
    \mu_\vI\br{F}
    \mu_{\vu\boxplus \vI}\br{E\setminus\Omega}\abs{\vI}.
\end{equation}
 Since \(2^n I_0\times 4^n I_1\not\subset \Omega\), there is \(\vJ\in \I \) such that the following three statement holds:
\begin{equation}
    \vJ\not\subset\Omega,\quad \abs{\vJ}\eqsim \abs{2^n I_0\times 4^n I_1} = 2^{3n}\abs{\vI},\quad \chi_\vI \lesssim \chi_\vJ.
\end{equation}
As a result, we have:
\begin{equation}\label{eq_exc_set_dil_bd}
    \mu_\vI\br{F}\lesssim 2^{3n}\mu_\vJ\br{F}\leq 2^{3n} \nrm{\cM\1_F}_{L^\infty\br{\Omega^c}} \lesssim 2^{3n} \abs{F}/\abs{E}.
\end{equation}
To sum over intervals in \(\I_n:=\BR{
        \vI\in\I\::\:
        \vI\times\vomega\in\cP_n
    }\), we observe that the defining property of \(\cP_n\) forces that for any \(\vI=I_0\times I_1,\, \vI'=I'_0\times I'_1\in\I_n\), the following statement holds:
\begin{equation}
    \vI\subset \vI'\implies  2 I_0\times 4 I_1 \not \subset I_0'\times I_1'.
\end{equation}
This implies a Carleson packing condition:
\begin{equation}\label{eq_Car_pack}
    \sum_{
        \vI\in\I_n\::\:
        \vI\subset \vJ
    }\abs{\vI}\lesssim \abs{\vJ},\quad
    \vJ\in\I_n,
\end{equation}
which implies that the function \(\chi_n:=\sum_{\vI\in\I_n}\1_\vI\) has bounded BMO norm \(\nrm{\chi_n}_{\bmo_\I}\lesssim 1\) and satisfies the corresponding John-Nirenberg inequality for some universal constant \(C>0\):
\begin{equation}\label{eq_J_N_ineq}
    \abs{
        \vJ\cap \chi_n^{-1}\bR{\varsigma,\infty}
    }
    \lesssim e^{-C\varsigma} \abs{\vJ},\quad \vJ\in\I_n.
\end{equation}
We now have two scenarios:
\begin{itemize}
\item For \(n> 2\log_2\br{1+\abs{\vu}}\),
We start with the trivial bound
\begin{equation}
    \mu_{\vu\boxplus \vI}\br{E\setminus \Omega}\lesssim \nrm{\chi_{\vu\boxplus \vI}^N}_{L^\infty\br{\Omega^c}}
    \lesssim \sup_{\vx\in\Omega^c}\ang{\frac{\dist\br{ I_0+u_0\abs{I_0},x_0}}{\abs{I_0}}}^{-N}
    \ang{\frac{\dist\br{I_0+u_1\abs{I_1},x_1}}{\abs{I_1}}}^{-N}.
\end{equation}
From the defining property of \(\cP_n\) and the fact that \(n>2\log_2\br{1+\abs{\vu}}\), we dominate the above with
\begin{align*}
    \lesssim & \ang{\frac{\dist\br{ I_0+u_0\abs{I_0},\br{2^{n-1} I_0}^c}}{\abs{I_0}}}^{-N}
     \ang{\frac{\dist\br{I_1+u_1\abs{I_1},\br{4^{n-1} I_1}^c}}{\abs{I_1}}}^{-N}\\
     \lesssim & \ang{\frac{\dist\br{ 2^{n/2}I_0,\br{2^{n-1} I_0}^c}}{\abs{I_0}}}^{-N}
     \ang{\frac{\dist\br{4^{n/2}I_1,\br{4^{n-1} I_1}^c}}{\abs{I_1}}}^{-N}
     \lesssim 2^{-3nN}.\numberthis\label{eq_exc_set_far_tile_est}
\end{align*}
Combine \eqref{eq_exc_set_dil_bd} and \eqref{eq_exc_set_far_tile_est}, we derive:
\begin{equation}
    \sum_{\substack{\vomega\in\W\\
    : \vI\times\vomega\in\cP_n
    }}
    \abs{\Lambda^{\br{j}}_{\vI\times\vomega}\br{\1_F,\1_{E'}}}
    \lesssim
    2^{3n\br{1-N}}
    \nrm{\Lambda}_{\frac{1}{2}}\cE_\ast
    \br{\abs{F}/\abs{E}} \abs{\vI}.
\end{equation}
We now sum over \(\vI\in\I_n\) by arranging the sum in the following way:
\begin{align*}
    \nrm{\abs{\Lambda^{\br{j}}_P\br{\1_F,\1_{E'}}}}_{\ell^1\br{P\in\cP_n}}
    \lesssim &
        2^{3n\br{1-N}}
        \nrm{\Lambda}_{\frac{1}{2}}\cE_\ast
        \br{\abs{F}/\abs{E}}
        \sum_{\vI\in\I_n}\abs{\vI}\\
    \lesssim &
    2^{3n\br{1-N}}
    \nrm{\Lambda}_{\frac{1}{2}}\cE_\ast
    \br{\abs{F}/\abs{E}}
    \sum_{\vJ\in \sup\I_n}\sum_{\substack{\vI\in\I_n\\ \vI\subset \vJ}}\abs{\vI},
\end{align*}
where for \(\J\subset \I\), \(\sup\J\) denotes the collection of maximal intervals in \(\J\) under the usual inclusion partial ordering in \(\I\).
By the Carleson packing condition \eqref{eq_Car_pack}, we dominate the above with
\begin{equation*}
    \lesssim
    2^{3n\br{1-N}}
    \nrm{\Lambda}_{\frac{1}{2}}\cE_\ast
    \br{\abs{F}/\abs{E}}
    \sum_{\vJ\in \sup\I_n}\abs{\vJ}
    \lesssim
    2^{3n\br{1-N}}
    \nrm{\Lambda}_{\frac{1}{2}}\cE_\ast
    \br{\abs{F}/\abs{E}}
    \abs{\Omega}
    \lesssim
    2^{3n\br{1-N}}
    \nrm{\Lambda}_{\frac{1}{2}}\cE_\ast
    \abs{F}.
\end{equation*}
We thus conclude that:
\begin{equation}
    \sum_{n> 2\log_2\br{1+\abs{\vu}}}
    \nrm{\Lambda^{\br{j}}_P\br{\1_F,\1_{E'}}}_{\ell^1\br{P\in\cP_n}}
    \lesssim
    \nrm{\Lambda}_{\frac{1}{2}}\cE_\ast
    \abs{F}.
\end{equation}

\item For \(n\leq 2\log_2\br{1+\abs{\vu}}\),
since we do not have an analog of \eqref{eq_exc_set_far_tile_est} to counter the exponential growth from \eqref{eq_exc_set_dil_bd}, we remedy this by applying the John-Nirenberg inequality on \(\chi_n\). To be precise, we consider the set \(\Omega_n:=\chi_n^{-1}\bR{\frac{3n \log 2}{C},\infty}\). On the one hand, we have the estimate
\begin{equation}\label{eq_J_H_hieght_bd}
    \sum_{\substack{
        \vI\in\I_n\\
        \vI\not\subset \Omega_n
    }}\1_\vI
    \lesssim n
    \quad
    \text{and thus},\quad
    \sum_{\substack{
        \vI\in\I_n\\
        \vI\not\subset \Omega_n
    }}\chi^N_\vI
    \lesssim n.
\end{equation}
On the other hand,
by John-Nirenberg inequality \eqref{eq_J_N_ineq}, we have:
\begin{equation}
    \abs{\Omega_n}=\sum_{\vJ\in \sup\I_n} \abs{\vJ\cap \Omega_n}\lesssim 2^{-3n}\sum_{\vJ\in \sup\I_n} \abs{\vJ} \leq 2^{-3n}\abs{\Omega}\leq 2^{-3n}\abs{E}.
\end{equation}
This suggests that we decompose the sum in the following manner:
\begin{align*}
    \nrm{\Lambda^{\br{j}}_P\br{\1_F,\1_{E'}}}_{\ell^1\br{P\in\cP_n}}
    \leq & 
    \nrm{\Lambda}_{\frac{1}{2}}\cE_\ast
    \bigg(
    \sum_{\substack{
        \vI\in\I_n\\
        \vI\not\subset \Omega_n
    }}
    +
    \sum_{\substack{
        \vI\in\I_n\\
        \vI\subset \Omega_n
    }}
    \bigg)
    \mu_\vI\br{F}
    \mu_{\vu\boxplus \vI}\br{E\setminus\Omega}\abs{\vI}\\
    \lesssim &
    \nrm{\Lambda}_{\frac{1}{2}}\cE_\ast
    \br{
    \bigg\Vert
    \sum_{\substack{
        \vI\in\I_n\\
        \vI\not\subset \Omega_n
    }}\chi_\vI^N \1_F
    \bigg\Vert_{L^1}
    +
    2^{3n}\br{\abs{F}/\abs{E}}
    \sum_{\substack{
        \vI\in\I_n\\
        \vI\subset \Omega_n
    }}
    \abs{\vI}
    }.
\end{align*}
We apply \eqref{eq_J_H_hieght_bd} for the first term and \eqref{eq_Car_pack} for the second term to dominate the above with
\begin{align*}
    \lesssim & \nrm{\Lambda}_{\frac{1}{2}}\cE_\ast
    \br{
    \nrm{n\1_F}_{L^1}+2^{3n}\br{\abs{F}/\abs{E}}
    \sum_{
    \vJ\in \I_n^\ast
    }
    \sum_{\substack{
        \vI\in \I_n\\
        \vI\subset \vJ
    }}
    \abs{\vI}}
    ,\quad
    \I_n^\ast:=\sup\BR{\vI\in\I_n\::\:\vI\subset \Omega_n}
    \\
    \lesssim & 
    \nrm{\Lambda}_{\frac{1}{2}}\cE_\ast
    \br{
    n\abs{F}
    +
    2^{3n}\br{\abs{F}/\abs{E}}\abs{\Omega_n}}
    \lesssim n\nrm{\Lambda}_{\frac{1}{2}}\cE_\ast
    \abs{F}.
\end{align*}
We thus conclude:
\begin{equation}
    \sum_{1\leq n\leq 2\log_2\br{1+\abs{\vu}}}
    \nrm{\Lambda^{\br{j}}_P\br{\1_F,\1_{E'}}}_{\ell^1\br{P\in\cP_n}}
    \lesssim 
    \nrm{\Lambda}_{\frac{1}{2}}\cE_\ast \abs{F}
    \sum_{1\leq n\leq 2\log_2\br{1+\abs{\vu}}} n\lesssim
    \nrm{\Lambda}_{\frac{1}{2}}\cE_\ast\cL^2\abs{F}.
\end{equation}
\end{itemize}
Putting together the two cases above, we conclude the proof of \textbf{Theorem \ref{thm_endpoint_1_infty}} as follows:
\begin{equation*}
    \nrm{\Lambda^{\br{j}}_P\br{\1_F,\1_{E'}}}_{\ell^1\br{P\in\P}}= 
    \sum_{n=0}^\infty
    \nrm{\Lambda^{\br{j}}_P\br{\1_F,\1_{E'}}}_{\ell^1\br{P\in\P_n}}
    \lesssim
        \br{\nrm{\Lambda}_{\frac{1}{2}}+\nrm{T_\nu}_2}
        \cE_\ast\cL^2
        \log\br{e+\frac{\cE_2^2\abs{E}}{\cE_\ast^2\abs{F}}}\abs{F}.
\end{equation*}
\end{proof}


\section{Proof of the key lemmas stated within the formal time-frequency analysis section}\label{sec_pf_key_lemmas}

In this section, we present the proofs of the four key lemmas stated in \textsc{Section \ref{sec_formal_tf_ana}}, that is, the proofs of \textbf{Lemmas \ref{lem_-_tree}}, \textbf{\ref{lem_gen_tree_c_half}},  \textbf{\ref{lem_u_shift_mass_sel}} and \textbf{\ref{lem_cE_2_bd}}. This is the most technically involved section of our paper, and thus, we invite the reader to treat the unfortunate technicalities below with patience.

\subsection{Prerequisites} In this very brief subsection, we address the following issue: the physical shift \(\vu\) in \(\vu\boxplus\vI\) breaks the order structure in the system of intervals \(\I\). We address this by introducing auxiliary supersets and the corresponding modified characteristic functions:
\begin{definition}
Recalling \eqref{eq_cU_cV_def}, we define:
    \begin{equation}
        n\br{\cU\boxplus \vI}:=
        \bigcup_{\vu_\ast\in\cU} n\br{\vu_\ast \boxplus \vI},\quad
        \chi^N_{\cU\boxplus \vI}:=
        \sum_{\vu_\ast\in\cU} \chi^N_{\vu_\ast \boxplus \vI},\quad
        \vI\in\I.
    \end{equation}
\end{definition}
\begin{observation}\label{obs_shifts_with_order}
    By filling in the \textbf{holes} between scales, we preserve the order structure:
    \begin{equation}
        \vI\subset \vJ\implies 
        \vI\subset n\br{\cU\boxplus\vI}\subset n\br{\cU\boxplus\vJ}
        \quad\text{and}\quad
        \chi^N_{\vu\boxplus \vI}
        \leq
        \chi^N_{\cU\boxplus\vI}\lesssim
        \chi^N_{\cU\boxplus\vJ}
        ,\quad
        \vI,\vJ\in\I,\quad n\geq 1.
    \end{equation}
\end{observation}

\subsection{Proof of Lemma \ref{lem_-_tree}: \texorpdfstring{\(\br{-,j}\)}{}-tree estimate}
We shall make a few reductions. Due to the positive nature of the left-hand side, we may assume without loss of generality that \(\cT\) is a finite \(\br{-,j}\)-tree. 
Moreover, we may sacrifice a \(O\br{1}\) loss, sparsify the scaling, and assume that the tree \(\cT\) has the following structure:
\begin{equation}\label{eq_tree_sparse}
    \cT=\bigsqcup_{k\in 10\Z}\cT_k,\quad
    \cT_k:=\cT\cap\P_k.
\end{equation}
Henceforth, \(k\) ranges within \(10\Z\) throughout the rest of the section.
Additionally, since \(j\) will be fixed throughout the argument, we will suppress the dependence whenever it does not affect the understanding.
Next, we aim to linearize the \(\ell^1\) norm. 
Consider uni-modulus coefficient \(\epsilon_{\vI\times\vomega}\in\C\)
such that: 
\begin{equation}\label{eq_remove_abs_w_coef}
    \abs{\Lambda^{\br{j}}_{\vI\times\vomega}\br{f,\1_E}}=\epsilon_{\vI\times\vomega}
    \Lambda^{\br{j}}_{\vI\times\vomega}\br{f,\1_E}=
    \ang{\cC_k \epsilon_{\vI\times\vomega} f^{\br{+,j}}_{\vI\times\vomega},\1_{E\cap A^{\br{-,j}}_\vomega}}.
\end{equation}
This suggests that we set \(f^\epsilon_{\vI\times\vomega}:=\epsilon_{\vI\times\vomega} f^{\br{+,j}}_{\vI\times\vomega}\) and define accordingly that \(f^\epsilon_\cP:=\sum_{P\in\cP}f^\epsilon_P\) for \(\cP\subset \P\).
We thus linearize the \(\ell^1\) norm:
\begin{equation}\label{eq_T_k_E_omega_def}
    \nrm{\Lambda^{\br{j}}_P\br{f,\1_E}}_{\ell^1\br{P\in\cT}}=
    \sum_{k\in 10\Z}
    \sum_{\vI\times\vomega \in\cT_k}\ang{\cC_k f^\epsilon_{\vI\times\vomega},\1_{E^-_\vomega}},\quad
    E^-_\vomega:=E\cap A^{\br{-,j}}_\vomega.
\end{equation}
\textbf{Lemma \ref{lem_-_tree}} can be reformulated as follows:
\begin{equation}\label{eq_the_tree_mod_sum}
    \sum_{k\in 10\Z}
    \sum_{\vI\times\vomega \in\cT_k}\ang{\cC_k f^\epsilon_{\vI\times\vomega},\1_{E^-_\vomega}}
    \lesssim
    \br{\nrm{\Lambda}_c+\nrm{T_\nu}_2}\cL\nrm{f^{\br{+,j}}}_{L^2\br{\cT}}M_\cT^{\frac{1}{2}+c}\br{E}\abs{\vI_\cT}^{\frac{1}{2}}.
\end{equation}
Due to the symmetry of the formulation, we may assume without loss of generality that \(P_\cT=\bR{0,1}^2\times\bR{-\lambda/2,\lambda/2}\times\bR{0,\lambda}\in\P_0\br{\lambda}\). 
We start with constructing a Whitney decomposition adapted to the tree \(\cT\). That is, we consider:
\begin{equation}\label{eq_J_def}
    \J:=
    \sup
    \bigcap_{\substack{
        \vI\times\vomega\in\cT\\
        \vu_\ast\in\cU
    }}
    \BR{
        \vJ\in\I\::\:
        \vu_\ast\boxplus \vI\not\subset 3\vJ
    },
\end{equation}
where the \(\sup\) denotes the collection of all maximal elements under the natural partial order in \(\ang{\I,\subset}\).
Since \(\J\subset \I\) is a partition of \(\R^2\), we have the following identity:
\begin{equation}
    \nrm{\Lambda^{\br{j}}_P\br{f,\1_E}}_{\ell^1\br{P\in\cT}}
    =
    \sum_{\vJ\in\J}
    \sum_{k\in 10\Z}
    \sum_{\vI\times\vomega \in\cT_k}\ang{\cC_k f^\epsilon_{\vI\times\vomega},\1_{E^-_\vomega\cap\vJ}}.
\end{equation}
This suggests that we perform our analysis relative to the location and scale of \(\vJ\in\J\). We set \(k_\vJ\in\Z\) be such that \(\vJ\in\I_{k_\vJ}\).
The above sum can be split into two:
\begin{equation}\label{eq_small_large_scale_decomp}
    =
    \sum_{\vJ\in\J}
    \sum_{k<k_\vJ+2}
    \sum_{\vI\times\vomega \in\cT_k}\ang{\cC_k f^\epsilon_{\vI\times\vomega},\1_{E^-_\vomega\cap\vJ}}
    +
    \sum_{\vJ\in\J}
    \sum_{k\geq k_\vJ+2}
    \sum_{\vI\times\vomega \in\cT_k}\ang{\cC_k f^\epsilon_{\vI\times\vomega},\1_{E^-_\vomega\cap\vJ}}.
\end{equation}
We analyze the two terms individually:
\subsubsection{Small scale contribution from \texorpdfstring{$k<k_\vJ+2$}{k<kJ+2}}\label{subsec_tree_small_scale_wide}
Recall \textbf{Definition \ref{def_sing_tile_est}}, we have the following estimate:
\begin{align*}
    \sum_{\vJ\in\J}
    \sum_{k<k_\vJ+2}
    \sum_{\vI\times\vomega \in\cT_k}&
        \abs{
            \ang{\cC_k f^\epsilon_{\vI\times\vomega},\1_{E^-_\vomega\cap \vJ}}
        }
    \leq 
    \nrm{\Lambda}_c
    \sum_{\vJ\in\J}
    \sum_{k<k_\vJ+2}
    \sum_{\vI\times\vomega\in\cT_k}
    \nrm{f^{\br{+,j}}}_{\vI\times\vomega}
    \mu_{\vu\boxplus\vI}^{\frac{1}{2}+c}\br{E^-_\vomega\cap \vJ}\abs{\vI}^{\frac{1}{2}}\\
    \leq &
    \nrm{\Lambda}_c
    \cM^c_\cT\br{E}
    \sum_{\vJ\in\J}
    \sum_{k<k_\vJ+2}
    \sum_{\vI\times\vomega\in\cT_k}
    \nrm{f^{\br{+,j}}}_{\vI\times\vomega}
    \mu_{\vu\boxplus\vI}^{\frac{1}{2}}\br{E^-_\vomega\cap \vJ}\abs{\vI}^{\frac{1}{2}}.
\end{align*}
A direct calculation shows the following:
\begin{observation}\label{obs_small_scale_meas_bd}
    \begin{align*}
        \mu_{\vu\boxplus \vI}\br{E^-_\vomega\cap \vJ}
        = &
        \nrm{\1_{E^-_\vomega\cap \vJ}\chi^{\frac{N}{2}}_{\vu\boxplus \vI}}_{L^1\br{d\mu_{\vu\boxplus\vI}}}
        \leq
        \nrm{\chi_{\vu\boxplus\vI}}_{L^\infty\br{\vJ}}^{\frac{N}{2}}
        \nrm{\1_{E^-_\vomega\cap \vJ}}_{L^1\br{d\mu_{\vu\boxplus\vI}}}\\
        \eqsim &
        \underbrace{\ang{\frac{\dist\br{u_0\boxplus I_0,J_0}}{\abs{I_0}}}^{-\frac{N}{2}}
        \ang{\frac{\dist\br{u_1\boxplus I_1,J_1}}{\abs{I_1}}}^{-\frac{N}{2}}}_{=:\chi^{\frac{N}{2}}_{\vu\boxplus \vI}\br{\vJ}}
        \mu_{\vu\boxplus\vI}\br{E^-_\vomega\cap\vJ}.
    \end{align*}
\end{observation}
We thus obtain:
\begin{align}\label{eq_z_tree_small_scale_useful_bd}
    &
    \sum_{\vJ\in\J}
    \sum_{k<k_\vJ+2}
    \sum_{\vI\times\vomega \in\cT_k}
        \abs{
            \ang{\cC_k f^\epsilon_{\vI\times\vomega},\1_{E^-_\vomega\cap \vJ}}
        }\nonumber\\
    \leq &
    \nrm{\Lambda}_c
    \cM^c_\cT\br{E}
    \sum_{\vJ\in\J}
    \sum_{k<k_\vJ+2}
    \sum_{\vI\times\vomega\in\cT_k}
    \nrm{f^{\br{+,j}}}_{\vI\times\vomega}
    \mu_{\vu\boxplus\vI}^{\frac{1}{2}}\br{E^-_\vomega\cap \vJ}
    \chi^{\frac{N}{4}}_{\vu\boxplus \vI}\br{\vJ}
    \abs{\vI}^{\frac{1}{2}}\\
    \leq &
    \nrm{\Lambda}_c
    \cM^c_\cT\br{E}
    \nrm{
        \nrm{f^{\br{+,j}}}_{\vI\times\vomega}
        \nrm{
            \mu_{\vu\boxplus\vI}^{\frac{1}{2}}\br{E^-_\vomega\cap \vJ}
        }_{\ell^2\br{\vJ\in\J}}
    }_{\ell^2\br{\vI\times\vomega\in\times\cT}}
    \nrm{
        \chi^{\frac{N}{4}}_{\vu\boxplus \vI}\br{\vJ}
        \abs{\vI}^{\frac{1}{2}}
    }_{\ell^2\br{\substack{\br{\vJ,\vI\times\vomega}\in\J\times\cT\\ \abs{\vI}<\abs{\vJ}}}}\nonumber\\
    \leq &
    \nrm{\Lambda}_c
    \nrm{f^{\br{+,j}}}_{L^2\br{\cT}}
    \cM^{\frac{1}{2}+c}_\cT\br{E}
    \br{
        \sum_{\vJ\in\J}
        \sum_{k<k_\vJ+2}
        \sum_{\vI\times\vomega \in\cT_k}
                \chi^{\frac{N}{2}}_{\vu\boxplus \vI}\br{\vJ}
                \abs{\vI}
    }^{\frac{1}{2}}.\nonumber
\end{align}
It remains to derive the following estimate:
\begin{equation}\label{eq_small_scale_vJ_bd}
    \sum_{\vJ\in\J}
    \sum_{k<k_\vJ+2}
    \sum_{\vI\times\vomega \in\cT_k}
            \chi^{\frac{N}{2}}_{\vu\boxplus \vI}\br{\vJ}
            \abs{\vI}
    \lesssim
    \abs{\cU\boxplus\vI_\cT}
    \lesssim
    \log\br{e+\abs{\vu}}
    \abs{\vI_\cT}
    \lesssim
    \cL
    \abs{\vI_\cT}.
\end{equation}
We shall perform the analysis for fixed \(\vJ\) and scale \(k<k_\vJ+2\) first:
\begin{equation*}
    \sum_{\vI\times\vomega \in\cT_k}
        \chi^{\frac{N}{2}}_{\vu\boxplus \vI}\br{\vJ}
        \abs{\vI}
    \leq
    2^{3k}
    \chi^{\frac{N}{4}}_{\cU\boxplus\vI_\cT}\br{\vJ}
    \sum_{\vI\times\vomega \in\cT_k}
    \ang{\frac{\dist\br{u_1\boxplus I_0,J_0}}{2^k}}^{-\frac{N}{4}}
    \ang{\frac{\dist\br{u_2\boxplus I_1,J_1}}{4^k}}^{-\frac{N}{4}}
    .
\end{equation*}
Since \(k<k_\vJ+2\), we have the following two relations:
\begin{equation*}
    2^k=\abs{I_0}<2^2\abs{J_0},\quad 4^k=\abs{I_1}<4^2\abs{J_1}.
\end{equation*}
Additionally, due to the Whitney structure of \(\J\), at least one of the following holds:
\begin{equation*}
    \dist\br{u_0\boxplus I_0,J_0}\gtrsim \abs{J_0}\text{, or }
    \dist\br{u_1\boxplus I_1,J_1}\gtrsim \abs{J_1}.
\end{equation*}
As a direct consequence, we thus have:
\begin{align*}
    &
    \sum_{\vI\times\vomega \in\cT_k}
    \ang{\frac{\dist\br{u_0\boxplus I_0,J_0}}{2^k}}^{-\frac{N}{4}}
    \ang{\frac{\dist\br{u_1\boxplus I_1,J_1}}{4^k}}^{-\frac{N}{4}}\\
    \lesssim &
    \sum_{m\in\N}
    \ang{
        \frac{\abs{J_0}}{2^k }+m
    }^{-\frac{N}{4}}
    \br{
        \frac{\abs{J_1}}{4^k}
        +
    \sum_{n\in\N}
    \ang{n}^{-\frac{N}{4}}
    }
    + 
    \br{
        \frac{\abs{J_0}}{2^k}
        +
    \sum_{m\in\N}
    \ang{m}^{-\frac{N}{4}}
    }
    \sum_{n\in\N}
    \ang{
        \frac{\abs{J_1}}{4^k }+n
    }^{-\frac{N}{4}}\\
    \lesssim &
    \sum_{m\in\N}
    \ang{
        2^{k_\vJ-k}+m
    }^{-\frac{N}{4}}
    \br{
        4^{k_\vJ-k}
        +
    \sum_{n\in\N}
    \ang{n}^{-\frac{N}{4}}
    }
    + 
    \br{
        2^{k_\vJ-k}
        +
    \sum_{m\in\N}
    \ang{m}^{-\frac{N}{4}}
    }
    \sum_{n\in\N}
    \ang{
        4^{k_\vJ-k}+n
    }^{-\frac{N}{4}}\\
    \lesssim &
    2^{-\frac{\br{k_\vJ-k}N}{4}}\cdot 4^{k_\vJ-k}
    +
    2^{k_\vJ-k}\cdot 4^{-\frac{\br{k_\vJ-k}N}{4}}
    \lesssim 1.
\end{align*}
As a result, we can sum over \(\vJ,k\) and conclude:
\begin{equation*}
    \sum_{\vJ\in\J}
    \sum_{k<k_\vJ+2}
    \sum_{\vI\times\vomega \in\cT_k}
            \chi^{\frac{N}{2}}_{\vu\boxplus \vI}\br{\vJ}
            \abs{\vI}
    \lesssim 
    \sum_{\vJ\in\J}
    \sum_{k<k_\vJ+2}
    2^{3k}
    \chi^{\frac{N}{4}}_{\cU\boxplus\vI_\cT}\br{\vJ}
    \lesssim
    \sum_{\vJ\in\J}
        \abs{\vJ}
        \chi^{\frac{N}{4}}_{\cU\boxplus\vI_\cT}\br{\vJ}
    \lesssim \abs{\cU\boxplus\vI_\cT}.
\end{equation*}
In combination, we obtain the desired estimate:
\begin{equation*}
    \sum_{\vJ\in\J}
    \sum_{k<k_\vJ+2}
    \sum_{\vI\times\vomega \in\cT_k}
        \abs{
            \ang{\cC_k f^\epsilon_{\vI\times\vomega},\1_{E^-_\vomega\cap \vJ}}
        }
    \lesssim 
    \nrm{\Lambda}_c\cL^{\frac{1}{2}}
    \nrm{f^{\br{+,j}}}_{L^2\br{\cT}}
    \cM^{\frac{1}{2}+c}_\cT\br{E}
    \abs{\vI_\cT}^{\frac{1}{2}}.
\end{equation*}

\subsubsection{Large scale contribution from \texorpdfstring{\(k\geq k_\vJ+2\)}{k >= k\_J+2}}\label{subsec_tree_large_scale_wide}
Guided by the heuristic provided by the argument of the usual Carleson operator, we expect the estimate to involve a form of maximal truncated singular integral operator. To perform the analysis, we recall the definition of the frequency projection \eqref{eq_LP_freq_proj}.
Due to the tree structure, \textbf{Observation \ref{obs_tree_freq_scale_rel}} suggests that we set:
\begin{equation}\label{eq_def_E_J_E_k}
    E_\vJ:=\vJ\cap\bigcup_{k\geq k_\vJ+2}E^-_k
    ,\quad
    E^-_k:=\left\{
    \begin{aligned}
        &E^-_{\vomega_{\cT_k}}&,\quad\cT_k\neq\varnothing,\\
        &\varnothing &,\quad \cT_k=\varnothing.
    \end{aligned}
    \right.
\end{equation}
Recall also the choice of Littlewood-Paley projection \eqref{eq_LP_freq_proj}.
The location of the tree top \(P_\cT\) and the sparsity assumption \eqref{eq_tree_sparse} guarantees that we have \(\pi_k f^\epsilon_\cT=\pi_k f^\epsilon_{\cT_k}=f^\epsilon_{\cT_k}\).
We thus rewrite the sum in terms of \(\pi_k\):
\begin{equation}\label{eq_E_J_tight}
    \sum_{\vJ\in\J}
    \sum_{k\geq k_\vJ+2}
        \sum_{\vI\times\vomega \in\cT_k}
            \ang{\cC_k f^\epsilon_{\vI\times\vomega},\1_{E^-_\vomega\cap \vJ}}
    = 
    \sum_{\vJ\in\J}
    \ang{
        \sum_{k=k_\bot\vee \br{k_\vJ+2}}^{k_\top}
            \cC_k \pi_k f^\epsilon_\cT,
        \1_{E_\vJ}
    }
    ,
\end{equation}
where the two measurable functions \(k_\top\) and \(k_\bot\) are defined below:
\begin{equation*}
    k_\top\br{\vx}:=\sup\Z_\cT\br{\vx},\quad
    k_\bot\br{\vx}:=\inf\Z_\cT\br{\vx},\quad
    \Z_\cT\br{\vx}
    :=
    \BR{
        k\in 10\Z\::\:
        \vx\in E^-_k
    }
\end{equation*}
with the convention that \(\sup\varnothing =-\infty\) and \(\inf\varnothing=\infty\).
Intuitively, our intention is to show that
\begin{equation}\label{eq_tree_Feff_heuri}
    \sum_{k=m}^{n}
        \cC_k \pi_k f^\epsilon_\cT\br{\vx}
    \approx 
    \abs{
    \sum_{k=m}^n
        \pi_k T_\nu f^\epsilon_\cT\br{\vx}
    }
    + 
    o\br{1}
    M_\nu T_\nu f^\epsilon_\cT\br{\vx}.
\end{equation}
To proceed, we observe that for \(\vx\in E_\vJ\), we have the estimate \(\abs{a\br{\vx}}\lesssim 4^{-k_\top\br{\vx}} \lambda\). As a direct consequence of \textbf{Assumption \ref{ass_phy_loc_meas}}, we have for \(\vtau\in\supp\nu_k\) the following estimate:
\begin{equation}\label{eq_exp_of_a_bd}
    \abs{\va\br{\vx}\cdot\vtau}\leq\abs{a\br{\vx}}\br{1+\abs{u_1}}4^k 
    \lesssim
    4^{k-k_\top\br{\vx}}
    \lambda
    \br{1+\abs{u_1}}.
\end{equation}
This suggests that we split the sum to separate at least the top \(O\br{\log\br{\lambda\br{1+\abs{u_1}}}}\) many scales of oscillatory nature. In fact, we will separate \(\cL_1:=10\cL\) many scales:
\begin{equation*}
    \ang{
        \sum_{k=k_\bot\vee \br{k_\vJ+2}}^{k_\top}
            \cC_k \pi_k f^\epsilon_\cT
        ,\1_{E_\vJ}
    }
    =
    \ang{
        \sum_{k=k_\top-\cL_1}^{k_\top}
            \1_{E^-_k}
            \cC_k \pi_k f^\epsilon_\cT
        ,\1_{E_\vJ}
    }
    +
    \ang{
        \sum_{k=k_\bot\vee \br{k_\vJ+2}}^{k_\top-\cL_1-1}
            \cC_k \pi_k f^\epsilon_\cT
        ,\1_{E_\vJ}
    }.
\end{equation*}
Moreover, we wish to obtain estimates almost constant in \(\vJ\). We shall utilize the averaging effect of \(\pi_k\) for \(k\) large enough in the appropriate sense. To be more precise, observe the following pointwise estimate:
\begin{equation}\label{eq_pi_to_max}
    \abs{\pi_k g\br{\vx}},
    \abs{\pi^k g\br{\vx}}
    \lesssim 
    \cA_k g\br{\vx}:=
    \int
    \frac{ \abs{g\br{\vx-\vz}} \lambda^2 2^{-3k} d\vz}{
        \ang{\lambda z_0/2^k}^N
        \ang{\lambda z_1/4^k}^N
    }
    .
\end{equation}
As soon as we know that \(\vx\in\vJ\) and that \(k\geq k_\vJ+\cL_1\), we can apply the pointwise estimate \(\cA_k g\br{\vx} \lesssim \inf_{\vJ} M^{\otimes 2}g\). This further motivates us to separate the lower \(\cL_1\) many scales as well. For now, we have:
\begin{align*}
    \sum_{\vJ\in\J}
    \ang{
        \sum_{k=k_\bot\vee \br{k_\vJ+2}}^{k_\top}
            \cC_k \pi_k f^\epsilon_\cT
        ,\1_{E_\vJ}
    }
    = &
    \sum_{\vJ\in\J}
    \ang{
        \br{
            \sum_{k=k_\bot\vee \br{k_\vJ+2}}^{\br{k_\bot\vee \br{k_\vJ+2}}+\cL_1}
            +
            \sum_{k=k_\top-\cL_1}^{k_\top}
        }
            \1_{E^-_k}
            \cC_k \pi_k f^\epsilon_\cT
        ,\1_{E_\vJ}
    }\numberthis \label{eq_tree_Log_scale}\\
    + &
    \sum_{\vJ\in\J}
    \ang{
        \sum_{k=\br{k_\bot\vee \br{k_\vJ+2}}+\cL_1+1}^{k_\top-\cL_1-1}
            \cC_k \pi_k f^\epsilon_\cT
        ,\1_{E_\vJ}
    }.
    \numberthis\label{eq_tree_MSIO_MM}
\end{align*}

\subsubsection{Treatment of \texorpdfstring{\eqref{eq_tree_Log_scale}}{}: control over \texorpdfstring{\(O\br{\cL}\)}{} many scales}\label{subsubsec_j_many_scale}
In this context, we use a square sum argument:
\begin{equation*}
    \br{
        \sum_{k=k_\bot\vee \br{k_\vJ+2}}^{\br{k_\bot\vee \br{k_\vJ+2}}+\cL_1}
        +
        \sum_{k=k_\top-\cL_1}^{k_\top}
    }
        \abs{
            \1_{E^-_k}
            \cC_k \pi_k f^\epsilon_\cT
        }
    \lesssim
    \cL^{\frac{1}{2}}
    \nrm{
        \1_{E^-_k} \cC_k \pi_k f^\epsilon_\cT
    }_{\ell^2\br{ k\in\Z }}.
\end{equation*}
As a result, we can control the corresponding inner product trivially via the Cauchy-Schwarz inequality:
\begin{align*}
&
\abs{
    \sum_{\vJ\in\J}
    \ang{
        \br{
            \sum_{k=k_\bot\vee \br{k_\vJ+2}}^{\br{k_\bot\vee \br{k_\vJ+2}}+\cL_1}
            +
            \sum_{k=k_\top-\cL_1}^{k_\top}
        }
        \cC_k \pi_k f^\epsilon_\cT,
        \1_{E_\vJ}
    }
}\\
\leq &
\cL^{\frac{1}{2}}
\sum_{\vJ\in\J}
\nrm{
    \nrm{
        \1_{E^-_k} \cC_k \pi_k f^\epsilon_\cT
    }_{L^2\br{J}}
}_{\ell^2\br{k}}
\abs{E_\vJ}^{\frac{1}{2}}
\leq 
\cL^{\frac{1}{2}}
\nrm{
    \nrm{
        \1_{E^-_k} \cC_k \pi_k f^\epsilon_\cT
    }_{L^2}
}_{\ell^2\br{k}}
\bigg\vert\bigsqcup_{\vJ\in\J} 
E_\vJ\bigg\vert^{\frac{1}{2}}.
\end{align*}
We recall the key fact that \(\pi_kf^\epsilon_\cT=
f^\epsilon_{ \cT_k }
=\sum_{P\in\cT_k}f^\epsilon_P\).
This allows us to compute the \(L^2\) norm via the bilinear bound given by \textbf{Definition \ref{def_sing_tile_est}}. To be more specific, for \(g\in L^2\)
\begin{align}\label{eq_lem_-_tree_large_scale_log_many_sqbd}
    \abs{\ang{\1_{E^-_k}\cC_k\pi_k f^\epsilon_\cT,g}}
    \leq &
    \sum_{\vI\times\vomega\in \cT_k}
    \abs{
        \ang{\1_{E^-_k}\cC_k      
            f^\epsilon_{\vI\times\vomega},
            g
        }
    }
    \leq
    \nrm{\Lambda}_c
    \sum_{\vI\times\vomega\in \cT_k}
    \mu_{\vu\boxplus\vI}^c\br{E^-_k}
    \nrm{f^{\br{+,j}}}_{\vI\times\vomega}
    \nrm{g}_{L^2\br{d\mu_{\vu\boxplus\vI}}}\abs{\vI}^{\frac{1}{2}}\nonumber\\
    \leq &
    \nrm{\Lambda}_c \cM^c_\cT\br{E}
    \nrm{
        \nrm{f^{\br{+,j}}}_P
    }_{\ell^2\br{P\in \cT_k}}
    \nrm{
        \nrm{g}_{L^2\br{d\mu_{\vu\boxplus\vI}}}\abs{\vI}^{\frac{1}{2}}
    }_{\ell^2\br{\vI\times\vomega\in \cT_k}}\nonumber\\
    \lesssim &
    \nrm{\Lambda}_c \cM^c_\cT\br{E}
    \nrm{f^{\br{+,j}}}_{L^2\br{\cT_k}}
    \nrm{g}_{L^2}.
\end{align}
In other words, we derive the estimate:
\begin{equation*}
    \nrm{
        \nrm{
            \1_{E^-_k} \cC_k \pi_k f^\epsilon_\cT
        }_{L^2}
    }_{\ell^2\br{k}}
    \lesssim 
    \nrm{\Lambda}_c \cM^c_\cT\br{E}
    \nrm{
        \nrm{f^{\br{+,j}}}_{L^2\br{\cT_k}}
    }_{\ell^2\br{k}}
    =
    \nrm{\Lambda}_c \cM^c_\cT\br{E}
    \nrm{f^{\br{+,j}}}_{L^2\br{\cT}}.
\end{equation*}
To complete the estimate for the \(O\br{\cL}\) many scales, it remains to show that
\begin{equation*}
    \bigg\vert\bigsqcup_{\vJ\in\J} E_\vJ\bigg\vert
    \lesssim \cL\cM_\cT\br{E}\abs{\vI_\cT}.
\end{equation*}
In fact, since \(\J\) is a partition, it suffices to prove the following localized result:
\begin{lemma}\label{lem_E_J_den_est_shift}
Let \(\cT\) be a tree, \(\vJ\in\J\) be defined as in \eqref{eq_J_def}, and \(E_\vJ\) be defined as in \eqref{eq_def_E_J_E_k}. We have:
    \begin{equation*}
        \abs{E_\vJ}\lesssim \cM_\cT\br{E}
        \abs{
            \vJ\cap 9\br{\cU\boxplus \vI_\cT}
        }.
    \end{equation*}
\end{lemma}
\begin{proof}[\textbf{Proof of} \textbf{Lemma \ref{lem_E_J_den_est_shift}}]
    Assume that \(E_\vJ\neq\varnothing\).
    There is a tile \( \vI\times\vomega\in\cT\) such that \(\abs{\vI}\geq 8^2 \abs{\vJ}\). As a direct consequence, \(8^2\abs{\vJ}\leq \abs{\vI_\cT}\).
    Moreover, due to the Whitney structure of \(\J\), we must have:
    \begin{equation*}
        \exists \vI_0\times\vomega_0\in\cT \text{ and }\vu_0\in\cU \text{ such that }
        \vu_0\boxplus \vI_0\subset 3\widehat{\vJ}.
    \end{equation*}
    On the other hand, by \textbf{Observation \ref{obs_shifts_with_order}}, we derive the following relation:
    \begin{equation*}
        \vu_0\boxplus \vI_0\subset \cU\boxplus \vI_0\subset
        \cU\boxplus \vI_\cT.
    \end{equation*}
    This shows that
    \begin{equation*}
        \vu_0\boxplus \vI_0\subset 3\widehat{\vJ}\cap \br{\cU\boxplus \vI_\cT} \neq \varnothing.
    \end{equation*}
    A simple calculation shows that \(\vJ\subset 9\br{\cU\boxplus \vI_\cT}\).
    We now choose the minimal tile \(\vI_1\times\vomega_1\in\P\) satisfying \(\vI_0\times\vomega_0\leq \vI_1\times\vomega_1\leq P_\cT\) and \(\abs{\vI_0}\leq 8^2\abs{\vJ}\leq \abs{\vI_1}\). On the one hand, we have
    \begin{equation*}
         \bigcup\BR{
            \vomega\in\W\::\:
            \vI\times\vomega\in\cT,\,\abs{\vI}\geq 8^2 \abs{\vJ}
         }\subset \vomega_1
    \end{equation*}
    and thus, \(E_\vJ\subset E\cap A_{\vomega_1}\cap \vJ\). On the other hand,
    by \textbf{Observation \ref{obs_shifts_with_order}}, we see that:
    \begin{equation*}
        \vu_0\boxplus \vI_0\subset
         3\widehat{\vJ}\cap \br{\cU\boxplus \vI_1}\neq \varnothing.
    \end{equation*}
    This implies the existence of \(\vu_1\in\cU\) such that \(3\widehat{\vJ}\cap 5\br{\vu_1\boxplus \vI_1}\neq \varnothing\). Thus, we have \(\1_\vJ\lesssim \chi_{\vu_1\boxplus \vI_1}^{N/2}\) and 
    \begin{equation*}
        \abs{E_\vJ} \leq\abs{E \cap A_{\vomega_1}\cap \vJ} \lesssim 
        \mu_{\vu_1\boxplus\vI_1}\br{E\cap A_{\vomega_1}}
        \abs{\vI_1}\eqsim
        \mu_{\vu_1\boxplus\vI_1}\br{E\cap A_{\vomega_1}}
        \abs{\vJ}.
    \end{equation*}
    Finally, since \(\vI_0\times\vomega_0\in\cT\) and \(\vI_0\times\vomega_0\leq \vI_1\times\vomega_1\in\P\), we can thus dominate the above with:
    \begin{equation*}
        \leq \cM_\cT\br{E}\abs{\vJ}
        = \cM_\cT\br{E}\abs{\vJ\cap 9\br{\cU\boxplus \vI_\cT}}.
    \end{equation*}
\end{proof}
By \textbf{Lemma \ref{lem_E_J_den_est_shift}}, we conclude the desired estimate:
\begin{equation*}
    \eqref{eq_tree_Log_scale}
    \lesssim 
    \nrm{\Lambda}_c\cL^{\frac{1}{2}}\nrm{f^{\br{+,j}}}_{L^2\br{\cT}}\cM^{\frac{1}{2}+c}_\cT\br{E}\abs{9\br{\cU\boxplus \vI_\cT}}^{\frac{1}{2}}
    \lesssim 
    \nrm{\Lambda}_c\cL\nrm{f^{\br{+,j}}}_{L^2\br{\cT}}\cM^{\frac{1}{2}+c}_\cT\br{E}\abs{\vI_\cT}^{\frac{1}{2}}.
\end{equation*}

\subsubsection{Treatment of \texorpdfstring{\eqref{eq_tree_MSIO_MM}}{}: control over the majority of the scales}
We now execute the idea from \eqref{eq_tree_Feff_heuri}. Fix \(\vJ\in\J\) and \(\vx\in E_\vJ\) for a brief moment, we split the sum into two terms:
\begin{align*}
    \sum_{k=\br{k_\bot\br{\vx}\vee \br{k_\vJ+2}}+\cL_1+1}^{k_\top\br{\vx}-\cL_1-1}
            \cC_k \pi_k f^\epsilon_\cT\br{\vx}
    = &
    \sum_{k=\br{k_\bot\br{\vx}\vee \br{k_\vJ+2}}+\cL_1+1}^{k_\top\br{\vx}-\cL_1-1}
    \int
        \pi_k f^\epsilon_\cT\br{\vx-\vtau}
    d\nu_k\br{\vtau}\numberthis\label{eq_tree_MSIO}\\
    +
    \sum_{k=\br{k_\bot\br{\vx}\vee \br{k_\vJ+2}}+\cL_1+1}^{k_\top\br{\vx}-\cL_1-1}
    &\int
        \pi_k f^\epsilon_\cT\br{\vx-\vtau}
        \br{e\br{\va\br{\vx}\cdot\vtau}-1}
    d\nu_k\br{\vtau}.\numberthis\label{eq_tree_MM}
\end{align*}
Starting with \eqref{eq_tree_MSIO}, we recall \textbf{Definition \ref{def_sio_max}} and observe the identity 
\begin{equation*}
    \sum_{k=m}^n
    \int
        \pi_k f^\epsilon_\cT\br{\vx-\vtau}
    d\nu_k\br{\vtau}= 
    \sum_{k=m}^n
    \pi_k T^{10\Z}_\nu f^\epsilon_\cT\br{\vx}
    =
    \pi^n T^{10\Z}_\nu f^\epsilon_\cT\br{\vx}-
    \pi^{m-1} T^{10\Z}_\nu f^\epsilon_\cT\br{\vx}.
\end{equation*}
where here we utilized condition \eqref{eq_tree_sparse} and the restriction of the summation in \(k\) within the set \(10\Z\).
In combination with relation \eqref{eq_pi_to_max} and the fact that \(k> k_\vJ+\cL_1\), we obtain the following estimate:
\begin{equation*}
    \abs{\eqref{eq_tree_MSIO}}\lesssim
    \sup_{k_\vJ+\cL_1< k <k_\top\br{\vx}-\cL_1}
    \cA_k T^{10\Z}_\nu f^\epsilon_\cT\br{\vx}
    \lesssim
    \inf_\vJ 
    M^{\otimes 2}
    T^{10\Z}_\nu
    f^\epsilon_\cT.
\end{equation*}
For \eqref{eq_tree_MM}, we apply a similar reasoning:
\begin{align*}
    \abs{\eqref{eq_tree_MM}}\leq &
    \sup_{\substack{
        m,n\in\Z\\
        k_\vJ+\cL_1 < m \leq n < k_\top\br{\vx}-\cL_1
    }}
    \sum_{k=m}^n
    \int
        \abs{
            \pi_k f^\epsilon_\cT\br{\vx-\vtau}
            \br{e\br{\va\br{\vx}\cdot\vtau}-1}
        }
    d\abs{\nu_k}\br{\vtau}\\
    \lesssim &
    \sup_{\substack{
        m,n\in\Z\\
        k_\vJ+\cL_1 < m \leq n < k_\top\br{\vx}-\cL_1
    }}
    \sum_{k=m}^n
    4^{k-k_\top\br{\vx}}\br{1+\abs{u_1}}\lambda
    \abs{\pi_k f^\epsilon_\cT }\ast d\abs{\nu_k}\br{\vx}
\end{align*}
Since \(k\leq k_\top\br{\vx}-\cL_1\), we have \(4^{k-k_\top\br{\vx}}\br{1+\abs{u_1}}\lambda \lesssim \nrm{T_\nu}_2/\nrm{M_\nu}_2\). We thus dominate the above with:
\begin{align*}
    \lesssim &
    \frac{\nrm{T_\nu}_2}{\nrm{M_\nu}_2}
    \cdot
    \sup_{
        k_\vJ+\cL_1 <k< k_\top\br{\vx}-\cL_1
    }
    \abs{\pi_k \circ \pi_k f^\epsilon_\cT }\ast d\abs{\nu_k}\br{\vx}\\
    \lesssim &
    \frac{\nrm{T_\nu}_2}{\nrm{M_\nu}_2}
    \cdot
    \sup_{
        k_\vJ+\cL_1 <k< k_\top\br{\vx}-\cL_1
    }
    \cA_k M_\nu f^\epsilon_\cT\br{\vx}
    \lesssim
    \frac{\nrm{T_\nu}_2}{\nrm{M_\nu}_2}
    \cdot
    \inf_\vJ
    M^{\otimes 2}M_\nu 
    f^\epsilon_\cT.
\end{align*}
In combination, we obtain for \(\vx\in E_\vJ\):
\begin{equation*}
    \abs{
        \sum_{k=\br{k_\bot\br{\vx}\vee \br{k_\vJ+2}}+\cL_1+1}^{k_\top\br{\vx}-\cL_1-1}
            \cC_k \pi_k f^\epsilon_\cT\br{\vx}
    }
    \lesssim
    \inf_\vJ 
    M^{\otimes 2}
    T_\nu
    f^\epsilon_\cT
    +
    \frac{\nrm{T^{10\Z}_\nu}_2}{\nrm{M_\nu}_2}
    \cdot
    \inf_\vJ
    M^{\otimes 2}M_\nu 
    f^\epsilon_\cT.
\end{equation*}
Using the above pointwise estimate, we deduce that
\begin{equation*}
    \abs{
        \eqref{eq_tree_MSIO_MM}
    }
    \lesssim
    \sum_{\vJ\in\J}
    \abs{E_\vJ}
    \br{
        \inf_\vJ 
        M^{\otimes 2}
        T^{10\Z}_\nu
        f^\epsilon_\cT
        +
        \frac{\nrm{T_\nu}_2}{\nrm{M_\nu}_2}
        \cdot
        \inf_\vJ
        M^{\otimes 2}M_\nu 
        f^\epsilon_\cT
    }.
\end{equation*}
Applying now \textbf{Lemma \ref{lem_E_J_den_est_shift}} we dominate the above expression by
\begin{align*}
    \lesssim &
    \cM_\cT\br{E}
    \sum_{\vJ\in\J}
    \abs{\vJ\cap 9\br{\cU\boxplus\vI_\cT}}
    \br{
        \inf_\vJ 
        M^{\otimes 2}
        T^{10\Z}_\nu
        f^\epsilon_\cT
        +
        \frac{\nrm{T_\nu}_2}{\nrm{M_\nu}_2}
        \cdot
        \inf_\vJ
        M^{\otimes 2}M_\nu 
        f^\epsilon_\cT
    }\\
    \lesssim &
    \cM_\cT\br{E}
    \br{
    \nrm{
        M^{\otimes 2}
        T^{10\Z}_\nu
        f^\epsilon_\cT
    }_{L^1\br{9\br{\cU\boxplus\vI_\cT}}}
    +
    \frac{\nrm{T_\nu}_2}{\nrm{M_\nu}_2}
    \nrm{
        M^{\otimes 2}M_\nu 
        f^\epsilon_\cT
    }_{L^1\br{9\br{\cU\boxplus\vI_\cT}}}
    }\\
    \lesssim &
    \cM_\cT\br{E}
    \br{
    \nrm{
        M^{\otimes 2}
        T^{10\Z}_\nu
        f^\epsilon_\cT
    }_{L^2}
    +
    \frac{\nrm{T_\nu}_2}{\nrm{M_\nu}_2}
    \nrm{
        M^{\otimes 2}M_\nu 
        f^\epsilon_\cT
    }_{L^2}
    }
    \abs{\cU\boxplus\vI_\cT}^{\frac{1}{2}}\\
    \lesssim &
    \nrm{T_\nu}_2
    \cL^{\frac{1}{2}}
    \nrm{
        f^\epsilon_\cT
    }_{L^2}
    \cM_\cT\br{E}
    \abs{\vI_\cT}^{\frac{1}{2}}.
\end{align*}
Lastly, by \textbf{Assumption \ref{ass_tf_proj_emb}}, we have \(\nrm{f^\epsilon_\cT}_{L^2}\lesssim \nrm{f^{\br{+,j}}}_{L^2\br{\cT}}\) and thus dominate the above with
\begin{equation*}
    \lesssim
    \nrm{T_\nu}_2
    \cL^{\frac{1}{2}}
    \nrm{f^{\br{+,j}}}_{L^2\br{\cT}}
    \cM_\cT\br{E}
    \abs{\vI_\cT}^{\frac{1}{2}}
    .
\end{equation*}
This completes the proof of \eqref{eq_the_tree_mod_sum} and thus of \textbf{Lemma \ref{lem_-_tree}}.\qed

\subsection{Proof of Lemma \ref{lem_gen_tree_c_half}: general tree estimate when \texorpdfstring{\(c=\frac{1}{2}\)}{}}
Given a tree \(\cT\subset \P\), we may again assume without loss of generality that \(\cT\) is finite and satisfies \eqref{eq_tree_sparse}. 
\subsubsection{Reduction to \texorpdfstring{\(\br{+,1}\)}{}-tree}
By definition, \(\cT\) can be split into disjoint union of a \(\br{-,j}\)-tree \(\cT_-\) and a \(\br{+,j}\)-tree \(\cT_+\) with the same top \(P_\cT\). However, since all the $j=0$ trees are $-$ trees, it suffices then to consider two cases: \emph{(i)} \(\cT\) is a \(\br{-,j}\)-tree for \(j\in\BR{0,1}\), or, \emph{(ii)} 
\(\cT\) is a \(\br{+,1}\)-tree. 

In the former, \(\br{-,j}\)-tree case, we may apply \textbf{Lemma \ref{lem_-_tree}} with \(c=\frac{1}{2}\) to obtain:
\begin{equation*}
    \nrm{\Lambda^{\br{j}}_P\br{f,\1_E}}_{\ell^1\br{P\in\cT}}
        \lesssim \br{\nrm{\Lambda}_{\frac{1}{2}}+\nrm{T_\nu}_2}\cL\nrm{f^{\br{+,j}}}_{L^2\br{\cT}}M_\cT\br{E}\abs{\vI_\cT}^{\frac{1}{2}}.
\end{equation*}
Using the fact that \(\nrm{f^{\br{+,j}}}_{L^2\br{\cT}}\leq \nrm{f^{\br{+,j}}}_{\bmo_2\br{\cT}}\abs{\vI_\cT}^\frac{1}{2}\), we dominate the above by
\begin{equation*}
    \leq \br{\nrm{\Lambda}_{\frac{1}{2}}+\nrm{T_\nu}_2}\cL\nrm{f^{\br{+,j}}}_{\bmo_2\br{\cT}}M_\cT\br{E}\abs{\vI_\cT},
\end{equation*}
which is the desired estimate.

Consequently, it remains to focus on the latter case, and thus, we assume from now on that \(\cT\) is a \(\br{+,1}\)-tree. We now recall the following key steps \eqref{eq_remove_abs_w_coef}, \eqref{eq_T_k_E_omega_def}, and \eqref{eq_J_def} in order to perform the analogous decomposition \eqref{eq_small_large_scale_decomp}.

\subsubsection{Small scale contribution from \texorpdfstring{\(k<k_\vJ+2\)}{}}\label{subsubsec_gen_tree_small_scale}
Recalling \textbf{Definition \ref{def_sing_tile_est}}, we have
\begin{align*}
    \sum_{\vJ\in\J}
    \sum_{k<k_\vJ+2}
    \sum_{\vI\times\vomega \in\cT_k}
        \abs{
            \ang{\cC_k f^\epsilon_{\vI\times\vomega},\1_{E^-_\vomega\cap \vJ}}
        }
    \leq &
    \nrm{\Lambda}_{\frac{1}{2}}
    \sum_{\vJ\in\J}
    \sum_{k<k_\vJ+2}
    \sum_{\vI\times\vomega\in\cT_k}
    \nrm{f^{\br{+,1}}}_{\vI\times\vomega}
    \mu_{\vu\boxplus\vI}\br{E^-_\vomega\cap \vJ}\abs{\vI}^{\frac{1}{2}}\\
    \leq &
    \nrm{\Lambda}_{\frac{1}{2}}
    \nrm{f^{\br{+,1}}}_{\bmo_2\br{\cT}}
    \sum_{\vJ\in\J}
    \sum_{k<k_\vJ+2}
    \sum_{\vI\times\vomega\in\cT_k}
    \mu_{\vu\boxplus\vI}\br{E^-_\vomega\cap \vJ}\abs{\vI}.
\end{align*}

Using \textbf{Observation \ref{obs_small_scale_meas_bd}} we dominate the above by
\begin{align*}
    \leq & 
    \nrm{\Lambda}_{\frac{1}{2}}
    \nrm{f^{\br{+,1}}}_{\bmo_2\br{\cT}}
    \sum_{\vJ\in\J}
    \sum_{k<k_\vJ+2}
    \sum_{\vI\times\vomega\in\cT_k}
    \chi_{\vu\boxplus \vI}^{\frac{N}{2}}\br{\vJ}
    \mu_{\vu\boxplus\vI}\br{E^-_\vomega\cap \vJ}\abs{\vI}\\
    \leq & 
    \nrm{\Lambda}_{\frac{1}{2}}
    \nrm{f^{\br{+,1}}}_{\bmo_2\br{\cT}}
    \cM_\cT\br{E}
    \sum_{\vJ\in\J}
        \sum_{k<k_\vJ+2}
            \sum_{\vI\times\vomega \in\cT_k}
                \chi^{\frac{N}{2}}_{\vu\boxplus \vI}\br{\vJ}
                \abs{\vI}
\end{align*}
Finally, recalling \eqref{eq_small_scale_vJ_bd}, we dominate the above by
\begin{equation*}
    \lesssim
    \nrm{\Lambda}_{\frac{1}{2}}
    \cL\nrm{f^{\br{+,1}}}_{\bmo_2\br{\cT}}M_\cT\br{E}\abs{\vI_\cT}.
\end{equation*}
This completes the estimate for small scales.

\subsubsection{Large scale contribution from \texorpdfstring{\(k\geq k_\vJ+2\)}{}}\label{subsubsec_gen_tree_large_scale}
As in Subsection \textsc{\ref{subsubsec_gen_tree_small_scale}}, we notice
\begin{equation*}
    \sum_{\vJ\in\J}
    \sum_{k\geq k_\vJ+2}
    \sum_{\vI\times\vomega \in\cT_k}
        \abs{
            \ang{\cC_k f^\epsilon_{\vI\times\vomega},\1_{E^-_\vomega\cap \vJ}}
        }
    \leq 
    \nrm{\Lambda}_{\frac{1}{2}}
    \nrm{f^{\br{+,1}}}_{\bmo_2\br{\cT}}
    \sum_{\vJ\in\J}
    \sum_{k\geq k_\vJ+2}
    \sum_{\vI\times\vomega\in\cT_k}
    \mu_{\vu\boxplus\vI}\br{E^-_\vomega\cap \vJ}\abs{\vI}.
\end{equation*}
Next, we focus on the innermost summation:
\begin{equation}\label{eq_lem_gen_tree_c_half_large_scale_sum_mu_E}
    \sum_{\vI\times\vomega\in\cT_k}
    \mu_{\vu\boxplus\vI}\br{E^-_\vomega\cap \vJ}\abs{\vI}
    \lesssim
    \abs{E^-_k\cap\vJ}
\end{equation}
Moreover, since \(\cT\) is a \(\br{+,1}\)-tree, the following sets are disjoint
\begin{equation*}
    \BR{\vomega^{\br{-,1}}\::\:\exists \vI\in\I,\,\vI\times\vomega\in\cT}=\BR{\vomega^{\br{-,1}}_{\cT_k}}_k
\end{equation*}
and thus the same holds for the collection of sets \(\BR{E^-_k}_k\). Using this and \textbf{Lemma \ref{lem_E_J_den_est_shift}}, we dominate the total sum as follows:
\begin{equation}\label{eq_lem_gen_tree_c_half_large_scale_sum_E}
    \sum_{\vJ\in\J}
    \sum_{k\geq k_\vJ+2}
    \abs{E^-_k\cap \vJ}
    =
    \sum_{\vJ\in\J}
    \abs{
    \vJ\cap\bigsqcup_{k\geq k_\vJ+2}E^-_k
    }
    =
    \sum_{\vJ\in\J}
    \abs{E_\vJ}
    \lesssim 
    \cL\,
    \cM_\cT\br{E}\abs{\vI_\cT}.
\end{equation}
Putting everything together, we deduce
\begin{equation*}
    \sum_{\vJ\in\J}
    \sum_{k\geq k_\vJ+2}
    \sum_{\vI\times\vomega \in\cT_k}
        \abs{
            \ang{\cC_k f^\epsilon_{\vI\times\vomega},\1_{E^-_\vomega\cap \vJ}}
        }
    \leq 
    \nrm{\Lambda}_{\frac{1}{2}}
    \cL
    \nrm{f^{\br{+,1}}}_{\bmo_2\br{\cT}}
    \cM_\cT\br{E}\abs{\vI_\cT}\,,
\end{equation*}
thus completing the proof of \textbf{Lemma \ref{lem_gen_tree_c_half}}. \qed

\subsection{Proof of Lemma \ref{lem_u_shift_mass_sel}: mass selection}
    For notation clarity, we write \(E_\vomega:=E\cap A_\vomega\).
    Recalling now the definition of \(\cM_{\br{\cdot}}\br{E}\), we collect the following tiles:
    \begin{equation*}
        \cP_\ast:=
        \bigcup_{P\in\cP}
        \bigcup_{
            \vu_\ast\in\cU
        }
        \BR{
            \vI\times\vomega\in\P\::\:
            \mu_{\vu_\ast \boxplus \vI}\br{E_\vomega}>\varsigma\text{ and }
            P\leq\vI\times\vomega
        }.
    \end{equation*}
    Observe the trivial bound \(\mu_{\vu_\ast\boxplus \vI}\br{E_\vomega}\lesssim \abs{E}/\abs{\vI}\). This suggests that the collection \(\cP_\ast\) is bounded under the partial ordering \(\leq\) in \(\P\). Moreover, the existence of \(P'\in\cP\) such that \(P'\leq P\) for \(P\in\cP_\ast\) confines the time-frequency location of \(P\). As a result, \(\cP_\ast\) is finite, and we can enumerate the maximal elements:
    \begin{equation*}
        \bigvee\cP_\ast:=\BR{
            P\in\cP_\ast\::\:
            \forall P'\in\cP_\ast,\,
            \br{
                P\leq P'\implies P=P'
            }
        }
        =\BR{
            \vI_i\times\vomega_i
        }_i.
    \end{equation*}
    We define trees iteratively:
    \begin{equation*}
        \cT_i:=\BR{
            P\in\cP\setminus
            \bigsqcup_{i'<i}\cT_{i'}\::\:
            P\leq \vI_i\times\vomega_i
        },\quad
        i\in\N.
    \end{equation*}
    By construction,
    \begin{equation*}
        \cM_{\cP\setminus \bigsqcup_i\cT_i}\br{E} \leq
        \sup_{P\in \cP\setminus \bigsqcup_i\cT_i}
        \sup_{
            \substack{
                \vI\times\vomega\in \P\setminus \cP_\ast\\
                P\leq \vI\times \vomega
            }
        }
        \max_{
            \vu_\ast\in\cU
        }
        \mu_{\vu_\ast \boxplus \vI}\br{E_\vomega}
        \leq \varsigma.
    \end{equation*}
    It remains to show \(\sum_i\abs{\vI_i}\lesssim\log\br{e+\abs{\vu}}\abs{E}/\varsigma\). By the defining property of \(\cP_\ast\), we have for all \(i\) a vector \(\vu_i\in\cU\) such that the following holds:
    \begin{equation}\label{eq_mass_heavy_cond_shift}
        \mu_{\vu_i\boxplus \vI_i}\br{E_{\vomega_i}}>\varsigma\text{ or equivalently, }
        \abs{\vI_i}<
        \frac{1}{\varsigma}
        \int_{E_{\vomega_i}}
        \chi_{\vu_i\boxplus \vI_i}^{N}\br{\vx}
        d\vx.
    \end{equation}
    On the other hand, a layer cake decomposition shows:
    \begin{equation*}
        \varsigma\abs{\vI_i}<
        \int_{E_{\vomega_i}}
        \chi_{\vu_i\boxplus \vI_i}^{N}\br{\vx}
        d\vx
        \lesssim 
        \sum_{n\in\N}
        n^{-N}
        \abs{E_{\vomega_i} \cap n\br{\vu_i\boxplus \vI_i}}
        \lesssim
        \sup_{n\in\N}
        n^{-\frac{N}{2}}
        \abs{E_{\vomega_i} \cap n\br{\vu_i\boxplus \vI_i}}.
    \end{equation*}
    Suggested by the pigeon-holing principle, we shall consider the following set:
    \begin{equation*}
        \cI_n:=\BR{i\::\: \varsigma\abs{\vI_i}\lesssim n^{-\frac{N}{2}}
        \abs{E_{\vomega_i} \cap n\br{\vu_i\boxplus \vI_i}}},\quad
        n\in\N.
    \end{equation*}
    By design, \(\bigcup_n\cI_n\) exhaust all index \(i\) and thus \(\sum_i \abs{\vI_i}\leq\sum_{n\in\N}\sum_{i\in\cI_n}\abs{\vI_i}\). We now fix \(n\in\N\) for a brief moment and estimate the contribution from \(\cI_n\). To salvage the lack of disjointness due to the physical shift \(\vu_i\) and \(n\) dilation of the interval, we employ a covering lemma type algorithm:
    \begin{enumerate}
        \item Initiate with \(\cJ_0:=\varnothing\).
        \item Given \(\BR{\cJ_k}_{k<m}\), we set \(\cK_m:=\cI_n\setminus\bigsqcup_{k<m}\cJ_k\).
        \begin{itemize}
            \item If \(\cK_m=\varnothing\), we terminate the algorithm.
            \item Else if \(\cK_m\neq \varnothing\), there is \(i_m\in \cK_m \) such that \(\abs{\vI_{i_m}}=\sup_{i\in \cK_m}\abs{\vI_i}\).
        \end{itemize} 
        \item Set \(\cJ_m:=\BR{i\in\cK_m\::\: n\br{\vu_i\boxplus \vI_i}\times\vomega_i\cap n\br{\vu_{i_m}\boxplus \vI_{i_m}}\times\vomega_{i_m}\neq\varnothing}\) and iterate.
    \end{enumerate}
    By design, the collection of data \(\BR{\br{\cJ_m,i_m}}_m\) enjoys the following properties:
    \begin{itemize}
        \item \(\BR{\cJ_m}_m\) partitions \(\cI_n\). In other words, \(\cI_n=\bigsqcup_m\cJ_m\).
        \item \(\BR{n\br{\vu_{i_m}\boxplus \vI_{i_m}}\times\vomega_{i_m}}_m\) are disjoint and thus, \(\BR{E_{\vomega_{i_m}}\cap n \br{\vu_{i_m}\boxplus \vI_{i_m}}}_m\) are disjoint.
        \item The estimate \(\sum_{i\in\cJ_m}\abs{\vI_i}\lesssim n^2\log\br{e+\abs{\vu}}\abs{\vI_{i_m}}\) holds.
    \end{itemize}
    To see the last item, we first observe that for \(i\in\cJ_m\),
    \begin{equation*}
        \abs{\vI_i}\leq\abs{\vI_{i_m}}\text{ and }n\br{\vu_i\boxplus \vI_i}\cap n\br{\vu_{i_m}\boxplus \vI_{i_m}}\neq \varnothing.
    \end{equation*}
    As a direct consequence, we have 
    \begin{equation*}
        \bigcup_{i\in\cJ_m}\vI_i\subset
        \bigcup_{i\in\cJ_m}
        n\br{\cU\boxplus \vI_i}
        \subset
        3n\br{\cU\boxplus \vI_{i_m}}.
    \end{equation*}
    The statement follows once we show the collection of intervals \(\BR{\vI_i}_{i\in\cJ_m}\) are disjoint. This is achieved by noticing for \(i\in\cJ_m\) the defining property implies:
    \begin{equation*}
        \abs{\vomega_i}\geq \abs{\vomega_{i_m}}
        \text{ and }\vomega_i\cap\vomega_{i_m}\neq \varnothing.
    \end{equation*}
    This implies \(\vomega_i\supset \vomega_{i_m}\) and thus, \(\vomega_{i_m}\subset \bigcap_{i\in\cJ_m}\vomega_i\neq\varnothing\). Yet, by the maximality of \(\bigvee\cP_\ast\), the collection of sets \(\BR{\vI_i\times\vomega_i}_i\) is an anti-chain under \(\leq\). This forces \(\BR{\vI_i}_{i\in\cJ_m}\) to be disjoint. Therefore, 
    \begin{equation*}
        \sum_{i\in\cJ_m}\abs{\vI_i}=
        \bigg\vert\bigsqcup_{i\in\cJ_m}\vI_i\bigg\vert\leq \abs{3n\br{\cU\boxplus \vI_{i_m}}}\eqsim n^2\log\br{e+\abs{\vu}}\abs{\vI_{i_m}}.
    \end{equation*}
    We now return to the estimate of the contribution from \(\cI_n\):
    \begin{equation}\label{eq_I_n_vv_to_log}
        \sum_{i\in\cI_n}
            \abs{\vI_i}
        =
        \sum_m\sum_{i\in\cJ_m}\abs{\vI_i}
        \lesssim
        \sum_m n^2\log\br{e+\abs{\vu}}\abs{\vI_{i_m}}.
    \end{equation}
    By the defining property of \(\cI_n\) and the disjointness among \(\BR{E_{\vomega_{i_m}}\cap n \br{\vu_{i_m}\boxplus \vI_{i_m}
    }}_m\), we obtain
    \begin{equation*}
        \eqref{eq_I_n_vv_to_log}
        \lesssim \frac{n^{2-\frac{N}{2}}
        \log\br{e+\abs{\vu}}}{\varsigma}\sum_m \abs{E_{\vomega_{i_m}}\cap n \br{\vu_{i_m}\boxplus \vI_{i_m}}}
        \leq \frac{n^{2-\frac{N}{2}}\log\br{e+\abs{\vu}}\abs{E}}{\varsigma}.
    \end{equation*}
    Lastly, summing over \(n\in\N\) we conclude our proof:
    \begin{equation*}
        \sum_i\abs{\vI_i}\leq \sum_{n\in\N}\sum_{i\in\cI_n}\abs{\vI_i}
        \lesssim \sum_{n\in\N}n^{2-\frac{N}{2}}\cdot
        \frac{\log\br{e+\abs{\vu}}\abs{E}}{\varsigma}\eqsim \frac{\log\br{e+\abs{\vu}}\abs{E}}{\varsigma}.\qed
    \end{equation*}

\subsection{Proof of Lemma \ref{lem_cE_2_bd}: energy selection}

Without loss of generality, we may assume the normalization condition \(\nrm{f^{\br{+,j}}}_{\bmo_2\br{\cP}}=1\).
We aim to extract a sequence of ``well-separated" \(\br{-,j}\)-trees \(\BR{\cT_i}_i\) whose contribution to the first estimate exceeds the desired threshold \(\nrm{f^{\br{+,j}}}_{\bmo_2\br{\cP}}/2=1/2\). Naturally, the actual collection \(\T\) of trees in the lemma will then be composed from the above-mentioned sequence of trees. 
With these, we are now ready to describe the tree selection algorithm:
    \begin{itemize}
        \item Initialize \(i=1\), \(\cP_0:=\cP\), and \(\T_0=\varnothing\).
        \item Consider the following collection of trees:
        \begin{equation*}
            \T^\ast_i:=\BR{
                \cT \subset \P_{i-1} \text{ be } (-,j)\text{-tree with top }P_\cT:=\vI_\cT\times \vomega_\cT\::\:
                \nrm{f^{\br{+,j}}}_{L^2\br{\cT}}>\abs{\vI_\cT}^{\frac{1}{2}}/2
            }.
        \end{equation*}
        
        \item If \(\T^\ast_i=\varnothing\), we terminate the algorithm. Else, we choose \(\cT_i\in\T^\ast_i\) with its top 
        having the ``highest" frequency location in the \(\xi_1\) coordinate among all tops of trees in \(\T^\ast_i\). To make the statement precise, we write for any tree \(\cT\) its top as:
        \begin{equation*}
            P_\cT:=
            \vI_\cT\times\vomega_\cT=
            \vI_\cT\times\omega_{\cT,0}\times\omega_{\cT,1}.
        \end{equation*}
        The tree \(\cT_i\in\T^\ast_i\) is then chosen to satisfy the following maximality condition
        \begin{equation}\label{eq_cE_2_rightmost}
            \inf\omega_{\cT_i,1} = 
                \sup 
                \BR{
                    \inf\omega_{\cT,1} \in\R\::\:
                    \cT\in\T^\ast_i
                }.
        \end{equation}
    \item We associate nine other trees with physically shifted tops:
    \begin{equation*}
        \cT_{i,\ve}:=\BR{
            P\in\cP_{i-1}\::\:
            P\leq 
            \br{\ve,\bm{0}}\boxplus P_{\cT_i}
            =\br{\ve\boxplus\vI_{\cT_i}}
                \times\vomega_{\cT_i}}
        ,\quad
        \ve\in\BR{-1,0,1}^2.
    \end{equation*}
    \item For the next iteration, we set:
    \begin{equation*}
        \cP_i:=\cP_{i-1}\setminus \bigsqcup_{\ve\in \BR{-1,0,1}^2} \cT_{i,\ve},\quad
        \T_i:=\T_{i-1}\sqcup\BR{\cT_{i,\ve}}_{\ve\in \BR{-1,0,1}^2},
    \end{equation*}
    increase \(i\) by one, and return to the second step.
    \end{itemize}
When the algorithm terminates at \(i=n\), we set \(\T:=\T_n\). By construction, 
\begin{equation*}
    \nrm{f^{\br{+,j}}}_{\bmo_2\br{\cP\setminus\T}}=\nrm{f^{\br{+,j}}}_{\bmo_2\br{\cP_n}}\leq 1/2= \nrm{f^{\br{+,j}}}_{\bmo_2\br{\cP}}/2.
\end{equation*}
On the other hand, due to the direct association between \(\cT_i\) and \(\cT_{i,\ve}\), the following trivial relation holds:
\begin{equation}\label{eq_cE_2_sum_I_T}
    \sum_{\cT\in\T}\abs{\vI_\cT}=
    \sum_i\sum_{\ve\in\BR{-1,0,1}^2}\abs{\vI_{\cT_{i,\ve}}}=
    3^2\sum_i\abs{\vI_{\cT_i}}.
\end{equation}
To complete the proof, it remains to properly control the right-hand side. 

We first notice that \(\nrm{f^{\br{+,j}}}_{L^2\br{\cT_i}}\eqsim \abs{\vI_{\cT_i}}^{\frac{1}{2}}\); indeed, by the definition of \(\T^\ast_i\) and the fact that \(\cT_i\in \T^\ast_i\), we have
\begin{equation*}
    1/2< \nrm{f^{\br{+,j}}}_{L^2\br{\cT_i}}/\abs{\vI_{\cT_i}}^{\frac{1}{2}}
    \leq
    \nrm{f^{\br{+,j}}}_{\bmo_2\br{\cP}}=1.
\end{equation*}
The \(j=0\) case is easy to treat since in this situation we may apply \textbf{Observation \ref{obs_total_L2_energy_bd}} to finish the proof:
\begin{equation*}
    \sum_{\cT\in\T}\abs{\vI_\cT}
    \eqsim 
    \sum_i\abs{\vI_{\cT_i}}
    \eqsim 
    \sum_i
        \nrm{f^{\br{+,0}}}^2_{L^2\br{\cT_i}}
    \leq
    \nrm{f^{\br{+,0}}}^2_{L^2\br{\P}}
    \lesssim\nrm{f}^2_{L^2}.
\end{equation*}
The \(j=1\) case is more delicate and is similar to the argument employed in the proof of the classical Carleson theorem, which
relies on a Bessel-type inequality. To be more specific, we intend to apply \textbf{Lemma \ref{lem_bessel}}, and in order to do so, we first have to verify its hypothesis. We start by observing that for any \(\vI\times\vomega\in\cT_i\) we have
\begin{equation*}
    \nrm{f^{\br{+,1}}}_{\vI\times\vomega}/\abs{\vI}^{\frac{1}{2}}\leq 
    \nrm{f^{\br{+,1}}}_{\bmo_2\br{\cT_i}}\leq
    \nrm{f^{\br{+,1}}}_{\bmo_2\br{\cP}}=1.
\end{equation*}
This shows \(\nrm{f^{\br{+,1}}}_{\vI\times\vomega}\lesssim \vI^{\frac{1}{2}}\). Since we have also shown that \(\nrm{f^{\br{+,j}}}_{L^2\br{\cT_i}}\eqsim \abs{\vI_{\cT_i}}^{\frac{1}{2}}\) for \(j\in\BR{0,1}\),
it remains to check whether \(\BR{\cT_i}_i\) form a sequence of strongly disjoint \(\br{-,1}\)-trees.

To verify the strong disjointness condition, we take \(\vI_m\times\vomega_m\in\cT_m\) and \(\vI_n\times\vomega_n\in\cT_n\) with \(\vomega^{\br{+,1}}_m\cap\vomega^{\br{+,1}}_n\neq \varnothing\).
We proceed in two steps. We first show that \(\vomega_m \subsetneq \vomega_n\) forces the condition \(m<n\). We then demonstrate that \(m<n\) by itself implies \(3\vI_{\cT_m}\cap \vI_n=\varnothing\).

\begin{proof}[\textbf{Proof of \(\vomega_m \subsetneq \vomega_n\) implies \(m<n\)}] 
    Recall that \(\cT_m,\cT_n\) are \(\br{-,1}\)-trees. We have \(\vomega^{\br{-,1}}_{\cT_n}\subset \vomega^{\br{-,1}}_n\) and \(\vomega^{\br{-,1}}_{\cT_m}\subset \vomega^{\br{-,1}}_m\), which gives the following implication:
    \begin{equation}\label{eq_omega_size_to_inf_ineq}
        \vomega^{\br{+,1}}_m\subsetneq
        \vomega^{\br{+,1}}_n
        \implies
        \inf\omega_{\cT_m,1}\geq \inf\omega_{m,1} \geq \sup \omega^{\br{-,1}}_n > \inf \omega_{\cT_n,1}.
    \end{equation}
    Indeed, since \(\vomega^{\br{+,1}}_m \cap \vomega^{\br{+,1}}_n\neq\varnothing\), the assumption \(\vomega_m \subsetneq \vomega_n\) implies \(\vomega^{\br{+,1}}_m \subsetneq \vomega^{\br{+,1}}_n\) due to the dyadic structure of the interval systems.
    By \eqref{eq_omega_size_to_inf_ineq}, the inequality \(\inf\omega_{\cT_m,1}>\inf\omega_{\cT_n,1}\) with the maximality assumption \eqref{eq_cE_2_rightmost} forces that \(m<n\).
\end{proof}

\begin{proof}[\textbf{Proof of \(m<n\) implies \(3\vI_{\cT_m}\cap \vI_n=\varnothing\)}]
Since \(m<n\) and (from the above argument) \(\vomega^{\br{+,1}}_m\subset\vomega^{\br{+,1}}_n\) we must have \(\vomega_{\cT_m} \subset \vomega_n\).
Yet, by construction:
\begin{equation*}
    P_n \not\leq \br{\ve,\vnull} \boxplus P_{\cT_m},\quad\forall\: \ve\in\BR{-1,0,1}^2.
\end{equation*}
This forces that
\begin{equation*}
    \br{\ve\boxplus\vI_{\cT_m}} \cap \vI_n=\varnothing,\quad
    \forall \:\ve\in\BR{-1,0,1}^2.
\end{equation*}
Namely, \(3\vI_{\cT_m}\cap\vI_n=\varnothing\).
\end{proof}

We thus proved that \(\BR{\cT_i}_i\) are strongly disjoint. Returning now to \eqref{eq_cE_2_sum_I_T} we apply \textbf{Lemma \ref{lem_bessel}} to conclude
\begin{equation*}
    \sum_{\cT\in\T}\abs{\vI_\cT}\eqsim
    \sum_i\abs{\vI_{\cT_i}}
    \lesssim \nrm{f}_{L^2}^2 =\br{\nrm{f}_{L^2}/\nrm{f^{\br{+,j}}}_{\bmo_2\br{\cP}}}^2.
\end{equation*}
This completes the proof of \textbf{Lemma \ref{lem_cE_2_bd}}

\section{The general dimensional setting: Reduction of Theorem \ref{thm_main} to Theorem \ref{thm_Hdim_HF_st_disc}}\label{sec_H_dim}

In the remainder of the paper, we focus on proving the general dimensional case $D\geq 1$ as stated in \textbf{Theorem \ref{thm_main}}. Our main goal in this section is the reduction of \textbf{Theorem \ref{thm_main}} to \textbf{Theorem \ref{thm_Hdim_HF_st_disc}}.

\subsection{Preliminaries}

Departing from  definition \eqref{def_C_ast_V_K} and the assumptions in \textbf{Theorem \ref{thm_main}}, we notice that via the rotation and dilation symmetry of the formulation, we may assume without loss of generality that
\begin{itemize}

\item  \(\vv_0=\ve_0:=\br{1,0,\dots,0}\); 

\item \(\V:=\ve_0^{\perp}=\BR{0}\times\R^D<\R^{D+1}\).
\end{itemize}
Next, we write:
\begin{equation}\label{red}
    \left\{
    \begin{aligned}
    \vx:= &\br{x_0,x_1,\dots,x_D},\\
    \vxi:= &\br{\xi_0,\xi_1,\dots,\xi_D},\\
    \va:=
    &
    \br{0,a_1,\dots,a_D},
    \end{aligned}
    \right.
    \quad
    \vX\br{\vt}:=\br{\vt,\abs{\vt}^2}=
    \br{t_0,t_1,\dots,t_{D-1},\sum_{j=0}^{D-1} t^2_j}.
\end{equation}
After the initial reduction presented in \textsc{Section \ref{subsec_kernel_decomp}} we obtain the linearized operator
\begin{equation*}
    \cC\cR^{\br{\va}}_{\BR{0}\times\R^D,K}f\br{\vx}:=
    \int
        \widehat{f}\br{\vxi}
        \sum_{k=\uk}^{\ok\br{\vx}}
            \widehat{\mu_k}
            \br{\vxi-\va\br{\vx}}
            e\br{\vx^\top\vxi}
    d\vxi,
\end{equation*}
and naturally focus on the analysis of the resulting symbol
\begin{equation*}
    \widehat{\mu_k}\br{\vxi}:
    =\int
        \overline{
            e\br{
                \vxi^\top \vX\br{\vt} 
            }
        }
        \rho_k\br{\abs{\vt}/\sqrt{D}}
        K\br{\vt}
    d\vt.
\end{equation*}
Next, for expository reasons, we introduce two positive definite diagonal matrices
\begin{equation}\label{eq_scale_settings_Ddim}
    \vS:=
    4\ve_D\ve_D^\top+
    2\sum_{j=0}^{D-1}\ve_j\ve_j^\top
    ,\quad
    \vR:=\sum_{j=0}^D R_j \ve_j\ve_j^\top
        \in M_{\br{D+1}\times \br{D+1}}\br{\R}\,,
\end{equation}
and, for vector \(\vz:=\br{z_0,\dots,z_m}\in\R^{m+1}\) and indices \(0\leq i<j\leq m\), the convention
\begin{equation}\label{eq_vec_hat}
    \widehat{\vz}_i:=\br{z_0,\dots,z_{i-1},z_{i+1},\dots,z_m}\in\R^m,\quad
    \widehat{\vz}_{i,j}:=
    \br{z_0,\dots,z_{i-1},z_{i+1},\dots,z_{j-1},z_{j+1},\dots,z_m}\in\R^{m-1}.
\end{equation}
The cone decomposition \eqref{eq_cone_dec} introduced in \textsc{Section \ref{subsec_symb_decomp}} has the following direct generalization:
\begin{equation*}
    \1_{\R^{D+1}\setminus\BR{0}}\br{\vxi}
    = 
    \int_{\R\setminus\BR{0}}\left\{
    c_0\,
    \phi\br{\frac{2\xi_0}{\lambda}-\frac{1}{2}}
    \varphi^{\otimes D}
    \br{\frac{\widehat{\vxi}_0}{\lambda}}
    + 
    \sum_{j=1}^D c_j\,\tarphi\br{\frac{\xi_0}{\lambda}}\tpsi\br{\frac{\xi_j}{\lambda}}
    \varphi^{\otimes D-1}\br{\frac{\widehat{\vxi}_{0,j}}{\lambda}}\right\}
    \frac{d\lambda}{\abs{\lambda}}.
\end{equation*}
Letting
\begin{equation}\label{eq_Hdim_chi_j_def}
    \chi_0\br{\vxi}:= 
    \phi\br{2\xi_0-\frac{1}{2}}
    \varphi^{\otimes D}
    \br{\widehat{\vxi}_0}
    ,\quad
    \chi_j\br{\vxi}:=
    \tarphi\br{\xi_0}
    \tpsi\br{\xi_j}
    \varphi^{\otimes D-1}\br{\widehat{\vxi}_{0,j}}
    ,
\end{equation}
we obtain the following symbol decomposition analog to \eqref{eq_symb_lamb_dec}:
\begin{equation}\label{eq_Hdim_symb_lamb_dec}
    \widehat{\mu}_k=
    \sum_{j=0}^D
    c_j
    \int_{\R\setminus\BR{0}}
        \widehat{\mu^{\br{j}}_{k,\lambda}}
    \frac{d\lambda}{\abs{\lambda}}
    ,\quad
    \widehat{\mu^{\br{j}}_{k,\lambda}}\br{\vxi}:=
    \widehat{\mu}_k\br{\vxi}
    \chi_j\br{\vS^k\vxi/\lambda}.
\end{equation}
Symmetric to \eqref{eq_C_V_decomp} and \eqref{eq_cC_j_lambda_def}, we now conclude
\begin{equation}\label{eq_Hdim_C_V_decomp}
    \cC\cR^{\br{\va}}_{\BR{0}\times\R^D} f\br{\vx}=
    \sum_{j=0}^D
    c_j
    \int_{\R\setminus\BR{0}}
        \cC\cR^{\br{\va}}_{\lambda,j} f\br{\vx}
    \frac{d\lambda}{\abs{\lambda}}\,,
\end{equation}
with each \(\cC\cR^{\br{\va}}_{\lambda,j}\) defined by
    $$\cC\cR^{\br{\va}}_{\lambda,j} f\br{\vx}:=
    \int \widehat{f}\br{\vxi}\left(\sum_{k=\uk}^{\ok\br{\vx}}\widehat{\mu^{\br{j}}_{k,\lambda}}\br{\vxi-\va\br{\vx}}\right)
    e\br{\vxi\cdot\vx}
    d\vxi.$$

\subsection{Reduction of Theorem \ref{thm_main} to Theorems \ref{prop_Hdim_4cases} and \ref{prop_Hdim_4cases1}}\label{Gendimred}

The proof of \textbf{Theorem \ref{thm_main}} follows a similar structure to the analog result in the planar setting. Indeed, proceeding as in \textsc{Section \ref{subsec_pf_main_planar}}, we decompose $\cC\cR^{\br{\va}}_{\BR{0}\times\R^D}$ into three sub-operators:
\begin{itemize}
\item a low-frequency component,

\item a high-frequency non-stationary component,

\item a high-frequency stationary component,
\end{itemize}
and then reduce our \textbf{Theorem \ref{thm_main}} to the following two results:

\begin{theorem}\label{prop_Hdim_4cases}
    There exists a universal constants \(\epsilon>0\) such that for any integer \(\uk\in\Z\), \(j\in\{0,\dots,D\}\), measurable functions \(\va:\R^{D+1}\to\BR{0}\times\R^D\) and \(\ok:\R^{D+1}\to\Z\cap\Br{\uk,\infty}\), and \(1<p<\infty\), the following hold:
    \begin{itemize}
        \item \textsf{Low-frequency} component---polynomial decay: given \(\abs{\lambda}\lesssim 1\) one has
        \begin{equation}\label{eq_prop_Hdim_LF}
            \nrm{
            \cC\cR^{\br{\va}}_{\lambda,j}f
            }_{L^p}
            \underset{p}{\lesssim} \abs{\lambda}^\epsilon  
            \nrm{f}_{L^p}.
        \end{equation}
        
        \item \textsf{High-frequency non-stationary} component---polynomial decay: given \(\abs{\lambda}\gtrsim 1\) one has
        \begin{equation}\label{eq_prop_Hdim_HF_nst}
                \nrm{
                \cC\cR^{\br{\va}}_{\lambda,D}f
                }_{L^p}
                \underset{p}{\lesssim}  \abs{\lambda}^{-\epsilon} \nrm{f}_{L^p}.
        \end{equation}

        \item \textsf{High-frequency stationary} component---tame log growing estimate: given \(\abs{\lambda}\gtrsim 1\) and  \(j<D\) one has
        \begin{equation}\label{eq_prop_Hdim_HF_log}
                \nrm{
                \cC\cR^{\br{\va}}_{\lambda,j}f
                }_{L^p}
                \underset{p}{\lesssim} \log^2\br{e+\abs{\lambda}} \nrm{f}_{L^p}.
        \end{equation}

        \end{itemize}
\end{theorem}

\begin{theorem}\label{prop_Hdim_4cases1}
   Under the same notations as above and assuming  \(\abs{\lambda}\gtrsim 1\) and  \(0\leq j<D\), the \textsf{high-frequency stationary} component obeys the following $L^2$-polynomial decay estimate: 
   \begin{equation}\label{eq_prop_Hdim_HF_st}
                \nrm{
                \cC\cR^{\br{\va}}_{\lambda,j}f
                }_{L^{2,\infty}} 
                \lesssim \abs{\lambda}^{-\epsilon} \nrm{f}_{L^2}.
   \end{equation}  
\end{theorem}

Once we assume \textbf{Theorems \ref{prop_Hdim_4cases}} and \textbf{\ref{prop_Hdim_4cases1}}, the proof of \textbf{Theorem \ref{thm_main}} follows arguments identical to those presented in \textsc{Section \ref{subsec_pf_main_planar}}. We thus omit the standard details.

It thus remains to prove the four estimates \eqref{eq_prop_Hdim_LF}, \eqref{eq_prop_Hdim_HF_nst}, \eqref{eq_prop_Hdim_HF_log}, and \eqref{eq_prop_Hdim_HF_st}. Due to the symmetry of the formulation, throughout the remainder of the paper, it suffices only to consider the case  \(\lambda >0\).

\subsection{Reduction of Theorem \ref{prop_Hdim_4cases} to Theorem \ref{thm_log_shift_est_D_dim}: shifted superposed Carleson operators}\label{ssco}
We essentially follow the same circle of ideas as in \textsc{Section \ref{subsec_class}}, that is, the three relations \eqref{eq_prop_Hdim_LF}, \eqref{eq_prop_Hdim_HF_nst}, and \eqref{eq_prop_Hdim_HF_log} will be proven via estimates on the corresponding symbols and an application of a high-dimensional generalization of \textbf{Theorem \ref{thm_log_shift_est}}, which we state below:
\begin{theorem}[\textsf{Control over the  shifted superposed Carleson operator}]\label{thm_log_shift_est_D_dim}
    Let \(j\in\{0,\dots,D\}\) and \(\vS\) be as in \eqref{eq_scale_settings_Ddim}, and \(\va,\uk,\ok\) be as in \textbf{Theorem \ref{prop_Hdim_4cases}}. For any sequence of coefficients \(\br{\varepsilon_k}_{k\in\Z} \in \ell^\infty\br{\C}\), vector \(\vu\in\R^{D+1}\), and positive definite diagonal matrix \(\vR\in M_{\br{D+1}\times\br{D+1}}\br{\R}\), the following operator\footnote{Recall \textbf{Definition \ref{def_sym_ops}}.}:
    \begin{equation*}
        \cC^{\br{\va,\varepsilon,\uk,\ok}}_{\vR,j,\vu}f\br{\vx}
        := 
        \sum_{k=\uk}^{\ok\br{\vx}}
            \varepsilon_k
        \int
        \widehat{f}\br{\vxi}
        \br{
             \Mod_{-\vu}
            \chi_j}
        \br{
            \vS^k \vR^{-1} \br{\vxi -\va\br{\vx}}
        }
        e\br{\vxi^\top\vx}
        d\vxi
    \end{equation*}
    satisfies the following estimate:
    \begin{equation*}
        \nrm{
            \cC^{\br{\va,\varepsilon,\uk,\ok}}_{\vR,j,\vu}f
        }_{L^p}\underset{p}{\lesssim}
        \log^2\br{e+\abs{\vu}}
        \nrm{\epsilon}_{\ell^\infty}
        \nrm{f}_{L^p},\quad
        p\in\br{1,\infty}.
    \end{equation*}
    Specifically, the implicit constant does not depend on the constant \(\uk\in\Z\), the matrix \(\vR\), and the measurable functions \(\va,\ok\). In fact, for this statement, the restriction \(\va:\R^{D+1}\to \BR{0}\times\R^D\) is not required.
\end{theorem}

\textbf{Theorem \ref{thm_log_shift_est_D_dim}} can be proven via similar arguments to those employed for proving \textbf{Theorem \ref{thm_log_shift_est}}.
Thus, we leave the standard adaptations to the interested reader.

With \textbf{Theorem \ref{thm_log_shift_est_D_dim}} verified, we decompose \(\cC\cR^{\br{\va}}_{\lambda,j}\) into a superposition of operators of the form \(\cC^{\br{\va,\varepsilon,\uk,\ok}}_{\vR,j,\vu}\) with suitable weights.
This will be achieved via a dual decomposition analogous to \eqref{eq_sym_4ier_sum} and \eqref{eq_sym_int_dec}. To be precise, we have
\begin{itemize}
\item \textsf{a frequency decomposition} which via a Fourier series argument relative to a well-chosen fixed lattice \(\Z^2:=\prod_{j=0}^D\theta_j\Z\) with \(\theta_j\eqsim 1\) gives
\begin{equation}\label{eq_Hdim_sym_4ier_sum}
    \widehat{\mu^{\br{j}}_{k,\lambda}}
    \br{\vxi}
    =
    \sum_{\vu\in \Z^2}
    m^{\br{j}}_{k,\lambda,\vu}\,
    \overline{e\br{
        \vu^\top\vS^k\vxi/\lambda
    }}\,
    \chi_j\br{\vS^k\vxi/\lambda}.
\end{equation}
\item \textsf{a spatial decomposition} based on the integral definition \eqref{eq_mu_k_def} in the form
\begin{equation}\label{eq_Hdim_sym_int_dec}
    \widehat{\mu^{\br{j}}_{k,\lambda}}\br{\vxi}
    =\int
        2^{kD} K\br{2^k \vt}\,
        \overline{e\br{
            \vX^\top\br{\vt}\vS^k \vxi
        }}\,
        \chi_j\br{\vS^k\vxi/\lambda}\,
    \rho\br{\abs{\vt}/\sqrt{D}}
    d\vt.
\end{equation}
\end{itemize}
The above two decompositions motivate the analog of \eqref{eq_cC_2ways_2_decomp} given by
\begin{equation}\label{eq_Hdim_cC_2ways_2_decomp}
    \cC\cR^{\br{\va}}_{\lambda,j}= \sum_{\vu\in \Z^2}
    \cC\cR^{\br{\va}}_{\lambda,j,\vu}= 
    \int
    \cC\cR^{\br{\va}}_{\lambda,j,\vt}\,\,
    \rho\br{\abs{\vt}/\sqrt{D}}
    d\vt,
\end{equation}
with the involved operators defined by
\begin{equation}\label{eq_Hdim_cC_4ier_dec}
    \cC\cR^{\br{\va}}_{\lambda,j,\vu}f\br{\vx}:= 
    \sum_{k=\uk}^{\ok\br{\vx}} 
    m^{\br{j}}_{k,\lambda,\vu}
    \cdot e\br{\va^\top\br{\vx}\vS^k\vu/\lambda}
    \cdot
    \int 
        \chi_j \br{\vS^k\vxi/\lambda}
        \widehat{f}\br{\vxi}
    e\br{\vxi^\top \br{\vx- \vS^k\vu/\lambda}}d\vxi,
\end{equation}
\begin{equation}\label{eq_Hdim_cC_int_dec}
    \cC\cR^{\br{\va}}_{\lambda,j,\vt}f\br{\vx}:= 
    \sum_{k=\uk}^{\ok\br{\vx}} 
    2^{kD}K\br{2^k \vt}
    \cdot e\br{\va^\top\br{\vx}\vS^k\vX\br{\vt}}
    \cdot
    \int 
        \chi_j \br{\vS^k\vxi/\lambda}
        \widehat{f}\br{\vxi}
    e\br{\vxi^\top \br{\vx- \vS^k\vX\br{\vt}}}d\vxi.
\end{equation}
By setting
\begin{equation*}
    \varepsilon^{\br{j}}_k:=
    m^{\br{j}}_{k,\lambda,\vu},\quad
    \varepsilon_k\br{\vt}:=2^{kD}K\br{2^k \vt},
\end{equation*}
and using $\vid:=\sum_{j=0}^D\ve_j\ve^\top_j\in M_{\br{D+1}\times\br{D+1}}\br{\R}$ we obtain the following two identities:
\begin{equation*}
    \cC\cR^{\br{\va}}_{\lambda,j,\vu}=
    \cC^{\br{\va,\varepsilon^{\br{j}},\uk,\ok}}_{\lambda \vid,j,\vu},\quad
    \cC\cR^{\br{\va}}_{\lambda,j,\vt}=\cC^{\br{\va,\varepsilon\br{\vt},\uk,\ok}}_{\lambda\vid,j,\lambda\vX\br{\vt}}\,.
\end{equation*}
As expected, we control the coefficient \(m^{\br{j}}_{k,\lambda,\vu}\) via the high-dimensional version of Mikhlin-type symbol condition derived in \textbf{Lemma \ref{lem_mikhlin_symb}} which takes the following form:
\begin{lemma}[\textsf{Mikhlin-type symbol condition}]\label{lem_Hdim_mikhlin_symb}
Let \(k\in\Z\) and \(j\in \{0,\dots,D\}\) and recall \eqref{eq_Hdim_chi_j_def} and \eqref{eq_Hdim_symb_lamb_dec}. For multi-indices \(\valpha\in \BR{0,\dots, \lfloor\frac{A-1}{D}\rfloor }^{D+1}\), the following hold:
    \begin{itemize}
        \item \textsf{Low-frequency} case: given \(\abs{\lambda}\lesssim 1\),
        \begin{equation}\label{eq_Hdim_LF_par}
            \abs{\partial^\valpha_\vxi
            \widehat{\mu^{\br{j}}_{k,\lambda}}\br{\vS^{-k}\vxi}
            }
            \lesssim
            \lambda
            \cdot
            \abs{\vxi}^{-\abs{\valpha}}.
        \end{equation}
        \item \textsf{High-frequency non-stationary} case: given \(\abs{\lambda}\gtrsim 1\),
        \begin{equation}\label{eq_Hdim_HF_nst_par}
            \abs{\partial^\valpha_\vxi
            \widehat{\mu^{\br{D}}_{k,\lambda}}\br{\vS^{-k}\vxi}
            }
            \lesssim
            \lambda^{-1}
            \cdot
            \abs{\vxi}^{-\abs{\valpha}}.
        \end{equation}
    \end{itemize}
\end{lemma}
The proof of \textbf{Lemma \ref{lem_Hdim_mikhlin_symb}} is similar to that of \textbf{Lemma \ref{lem_mikhlin_symb}} with some slight modifications. Specifically, when treating the second item, one needs to address some subtleties. Concretely, we write
\begin{equation*}
    \rho\br{\abs{\vt}/\sqrt{D}}=\sum_{j=0}^{D-1}
    \rho_j\br{\vt},\quad
    \supp\rho_j \subset 
    \BR{
        \vt\in\R^D
    \::\:
        \abs{\vt}\lesssim 1\leq 3 \abs{t_j}
    }\,,
\end{equation*}
and decompose the normalized symbol as follows
\begin{equation*}
    \partial^\valpha_\vxi
    \widehat{\mu_k}\br{\vS^{-k}\vxi}
    =
    \sum_{j=0}^{D-1}
    \int
        \overline{e\br{\vxi^\top\vX\br{\vt}}}\,
    2^{kD}\, K\br{2^k \vt}\,
    \vt^{\tilde{\valpha}}\rho_j\br{\vt}
    d\vt=:
    \sum_{j=0}^{D-1}I_j,
\end{equation*}
where in the above, we set the multi-index
\begin{equation*}
    \tilde{\valpha}:=\br{\alpha_0+2\alpha_D,\dots,\alpha_{D-1}+2\alpha_D}.
\end{equation*}
Let now \(\vxi\in\supp \left(\Dil^\infty_\lambda\chi_D\right)\) and \(j\in \{0,\dots,D-1\}\). By \eqref{eq_tar_tpsi_supp} and \eqref{eq_Hdim_chi_j_def}
we deduce the following conditions:
\begin{equation*}
    \abs{\xi_j}\leq 0.26\lambda,\quad
    \xi_D\geq 0.49\lambda.
\end{equation*}
As a result, we may now integrate by parts each individual expression \(I_j\):
\begin{equation}\label{eq_Hdim_symb_by_part}
    \abs{I_j}
    \underset{\valpha}{\eqsim}
    \abs{
    \int
    \overline{e\br{\vxi^\top\vX\br{\vt}}}
        \br{
            \frac{\partial}{\partial t_j}
            \cdot
            \frac{1}{\xi_j+2\xi_D t_j}\cdot
        }^{1+\abs{\valpha}}
        2^{kD} K\br{2^k \vt}
    \vt^{\tilde{\valpha}}\rho_j\br{\vt}
    d\vt
    },
\end{equation}
and using \(\abs{\xi_j+2\xi_D t_j}\geq 2\abs{\xi_D} \cdot \abs{t_j}-\abs{\xi_j}\geq \br{\frac{2\times 0.49}{3}-0.26}\lambda \eqsim\lambda\) we conclude
\begin{equation*}
    \abs{\partial^\valpha_\vxi
    \widehat{\mu_k}\br{\vS^{-k}\vxi}}
    \leq
    \sum_{j=0}^{D-1}
    \abs{I_j}
    \underset{\valpha}{\lesssim}
    \sum_{j=0}^{D-1}
    \lambda^{-1-\abs{\valpha}}
    \underset{\valpha}{\eqsim}
    \lambda^{-1}
    \abs{\vxi}^{-\abs{\valpha}}.
\end{equation*}
\subsection{Reduction of Theorem \ref{prop_Hdim_4cases1} to Theorem \ref{thm_Hdim_HF_st_disc}: non-classical phenomenon}\label{redmain}
To perform the discretization process shown in \textsc{Section \ref{subsec_Freq_disc}},
we first extend \textbf{Definition \ref{def_W_n_Wk}}:
\begin{equation}\label{eq_Hdim_W_n_Wk}
    \W:=\W\br{\lambda}:=\bigsqcup_{k\in\Z}\W_k,\quad
    \W_k:=\W_k\br{\lambda}:=
    \BR{
        \br{0,\vw}\boxplus 
        \vS^{-k}\br{\bR{-\lambda/2,\lambda/2}\times\bR{0,\lambda}^D}
        \::\:
        \vw\in \Z^D
    }.
\end{equation}
However, due to the presence of \(\tpsi\br{\xi_j}\) in \(\chi_j\br{\vxi}\), we have to properly rephrase \eqref{eq_chi_01_nice} for \(0<j<D\). To replace \(\tpsi\br{z}\) with \(\psi\br{z}:=\varphi\br{2z-1}\), we invoke again a Fourier series trick:
\begin{equation}\label{eq_tspi_2_psi}
    \tpsi=\sum_{u\in\Z}
    m_u \Mod_{-u}\psi,\quad
    \abs{m_u}\underset{N}{\lesssim}\ang{u}^{-N}.
\end{equation}
Consequently, the analog of
\(\cC^{\br{\va}}_{\Delta,\lambda,0}\) for \(j<D\) shall include a shifting parameter \(\vu\in\R^{D+1}\) as given by
\begin{equation*}
    \cC^{\br{\va}}_{\Delta,\vu,\lambda,j}
    f\br{\vx}:=
    \sum_{k=\uk}^{\ok\br{\vx}}
        \sum_{\vomega\in\W_k\br{\lambda}} 
        \1_{\vomega^{\br{-,j}}}\br{\va\br{\vx}}
        \int
            \widehat{\mu_k}\br{\vxi-\va\br{\vx}}
        \phi_{\vomega^{\br{+,j}}}\br{\vxi}
        \widehat{f}\br{\vxi}
        e\br{\br{\vx-\vS^k\vu}^\top\vxi}
        d\vxi.
\end{equation*}
The original operator \(\cC\cR^{\br{\va}}_{\lambda,j}\) can thus be recovered via the following formulas
\begin{equation*}
    \cC\cR^{\br{\va}}_{\lambda,0}f\br{\vx}:=
    c_0
    \fint_{\lambda \vS^{-\uk}\br{\BR{0}\times\bR{0,1}^D}}
        \Mod_\vzeta
        \cC^{\br{\va-\vzeta}}
        _{\Delta,\vnull,\lambda,j}
        \Mod_{-\vzeta}
            f\br{\vx}
    d\cH^D\br{\zeta};
\end{equation*}
\begin{equation*}
    \cC\cR^{\br{\va}}_{\lambda,j}f\br{\vx}:=
    c_1
    \sum_{u\in\Z}m_u
    \fint_{\lambda \vS^{-\uk}\br{\BR{0}\times\bR{0,1}^D}}
        \Mod_\vzeta
        \cC^{\br{\va-\vzeta}}
        _{\Delta,u\ve_j/\lambda,\lambda,j}
        \Mod_{-\vzeta}
            f\br{\vx}
    d\cH^D\br{\zeta},\quad 0<j<D.
\end{equation*}
Due to the fast decay of the coefficient \(m_u\) in \eqref{eq_tspi_2_psi}, \eqref{eq_prop_Hdim_HF_st} follows once we prove the following statement:
\begin{theorem}[\textsf{Model sum estimate for stationary-phase contribution with shift}]\label{thm_Hdim_HF_st_disc}
There exists \(\epsilon>0\) such that for any  \(\lambda\gtrsim 1\), $0\leq j<D$, \(\uk\in\Z\), \(\vu\in\R^{D+1}\) and any measurable functions \(f\in L^2\br{\R^2}\) and \(\va:\R^{D+1}\to\BR{0}\times\R^D\), \(\ok:\R^{D+1}\to\Z\cap\Br{\uk,\infty}\), the following estimate holds:
    \begin{equation*}
        \nrm{
            \cC\cR^{\br{\va}}
            _{\Delta,\vu,\lambda,j}
            f
        }_{L^{2,\infty}}
        \lesssim
        \lambda^{-\epsilon}\log^2\br{e+\lambda\br{1+\abs{\vu}}}
        \nrm{f}_{L^2}.
    \end{equation*}
\end{theorem}
The rest of the paper will be dedicated to the proof of \textbf{Theorem \ref{thm_Hdim_HF_st_disc}}. 
\medskip

Recalling now \textbf{Definition \ref{def_Freq_proj} and \ref{def_A_omega}} we reformulate the defining formula for \(\cC^{\br{\va}}_{\Delta,\vu,\lambda,j}\) by
\begin{equation}\label{eq_def_Hdim_cC_Delta_lambda_j}
    \cC^{\br{\va}}_{\Delta,\vu,\lambda,j} f\br{\vx}:=
    \sum_{k\in\Z}
        \sum_{\vomega\in\W_k}
            \1_{A^{\br{-,j}}_\vomega}\br{\vx}
            \int
                \pi^{\br{+,j}}_\vomega f\br{\vx-\vtau}e\br{\va\br{\vx}^\top \vtau}
            d\mu_k\br{\vtau-\vS^k \vu}.
\end{equation}
Thus, in the language introduced in \textbf{Definition \ref{def_cC_nu}}, we have established
\begin{equation*}
    \cC_{\lambda,j}:=\cC^{\br{\va}}_{\nu,\lambda,j}=\cC^{\br{\va}}_{\Delta,\vu,\lambda,j}
    ,\quad
    d\nu_k\br{\vtau}:=d\mu_k\br{\vtau-\vS^k \vu}.
\end{equation*}

At this stage, we are ready to perform the overall scheme presented in \textsc{Section \ref{sec_formal_tf_ana}}.
The first step is to perform a Gabor/wave-packet decomposition analog to \eqref{eq_planar_tf_proj} for the frequency truncated function \(\pi^{\br{+,j}}_\vomega f\) with respect to the following system of $(D+1)-$intervals:
\begin{equation}\label{spinterv}
    \I:=\bigsqcup_{k\in\Z}\I_k,\quad 
    \I_k:=\BR{
        \valpha\boxplus \vS^k \Br{0,1}^{D+1}
        \::\:
        \valpha\in\Z^{D+1}
    }.
\end{equation}
\medskip

Once at this point, we mention the following: 

\begin{remark}[\textsf{Dual high-resolution analysis within the low-resolution analysis of the (HS) component}]\label{dualhres}
In contrast with the one-dimensional case where the treatment of the (HS) component involves a single high-resolution analysis dictated by the LGC-method, in the higher-dimensional setting, we will need a \textsf{dual high-resolution analysis} in order to address new issues emerging in the \(0<j<D\) regimes for the (HS) component.\footnote{For more on these, see the introductory paragraph of the next section.} This requires a redesign of the wave-packet decomposition, thus motivating the following: 
\end{remark}

\begin{observation}
Given \(\phi_1,\phi_2\in \mathcal{S}\br{\R^{D+1}}\) with \(\nrm{\phi_k}_{L^2}=1\), define the following functions:
    \begin{equation*}
        \chi^\Phi:=
        \1_{\Br{1/2,1/2}^{D+1}}\ast
        \abs{\phi_1}^2\ast\abs{\phi_2}^2,\quad
        \chi^{\phi_1}:=
        \1_{\Br{1/2,1/2}^{D+1}}\ast
        \abs{\phi_2}^2,\quad
        \chi^{\phi_2}:=
        \1_{\Br{1/2,1/2}^{D+1}}\ast
        \abs{\phi_1}^2.
    \end{equation*}
By design, the following identities hold:
\begin{equation}\label{eq_Gabor_part_unity}
    \chi^\Phi=\abs{\phi_k}^2\ast \chi^{\phi_k},\quad
    1=\sum_{\vz\in\Z^{D+1}}\Tr_\vz \chi^\Phi=\sum_{\vz\in\Z^{D+1}}\Tr_\vz\chi^{\phi_k}
    ,\quad k\in\{1,2\}
    .
\end{equation}
Moreover,  by \textbf{Proposition \ref{prop_Gabor_dec_conti}}, for any \(g\in L^2\br{\R^{D+1}}\) we have
\begin{align*}
    \chi^\Phi g= 
    \br{\chi^{\phi_k}\ast\abs{\phi_k}^2} g
    = & \int \ang{g, \Mod_\vxi \Tr_\vx \phi_k} \Mod_\vxi\Tr_\vx \phi_k \cdot \chi^{\phi_k}\br{\vx}d\vx d\vxi\\
    = & \int \ang{g, \Tr_\vx \Mod_\vxi \phi_k} \Tr_\vx \Mod_\vxi \phi_k \cdot \chi^{\phi_k}\br{\vx}d\vx d\vxi,\quad
    k\in\{1,2\}.
\end{align*}
Deduce this way, that \(\chi^\Phi g\) can be represented in two different ways---and at two distinct scales---as a weighted Gabor superposition physically localized around the unit interval \(\bR{-1/2,1/2}^{D+1}\). Along the way, we notice that one may naturally generalize the above construction to allow \(m>2\) many representations. 
\end{observation}
We may apply the above observation to \(\pi^{\br{+,j}}_\vomega f\) to physically localize the function while preserving the frequency localization. We choose the following two functions:
\begin{equation*}
    \phi_1:=\varphi^{\br{C}}:=
    \br{\widecheck{\phi}/\nrm{\phi}_{L^2}}^{\otimes (D+1)},\quad \phi_2:=\varphi^{\br{L}}:=\Dil^2_{1/\sqrt{\lambda}}\phi_1
\end{equation*}
and, for convenience, write \(\chi^{\br{\ast}}:=\chi^{\varphi^{\br{\ast}}}\) for \(\ast\in\{C,L\}\).
 
We are now ready to implement a time-frequency localization adapted to the following system of tiles:
\begin{equation}\label{superHeis}
    \P:=\P\br{\lambda}:=\bigsqcup_{k\in\Z}\P_k,\quad
    \P_k:=
    \P_k\br{\lambda}:=\BR{
        \vI\times\vomega\::\:
        \br{\vI,\vomega}\in\I_k\times \W_k\br{\lambda}
    }.
\end{equation}
For this, for each \(P:=\vI\times\vomega\in\P\), we define
\begin{equation}\label{eq_def_tf_proj_dual_scale}
    f^{\br{+,j}}_P:=f^{\br{+,j}}_{\vI\times\vomega}:=\chi^\Phi_\vI \pi^{\br{+,j}}_\vomega f.
\end{equation}
and notice that by \eqref{eq_Gabor_part_unity} and rescaling, we have \eqref{eq_time_decomp} as a direct consequence. 

More importantly, for each tile \(\vI\times\vomega\in\P_k\), we have the following dual-scale identity:
\begin{equation}\label{eq_Hdim_wp_CL_decomp}
    f^{\br{+,j}}_{\vI\times\vomega}=
    \int
    \int_{\R^{D+1}} \ang{\pi^{\br{+,j}}_\vomega f, \varphi^{\br{\ast}}_{k,\vx,\vxi}} \varphi^{\br{\ast}}_{k,\vx,\vxi} \chi^{\br{\ast}}_{\vI}\br{\vx}d\vx d\vxi,\quad
    \varphi^{\br{\ast}}_{k,\vx,\vxi}:=\Tr_\vx\Mod_{\vxi}\Dil^2_{\vS^k}
    \varphi^{\br{\ast}},\quad 
    \ast\in\{C,L\}.
\end{equation}
Correspondingly, we introduce the analog of \eqref{eq_wp_nrm}:
\begin{equation}\label{eq_Hdim_wp_nrm}
    \nrm{f^{\br{+,j}}}_{\vI\times\vomega}:=
    \nrm{\pi^{\br{+,j}}_\vomega f\br{\vx}}_{L^2\br{\chi^\Phi_\vI\br{\vx}d\vx}}=
    \nrm{
        \ang{\pi^{\br{+,j}}_\vomega f, \varphi^{\br{\ast}}_{k,\vx,\vxi}}
    }_{L^2\br{\chi^{\br{\ast}}_\vI\br{\vx}d\vx d\xi}},\quad \ast\in\{C,L\}.
\end{equation}
In combination with \textbf{Definition \ref{def_sing_tile_bili}}, we derive the following model sum representation:
\begin{equation*}
    \abs{\ang{\cC^{\br{\va}}_{\Delta,\vu,\lambda,j}f,g}}
    =
    \abs{\ang{\cC_{\lambda,j}f,g}}
    \leq
    \sum_{k\in\Z}
    \sum_{\vI\times\vomega\in\P_k}
        \abs{\ang{
            \cC_k 
                f^{\br{+,j}}_{\vI\times\vomega},
            \1_{A^{\br{-,j}}_\vomega}g
        }}
    =\nrm{
        \Lambda^{\br{j}}_P\br{f,g}
    }_{\ell^1\br{P\in\P}}
\end{equation*} 
with the same convention as the one in \textbf{Definition \ref{def_cC_nu}}, that is
\begin{equation*}
    \cC_k f\br{\vx}:=\cC_{\nu_k}^{\br{\va}} f\br{\vx}:=
    \int f\br{\vx-\vtau}e\br{\va\br{\vx}^\top \vtau} d\nu_k\br{\vtau},\quad 
    d\nu_k\br{\vtau}:=
    d\mu_k\br{\vtau-\vS^k\vu}.
\end{equation*}
Now \textbf{Theorem \ref{thm_Hdim_HF_st_disc}} follows once we show that for all finite measurable set \(E\subset \R^{D+1}\) the following holds:
\begin{equation}\label{eq_thm_Hdim_HF_st_mod}
    \nrm{\Lambda^{\br{j}}_P\br{f,\1_E}}_{\ell^1\br{P\in\P\br{\lambda}}}\lesssim
    \lambda^{-\epsilon} 
    \log^2\br{e+\lambda\br{1+\abs{\vu}}}
    \nrm{f}_{L^2}\abs{E}^{\frac{1}{2}}.
\end{equation}

\section{Proof of Theorem \ref{thm_Hdim_HF_st_disc} subject to Proposition \ref{thm_Hdim_sing_tile_ests} and Lemma \ref{lem_Hdim_energy_sel}}\label{sec_pf_Hdim_HF_st_mod}

In this section, our main focus will be the proof of \eqref{eq_thm_Hdim_HF_st_mod} assuming that both \textbf{Proposition \ref{thm_Hdim_sing_tile_ests}} and \textbf{Lemma \ref{lem_Hdim_energy_sel}} hold. As a first comment, we notice that due to the permutation symmetry among the variables \(x_1,\dots,x_{D-1}\), it is enough to limit our discussion to only two cases represented by \(j=0\) and \(j=1\).
\smallskip

In what follows, we elaborate on
\begin{itemize}
\item several arguments in the spirit of those in \textsc{Section \ref{sec_pf_HF_st_disc}}, especially when treating the  \(j=0\) case;
\item the required modifications for treating the more difficult \(j=1\) case.
\end{itemize}

More concretely, departing from \textbf{Assumption \ref{ass_tf_proj_emb}}, one may easily check that the setting \eqref{eq_Hdim_wp_nrm}
produces the desired properties by a line-by-line adaptation of the argument presented in \textsc{Section \ref{subsec_ver_tf_proj_emb}} with just trivial modifications, such as replacing \(\br{\cdot}^{\br{+,0}}\) by \(\br{\cdot}^{\br{+,j}}\) and \(\1_\vI\) by \(\chi^{\br{\ast}}_\vI\). 

However, due to a subtlety that only arises in the \(j=1\) case, a refined version of \eqref{eq_ass_wp_t_in_prod_V} for the \(j=1\) case is needed in order to preserve the high-resolution frequency structure. Moreover, besides the single tile estimate analogous to \eqref{eq_sing_tile_est} with \(\nrm{\Lambda}_c \lesssim \lambda^{-\sigma}\) for some small universal constants \(c,\sigma>0\), we will need to investigate a special case where the estimate with \(c=1/2\) produces \(\lambda^{-\sigma}\) decay. In order to perform this analysis, we will introduce frequency sub-systems, define a generalization of energy \textbf{Definition \ref{def_L2_BMO_energy}}, and prove a more refined version of the single tree estimate generalizing both \textbf{Lemma \ref{lem_-_tree} and \ref{lem_gen_tree_c_half}}.

The final step in this required higher-dimensional adaptation will be presented in \textsc{Section \ref{sec_pf_gen_energy_est}}: therein, in order to complete the proof of \eqref{eq_thm_Hdim_HF_st_mod} and hence \textbf{Theorem \ref{thm_Hdim_HF_st_disc}}, we will have to modify \textbf{Definition \ref{def_out_L2_infty_energy_emb}} and develop an energy selection algorithm for \(\P:=\P\br{\lambda}\) tile-sets with volume \(V_\P=\lambda^{D+1}\gg 1\).

\subsection{Single tile estimate: statement}

In this subsection, we only state the refined single tile estimate that we need for the general case $D\geq 1$:

\begin{proposition}\label{thm_Hdim_sing_tile_ests}
    There are universal constant \(c,\sigma>0\) such that for \(j\in\{0,1\}\), \(k\in\Z\), and \(\vI\times\vomega\in\P_k:=\P_k\br{\lambda}\), we have the following single tile estimate:
    \begin{equation}\label{eq_thm_Hdim_sing_tile_est_c}
        \abs{\Lambda_{\vI\times\vomega}^{\br{j}}\br{f,g}}\lesssim 
        \lambda^{-\sigma}
        \mu_{\vu\boxplus\vI}^c\br{A^{\br{-,j}}_\vomega\cap\supp g}
        \nrm{f^{\br{+,j}}}_{\vI\times\vomega}
        \nrm{\1_{A^{\br{-,j}}_\vomega}g}_{L^2\br{d\mu_{\vu\boxplus \vI}}}\abs{\vI}^{\frac{1}{2}}.
    \end{equation}
    Additionally, when \(j=1\), and the function \(f\) satisfies the frequency assumptions:
    \begin{equation}\label{eq_thm_Hdim_sing_tile_freq_cond}
        \dia\br{\pi_{\R\ve_0\oplus\R\ve_1} \supp\widehat{f}} \lesssim 2^{-k},
    \end{equation}
    the following stronger estimate holds:
    \begin{equation}\label{eq_thm_Hdim_sing_tile_est_half}
        \abs{\Lambda_{\vI\times\vomega}^{\br{1}}\br{f,g}}\lesssim 
        \lambda^{-1/2}
        \mu_{\vu\boxplus\vI}^{1/2}\br{A^{\br{-,1}}_\vomega\cap\supp g}
        \nrm{f^{\br{+,1}}}_{\vI\times\vomega}
        \nrm{\1_{A^{\br{-,1}}_\vomega}g}_{L^2\br{d\mu_{\vu\boxplus \vI}}}\abs{\vI}^{\frac{1}{2}}.
    \end{equation}
\end{proposition}

The proof of the above theorem will be left for \textsc{Section \ref{singltilesestimgen}}.

\subsection{Proof of \texorpdfstring{\eqref{eq_thm_Hdim_HF_st_mod}}{}: the \texorpdfstring{\(j=0\)}{} case}\label{subsec_pf_Hdim_HF_st_mod_zero_case}
In this subsection, we derive the high-dimensional analogs of \textbf{Lemma \ref{lem_-_tree} and \ref{lem_u_shift_mass_sel}} and adapt the argument in \textsc{Section \ref{subsec_mass_sel_ez_22_bd}} in order to prove \eqref{eq_thm_Hdim_HF_st_mod} in the \(j=0\) case.

We first recall \textbf{Definition \ref{def_tile_leqs}, \ref{def_tree}, and \ref{def_L2_BMO_energy}} and in the desire to give sense to \textbf{Definitions \ref{def_shift_mass} and \ref{def_sio_max}} in the \(D\geq 1\) regime, we generalize \eqref{eq_cU_cV_def} to
\begin{equation}\label{eq_Hdim_cU_def}
    \cU:=\BR{
        \vS^{-k}\vu\in\R^{D+1}\setminus\mr{-1,1}^{D+1}\::\:
        k\in\N\sqcup\BR{0}
    }\sqcup\BR{-1,0,1}^{D+1}
\end{equation}
and extend \eqref{eq_LP_freq_proj} to
\begin{equation}\label{eq_Hdim_LP_freq_proj}
    \widehat{\pi^k f}\br{\vxi}
    :=
    \varphi^{\otimes 2}\br{
        \vS^k
        \cdot
        \frac{\vxi}{2\lambda}
    }
    \widehat{f}\br{\vxi}
    ,\quad
    \pi_k:=\pi^k-\pi^{k+5}
\end{equation}
With the above modifications, the arguments presented in \textsc{Section \ref{sec_pf_key_lemmas}} yield the desired high-dimensional statement for \textbf{Lemma \ref{lem_-_tree}}. That is, for a \(\br{-,j}\)-tree \(\cT\subset \P\br{\lambda}\), the following estimate holds:
\begin{equation*}
    \nrm{\Lambda^{\br{j}}_P\br{f,\1_E}}_{\ell^1\br{\cT}}
    \lesssim \br{\nrm{\Lambda}_c+\nrm{T_\nu}_2}\log\br{e+\nrm{M_\nu}_2/\nrm{T_\nu}_2+\lambda\br{1+\abs{\vu}}}\nrm{f^{\br{+,j}}}_{L^2\br{\cT}}\cM^{\frac{1}{2}+c}_\cT\br{E}\abs{\vI_\cT}^{\frac{1}{2}}.
\end{equation*}
A line-by-line adaptation of the argument presented in \textsc{Section \ref{subsec_est_T_M_4_st_pt}} produces the following two estimates:
\begin{equation*}
    \nrm{T_\nu}_2\lesssim\lambda^{-\frac{1}{2}},\quad
    \nrm{M_\nu}_2\lesssim \log^{\frac{1}{2}}\br{e+\lambda\abs{\vu}}.
\end{equation*}
The \(\log\) factor comes from the use of the shifted maximal operator in the spirit of \eqref{eq_M_nu_st_pw_est}:
\begin{equation*}
    M_\nu f\br{\vx} \lesssim M_\para \br{\bigotimes_{k=0}^D M_{\lambda u_k}} f\br{\vx},\quad
    M_u g\br{z}:=\sup_{k\in\Z}
    \abs{g}\ast\Dil^1_{2^k}\Tr_u \chi\br{z}
\end{equation*}
where here  we used the fact that \(\nrm{M_u g}_{L^2\br{\R}}\lesssim \log^{\frac{1}{2}}\br{e+\abs{u}}\nrm{g}_{L^2\br{\R}}\).
Recalling our assumption, we appeal to \eqref{eq_thm_Hdim_sing_tile_est_half} in \textbf{Proposition \ref{thm_Hdim_sing_tile_ests}}, and, by simply estimating \(\nrm{\Lambda}_c\), we obtain:
\begin{lemma}[\(-,j\)-\textsf{tree estimate: the high-dimensional setting}]\label{lem_Hdim_-_tree}
There are universal constants \(\sigma,c>0\) such that for any \(\br{-,j}\)-tree \(\cT\subset \P\br{\lambda}\), we have the following estimate:
\begin{equation}\label{eq_lem_Hdim_-_tree}
    \nrm{\Lambda^{\br{j}}_P\br{f,\1_E}}_{\ell^1\br{\cT}}
    \lesssim 
    \lambda^{-\sigma}\log\br{e+\lambda\br{1+\abs{\vu}}}\nrm{f^{\br{+,j}}}_{L^2\br{\cT}}\cM^{\frac{1}{2}+c}_\cT\br{E}\abs{\vI_\cT}^{\frac{1}{2}}.
\end{equation}
\end{lemma}

Similarly, the arguments presented in \textsc{Section \ref{sec_pf_key_lemmas}} yield the high-dimensional analog of \textbf{Lemma \ref{lem_u_shift_mass_sel}}:
\begin{lemma}[\textsf{Mass selection: the high-dimensional setting}]\label{lem_Hdim_u_shift_mass_sel}
    For all \(\varsigma>0\), \(\cP\subset \P\) finite, and \(E\subset \R^{D+1}\) measurable, there is a finite list of disjoint trees \(\BR{\cT_i}_i\) with \(\cT_i\subset \cP\) such that the following hold:
    \begin{equation}\label{eq_lem_Hdim_u_shift_mass_sel}
        \cM_{\cP\setminus\bigsqcup_i\cT_i}\br{E}\leq \varsigma \qquad\textrm{and}\qquad
        \sum_{i}\abs{\vI_{\cT_i}}
        \lesssim \frac{\log\br{e+\abs{\vu}}\abs{E}}{\varsigma}.
    \end{equation}
\end{lemma}
Given \textbf{Lemmas \ref{lem_Hdim_-_tree}} and \textbf{\ref{lem_Hdim_u_shift_mass_sel}}, we may adopt the treatment presented in \textsc{Section \ref{subsec_mass_sel_ez_22_bd}} in order to prove \eqref{eq_thm_Hdim_HF_st_mod} when \(j=0\).
To elaborate, we use again the key fact that \textbf{all trees are \(\br{-,0}\)-trees} to apply \eqref{eq_lem_Hdim_-_tree} on all of the trees selected in \textbf{Lemma \ref{lem_Hdim_u_shift_mass_sel}}. Finally, \textbf{Observation \ref{obs_total_L2_energy_bd}} provides the key square summability, which enables a simplified argument and implies the desired estimate
\begin{equation*}
    \nrm{\Lambda^{\br{0}}_P\br{f,\1_E}}_{\ell^1\br{P\in\P\br{\lambda}}}\lesssim
    \lambda^{-\epsilon} 
    \log^2\br{e+\lambda\br{1+\abs{\vu}}}
    \nrm{f}_{L^2}\abs{E}^{\frac{1}{2}}.
\end{equation*}
This concludes the treatment for the \(j=0\) case.

\subsection{Refinement of \texorpdfstring{\eqref{eq_ass_wp_t_in_prod_V}}{} for the \texorpdfstring{\(j=1\)}{} case: preserving the high-resolution frequency structure}

\begin{proposition}\label{prop_Hdim_wp_in_prod_psi_VP}
Given \(f,g\in L^2\br{\R^{D+1}}\) and \(P_i:=\vI_i\times\vomega_i\in\P_{k_i}\) for \(i=1,2\), let \(f^{\br{+,1}}_{P_1},g^{\br{+,1}}_{P_2}\) be defined as in \eqref{eq_def_tf_proj_dual_scale} and \(k_2\leq k_1\). For any \(\Psi\in \mathcal{S}'\br{\R^{D+1}}\) satisfying \(\abs{\Psi} \ast \chi^N_{\vI_1} \underset{N}{\lesssim} \chi^N_{\vI_1}\),
the following estimate holds
    \begin{equation}\label{eq_prop_Hdim_wp_in_prod_psi_VP}
        \abs{
            \ang{
                \Psi\ast f^{\br{+,1}}_{P_1},
                g^{\br{+,1}}_{P_2}
            }
        }
        \underset{N}{\lesssim} 
        \nrm{f^{\br{+,1}}}_{P_1}
        \nrm{g^{\br{+,1}}}_{P_2}
        \min\br{1,\sqrt{\frac{\lambda^{D+1}\abs{\vI_2}}{\abs{\vI_1}}}}
        \nrm{\chi_{\vI_1}^N}_{L^\infty\br{\vI_2}}.
    \end{equation}
\end{proposition}
\begin{remark}
    Typically, \(\widehat{\Psi}\) will be a smooth function adapted to a $(D+1)-$interval of dimension greater than \(\br{2^{-k_1}}^D\times 4^{-k_1}\). Notably, the condition \(\abs{\Psi}\ast\chi^N_\vI\lesssim\chi^N_\vI\) for some \(\vI\in\I_k\) implies 
    \begin{equation}\label{eq_Psi_preserv_phy}
    \int
        \abs{\Psi}\br{\vz}\ang{\vS^{-k}\br{\vx-\vz}}^{-N}_\otimes 
    d\vz
    \underset{N}{\lesssim}
        \ang{\vS^{-k}\vx}^{-N}_\otimes 
    \end{equation}
    since translation commutes with convolution: \(\Tr_\vu \br{f\ast g}=f\ast \Tr_\vu g\). 
\end{remark}

\begin{proof}[\textbf{Proof of Proposition \ref{prop_Hdim_wp_in_prod_psi_VP}}]
    Let \(f,g,\Psi\) and \(P_i:=\vI_i\times\vomega_i\in\P_{k_i}\) be as provided by our hypothesis. We will first show the second estimate
    \begin{equation}\label{eq_Hdim_wp_Psi_est_w_lambda}
        \abs{
            \ang{
                \Psi\ast f^{\br{+,1}}_{P_1},
                g^{\br{+,1}}_{P_2}
            }
        }
        \underset{N}{\lesssim} 
        \nrm{f^{\br{+,1}}}_{P_1}
        \nrm{g^{\br{+,1}}}_{P_2}
        \sqrt{\frac{\lambda^{D+1}\abs{\vI_2}}{\abs{\vI_1}}}
        \nrm{\chi_{\vI_1}^N}_{L^\infty\br{\vI_2}}.
    \end{equation}
    Recalling \eqref{eq_Hdim_wp_CL_decomp} we slightly modify the argument in \eqref{eq_gaba_wp_inner_prod} in order to obtain the following:
    \begin{align*}
    \abs{\ang{
        \Psi\ast
        f^{\br{+,1}}_{P_1},
        g^{\br{+,1}}_{P_2}
    }}
    \leq 
    \int &
    \abs{\ang{\pi^{\br{+,1}}_{\vomega_1}f,\varphi^{\br{C}}_{k_1,\vx_1,\vxi_1}}}
    \chi^{\br{C}}_{\vI_1}\br{\vx_1}\1_{\vomega^{\br{+,0}}_1}\br{\vxi_1}\\
    \cdot &
    \abs{\ang{\pi^{\br{+,1}}_{\vomega_2}g,
    \varphi^{\br{C}}_{k_2,\vx_2,\vxi_2}}}
    \chi^{\br{C}}_{\vI_2}\br{\vx_2}\1_{\vomega^{\br{+,0}}_2}\br{\vxi_2}\\
     \cdot &
     \abs{\ang{
         \Psi\ast
         \varphi^{\br{C}}_{k_1,\vx_1,\vxi_1},
         \varphi^{\br{C}}_{k_2,\vx_2,\vxi_2}
     }}
     d\vx_1 d\vxi_1 d\vx_2 d\vxi_2.
\end{align*}
By Cauchy-Schwarz inequality and the trivial bound \(\chi^{\br{C}}_{\vI_k}\underset{N}{\lesssim} \chi_{\vI_k}^N\), we dominate the above by
\begin{equation*}
    \underset{N}{\lesssim}
    \nrm{f^{\br{+,1}}}_{P_1}
    \nrm{g^{\br{+,1}}}_{P_2}
    \nrm{
        \ang{
            \Psi\ast
            \varphi^{\br{C}}_{k_1,\vx_1,\vxi_1},
            \varphi^{\br{C}}_{k_2,\vx_2,\vxi_2}
        }
        \prod_{k=0,1}
        \chi_{\vI_k}^N\br{\vx_k}
        \1_{\vomega^{\br{+,0}}_k}\br{\vxi_k}
    }_{L^2\br{d\vx_1 d\vxi_1
    d\vx_2 d\vxi_2}}
    .
\end{equation*}
Note that the only key difference between the above expression and \eqref{eq_ass_wp_pf_f_g_phi} is the existence of the convolution with \(\Psi\) within the inner product term. However, since convolution with \(\Psi\) does not affect either the frequency or the physical localization, we have
\begin{itemize}
    \item \textsf{Preservation of the frequency localization:}
    \begin{equation*}
        \supp\br{
            \Psi\ast
            \varphi^{\br{C}}_{k_1,\vx_1,\vxi_1}
        }^{\wedge}
        =\supp
            \widehat{\Psi}
            \cdot
            \widehat{\varphi^{\br{C}}_{k_1,\vx_1,\vxi_1}}
        \subset 
        \supp
        \widehat{\varphi^{\br{C}}_{k_1,\vx_1,\vxi_1}}.
    \end{equation*}
    \item \textsf{Preservation of the physical localization:}
    \begin{align*}
        \abs{
            \Psi\ast
            \varphi^{\br{C}}_{k_1,\vx_1,\vxi_1}
        }
        \br{\vx}
        \underset{N}{\lesssim} &
        \int
            \abs{\Psi}\br{\vz}
            \sqrt{
            \det\br{\vS^{-k_1}}}
            \ang{\vS^{-k_1}\br{\vx-\vz-\vx_1}}^{-N}_\otimes
        d\vz\\
        \text{by \eqref{eq_Psi_preserv_phy},}\quad
        \underset{N}{\lesssim} &
        \sqrt{
        \det\br{\vS^{-k_1}}}
        \ang{\vS^{-k_1}\br{\vx-\vx_1}}^{-N}_\otimes
        .
    \end{align*}
\end{itemize}
    With these we can follow line-by-line the rest of the argument presented in \textsc{Section \ref{subsec_ver_tf_proj_emb}} to deduce
    \begin{equation*}
        \nrm{
            \ang{
                \Psi\ast
                \varphi^{\br{C}}_{k_1,\vx_1,\vxi_1},
                \varphi^{\br{C}}_{k_2,\vx_2,\vxi_2}
            }
            \prod_{k=0,1}
            \chi_{\vI_k}^N\br{\vx_k}
            \1_{\vomega^{\br{+,0}}_k}\br{\vxi_k}
        }_{L^2\br{d\vx_1 d\vxi_1
        d\vx_2 d\vxi_2}}
        \underset{N}{\lesssim}
        \sqrt{\frac{\lambda^{D+1}\abs{\vI_2}}{\abs{\vI_1}}}
        \nrm{\chi_{\vI_1}^N}_{L^\infty\br{\vI_2}}
    \end{equation*}
    and conclude \eqref{eq_Hdim_wp_Psi_est_w_lambda}.
    
    It now remains to show the first estimate, that is
    \begin{equation}\label{eq_Hdim_wp_Psi_est_wo_lambda}
        \abs{
            \ang{
                \Psi\ast f^{\br{+,1}}_{P_1},
                g^{\br{+,1}}_{P_2}
            }
        }
        \underset{N}{\lesssim} 
        \nrm{f^{\br{+,1}}}_{P_1}
        \nrm{g^{\br{+,1}}}_{P_2}
        \nrm{\chi_{\vI_1}^N}_{L^\infty\br{\vI_2}}.
    \end{equation}
    Suggested by the treatment of the analogous statement presented in \textsc{Section \ref{subsec_ver_tf_proj_emb}}, we shall now utilize the expression \eqref{eq_def_tf_proj_dual_scale} and write:
    \begin{align*}
        \ang{
            \Psi\ast f^{\br{+,1}}_{P_1},
            g^{\br{+,1}}_{P_2}
        }
        =&\ang{
            \Psi\ast\br{\chi^\Phi_{\vI_1}\pi^{\br{+,1}}_{\vomega_1} f},
            \chi^\Phi_{\vI_2}\pi^{\br{+,1}}_{\vomega_2} g
        }\\
        =&
        \int
            \Psi\br{\vx-\vy}\chi^\Phi_{\vI_1}\br{\vy}\pi^{\br{+,1}}_{\vomega_1} f\br{\vy}
            \chi^\Phi_{\vI_2}\br{\vx}\overline{\pi^{\br{+,1}}_{\vomega_2} g}\br{\vx}
        d\vx d\vy.
    \end{align*}
With some application of trivial estimates, we obtain:
\begin{equation}\label{eq_Hdim_wp_Psi_est_wo_lambda_exp}
    \abs{
        \ang{
            \Psi\ast f^{\br{+,1}}_{P_1},
            g^{\br{+,1}}_{P_2}
        }
    }
    \underset{N}{\lesssim}
    \int
        \chi^N_{\vI_2}\br{\vx}
        \abs{\Psi}\br{\vx-\vy}\chi^N_{\vI_1}\br{\vy}
        \abs{\sqrt{\chi^\Phi_{\vI_1}\br{\vy}}\pi^{\br{+,1}}_{\vomega_1} f\br{\vy}}
        \abs{\sqrt{\chi^\Phi_{\vI_2}\br{\vx}}\pi^{\br{+,1}}_{\vomega_2} g}\br{\vx}
    d\vx d\vy.
\end{equation}
This suggests we consider the following positive operator:
\begin{equation*}
    Tf\br{\vx}:=
    \int
        \chi^N_{\vI_2}\br{\vx}
        \abs{\Psi}\br{\vx-\vy}\chi^N_{\vI_1}\br{\vy}
        f\br{\vy}
    d\vy.
\end{equation*}
Once we show the following operator bound:
\begin{equation}\label{eq_Hdim_wp_Psi_est_wo_lambda_op_bd}
    \nrm{T f}_{L^2}
    \underset{N}{\lesssim}
    \nrm{\chi_{\vI_1}^N}_{L^\infty\br{\vI_2}}
    \nrm{f}_{L^2},
\end{equation}
we recover \eqref{eq_Hdim_wp_Psi_est_wo_lambda} by rewriting \eqref{eq_Hdim_wp_Psi_est_wo_lambda_exp} and invoking \eqref{eq_Hdim_wp_Psi_est_wo_lambda_op_bd} and \eqref{eq_Hdim_wp_nrm}:
\begin{equation*}
    \abs{
        \ang{
            \Psi\ast f^{\br{+,1}}_{P_1},
            g^{\br{+,1}}_{P_2}
        }
    }
    \underset{N}{\lesssim}
    \ang{
        T
        \abs{\sqrt{\chi^\Phi_{\vI_1}}\pi^{\br{+,1}}_{\vomega_1} f},
        \abs{\sqrt{\chi^\Phi_{\vI_2}}\pi^{\br{+,1}}_{\vomega_2} g}
    }
    \underset{N}{\lesssim} 
    \nrm{f^{\br{+,1}}}_{P_1}
    \nrm{g^{\br{+,1}}}_{P_2}
    \nrm{\chi_{\vI_1}^N}_{L^\infty\br{\vI_2}}.
\end{equation*}
To prove the operator bound \eqref{eq_Hdim_wp_Psi_est_wo_lambda_op_bd}, we first observe that the assumption \(\abs{\Psi}\ast\chi^N_{\vI_1}\underset{N}{\lesssim}\chi^N_{\vI_1}\) implies:
\begin{equation}\label{eq_Hdim_wp_Psi_est_wo_lambda_op_bd_Linfty}
    \nrm{Tf}_{L^\infty}
    \underset{N}{\lesssim}
    \nrm{\chi^N_{\vI_2}\chi^N_{\vI_1}}_{L^\infty}\nrm{f}_{L^\infty}
    \underset{N}{\eqsim}
    \nrm{\chi^N_{\vI_1}}_{L^\infty\br{\vI_2}}
    \nrm{f}_{L^\infty}.
\end{equation}
Next, we rewrite the defining formula for \(Tf\):
\begin{equation*}
    Tf\br{\vx}:=
    \int
        \chi^N_{\vI_2}\br{\vx}
        \abs{\Psi}\br{\vy}\chi^N_{\vI_1}\br{\vx-\vy}
        f\br{\vx-\vy}
    d\vy.
\end{equation*}
Via the trivial estimate \(\chi_{\vI_2}\br{\vx}\chi_{\vI_1}\br{\vx-\vy}\lesssim\chi_{\vI_1}\br{\vc_{\vI_2}-\vy}\) and the assumption \(\abs{\Psi}\ast\chi^N_{\vI_1}\underset{N}{\lesssim}\chi^N_{\vI_1}\), we obtain:
\begin{equation}\label{eq_Hdim_wp_Psi_est_wo_lambda_op_bd_L1}
    \nrm{Tf}_{L^1}
    \underset{N}{\lesssim}
    \int
        \abs{\Psi}\br{\vy}\chi^N_{\vI_1}\br{\vc_{\vI_2}-\vy}
        \nrm{
            f\br{\vx-\vy}
        }_{L^1\br{d\vx}}
    d\vy
    \underset{N}{\lesssim}
    \chi^N_\vI\br{\vc_{\vI_2}}\nrm{f}_{L^1}
    \leq
    \nrm{\chi^N_{\vI_1}}_{L^\infty\br{\vI_2}}
    \nrm{f}_{L^1}.
\end{equation}
We conclude \eqref{eq_Hdim_wp_Psi_est_wo_lambda_op_bd} by interpolating between \eqref{eq_Hdim_wp_Psi_est_wo_lambda_op_bd_Linfty} and \eqref{eq_Hdim_wp_Psi_est_wo_lambda_op_bd_L1}.

\end{proof}

\subsection{Specific issues in the \texorpdfstring{\(j=1\)}{} case: frequency sub-systems and scale sparsification}\label{subsec_Hdim_issue_j=1}
In this section, we discuss the various issues arising when trying to mimic the $j=0$ argument in \textsc{Section \ref{subsec_pf_Hdim_HF_st_mod_zero_case}} as well as means to identify a remedy, based on the notion of frequency sub-systems. 

We start by briefly justifying the need for \textbf{Proposition \ref{prop_Hdim_wp_in_prod_psi_VP}} and explaining the motivation for considering fine-scale frequency structures: The origin of our difficulties rests in the fact that the resulting trees in the selection algorithm of \textbf{Lemma \ref{lem_Hdim_u_shift_mass_sel}} are not necessarily \(\br{-,1}\)-trees. This suggests a generalization of \textbf{Lemma \ref{lem_gen_tree_c_half}} which incorporates \eqref{eq_sing_tile_est} with \(c=\frac{1}{2}\)
for \(\br{+,1}\)-trees corresponding to the situation when \(\BR{f^{\br{+,1}}_P}_{P\in\cT}\) have overlapping time-frequency supports. 
\smallskip

\emph{However,  in the higher-dimensional setting, direct calculations show that} 
\begin{equation}\label{blowup}
\nrm{\Lambda}_{\frac{1}{2}}=\nrm{\Lambda}_{\frac{1}{2},\lambda,1}\to\infty\quad \textit{as}\quad \lambda\to\infty\,.
\end{equation}

The growth of \(\nrm{\Lambda}_{\frac{1}{2},\lambda,1}\) removes the possibility of proving \eqref{eq_thm_Hdim_HF_st_mod} via a naive generalization of \textbf{Lemma \ref{lem_gen_tree_c_half}}. Thus, when dealing with \((+,1)\)-trees, the remedy will be to replace \eqref{eq_sing_tile_est} by \eqref{eq_thm_Hdim_sing_tile_est_half} subject to condition \eqref{eq_thm_Hdim_sing_tile_freq_cond}. As a consequence, for the proper generalization of \textbf{Definition \ref{def_L2_BMO_energy}} and \textbf{Lemma \ref{lem_gen_tree_c_half}}, requirement \eqref{eq_thm_Hdim_sing_tile_freq_cond} invites us to consider a fine-scale frequency decomposition in the \(0\)-th and \(1\)-st coordinates.
\smallskip

For expository reasons, we let \(k_\lambda\in\N\) be such that \(2^{k_\lambda-1}< \lambda \leq 2^{k_\lambda}\) and set:
\begin{equation}\label{eq_one_lambda_and_slab_mat}
    1_\lambda:=\lambda /2^{k_\lambda}\in\bR{1/2,1},\quad
    \fS:=\fS_\lambda:=2^{-k_\lambda} \ve_0\ve^\top_0+\sum_{i=1}^D\ve_i\ve_i^\top\in M_{\br{D+1}\times\br{D+1}}\br{\R}.
\end{equation}
As a natural extension of \eqref{eq_Hdim_W_n_Wk}, we define for each \(k\in\Z\) the following system of intervals:
\begin{definition}[\textsf{Frequency sub-system}]\label{def_Hdim_freq_sub}
\begin{equation}\label{eq_slabs}
    \fS\br{\W_k}:=
    \BR{
        \fS\vomega\::\:\vomega\in\W_k
    }=
    \BR{
        \br{0,\vw}\boxplus \vS^{-k}\bR{-1_\lambda/2,1_\lambda/2}\times\bR{0,\lambda}^D
        \::\:
        \vw\in\Z^D
    },
\end{equation}
\begin{equation}\label{eq_sticks}
    \fs\br{\W_k}:=\BR{
        \br{0,\vw}\boxplus \vS^{-k}\bR{-1_\lambda/2,1_\lambda/2}\times\bR{0,1_\lambda}\times\bR{0,\lambda}^{D-1}
        \::\:
        \vw\in\Z^D
    }.
\end{equation}
\end{definition}
\begin{observation}
    For \(\fs\vomega\in\fs\br{\W_k}\), there is a unique \(\vomega\in\W_k\) such that \(\fs\vomega\subset \vomega\) and also \(\fs\vomega\subset \fS\vomega\).
\end{observation}
Take now an interval \(\vomega:=\br{0,\vw}\boxplus \vS^{-k}\bR{-\lambda/2,\lambda/2}\times\bR{0,\lambda}^D\in\W_k\), the interval \(\fS\vomega\in\fS\br{\W_k}\) can be subdivided into \(2^{k_\lambda}\) many sub-intervals in \(\fs\br{\W_k}\) of the following form:
\begin{equation}\label{eq_stick_v}
        \fs_v\vomega:=\br{0,2^{k_\lambda} w_1+v,\widehat{\vw}_1}\boxplus
        \vS^{-k}\bR{-1_\lambda/2,1_\lambda/2}\times\bR{0,1_\lambda}\times\bR{0,\lambda}^{D-1},\quad
        v\in\Z.
\end{equation}
Specifically, the relevant collection of intervals are the following:

\begin{definition}[\textsf{Frequency sub-system subject to} \(\vomega\) \textsf{and} \(\vomega^{\br{+,1}}\)]\label{def_Hdim_freq_sub_w}
\begin{equation}\label{eq_stick_w}
    \fs\br{\vomega}:=\BR{
        \fs\vomega\in\fs\br{\W_k}\::\: \fs\vomega\subset \vomega
    }=\BR{\fs_v\vomega\in\fs\br{\W_k}\::\: 0\leq v<2^{k_\lambda}}
\end{equation}
\begin{equation}\label{eq_stick_w_uphalf}
    \supset
    \fs^{\br{+,1}}\br{\vomega}:=
    \BR{
        \fs\vomega\in\fs\br{\vomega}\::\: \fs\vomega\subset \vomega^{\br{+,1}}
    }
    =
    \BR{\fs_v\vomega\in\fs\br{\W_k}\::\: 2^{k_\lambda-1}-1\leq v<2^{k_\lambda}}.
\end{equation}
Notably, we have two identities:
\begin{equation}\label{eq_stick_2_slab}
    \fS\vomega=\bigsqcup\fs\br{\vomega}=
    \bigsqcup_{v=0}^{2^{k_\lambda}-1}
        \fs_v\vomega,\quad
    \fS\vomega^{\br{+,1}}=\bigsqcup\fs^{\br{+,1}}\br{\vomega}=
    \bigsqcup_{v=2^{k_\lambda-1}-1}^{2^{k_\lambda}-1}
        \fs_v\vomega.
\end{equation}
\end{definition}
%
%
With the above definitions, we make the following observation:\footnote{For convenience, throughout this section we will suppress the dependence on \(j\) whenever it does not affect our reading.}
%
%
\begin{observation}[\textsf{Nested structures}]\label{obs_Hdim_freq_sub_nest}
    Let \(l\geq k+k_\lambda\) and \(\vomega'\in\W_l\). Within the unique \(\vomega\in\W_k\) with \(\vomega'\subset \vomega\), there is a unique \(\fs\vomega\in\fs\br{\vomega}\) such that \(\vomega'\subset \fs\vomega\). Moreover, if \(\vomega'\subset \vomega^{\br{+,1}}\), we must also have \(\fs\vomega\in\fs^{\br{+,1}}\br{\vomega}\).
\end{observation}
This strengthens \textbf{Observation \ref{obs_tree_freq_scale_rel}} and suggests pigeonholing on scales to utilize the nested structures.

\begin{observation}[\textsf{Scale sparsification}]\label{obs_Hdim_P_sparse}
Consider the following system of tiles:
\begin{equation}\label{eq_Hdim_P_sparse}
    \P:=\bigsqcup_{l=0}^{k_\lambda}
    \P_{\equiv l},\quad
    \P_{\equiv l}:=
    \bigsqcup_{k\in l+k_\lambda\Z} \P_k
\end{equation}
It suffices to prove \eqref{eq_thm_Hdim_HF_st_mod} with \(\P\) replaced by \(\P_{\equiv l}\) since \(k_\lambda \eqsim \log\br{e+\lambda}\) and thus, the pigeonholing process only contributes a negligible \(\log\br{e+\lambda}\) loss.    
\end{observation}

Henceforth, we fix \(0\leq l<k_\lambda\), impose throughout the rest of the section that \(k\in l+k_\lambda\Z\), and make the following abuse of notations:
\begin{equation}\label{eq_Hdim_sparse_setting}
    \I:=\bigsqcup_{k\in l+k_\lambda\Z} \I_k,\quad
    \W:=\bigsqcup_{k\in l+k_\lambda\Z} \W_k,\quad
    \P:=\bigsqcup_{k\in l+k_\lambda\Z} \P_k.
\end{equation}

We finish this subsection with a last comment on the needed modifications for the treatment of the \(j=1\) case: Instead of the total \(L^2\) energy bound \eqref{eq_obs_total_L2_energy_bd} in \textbf{Observation \ref{obs_total_L2_energy_bd}} a selection algorithm in the spirit of \eqref{eq_def_out_L2_infty_energy_emb}  will be required for a proper generalization of \(\nrm{f^{\br{+,1}}}_{\bmo_2\br{\cP}}\) in \textbf{Definition \ref{def_L2_BMO_energy}}.

\subsection{Single tree estimate for the \texorpdfstring{\(j=1\)}{} case}\label{sgtreejone}

We first elaborate on some finer tree structure notions when considering sub-systems and sparse scales:
\begin{definition}\label{def_Hdim_tree_anchor_to_subsys}
    Let \(\cT\subset\P\) be a tree with top \(\vI_\cT\times\vomega_\cT\in\P\), and let \(\vc_\cT\in\R^{D+1}\) denote the center of \(\vomega^{\br{+,1}}_\cT\).
    For \(\vI\times\vomega\in \cT\), we denote \(v_\cT\) the unique integer  such that \(\fs_\cT\vomega:=\fs_{v_\cT}\vomega\in\fs\br{\vomega}\) is the unique sub-system interval satisfying \(\vc_\cT\in\fs_\cT\vomega
    \subset \vomega\).
\end{definition}
Given a tree \(\cT\subset \P\) with top \(\vI_\cT\times\vomega_\cT\in\P\) and two tiles \(\vI\times\vomega, \vI'\times\vomega'\in\cT\),
\textbf{Observation \ref{obs_Hdim_freq_sub_nest}} forces that:
\begin{equation}\label{eq_obs_Hdim_freq_sub_nest_w_tree}
    \vomega_\cT\subset \vomega'\subsetneq \vomega\implies
    \vomega'\subset \fs_\cT\vomega.
\end{equation}
In other words, whenever \(\vI\subsetneq \vI' \), the worst frequency overlap for two tiles is controlled by \(\fs_\cT\vomega\). This suggests that we split \(f^{\br{+,1}}_{\vI\times\vomega}\) for all \(\vI\times\vomega\in\cT\) into the outer parts with frequency support on \(\vomega^{\br{+,1}}\setminus\fs_\cT\vomega\) and the central parts with frequency support on \(\fs_\cT\vomega\). The former can then be treated via an orthogonality argument, while the latter would satisfy \eqref{eq_thm_Hdim_sing_tile_freq_cond} and enable \eqref{eq_thm_Hdim_sing_tile_est_half} in \textbf{Proposition \ref{thm_Hdim_sing_tile_ests}} to address the lack of orthogonality condition. 
To elaborate, we define:

\begin{definition}[\textsf{Frequency projection for a sub-system}]\label{def_Hdim_freq_proj_sub_sys}
Given \(k\) and \(\fs\vomega\in\fs\br{\W_k}\), we set
    \begin{equation}\label{eq_Hdim_freq_proj_for_subsys}
        \widehat{\fp_k f}\br{\vxi}:=
        \phi\br{2^k\xi_0/1_\lambda}\widehat{f}\br{\vxi},\quad
        \widehat{\fp_{\fs\vomega} f}
        :=\br{\Dil^\infty_2\varphi}_{\fs\vomega}
        \widehat{\fp_k f}
        ,\quad
        \fp_{\setminus\fs\vomega}:=\id-\fp_{\fs\vomega}.
    \end{equation}
\end{definition}

The exact definition of outer parts and central parts will be the following:

\begin{definition}[\textsf{Outer and central parts of a tree}]\label{def_Hdim_cen_ext}
    Given a tree \(\cT\subset P\) and a tile \(\vI\times\vomega\in\cT\), we set:
    \begin{equation}\label{eq_cent_n_out}
        f^{\br{+,1}}_{\vI\times\vomega}=
        \br{
            \fp_{\setminus\fs_\cT\vomega}f
        }^{\br{+,1}}_{\vI\times\vomega}
        +
        \br{
            \fp_{\fs_\cT\vomega}f
        }^{\br{+,1}}_{\vI\times\vomega}
        =:f^{\cT_{out}}_{\vI\times\vomega}+f^{\cT_{cen}}_{\vI\times\vomega}.
    \end{equation}
    Correspondingly, we introduce the following modifications of \eqref{eq_Hdim_wp_nrm}:
    \begin{equation}
    \label{eq_cent_n_out_nrm}
        \nrm{f^{\cT_{out}}}
        _{\vI\times\vomega}:=
        \nrm{
            \br{
                \fp_{\setminus
                \fs_\cT\vomega}
                f
            }^{\br{+,1}}
        }_{\vI\times\vomega},\quad
        \nrm{f^{cen}}
        _{\vI\times\vomega}
        :=
        \sup
        _{\fs\vomega\in\fs^{\br{+,1}}\br{\vomega}}
        \nrm{
            \br{
                \fp_{\fs\vomega}f
            }^{\br{+,1}}
        }_{\vI\times\vomega}.
    \end{equation}
    Notice that the latter does not require a reference to a tree.
\end{definition}

At this stage, we are ready to introduce the suitable replacement for \textbf{Definition \ref{def_L2_BMO_energy}}:

\begin{definition}[\textsf{The generalized energy}]\label{def_Hdim_energy}
    Let \(\cT\subset \P\) be a tree, we define for \(f\in L^2\br{\R^{D+1}}\) the following  modified \(L^2\) energy:
    \begin{equation*}
        \nrm{f^{out}}_{L^2\br{\cT}}:=
        \nrm{
            \nrm{f^{\cT_{out}}}_{\vI\times\vomega}
        }_{\ell^2\br{P\in\cT}}.
    \end{equation*}
    As the replacement of the BMO quantity, we define:
    \begin{equation*}
        \nrm{f}_{E\br{\cP}}:=
        \max\br{
            \sup_{\substack{
                \cT\subset \cP\\
                \cT\text{ be a tree}
            }}
            \nrm{f^{out}}_{L^2\br{\cT}}
            /\abs{\vI_\cT}^{\frac{1}{2}},
            \sup_{\vI\times\vomega\in\cP}
                \nrm{f^{cen}}_{\vI\times\vomega}/\abs{\vI}^{\frac{1}{2}}
        }.
    \end{equation*}
\end{definition}

Finally, we state the appropriate generalization of \textbf{Lemma \ref{lem_gen_tree_c_half}}:

\begin{lemma}[\textsf{General tree estimate: the high-dimensional setting}]\label{lem_Hdim_gen_tree}
There are universal constants \(\sigma,c>0\) such that for a tree \(\cT\subset \P\), the following estimate holds:
    \begin{equation}\label{eq_lem_Hdim_gen_tree}
        \nrm{\Lambda^{\br{1}}_P\br{f,\1_E}}_{\ell^1\br{P\in\cT}}
        \lesssim \lambda^{-\sigma}\log\br{e+\lambda\br{1+\abs{\vu}}}\nrm{f}_{E\br{\cT}}M^{c+\frac{1}{2}}_\cT\br{E}\abs{\vI_\cT}.
    \end{equation}
\end{lemma}
 The rest of the subsection will be dedicated to the proof of the above lemma.

\subsubsection{Proof of Lemma \ref{lem_Hdim_gen_tree}: the tree top}
Let now \(\cT\) be a tree with top \(P_\cT:=\vI_\cT\times\vomega_\cT\in\P\). We claim that it suffices to consider the \textbf{topless} situation: \(P_\cT\notin\cT\). 

Assume the alternative: \(P_\cT\in\cT\). We have the following trivial fact:
\begin{equation}\label{eq_Hdim_gen_tree_w_top}
    \nrm{\Lambda^{\br{1}}_P\br{f,\1_E}}_{\ell^1\br{P\in\cT}}
    =
    \nrm{\Lambda^{\br{1}}_P\br{f,\1_E}}_{\ell^1\br{P\in\cT\setminus\BR{P_\cT}}}
    +
    \abs{\Lambda^{\br{1}}_{P_\cT}\br{f,\1_E}}.
\end{equation}

To estimate the last term, we introduce the notion of trivial trees:

\begin{definition}[\textsf{Trivial trees}]
    A tile \(P\in\P\) can be viewed as a tree and specifically a \(\br{\pm,j}\)-tree for any choice of sign \(\pm\) and \(j\), by setting \(\cT_P:=\BR{P}\) and \(P_\cT:=P\). 
\end{definition}

\begin{observation}\label{obs_Hdim_tile_leq_engery}
    Let \(P:=\vI\times\vomega\in\P\), the following trivial inequality holds:
    \begin{equation*}
        \nrm{f^{\br{+,1}}}_P/\abs{\vI}^{\frac{1}{2}}\leq
        \nrm{f^{\cT_{P,out}}}_P/\abs{\vI}^{\frac{1}{2}}+
        \nrm{f^{cen}}_P/\abs{\vI}^{\frac{1}{2}}
        \lesssim 
        \nrm{f}_{E\br{\BR{P}}}.
    \end{equation*}
\end{observation}

As a direct consequence, \textbf{Proposition \ref{thm_Hdim_sing_tile_ests}} implies:
\begin{align*}
    \abs{\Lambda_{P_\cT}^{\br{1}}\br{f,\1_E}}
    \lesssim &
        \lambda^{-\sigma}
        \nrm{f^{\br{+,1}}}_{P_\cT}
        \mu_{\vu\boxplus\vI_\cT}^{c+\frac{1}{2}}\br{E\cap A^{\br{-,1}}_{\vomega_\cT}}
        \abs{\vI_\cT}^{\frac{1}{2}}\\
    \lesssim & \lambda^{-\sigma}
        \nrm{f}_{E\br{\BR{P_\cT}}}
        M^{c+\frac{1}{2}}_{\BR{P_\cT}}\br{E}\abs{\vI_\cT}
    \leq \lambda^{-\sigma}
    \nrm{f}_{E\br{\cT}}
    M^{c+\frac{1}{2}}_{\cT}\br{E}\abs{\vI_\cT}.
\end{align*}
This addresses the last term in \eqref{eq_Hdim_gen_tree_w_top}. Henceforth, it suffices to consider the topless situation: \(\cT=\cT\setminus\BR{P_\cT}\).
For the remaining discussion, the tree \(\cT\) will be realized as the disjoint union of a topless \(\br{-,1}\)-tree and a topless \(\br{+,1}\)-tree. The analysis of the two cases will be carried out separately.

\subsubsection{Proof of Lemma \ref{lem_Hdim_gen_tree}: the topless \texorpdfstring{\(\br{-,1}\)}{}-tree}
Let \(\cT\) be a topless \(\br{-,1}\)-tree. For any \(P:=\vI\times\vomega\in \cT\), we observe the topless assumption implies \(\vI\subsetneq \vI_\cT\) and by \eqref{eq_obs_Hdim_freq_sub_nest_w_tree}:
\begin{equation}
    \vomega_\cT\subset \fs_\cT\vomega,\quad
    \vomega^{\br{-,1}}_\cT\subset \vomega_\cT \subset \vomega^{\br{-,1}}.
\end{equation}
As a result, \(\fs_\cT\vomega\cap\vomega^{\br{+,1}}=\varnothing\) and the following two pairs of notions match:
\begin{equation*}
    f^{\cT_{out}}_{\vI\times\vomega}=f^{\br{+,1}}_{\vI\times\vomega},\quad
    \nrm{f^{\cT_{out}}}_{\vI\times\vomega}=
    \nrm{f^{\br{+,1}}}_{\vI\times\vomega}.
\end{equation*}
Specifically, we have:
\begin{equation*}
    \nrm{f^{\br{+,1}}}_{L^2\br{\cT}}/\abs{\vI_\cT}^{\frac{1}{2}}=\nrm{f^{out}}_{L^2\br{\cT}}/\abs{\vI_\cT}^{\frac{1}{2}}\leq \nrm{f}_{E\br{\cT}}.
\end{equation*}
As a direct consequence of the above inequality and \textbf{Lemma \ref{lem_Hdim_-_tree}}, we obtain the desired estimate:
\begin{align*}
    \nrm{\Lambda^{\br{1}}_P\br{f,\1_E}}_{\ell^1\br{\cT}}
    \lesssim &
    \lambda^{-\sigma}\log\br{e+\lambda\br{1+\abs{\vu}}}\nrm{f^{\br{+,1}}}_{L^2\br{\cT}}\cM^{\frac{1}{2}+c}_\cT\br{E}\abs{\vI_\cT}^{\frac{1}{2}}\\
    \leq &
    \lambda^{-\sigma}\log\br{e+\lambda\br{1+\abs{\vu}}}\nrm{f}_{E\br{\cT}}\cM^{\frac{1}{2}+c}_\cT\br{E}\abs{\vI_\cT}.
\end{align*}

\subsubsection{Proof of Lemma \ref{lem_Hdim_gen_tree}: the topless \texorpdfstring{\(\br{+,1}\)}{}-tree}
We recall the linearization trick \eqref{eq_remove_abs_w_coef}, notation \eqref{eq_T_k_E_omega_def}, and the Whitney decomposition \eqref{eq_J_def} adapted to the tree \(\cT\). With these, we reformulate \eqref{eq_lem_Hdim_gen_tree} as
\begin{equation*}
    \sum_{\vJ\in\J}\sum_{k\in\Z}\sum_{\vI\times\vomega\in\cT_k}
    \ang{\cC_k f^\epsilon_{\vI\times\vomega},\1_{E^-_\vomega\cap\vJ}}
    \lesssim \lambda^{-\sigma}\log\br{e+\lambda\br{1+\abs{\vu}}}\nrm{f}_{E\br{\cT}}M^{c+\frac{1}{2}}_\cT\br{E}\abs{\vI_\cT}.
\end{equation*}
Let again \(k_\vJ\) be the \(k\in\Z\) such that \(\vJ\in\I_k\) and split the left-hand side of the above expression into the small-scale and the large-scale contributions as in the proof of \textbf{Lemma \ref{lem_-_tree} and \ref{lem_gen_tree_c_half}}:
\begin{equation*}
    \sum_{\vJ\in\J}\sum_{k\in\Z}\br{\cdots}=
    \sum_{\vJ\in\J}\sum_{k<k_\vJ+2}\br{\cdots}+
    \sum_{\vJ\in\J}\sum_{k\geq k_\vJ+2}\br{\cdots}.
\end{equation*}
For the small-scale contribution when \(k<k_\vJ+2\),
the desired estimate follows almost line-by-line the argument presented in \textsc{Section \ref{subsubsec_gen_tree_small_scale}} with the following list of modifications:
\begin{itemize}
    \item Apply \eqref{eq_sing_tile_est} with the \(c>0\) shown in \textbf{Proposition \ref{thm_Hdim_sing_tile_ests}} instead of \(c=\frac{1}{2}\).
    \item Replace \(\mu_{\vu\boxplus\vI}\br{E^-_\vomega\cap\vJ}\) with \(\mu^{c+\frac{1}{2}}_{\vu\boxplus\vI}\br{E^-_\vomega\cap\vJ}\).
    \item Replace the trivial bound \(\nrm{f^{\br{+,1}}}_{\vI\times\vomega
    }/\abs{\vI}^{\frac{1}{2}}\leq \nrm{f^{\br{+,1}}}_{\bmo_2\br{\cT}}\) with \textbf{Observation \ref{obs_Hdim_tile_leq_engery}}.
    \item Replaced the large constant \(\frac{N}{2}\) with \(\frac{N}{2}\cdot\br{c+\frac{1}{2}}\).
\end{itemize}
In combination, we derive:
\begin{equation}
    \sum_{\vJ\in\J}\sum_{k<k_\vJ+2}\sum_{\vI\times\vomega\in\cT_k}
    \abs{
        \ang{\cC_k f^\epsilon_{\vI\times\vomega},\1_{E^-_\vomega\cap\vJ}}
    }
    \lesssim  \nrm{\Lambda}_c \cL\nrm{f}_{E\br{\cT}}M^{c+\frac{1}{2}}_\cT\br{E}\abs{\vI_\cT},
\end{equation}
which is bounded by the right-hand side of \eqref{eq_lem_Hdim_gen_tree}. 

For the large-scale contribution when \(k\geq k_\vJ+2\),
we recall \eqref{eq_cent_n_out} and decompose the expression \(f^\epsilon_P\):
\begin{equation*}
    f^\epsilon_P:=\epsilon_Pf^{\br{+,1}}_P=\epsilon_P f^{\cT_{cen}}_P+\epsilon f^{\cT_{out}}_P=: f^{\epsilon,cen}_P + f^{\epsilon,out}_P,\quad P\in\cT.
\end{equation*}
It suffices to show for \(\ast\in\BR{cen,out}\) the following:
\begin{equation}\label{eq_lem_Hdim_gen_tree_large_scale_split}
    \abs{
        \sum_{\vJ\in\J}\sum_{k\geq k_\vJ+2}\sum_{\vI\times\vomega\in\cT_k}
            \ang{\cC_k f^{\epsilon,\ast}_{\vI\times\vomega},\1_{E^-_\vomega\cap\vJ}}
    }
    \lesssim \lambda^{-\sigma} \log\br{e+\lambda\br{1+\abs{\vu}}}\nrm{f}_{E\br{\cT}}M^{c+\frac{1}{2}}_\cT\br{E}\abs{\vI_\cT}.
\end{equation}

For the \(\ast=cen\) case, the argument is similar to that presented in \textsc{Section \ref{subsubsec_gen_tree_large_scale}}. However, we do not apply \eqref{eq_sing_tile_est} with \(c=\frac{1}{2}\) directly. Instead, we observe that for \(P\in\cT_k:=\cT\cap \P_k\) relation \eqref{eq_thm_Hdim_sing_tile_freq_cond} holds:
\begin{equation*}
    f^{\epsilon,cen}_P:=\epsilon_P f^{\cT_{cen}}_P:=\br{
        \epsilon_P
            \fp_{\fs_\cT\vomega_{\cT_k}}f
        }^{\br{+,1}}_P,\quad
    \dia\br{\pi_{\R\ve_0\oplus\R\ve_1}\supp \widehat{\epsilon_P\fp_{\fs_\cT\vomega_{\cT_k}}f }}
    \lesssim 2^{-k}.
\end{equation*}
We thus may apply \eqref{eq_thm_Hdim_sing_tile_est_half} in \textbf{Proposition \ref{thm_Hdim_sing_tile_ests}} instead of \eqref{eq_sing_tile_est} with \(c=\frac{1}{2}\) to gain an extra decaying factor:
\begin{align*}
    &
    \sum_{\vJ\in\J}
    \sum_{k\geq k_\vJ+2}
    \sum_{\vI\times\vomega \in\cT_k}
        \abs{
            \ang{\cC_k f^{\epsilon,cen}_{\vI\times\vomega},\1_{E^-_\vomega\cap\vJ}}
        }\\
    \lesssim &
    \lambda^{-\frac{1}{2}}
    \sum_{\vJ\in\J}
    \sum_{k\geq k_\vJ+2}
    \sum_{\vI\times\vomega\in\cT_k}
    \nrm{\br{\epsilon_P\fp_{\fs_\cT\vomega_{\cT_k}}f}^{\br{+,1}}}_{\vI\times\vomega}
    \mu_{\vu\boxplus\vI}\br{E^-_\vomega\cap\vJ}\abs{\vI}^{\frac{1}{2}}\\
    \leq &
    \lambda^{-\frac{1}{2}}
    \sum_{\vJ\in\J}
    \sum_{k\geq k_\vJ+2}
    \sum_{\vI\times\vomega\in\cT_k}
    \nrm{f^{cen}}_{\vI\times\vomega}
    \mu_{\vu\boxplus\vI}\br{E^-_\vomega\cap\vJ}\abs{\vI}^{\frac{1}{2}}\\
    \leq &
    \lambda^{-\frac{1}{2}}
    \nrm{f}_{E\br{\cT}}
    \sum_{\vJ\in\J}
    \sum_{k\geq k_\vJ+2}
    \sum_{\vI\times\vomega\in\cT_k}
    \mu_{\vu\boxplus\vI}\br{E^-_\vomega\cap\vJ}\abs{\vI}.
\end{align*}
The rest of the argument follows almost word-for-word the presentation in \textsc{Section \ref{subsubsec_gen_tree_large_scale}}. Namely, combining \eqref{eq_lem_gen_tree_c_half_large_scale_sum_mu_E} and \eqref{eq_lem_gen_tree_c_half_large_scale_sum_E} yields:
\begin{equation*}
    \sum_{\vJ\in\J}
    \sum_{k\geq k_\vJ+2}
    \sum_{\vI\times\vomega \in\cT_k}
        \abs{
            \ang{\cC_k f^{\epsilon,cen}_{\vI\times\vomega},\1_{E^-_\vomega\cap\vJ}}
        }
    \leq 
    \lambda^{-\frac{1}{2}}
    \cL
    \nrm{f}_{E\br{\cT}}
    \cM_\cT\br{E}\abs{\vI_\cT}.
\end{equation*}
This shows the central part contribution does not exceed the right-hand side of \eqref{eq_lem_Hdim_gen_tree}. 

As for the \(\ast=out \) case, we recall \eqref{eq_def_E_J_E_k} and rearrange the left-hand side of \eqref{eq_lem_Hdim_gen_tree_large_scale_split}:
\begin{equation*}
    \abs{
        \sum_{\vJ\in\J}\sum_{k\geq k_\vJ+2}
            \ang{\cC_k
            f^{\epsilon,out}_{\cT_k},
            \1_{E^-_k\cap\vJ}
        }
    },\quad
    f^{\epsilon,out}_{\cT_k}
    :=
    \sum_{P\in\cT_k}
        f^{\epsilon,out}_P.
\end{equation*}
By the Cauchy-Schwarz inequality, we dominate the above by
\begin{align*}
    &
    \sum_{\vJ\in\J}\sum_{k\geq k_\vJ+2}
    \nrm{
        \1_{E^-_k}
        \cC_k
        f^{\epsilon,out}_{\cT_k}
    }_{L^2\br{\vJ}}
    \abs{E^-_k\cap\vJ}^{\frac{1}{2}}
    \leq 
    \sum_{\vJ\in\J}
    \nrm{
        \nrm{
            \1_{E^-_k}
            \cC_k
            f^{\epsilon,out}_{\cT_k}
        }_{L^2\br{\vJ}}
    }_{\ell^2\br{k}}
    \abs{E_\vJ}^{\frac{1}{2}}\\
    = &
    \sum_{\vJ\in\J}
    \nrm{
        \nrm{
            \1_{E^-_k}
            \cC_k
            f^{\epsilon,out}_{\cT_k}
        }_{\ell^2\br{k}}
    }_{L^2\br{\vJ}}
    \abs{E_\vJ}^{\frac{1}{2}}
    \leq 
    \nrm{
        \nrm{
            \1_{E^-_k}
            \cC_k
            f^{\epsilon,out}_{\cT_k}
        }_{\ell^2\br{k}}
    }_{L^2}
    \Big\vert
        \bigsqcup_{\vJ\in\J}
            E_\vJ
    \Big\vert^{\frac{1}{2}}.
\end{align*}
By a line-by-line adaptation of the argument in \eqref{eq_lem_-_tree_large_scale_log_many_sqbd} with all the \(f\)-related terms replaced by the outer parts analog, we show the following estimate:
\begin{equation*}
    \nrm{
        \nrm{
            \1_{E^-_k}
            \cC_k
            f^{\epsilon,out}_{\cT_k}
        }_{\ell^2\br{k}}
    }_{L^2}
    \lesssim
    \nrm{\Lambda}_c
    \cM^c_\cT\br{E}
    \nrm{f^{out}}_{L^2\br{\cT}}.
\end{equation*}
Combined with \textbf{Lemma \ref{lem_E_J_den_est_shift}}, we derive
\begin{align*}
    \abs{
        \sum_{\vJ\in\J}\sum_{k\geq k_\vJ+2}
            \ang{\cC_k
            f^{\epsilon,out}_{\cT_k},
            \1_{E^-_k\cap\vJ}
        }
    }
    \lesssim &
    \nrm{\Lambda}_c
    \cL^{\frac{1}{2}}
    \nrm{f^{out}}_{L^2\br{\cT}}
    \cM^{c+\frac{1}{2}}_\cT\br{E}
    \abs{\vI_\cT}^{\frac{1}{2}}\\
    \lesssim &
    \nrm{\Lambda}_c
    \cL^{\frac{1}{2}}
    \nrm{f}_{E\br{\cT}}
    \cM^{c+\frac{1}{2}}_\cT\br{E}
    \abs{\vI_\cT}\,,
\end{align*}
thus showing that the outer part contribution is under control. This concludes the proof of \textbf{Lemma \ref{lem_Hdim_gen_tree}}.\qed

\subsection{Proof of \texorpdfstring{\eqref{eq_thm_Hdim_HF_st_mod}}{}: the \texorpdfstring{\(j=1\)}{} case}
As is discussed in \textsc{Section \ref{subsec_Hdim_issue_j=1}}, to addess the lack of the total \(L^2\) energy bound, we shall develop a suitable replacement of \eqref{eq_def_out_L2_infty_energy_emb}:

\begin{lemma}[\textsf{Energy selection for} \(\lambda\geq 1\)]\label{lem_Hdim_energy_sel}
    Given \(f\in L^2\br{\R^{D+1}}\) and a finite \(\cP\subset \P:=\P\br{\lambda}\), for any \(\epsilon,\varsigma>0\), there is a finite list of disjoint trees \(\BR{\cT_i}_i\) in \(\cP\) such that the following two estimates hold:
    \begin{equation*}
        \nrm{f}_{E\br{\cP\setminus\bigsqcup_i\cT_i}}\leq \varsigma \qquad\textrm{and}\qquad
        \sum_i\abs{\vI_{\cT_i}}\underset{\epsilon}{\lesssim} \lambda^\epsilon \cdot \frac{\nrm{f}_{L^2}^2}{\varsigma^2}.
    \end{equation*}
\end{lemma}
In what follows, we assume that the above lemma holds--- its proof will be delivered in \textsc{Section \ref{sec_pf_gen_energy_est}}---in order to complete the proof of \eqref{eq_thm_Hdim_HF_st_mod} in the \(j=1\) case. Let \(f\in L^2\br{\R^{D+1}}\) and \(E\subset \R^{D+1}\) be measurable. Due to the positivity of the statement \eqref{eq_thm_Hdim_HF_st_mod}, it suffices to show for any finite \(\cP\subset \P\) the following holds\footnote{In the expression below
$\sigma$ is the same parameter as the one appearing in \eqref{eq_lem_Hdim_gen_tree}.}:
\begin{equation}\label{eq_prop_Hdim_HF_st_sparse_truncate}
    \nrm{\Lambda^{\br{1}}_P\br{f,\1_E}}_{\ell^1\br{P\in\cP}}
    \lesssim
    \lambda^{-\sigma/2} 
    \log^2\br{e+\lambda\br{1+\abs{\vu}}}
    \nrm{f}_{L^2}\abs{E}^{\frac{1}{2}}.
\end{equation}
Without loss of generality, we make the following harmless assumptions:
\begin{equation}\label{eq_nonvanishing}
    0<\abs{E}<\infty,\quad
    \Lambda^{\br{1}}_P\br{f,\1_E}\neq 0,\quad \forall P\in\cP.
\end{equation}
We initiate the argument by decomposing \(\cP\) into collection of trees \(\BR{\cT_{n,i}}_{n,i}\) in the following manner:
\begin{enumerate}
    \item Set \(\cP_{n_0}:=\cP\) for \(n_0\in\Z\) being the largest number such that \(\nrm{f}_{E\br{\cP}}\leq 2^{-n_0/2}\nrm{f}_{L^2}\).
    \item Given \(\cP_{n-1}\), we apply \textbf{Lemma \ref{lem_Hdim_u_shift_mass_sel}} and the trivial bound to obtain a disjoint collection of trees \(\BR{\cT_{n,i}}_{i=1}^m\) in \(\cP_{n-1}\) such that the following two estimates hold:
    \begin{equation*}
        \cM_{\cP_{n-1}\setminus\bigsqcup_i\cT_{n,i}}\br{E}\lesssim \min\br{1, 2^{-n}\log\br{e+\abs{\vu}}\abs{E}},\quad
        \sum_{i=1}^m\abs{\vI_{\cT_{n,i}}}\lesssim 2^n.
    \end{equation*}
    \item By \textbf{Lemma \ref{lem_Hdim_energy_sel}}, there is a disjoint collection of trees \(\BR{\cT_{n,i}}_{i=m+1}^M\) in \(\cP_{n-1}\setminus\bigsqcup_{i=1}^m\cT_{n,i}\) such that the following two estimates hold:
    \begin{equation*}
        \nrm{f}_{E\br{
            \cP_{n-1}\setminus\bigsqcup_{i=1}^M\cT_{n,i}
        }}
        \underset{\epsilon}{\lesssim} \lambda^{\epsilon/2} \cdot 2^{-n/2}\nrm{f}_{L^2},\quad
        \sum_{i=m+1}^M\abs{\vI_{\cT_{n,i}}}\lesssim 2^n.
    \end{equation*}
    \item Set \(\cP_n:=\cP_{n-1}\setminus\bigsqcup_{i=1}^M\cT_{n,i}\), increase \(n\) by one, and go to step (2).
\end{enumerate}
By \eqref{eq_nonvanishing}, the algorithm will exhaust the collection of tiles. We thus obtain \(\cP=\bigsqcup_{n,i}\cT_{n,i}\) and
\begin{equation}\label{eq_Hdim_energy_mass_bd_in_proof}
    \nrm{f}_{E\br{\cT_{n,i}}}
    \underset{\epsilon}{\lesssim} \lambda^{\epsilon/2} \cdot 2^{-n/2}\nrm{f}_{L^2},\quad
    \cM_{\cT_{n,i}}\br{E}\lesssim
    \log\br{e+\abs{\vu}}\min\br{1,2^{-n}\abs{E}}
    .
\end{equation}
We now apply the decomposition to rewrite the model sum:
\begin{equation*}
    \nrm{\Lambda^{\br{1}}_P\br{f,\1_E}}_{\ell^1\br{P\in\cP}}=\sum_n \sum_i \nrm{\Lambda^{\br{1}}_P\br{f,\1_E}}_{\ell^1\br{P\in\cT_{n,i}}}.
\end{equation*}
By \textbf{Lemma \ref{lem_Hdim_gen_tree}}, we dominate the above by
\begin{equation*}
  \lambda^{-\sigma}\log\br{e+\lambda\br{1+\abs{\vu}}}
    \sum_{n,i}
        \nrm{f}_{E\br{\cT_{n,i}}}
        \cM^{c+\frac{1}{2}}_{\cT_{n,i}}\br{E}
        \abs{\vI_{\cT_{n,i}}}.
\end{equation*}
By \eqref{eq_Hdim_energy_mass_bd_in_proof}, we further dominate the above by
\begin{align*}
    &
    \lambda^{-\sigma}\log\br{e+\lambda\br{1+\abs{\vu}}}
    \lambda^{\epsilon/2}
    \sum_n
        2^{-n/2}\nrm{f}_{L^2}
        \br{\log\br{e+\abs{\vu}}\min\br{1,2^{-n}\abs{E}}}^{c+\frac{1}{2}}
        \sum_i
            \abs{\vI_{n,i}}\\
    \lesssim &
    \lambda^{\epsilon/2-\sigma}
    \log^2\br{e+\lambda\br{1+\abs{\vu}}}
    \sum_n
        2^{n/2}
        \nrm{f}_{L^2}
        \min\br{1,2^{-n\br{c+\frac{1}{2}}}\abs{E}^{c+\frac{1}{2}}}\\
    = &
    \lambda^{\epsilon/2-\sigma}
    \log^2\br{e+\lambda\br{1+\abs{\vu}}}
    \nrm{f}_{L^2}
    \abs{E}^{\frac{1}{2}}
    \sum_n
        \min\br{\br{\frac{2^n}{\abs{E}}}^{1/2},\br{\frac{\abs{E}}{2^n}}^c}\\
    \underset{c}{\eqsim} &
    \lambda^{\epsilon/2-\sigma}
    \log^2\br{e+\lambda\br{1+\abs{\vu}}}
    \nrm{f}_{L^2}
    \abs{E}^{\frac{1}{2}}.
\end{align*}
Taking now \(\epsilon=\sigma\) we prove \eqref{eq_prop_Hdim_HF_st_sparse_truncate} and thus \eqref{eq_thm_Hdim_HF_st_mod}.

\section{The energy selection algorithm for \texorpdfstring{\(\lambda\gtrsim 1\)}{}: Proof of Lemma \ref{lem_Hdim_energy_sel}}\label{sec_pf_gen_energy_est}

The proof of \textbf{Lemma \ref{lem_Hdim_energy_sel}} follows the strategy presented in \textsc{Section \ref{subsec_ver_tf_proj_emb}}. That is, we will first normalize the above formulation, design a suitable algorithm, and reduce the final estimate to a Bessel-type inequality similar to \textbf{Lemma \ref{lem_bessel}}. 

However, several standard notions appearing in the original argument need to be modified. Specifically, since the definition of \(\nrm{f}_{E\br{\cP}}\) involves frequency projections onto frequency sub-systems as in \textbf{Definition \ref{def_Hdim_cen_ext}}, the notion of the strong disjointness should be adapted to the individual frequency sub-tile location. In other words, the frequency support should be described by unions of frequency sub-intervals \(\fs\vomega\). In order to formalize this, we introduce the following 
\begin{definition}[\textsf{Frequency sub-tile assignment}]\label{def_Hdim_freq_asign}
    Given \(\cP\subset \P\),
    the following map 
    \begin{equation*}
        \cP\to 2^{\R^{D+1}}\times
        L^1\br{\R^{D+1}}:
        P\mapsto\br{\Omega_P,\Psi_P}    
    \end{equation*}
     is a frequency sub-tile assignment on \(\cP\) if for all \(P:=\vI\times\vomega\in\cP\) there is \(\cW_P\subset \fs^{\br{+,1}}\br{\vomega}\) such that 
     \begin{equation}\label{eq_def_Hdim_freq_asign_freq}
         \abs{\Psi_P}\ast\chi^N_\vI\underset{N}{\lesssim}\chi^N_\vI
         ,\quad
         \varnothing\neq \supp\widehat{\Psi_P}\subset \Omega_P:=\bigsqcup\cW_P\subset \vomega^{\br{+,1}}.
     \end{equation}
\end{definition}

\begin{observation}\label{obs_Hdim_subfreq_inclus}
    Fix \(\br{\Omega_{\br{\cdot}},\Psi_{\br{\cdot}}}\) a frequency sub-tile assignment on \(\cP\subset \P\). For any \(P:=\vI\times\vomega,P':=\vI'\times\vomega'\in\cP\) with \(\Omega_P\cap\Omega_{P'}\neq\varnothing\), if \(\abs{\vomega}<\abs{\vomega'}\), there is \(\fs\vomega'\in\cW_{P'}\) such that the following relations hold:
    \begin{equation*}
        \Omega_P\subset \vomega^{\br{+,1}}\subset 
        \fs\vomega' \subset \Omega_{P'}\subset \vomega^{\prime \br{+,1}}
    \end{equation*}
    In other words, the collection \(\BR{\Omega_P}_{P\in\cP}\) has a behavior resembling the usual dyadic system. 
\end{observation}

We now adapt the notion of strong disjointness to a frequency sub-tile assignment:
\begin{definition}[\(C\)-\textsf{strong disjointness}]
    Given a large constant \(C\gg 1\), a sequence of trees \(\BR{\cT_i}_i\) in \(\P\) and a frequency sub-tile assignment \(\br{\Omega_{\br{\cdot}},\Psi_{\br{\cdot}}}\) on \(\bigcup_i\cT_i\), we say that the sequence of trees \(\BR{\cT_i}_i\) is \(C\)-strongly disjoint with respect to \(\br{\Omega_{\br{\cdot}},\Psi_{\br{\cdot}}}\) if for all pairs of tiles \(\br{P,P'}:=\br{\vI\times\vomega,\vI'\times\vomega'} \in\cT_m\times\cT_n\) with \(\Omega_{P}\cap\Omega_{P'}\neq \varnothing\) the following implication holds:
    \begin{equation}\label{strdisj}
        \br{
            \abs{\vomega}<\abs{\vomega'}\text{ or }
            m<n
        }\implies 
        C \vI_{\cT_m}\cap \vI'=\varnothing.
    \end{equation}
\end{definition}
We are now ready to introduce a Bessel-type inequality that generalizes \textbf{Lemma \ref{lem_bessel}}:
\begin{lemma}[\textsf{Generalized Bessel-type inequality}]\label{lem_Hdim_bessel}
    Let \(\lambda\gg 1\gg \epsilon>0\) and \(f\in L^2\br{\R^{D+1}}\). Given a sequence of trees \(\BR{\cT_i}_i\) in \(\P:=\P\br{\lambda}\) and a frequency sub-tile assignment \(\br{\Omega_{\br{\cdot}},\Psi_{\br{\cdot}}}\) on \(\bigcup_i\cT_i\),
    if the sequence of trees \(\BR{\cT_i}_i\) is \(\lambda^\epsilon\)-strongly disjoint with respect to \(\br{\Omega_{\br{\cdot}},\Psi_{\br{\cdot}}}\) and satisfies the two estimates:
    \begin{equation}\label{eq_lem_Hdim_bessel_tile_ass}
        \nrm{\br{\Psi_P\ast f}^{\br{+,1}}}_P\lesssim \abs{\vI}^{\frac{1}{2}},\quad \forall P:=\vI\times\vomega\in\bigcup_i\cT_i;
    \end{equation}
    \begin{equation}\label{eq_lem_Hdim_bessel_tree_ass}
        \nrm{\nrm{\br{\Psi_P\ast f}^{\br{+,1}}}_P}_{\ell^2\br{P\in\cT_i}}\eqsim \abs{\vI_{\cT_i}}^{\frac{1}{2}},\quad \forall i,
    \end{equation}
    the following estimate holds:
    \begin{equation}\label{eq_lem_Hdim_bessel_concl}
        \sum_i\abs{\vI_{\cT_i}}
        \underset{\epsilon}{\lesssim} \nrm{f}_{L^2}^2\,.
    \end{equation}
\end{lemma}
We will prove \textbf{Lemma \ref{lem_Hdim_energy_sel}} assuming \textbf{Lemma \ref{lem_Hdim_bessel}}.
The proof of \textbf{Lemma \ref{lem_Hdim_bessel}} will be postponed to the end of the section.

Let \(f\in L^2\br{\R^{D+1}}\) and \(\cP\subset \P\)  finite be given. It suffices to show that for any small \(\epsilon>0\) there is a collection of disjoint trees \(\T\) in \(\cP\) such that the following two estimates hold:
\begin{equation*}
    \nrm{f}_{E\br{\cP\setminus\bigsqcup\T}}\leq
    \nrm{f}_{E\br{\cP}}/2,\quad
    \sum_{\cT\in\T}\abs{\vI_{\cT}}\underset{\epsilon}{\lesssim} \lambda^{\epsilon\br{D+1}}
    \br{
        \nrm{f}_{L^2}/\nrm{f}_{E\br{\cP}}
    }^2.
\end{equation*}
Without loss of generality, we impose the normalization \(\nrm{f}_{E\br{\cP}}=1\).
We aim to extract sequences of \(\lambda^\epsilon\)-strongly disjoint trees \(\BR{\cT_i}_i\) with respect to certain frequency sub-tile assignment whose contribution corresponding to the first estimate exceeds the threshold. The actual collection of trees \(\T\) will be associated to \(\BR{\cT_i}_i\). To proceed, recalling \textbf{Definition \ref{def_Hdim_energy}}, we deduce
\begin{equation*}
    \nrm{f}_{E\br{\cP\setminus\bigsqcup\T}}\leq \frac{1}{2} \iff
    \sup_{\substack{
        \cT\subset \cP\setminus \bigsqcup\T\\
        \cT\text{ a tree}
    }}\nrm{f^{out}}_{L^2\br{\cT}}/\abs{\vI_\cT}^{\frac{1}{2}}
    \leq \frac{1}{2} \quad \text{and}\quad
    \sup_{\vI\times\vomega\in \cP\setminus \bigsqcup\T}
    \nrm{
        f^{cen}
    }_{\vI\times\vomega}/\abs{\vI}^{\frac{1}{2}}
    \leq \frac{1}{2}.
\end{equation*}
We will address the two components separately.
\subsection{The central-part algorithm}
Starting with the central parts,
\begin{itemize}
    \item Initialize \(i=1\), \(\cP_0:=\cP\), and \(\T_0=\varnothing\).
        \item If \(\cP_{i-1}=\varnothing\) or \(\sup_{\vI\times\vomega\in \cP_{i-1}}
    \nrm{
        f^{cen}
    }_{\vI\times\vomega}/\abs{\vI}^{\frac{1}{2}}\leq\frac{1}{2}\),
        we terminate the algorithm. Otherwise, the set
        \begin{equation*}
        \cP^\ast_i:=
        \BR{
            \vI\times\vomega\in\cP_{i-1}
        \::\:
            \nrm{
                f^{cen}
            }_{\vI\times\vomega}/\abs{\vI}^{\frac{1}{2}}
            > \frac{1}{2}
        }
        \end{equation*}
        is non-empty. Due to the finiteness of \(\cP^\ast_i\subset\cP_{i-1}\subset \cP\), there is \( P_i:=\vI_i\times\vomega_i\in \cP^\ast_i\) with maximal physical interval size, that is
        \begin{equation}\label{eq_lem_Hdim_sel_I_longest}
            \sup_{\vI\times\vomega\in\cP^\ast_i} \abs{\vI}\leq\abs{\vI_i}.
        \end{equation}
    \item Given such \(P_i\in \cP^\ast_i\), we construct trees of the following form:
        \begin{equation*}
            \cT_{i,\ve}:=\BR{
                P\in\cP_{i-1}\::\:
                P\leq 
                \br{\ve,\vnull}\boxplus P_i
                =\br{\ve\boxplus \vI_i}\times\vomega_i}
            ,\quad
            \ve\in \vZ_\epsilon:=\br{\Z\cap
            \mr{-2\lambda^\epsilon,2\lambda^\epsilon}}^{D+1},
        \end{equation*}
        where we assign \(\br{\ve,\vnull}\boxplus P_i\) to be the tree top of \(\cT_{i,\ve}\).
    \item For the next iteration, we set:
    \begin{equation*}
        \T_i:=\T_{i-1}\sqcup\BR{\cT_{i,\ve}}_{\ve\in \vZ_\epsilon},\quad
        \cP_i:=\cP_{i-1}\setminus \bigsqcup_{\ve\in\vZ_\epsilon} \cT_{i,\ve}=\cP\setminus\bigsqcup \T_i,
    \end{equation*}
    increase \(i\) by one, and return to the second step.
    \end{itemize}
When the algorithm terminates at some \(i=n\), we set
 \begin{equation}\label{centrpart}
 \T^{cen}:=\T_n\,. 
 \end{equation}
 
By construction, 
\begin{equation*}
    \sup_{\vI\times\vomega\in \cP\setminus \bigsqcup\T^{cen}}
    \nrm{
        f^{cen}
    }_{\vI\times\vomega
    }/\abs{\vI}^{\frac{1}{2}}
    \leq \frac{1}{2}.
\end{equation*}
To obtain control over the sum of the tree top physical interval:
\begin{equation*}
    \sum_{\cT\in\T^{cen}}\abs{\vI_\cT}\underset{\epsilon}{\lesssim }\lambda^{\epsilon\br{D+1}} \nrm{f}_{L^2}^2,
\end{equation*}
we invoke the relation between \(\cT_{i,\ve}\) and \(P_i:=\vI_i\times\vomega_i\):
\begin{equation*}
    \sum_{\cT\in\T^{cen}}\abs{\vI_\cT}=
    \sum_i\sum_{\ve\in\vZ_\epsilon}\abs{\vI_{\cT_{i,\ve}}}= \#\vZ_\epsilon \sum_i\abs{\vI_{\cT_{i,\vnull}}}\eqsim
    \lambda^{\epsilon\br{D+1}}
    \sum_i\abs{\vI_i}.
\end{equation*}
It thus suffices to show:
\begin{equation}\label{eq_Hdim_energy_sel_cen_I_sum}
    \sum_i\abs{\vI_i}\underset{\epsilon}{\lesssim}\nrm{f}_{L^2}^2.
\end{equation}
We will now set up a frequency sub-tile assignment for this specific setting and apply the generalized Bessel inequality \textbf{Lemma \ref{lem_Hdim_bessel}} to prove the above estimate. 

Observe that, by \textbf{Definition \ref{def_Hdim_cen_ext}}, for \(P_i:=\vI_i\times\vomega_i\in\P\) there is \(\fs_{v_i}\vomega_i\in\fs^{\br{+,1}}\br{\vomega_i}\) such that
\begin{equation*}
    \frac{1}{2}<\nrm{\fp_{\fs_{v_i}\vomega_i}f}_{P_i}/\abs{\vI_i}^{\frac{1}{2}}
    \leq
    \nrm{f^{cen}}_{P_i}/\abs{\vI_i}^{\frac{1}{2}}
    \leq\nrm{f}_{E\br{\cP}}=1
\end{equation*}
This suggests that we consider trivial trees \(\BR{\cT_i}_i:=\BR{\cT_{P_i}}:=\BR{\BR{P_i}}_i\) and set the frequency sub-tile assignment for \(P_i\in \P_k\) the following set \(\Omega_{P_i}\) and function \(\Psi_{P_i}\):
\begin{equation}\label{eq_lem_Hdim_sel_cen_assign}
    \Omega_{P_i}:=\br{\fs_{v_i-1}\vomega_i\sqcup\fs_{v_i}\vomega_i\sqcup \fs_{v_i+1}\vomega_i}\cap\vomega^{\br{+,1}}_i,\quad
    \widehat{\Psi_{P_i}}\br{\vxi}:= 
    \br{\Dil^\infty_2\varphi}_{\fs_{v_i}\vomega_i}\br{\vxi}\phi\br{2^k\xi_0/\1_\lambda}.
\end{equation}
By the Heisenberg principle, one easily verifies the first condition in \eqref{eq_def_Hdim_freq_asign_freq}. Since we have:
\begin{equation*}
    \nrm{f^{cen}}_{P_i}\eqsim
    \nrm{\br{\fp_{\fs_{v_i}\vomega_i}f}^{\br{+,1}}}_{P_i}=
    \nrm{\br{\Psi_{P_i}\ast f}^{\br{+,1}}}_{P_i}\eqsim \abs{\vI_i}^{\frac{1}{2}},
\end{equation*}
once we demonstrate the sequence of trees \(\BR{\cT_i}_i\) is \(\lambda^\epsilon\)-strongly disjoint with respect to \eqref{eq_lem_Hdim_sel_cen_assign}, we may apply \textbf{Lemma \ref{lem_Hdim_bessel}} and conclude \eqref{eq_Hdim_energy_sel_cen_I_sum}. To show the \(\lambda^\epsilon\)-strongly disjointness assumption, we take a pair \(P_m,P_n\), with \(\Omega_{P_m}\cap\Omega_{P_n}\neq \varnothing\) and verify the validity of \eqref{strdisj}.

 If \(m<n\), the above algorithm implies that \(P_n\notin \sqcup\cP^\ast_m\). In particular, we have:
\begin{equation}\label{eq_lem_Hdim_sel_cen_PlessP}
    P_n \not\leq \br{\ve,\vnull}\boxplus P_m
    =\br{\ve\boxplus \vI_m}\times\vomega_m,\quad \forall \ve\in\vZ_\epsilon.
\end{equation}
Additionally, due to the maximality assumption \eqref{eq_lem_Hdim_sel_I_longest}, we also have \(\abs{\vI_n}\leq \abs{\vI_m}\). By the fact that
\begin{equation*}
    \varnothing\neq\Omega_{P_m}\cap\Omega_{P_n}\subset 
    \vomega^{\br{+,1}}_m \cap\vomega^{\br{+,1}}_n=\vomega^{\br{+,1}}_m\subset \vomega^{\br{+,1}}_n,
\end{equation*}
we conclude that \eqref{eq_lem_Hdim_sel_cen_PlessP} can only happen when:
\begin{equation*}
    \lambda^\epsilon\vI_m\cap \vI_n\subset 
    \vI_n\cap \bigsqcup_{\ve\in\vZ_\epsilon}\ve\boxplus \vI_m=\varnothing.
\end{equation*}
Alternatively, if \(\abs{\vomega_m}<\abs{\vomega_n}\), we have \(\abs{\vI_m}>\abs{\vI_n}\).
By the maximality assumption \eqref{eq_lem_Hdim_sel_I_longest}, we must have \(m<n\) and thus \(\lambda^\epsilon\vI_m\cap \vI_n=\varnothing\) by the previous discussion. This demonstrates that the sequence of trivial trees \(\BR{\cT_i}_i:=\BR{\BR{P_i}}_i\) is \(\lambda^\epsilon\)-strongly disjoint with respect to \eqref{eq_lem_Hdim_sel_cen_assign} completing the central part treatment. 

\subsection{The outer-part algorithm}
We now focus on treating the outer parts. For the sake of discussion, we set\footnote{Recall \eqref{centrpart}.} \(\cP^{cen}:=\cP\setminus\bigsqcup\T^{cen}\). We aim to extract a finite collection of disjoint trees \(\T^{out}\) out of \(\cP^{cen}\) such that:
\begin{equation*}
    \sup_{\substack{
        \cT\subset \cP^{cen}\setminus\bigsqcup\T^{out}\\
        \cT\text{ be a tree}
    }}
    \nrm{f^{out}}_{L^2\br{\cT}}/\abs{\vI_\cT}^{\frac{1}{2}}\leq \frac{1}{2}\qquad\textrm{and}\qquad
    \sum_{\cT\in\T^{out}}\abs{\vI_\cT}
    \lesssim \lambda^{\epsilon\br{D+1}}\nrm{f}^2_{L^2}.
\end{equation*}
For the sake of analysis, we will further split the outer part contributions into three sub-components. 
To motivate this latter decomposition, we shall first study the related expression involved in \textbf{Definition \ref{def_Hdim_cen_ext}}. Given a tree \(\cT\subset \P\) and a tile \(\vI\times\vomega\in\cT_k:=\cT\cap\P_k\), recall that:
\begin{equation*}
    f^{\cT_{out}}_{\vI\times\vomega}:=
    \br{
        \fp_{\setminus\fs_\cT\vomega}f
    }^{\br{+,1}}_{\vI\times\vomega}
    :=
    \chi^{\Phi}_\vI\cdot
    \pi^{\br{+,1}}_\vomega
    \fp_{\setminus\fs_\cT\vomega}f
    .
\end{equation*}
For the sake of analysis, we insert an additional projection and identify the above with:
\begin{equation*}
    =
    \chi^{\Phi}_\vI\cdot
    \pi^{\br{+,1}}_\vomega p^{\br{+,1}}_\vomega
    \fp_{\setminus\fs_\cT\vomega}f,\quad
    \widehat{p^{\br{+,1}}_\vomega g}:=\br{\Dil^{\infty}_2\phi}_{\vomega^{\br{+,1}}} \widehat{g}.
\end{equation*}
The decomposition of the outer part will be executed at the level of the expression \(p^{\br{+,1}}_\vomega\fp_{\setminus\fs_\cT\vomega}f\). Note that
\begin{equation*}
    \supp\br{p^{\br{+,1}}_\vomega\fp_{\setminus\fs_\cT\vomega}f}^{\wedge}
    \subset \vomega^{\br{+,1}}\setminus\fs_\cT\vomega.
\end{equation*}
Naively, we would like to decompose the above function so that the corresponding frequency support is subdivided into the three pieces shown below:
\begin{align*}
    \vomega^{\br{+,1}}\setminus\fs_\cT\vomega
    =&
    \br{
        \vomega^{\br{+,1}}
        \setminus\fS\vomega^{\br{+,1}}
    }
    \sqcup
    \br{
        \fS\vomega^{\br{+,1}}\setminus
        \fs_\cT\vomega
    }
    \\
    =&
    \br{
    \vomega^{\br{+,1}}\setminus\fS\vomega^{\br{+,1}}
    }
    \sqcup
    \br{
        \bigsqcup_{v_\cT<v<2^{k_\lambda}}
        \fs_v\vomega
    }
    \sqcup
    \br{
        \bigsqcup_{2^{k_\lambda-1}-1\leq v<v_\cT}
        \fs_v\vomega
    }.
\end{align*}
Guided by the above, we perform the actual decomposition on the function:
\begin{equation*}
    p^{\br{+,1}}_\vomega
    \fp_{\setminus\fs_\cT\vomega}f
    :=
    p^{\br{+,1}}_\vomega
    \br{\id-\fp_{\fs_\cT\vomega}}
    f.
\end{equation*}
Recalling \eqref{eq_Hdim_freq_proj_for_subsys}, we insert an additional projection \(\fp_k\) and after a telescoping we equate the above with
\begin{equation}\label{eq_Hdim_out_split_func_ini}
=
    p^{\br{+,1}}_\vomega
    \br{\id-\fp_k}
    f
    +
    p^{\br{+,1}}_\vomega
    \br{
        \fp_k-
        \fp_{\fs_\cT\vomega}
    }
    f.
\end{equation}
Since \(\vomega\in\W_k\), we have for the first term in \eqref{eq_Hdim_out_split_func_ini} satisfies
\begin{equation*}
    \supp\br{p^{\br{+,1}}_\vomega
    \br{\id-\fp_k}
    f}^{\wedge}
    \subset
    \BR{
        \vxi\in\vomega^{\br{+,1}}
    \::\:
        \abs{\xi_0}\gtrsim 2^k
    }
    \approx
    \vomega^{\br{+,1}}\setminus\fS\vomega^{\br{+,1}}\,,
\end{equation*}
while the second term in \eqref{eq_Hdim_out_split_func_ini} verifies
\begin{equation}\label{eq_Hdim_out_disconnect_comp}
    \supp
    \br{
    p^{\br{+,1}}_\vomega
    \br{
        \fp_k-
        \fp_{\fs_\cT\vomega}
    }
    f}^{\wedge}
    \subset 
    \br{
        \bigsqcup_{v_\cT<v<2^{k_\lambda}}
            \fs_v\vomega
    }
    \sqcup
    \br{
        \bigsqcup_{2^{k_\lambda-1}-1\leq v<v_\cT}\fs_v\vomega
    }.
\end{equation}
To separate the two disconnected components in \eqref{eq_Hdim_out_disconnect_comp}, we introduce the following 
\begin{definition}[\textsf{Auxiliary notions for frequency projections}]\label{def_Hdim_freq_proj_sub_sys_aux}
    Recalling \textbf{Definition \ref{def_Hdim_freq_proj_sub_sys}}, for
    given \(\vomega\in\W_k\) and \(\fs_v\vomega\in\fs\br{\W_k}\), we define 
    \begin{equation*}
        \widehat{\fp_{\fs^+_v\vomega}f}\br{\vxi}
        :=
        \1_{\bigsqcup_{v'>v}\fs_{v'}\vomega}\br{\vxi}
        \widehat{\br{\fp_k-\fp_{\fs_v\vomega}}f}\br{\vxi},
    \end{equation*}
    \begin{equation*}
        \widehat{\fp_{\fs^-_v\vomega}f}\br{\vxi}
        :=
        \1_{\bigsqcup_{v'<v}\fs_{v'}\vomega}\br{\vxi}
        \widehat{\br{\fp_k-\fp_{\fs_v\vomega}}f}\br{\vxi}.
    \end{equation*}
\end{definition}
With the above notions, we can now write
\begin{equation*}
    p^{\br{+,1}}_\vomega
    \fp_{\setminus\fs_\cT\vomega}f=
    p^{\br{+,1}}_\vomega
    \br{\id-\fp_k}
    f
    +
    p^{\br{+,1}}_\vomega
    \fp_{\fs^+_{v_\cT}\vomega}
    f
    +
    p^{\br{+,1}}_\vomega
    \fp_{\fs^-_{v_\cT}\vomega}
    ,
\end{equation*}
where by design, the last two terms satisfy:
\begin{equation*}
    \supp 
    \br{
        p^{\br{+,1}}_\vomega
        \fp_{\fs^+_{v_\cT}\vomega}
        f
    }^{\wedge}
    \subset
    \bigsqcup_{v_\cT<v<2^{k_\lambda}}
        \fs_v\vomega,\quad
    \supp 
    \br{
        p^{\br{+,1}}_\vomega
        \fp_{\fs^+_{v_\cT}\vomega}
        f
    }^{\wedge}
    \subset
    \bigsqcup_{2^{k_\lambda-1}-1\leq v<v_\cT}\fs_v\vomega.
\end{equation*}
We now decompose the outer part accordingly.

\begin{definition}[\textsf{Auxiliary notions for outer parts}]\label{def_Hdim_aux_out}
    Given a tree \(\cT\subset \P\), a collection \(\cP\subset \P\), and a tile \(P\equiv\vI\times\vomega\in\cT\) with \(P\in\P_k\), we decompose \(f^{\cT_{out}}_{\vI\times\vomega}\) into three sub-components:
    \begin{equation*}
        f^{\cT_{out}}_{\vI\times\vomega}=f^{out,0}_{\vI\times\vomega}
        +
        f^{\cT_{out,+}}_{\vI\times\vomega}
        +
        f^{\cT_{out,-}}_{\vI\times\vomega}
    \end{equation*}
    with each of them accompanied by a norm quantity as shown below:
    \begin{equation*}
        f^{out,0}_P:=
        \br{\br{\id-\fp_k}f}^{\br{+,1}}_P,\quad
        \nrm{f^{out,0}}_P:=\nrm{\br{\br{\id-\fp_k}f}^{\br{+,1}}}_P;
    \end{equation*}
    \begin{equation*}
        f^{\cT_{out},\pm}_{\vI\times\vomega}:=
        \br{
            \fp_{\fs^\pm_{v_\cT}\vomega}
        f}^{\br{+,1}}_{\vI\times\vomega},\quad
        \nrm{f^{\cT_{out,\pm}}}_{\vI\times\vomega}:=\nrm{
            \br{
                \fp_{\fs^\pm_{v_\cT}\vomega}f
            }^{\br{+,1}}
        }_{\vI\times\vomega}.
    \end{equation*}
    By extension, we have the following modified notion of \(L^2\) energies:
    \begin{equation*}
        \nrm{f^{out,0}}_{L^2\br{\cP}}:=
        \nrm{
        \nrm{f^{out,0}}_P
        }_{\ell^2\br{P\in\cP}},\quad
        \nrm{f^{out,\pm}}:=
        \nrm{
            \nrm{
                f^{\cT_{out,\pm}}
            }_P
        }_{\ell^2\br{P\in\cT}}.
    \end{equation*}
\end{definition}
Note that by construction, we have the following trivial estimate:
\begin{equation*}
    \nrm{f^{out}}_P
    \leq
    \nrm{f^{out,0}}_P
    +\nrm{f^{\cT_{out,+}}}_P
    +\nrm{f^{\cT_{out,-}}}_P
    ,\quad P\in\cT;
\end{equation*}
\begin{equation*}
    \nrm{f^{out}}_{L^2\br{\cT}}\leq \sum_{\ast \in\BR{0,+,-}}\nrm{f^{out,\ast}}_{L^2\br{\cT}}.
\end{equation*}
As a result, the problem reduces to extracting a finite collection of disjoint trees \(\T^{out}=\T^{out,0}\sqcup \T^{out,+}\sqcup \T^{out,-}\) out of \(\cP^{cen}\) such that for \(\ast\in\BR{0,+,-}\) the following holds:
\begin{equation*}
    \sup_{\substack{
        \cT\subset \cP^{cen}\setminus\bigsqcup\T^{out}\\
        \cT\text{ a tree}
    }}
    \nrm{
        f^{out,\ast}
    }_{L^2\br{\cT}}
    /\abs{\vI_{\cT}}^{\frac{1}{2}}\leq \frac{1}{6},\quad
    \sum_{\cT\in\T^{out,\ast}}\abs{\vI_\cT}\underset{\epsilon}{\lesssim}
    \lambda^{\epsilon\br{D+1}}
    \nrm{f}^2_{L^2}.
\end{equation*}

\subsubsection{The (out,\texorpdfstring{\(0\)}{}) contribution}
We first address the \(\ast=0\) case. 
We claim that a version of total \(L^2\) energy bound holds:
\begin{equation}\label{eq_Hdim_L2_out0_total}
    \nrm{f^{out,0}}_{L^2\br{\P}}
    \lesssim
    \nrm{f}_{L^2}.
\end{equation}
Indeed, from \textbf{Definition \ref{def_Hdim_aux_out}} we have
\begin{equation*}
    \nrm{f^{out,0}}_{L^2\br{\P}}
    :=
    \nrm{
        \nrm{
            f^{out,0}
        }_P
    }_{\ell^2\br{P\in\P}}
    :=
    \nrm{
        \nrm{
            \nrm{
                \br{\br{\id-\fp_k} f}^{\br{+,1}}
            }_{P}
        }_{\ell^2\br{P\in\P_k}}
    }_{\ell^2\br{k\in l+k_\lambda\Z}}
\end{equation*}
By \eqref{eq_ass_sum_I}, we further dominate the above by
\begin{equation*}
 \nrm{
        \nrm{
            \nrm{
                \pi_{\vomega^{\br{+,1}}}
                \br{\id-\fp_k}
                f
            }_{L^2}
        }_{\ell^2\br{\vomega\in\W_k}}
    }_{\ell^2\br{k\in l+k_\lambda\Z}}
    \leq
    \nrm{
        \nrm{
            \br{\fp_{k-2k_\lambda}-\fp_k}
            f
        }_{L^2}
    }_{\ell^2\br{k\in l+k_\lambda\Z}}
    \lesssim\nrm{f}_{L^2}.
\end{equation*}

Given the total \(L^2\) energy bound \eqref{eq_Hdim_L2_out0_total}, we now choose via a greedy algorithm a maximal disjoint collection of trees \(\T^{out,0}\) out of \(\cP^{cen}\) which satisfies the following condition:
\begin{equation*}
    \nrm{f^{out,0}}_{L^2\br{\cT}}/\abs{\vI_\cT}^{\frac{1}{2}}>\frac{1}{6},\quad
    \forall \cT\in\T^{out,0}.
\end{equation*}
On the one hand, from the above construction, we trivially have
\begin{equation*}
    \sup_{\substack{
        \cT\subset \cP^{cen}\setminus\bigsqcup\T^{out,0}\\
        \cT\text{ a tree}
    }}
    \nrm{f^{out,0}}_{L^2\br{\cT}}/\abs{\vI_\cT}^{\frac{1}{2}}\leq\frac{1}{6}.
\end{equation*}
On the other hand, from Chebyshev's inequality, we get
\begin{equation*}
    \sum_{\cT\in\T^{out,0}}\abs{\vI_\cT}
    <\frac{1}{36}
    \sum_{\cT\in\T^{out,0}}
    \nrm{f^{out,0}}_{L^2\br{\cT}}^2
    \lesssim
    \nrm{
        f^{out,0}
    }^2_{L^2\br{\bigsqcup\T^{out,0}}}
    \leq
    \nrm{f^{out,0}}^2_{L^2\br{\P}}
    \lesssim
    \nrm{f}^2_{L^2}.
\end{equation*}
This completes the treatment of the \(\ast=0\) case. We now set \(\cP^{out,0}:=\cP^{cen}\setminus\bigsqcup \T^{out,0}\).

\subsubsection{The (out,\texorpdfstring{\(\pm\)}{}) contribution}
Treatment of the \(\ast=+\) and \(\ast=-\) cases are almost identical. We will demonstrate the argument for the \(\ast=+\) case and only make some brief comments on the required modifications for the \(\ast=-\) case. To reiterate, we aim to find disjoint collection of trees \(\T^{out,+}\) in \(\cP^{out,0}\) such that:
\begin{equation*}
    \sup_{\substack{
        \cT\subset \cP^{out,0}\setminus\bigsqcup\T^{out,+}\\
        \cT\text{ be a tree}
    }}
    \nrm{
        f^{out,+}
    }_{L^2\br{\cT}}
    /\abs{\vI_{\cT}}^{\frac{1}{2}}\leq \frac{1}{6},\quad
    \sum_{\cT\in\T^{out,+}}\abs{\vI_\cT}\underset{\epsilon}{\lesssim}
    \lambda^{\epsilon\br{D+1}}
    \nrm{f}^2_{L^2}.
\end{equation*}
Our selection of the trees is a result of the following algorithm:
\begin{itemize}
    \item Initialize \(i=1\), \(\cP_0:=\cP^{out,0}\), and \(\T_0=\varnothing\).
    \item Consider the following set:
        \begin{equation*}
        \T^\ast_i:=
        \BR{
            \cT\subset \cP_{i-1}
        \::\:
            \cT\text{ a tree},\quad
            \nrm{
                f^{out,+}
            }_{L^2\br{\cT}}
            /\abs{\vI_{\cT}}^{\frac{1}{2}}> \frac{1}{6}
        }
        \end{equation*}
        If \(\T^\ast_i=\varnothing\), we terminate the algorithm. Otherwise, we choose among all trees in \(\T^\ast_i\) the one with the top having its \(1\)st frequency component being the "rightmost". To be more precise, since the set \(\T^\ast_i\) is finite, there is a tree \(\cT_i\in\T^\ast_i\) with top \(P_{\cT_i}:=\vI_{\cT_i}\times\vomega_{\cT_i}\in\P\) such that:
        \begin{equation}\label{eq_lem_Hdim_sel_rightmost}
            \cT\in\T^\ast_i\implies
            \inf\pi_{\R\ve_1}\vomega^{\br{+,1}}_\cT \leq \inf\pi_{\R\ve_1}\vomega^{\br{+,1}}_{\cT_i}.
        \end{equation}
        
    \item Given such \(\cT_i\in \T^\ast_i\), we construct trees of the following form:
        \begin{equation*}
            \cT_{i,\ve}:=\BR{
                P\in\cP_{i-1}\::\:
                P\leq 
                \br{\ve,\vnull}\boxplus P_{\cT_i}}
            ,\quad
            \ve\in \vZ_\epsilon:=\br{\Z\cap
            \mr{-2\lambda^\epsilon,2\lambda^\epsilon}}^{D+1},
        \end{equation*}
        where we assign \(\br{\ve,\vnull}\boxplus P_{\cT_i}\) to be the tree top of \(\cT_{i,\ve}\).
    \item For the next iteration, we set:
    \begin{equation*}
        \T_i:=\T_{i-1}\sqcup\BR{\cT_{i,\ve}}_{\ve\in \vZ_\epsilon},\quad
        \cP_i:=\cP_{i-1}\setminus \bigsqcup_{\ve\in\vZ_\epsilon} \cT_{i,\ve}=\cP\setminus\bigsqcup \T_i,
    \end{equation*}
    increase \(i\) by one, and return to the second step.
\end{itemize}
When the algorithm terminates at some \(i=n\), we set \(\T^{out,+}:=\T_n\). By construction,
\begin{equation*}
    \sup_{\substack{
        \cT\subset \cP^{out,0}\setminus\bigsqcup\T^{out,+}\\
        \cT\text{ a tree}
    }}
    \nrm{
        f^{out,+}
    }_{L^2\br{\cT}}
    /\abs{\vI_{\cT}}^{\frac{1}{2}}\leq \frac{1}{6}.
\end{equation*}
It remains to prove the following estimate:
\begin{equation*}
    \sum_{\cT\in\T^{out,+}}\abs{\vI_\cT}\underset{\epsilon}{\lesssim }\lambda^{\epsilon\br{D+1}} \nrm{f}_{L^2}^2.
\end{equation*}
By invoking the relation between \(\cT_{i,\ve}\) and \(P_{\cT_i}\), we have
\begin{equation*}
    \sum_{\cT\in\T^{out,+}}\abs{\vI_\cT}=
    \sum_i\sum_{\ve\in\vZ_\epsilon}\abs{\vI_{\cT_{i,\ve}}}= \#\vZ_\epsilon \sum_i\abs{\vI_{\cT_i}}\eqsim
    \lambda^{\epsilon\br{D+1}}
    \sum_i\abs{\vI_{\cT_i}}.
\end{equation*}
It thus suffices to show
\begin{equation}\label{eq_Hdim_energy_sel_out_I_sum}
    \sum_i\abs{\vI_{\cT_i}}\underset{\epsilon}{\lesssim}\nrm{f}_{L^2}^2.
\end{equation}
We will again set up a frequency sub-tile assignment on \(\bigsqcup_i \cT_i\), show the \(\BR{\cT_i}_i\) is \(\lambda^\epsilon\)-strongly disjoint with respect to the assignment, and apply the generalized Bessel-type inequality to prove the above estimate. 

Given \(P:=\vI\times\vomega\in \bigsqcup_i\cT_i\), there is a unique \(i\) such that \(P\in\cT_i\) by the algorithm. Let \(k\) be such that \(P\in\cT\cap\P_k\). We now recall \eqref{eq_Hdim_out_disconnect_comp} and \textbf{Definitions \ref{def_Hdim_freq_proj_sub_sys_aux}} and \textbf{\ref{def_Hdim_aux_out}} and assign:
\begin{equation}\label{eq_lem_Hdim_sel_out_assign}
    \Omega_P:=
    \bigsqcup_{v_{\cT_i}<v<2^{k_\lambda}}\fs_v\vomega
    ,\quad
    \Psi_P:=
        \fp_{\fs^+_{v_{\cT_i}}\vomega}
        \br{\br{\Dil^\infty_2\phi}_{\vomega^{\br{+,1}}}}^\vee
        .
\end{equation}
By design, \(\supp\widehat{\Psi_P}\subset \Omega_P\), while the fact that \(\Psi_P\ast\chi^N_\vI\underset{N}{\lesssim}\chi^N_\vI\) follows from the Heisenberg principle. This justifies that the above is a frequency sub-tile assignment. 
To show the \(\lambda^\epsilon\)-strong disjointness assumption with respect to the above assignment, we take \(P_m:=\vI_m\times\vomega_m\in\cT_m\) and \(P_n:=\vI_n\times\vomega_n\in\cT_n\) with \(\Omega_{P_m}\cap\Omega_{P_n}\neq\varnothing\).
We proceed in two steps. We first show that \(\abs{\vomega_m}<\abs{\vomega_n}\) forces the condition \(m<n\). We then demonstrate that \(m<n\) by itself implies \(\lambda^\epsilon\vI_{\cT_m}\cap \vI_n=\varnothing\).

\begin{proof}[\textbf{Proof of \(\abs{\vomega_m}<\abs{\vomega_n}\) implies \(m<n\)}] 
Recall that
\begin{equation}\label{eq_Hdim_out_sel_intersect}
    \varnothing\neq \Omega_{P_m}\cap\Omega_{P_n}\subset \vomega^{\br{+,1}}_m\cap\vomega^{\br{+,1}}_n.
\end{equation}
By the dyadic structure, we must have
\begin{equation*}
    \vomega^{\br{+,1}}_m=\vomega^{\br{+,1}}_m\cap\vomega^{\br{+,1}}_n\subsetneq\vomega^{\br{+,1}}_n.
\end{equation*}
By \textbf{Observation \ref{obs_Hdim_subfreq_inclus}}, the above relation can be strengthen to the following:
\begin{equation*}
    \vomega^{\br{+,1}}_{\cT_m}\subset\vomega^{\br{+,1}}_m\subset\Omega_{P_n}:=\bigsqcup_{v_{\cT_n}<v<2^{k_\lambda}}\fs_v\vomega_n.
\end{equation*}
As a direct consequence, we obtain the following chain of inequalities:
\begin{equation*}
    \inf \pi_{\R\ve_1}\vomega^{\br{+,1}}_{\cT_m} \geq
    \inf\pi_{\R\ve_1}\Omega_{P_n}
    =
    \inf\pi_{\R\ve_1}\fs_{v_{\cT_n}+1}\vomega_n
    >
    \inf\pi_{\R\ve_1}\fs_{v_{\cT_n}}\vomega_n
    \geq
    \inf \pi_{\R\ve_1}\vomega^{\br{+,1}}_{\cT_n}.
\end{equation*}
By the maximality assumption \eqref{eq_lem_Hdim_sel_rightmost}, we must have \(m<n\).
\end{proof}

\begin{proof}[\textbf{Proof of \(m<n\) implies \(\lambda^\epsilon\vI_{\cT_m}\cap \vI_n=\varnothing\)}] 
Assume \(m<n\). By the statement proven above, we have \(\abs{\vomega_m}\leq\abs{\vomega_n}\). Combined with \eqref{eq_Hdim_out_sel_intersect}, we deduce that \( \vomega_{\cT_m}\subset \vomega_m\subset \vomega_n\).
However, the algorithm forces that
\begin{equation}\label{eq_lem_Hdim_sel_out_PlessP}
    P_n\not\leq \br{\ve,\vnull}\boxplus P_{\cT_m},\quad \forall \ve\in\vZ_\epsilon.
\end{equation}
This can only happen when:
\begin{equation*}
    \br{\ve\boxplus\vI_{\cT_m}}\cap\vI_n=\varnothing,\quad \forall \ve\in\vZ_\epsilon,
\end{equation*}
which implies \(\lambda^\epsilon \vI_{\cT_m}\cap\vI_n=\varnothing\).
\end{proof}
This demonstrates that the sequence of trees \(\BR{\cT_i}_i\) is \(\lambda^\epsilon\)-strongly disjoint with respect to \eqref{eq_lem_Hdim_sel_out_assign}. It remains to check the two inequalities required for the application of the generalized Bessel-type inequality \textbf{Lemma \ref{lem_Hdim_bessel}}.
On the one hand, we notice that for such \(P:=\vI\times\vomega\in \cT_i\), the condition \(\nrm{f}_{E\br{\cT_i}}\leq\nrm{f}_{E\br{\cP}}=1\) gives the following relations:
\begin{equation*}
    \nrm{f^{\cT_{out,+}}}_P=\nrm{\br{\Psi_P\ast f}^{\br{+,1}}}_P\leq \abs{\vI}^{\frac{1}{2}}.
\end{equation*}
On the other hand, by the above identity and the construction of \(\cT_i\in\T^\ast_i\), we have
\begin{equation*}
    \nrm{f^{out,+}}_{L^2\br{\cT_i}}=\nrm{
        \nrm{
            \br{\Psi_P\ast f}^{\br{+,1}}
        }_P
    }_{\ell^2\br{P\in\cT_i}}
    \eqsim \abs{\vI_{\cT_i}}^{\frac{1}{2}}.
\end{equation*}
As a consequence, we can apply the generalized Bessel inequality \textbf{Lemma \ref{lem_Hdim_bessel}}, thus proving \eqref{eq_Hdim_energy_sel_out_I_sum}.

The treatment of the \(\ast=-\) case is essentially the mirror image of the \(\ast=+\) case. The two key modifications are replacing \eqref{eq_lem_Hdim_sel_rightmost} with the following requirement
\begin{equation*}
    \cT\in\T^\ast_i\implies
    \sup\pi_{\R\ve_1}\vomega^{\br{+,1}}_\cT \geq \sup\pi_{\R\ve_1}\vomega^{\br{+,1}}_{\cT_i}
\end{equation*}
in the algorithm and replacing the assignment \eqref{eq_lem_Hdim_sel_out_assign} with
\begin{equation*}
    \Omega_P:=
    \bigsqcup_{2^{k_\lambda-1}-1\leq v< v_{\cT_i}}\fs_v\vomega
    ,\quad
    \Psi_P:=
        \fp_{\fs^-_{v_{\cT_i}}\vomega}
        \br{\br{\Dil^\infty_2\phi}_{\vomega^{\br{+,1}}}}^\vee
        .
\end{equation*}
Once the above two adjustments are made, the argument for \(\ast=-\) case follows identically. 

In summary, by setting \(\T:=\T^{cen}\sqcup\T^{out,0}\sqcup\T^{out,+}\sqcup\T^{out,-}\), we complete the proof of \textbf{Lemma \ref{lem_Hdim_energy_sel}}.

\subsection{Proof of Lemma \ref{lem_Hdim_bessel}: the generalized Bessel-type inequality}

\begin{proof}[\textbf{Proof of Lemma \ref{lem_Hdim_bessel}}]
    Starting with the left-hand side of \eqref{eq_lem_Hdim_bessel_concl} we apply \eqref{eq_lem_Hdim_bessel_tree_ass} to deduce the following:
    \begin{equation*}
        \sum_i\abs{\vI_{\cT_i}}\eqsim 
        \sum_i \nrm{\nrm{\br{\Psi_P\ast f}^{\br{+,1}}}_P}_{\ell^2\br{P\in\cT_i}}^2=
        \sum_i\sum_{P\in\cT_i}\nrm{\br{\Psi_P\ast f}^{\br{+,1}}}_P^2.
    \end{equation*}
    By \eqref{eq_Hdim_wp_nrm}, we reformulate and equate the above with
    \begin{equation*}
        =
        \sum_i\sum_{\vI\times\vomega\in\cT_i}
        \nrm{\pi_{\vomega^{\br{+,1}}} \br{\Psi_{\vI\times\vomega}\ast f} \br{\vx} }_{L^2\br{\chi^\Phi_\vI\br{\vx}d\vx}}^2
        =
        \sum_i\sum_{P\in\cT_i}
        \ang{\Psi^\ast_{\vI\times\vomega}\ast\pi_{\vomega^{\br{+,1}}}\br{\Psi_{\vI\times\vomega}\ast f}^{\br{+,1}}_{\vI\times\vomega},f},
    \end{equation*}
    where, to simplify the notation, we set
    \begin{equation*}
        F\br{\vI\times\vomega}:=
        \Psi^\ast_{\vI\times\vomega}\ast
        \pi_{\vomega^{\br{+,1}}}
        \br{
            \Psi_{\vI\times\vomega}\ast f
        }^{\br{+,1}}_{\vI\times\vomega}, \quad g^\ast\br{\vz}:=\overline{g\br{-\vz}}
    \end{equation*}
    for \(\vI\times\vomega \in\cP:=\bigcup_i\cT_i\). We thus have, by Cauchy-Schwarz inequality, the following estimate:
    \begin{equation*}
        \sum_i\abs{\vI_{\cT_i}}\eqsim
        \sum_i\sum_{P\in\cT_i}\ang{F\br{P},f}
        =
        \ang{\sum_i\sum_{P\in\cT_i}F\br{P},f}
        \leq \nrm{\sum_i\sum_{P\in\cT_i} F\br{P}}_{L^2}\nrm{f}_{L^2}.
    \end{equation*}
    Thus, to show \eqref{eq_lem_Hdim_bessel_concl}, it suffices to prove the following:
    \begin{equation}\label{eq_lem_Hdim_bessel_reduced}
        \nrm{\sum_i\sum_{P\in\cT_i} F\br{P}}_{L^2}^2 \underset{\epsilon}{\lesssim } \sum_i \abs{\vI_{\cT_i}}.
    \end{equation}
    We expand the left-hand side of the above statement
    \begin{equation*}
        \nrm{\sum_i\sum_{P\in\cT_i} F\br{P}}_{L^2}^2
        =
        \sum_{m,n}
        \sum_{\substack{
            P_m\in\cT_m,\\
            P_n\in\cT_n
        }}
        \ang{F\br{P_m},F\br{P_n}}
    \end{equation*}
    and equate the above with the following initial decomposition:
    \begin{equation*}
        =
        \sum_{m,n}
        \sum_{\substack{
            P_m:=\vI_m\times\vomega_m\in\cT_m,\\
            P_n:=\vI_n\times\vomega_n\in\cT_n\\
            \abs{\vomega_m}=\abs{\vomega_n}
        }}
        \ang{F\br{P_m},F\br{P_n}}
        +2\Re
        \br{
        \sum_{m,n}
        \sum_{\substack{
            P_m:=\vI_m\times\vomega_m\in\cT_m,\\
            P_n:=\vI_n\times\vomega_n\in\cT_n\\
            \abs{\vomega_m}<\abs{\vomega_n}
        }}
        \ang{F\br{P_m},F\br{P_n}}
        }=:
        S_1+
        2\Re S_2.
    \end{equation*}
    For \(S_1\), we make the further decomposition
    \begin{equation*}
        S_1=
        \sum_i
        \sum_{\substack{
            P:=\vI\times\vomega\in\cT_i,\\
            P':=\vI'\times\vomega'\in\cT_i\\
            \abs{\vomega}=\abs{\vomega'}
        }}
        \ang{F\br{P},F\br{P'}}
        +
        2\Re\br{
        \sum_{m<n}
        \sum_{\substack{
            P_m:=\vI_m\times\vomega_m\in\cT_m,\\
            P_n:=\vI_n\times\vomega_n\in\cT_n\\
            \abs{\vomega_m}=\abs{\vomega_n}
        }}
        \ang{F\br{P_m},F\br{P_n}}
        }=:
        S_3+
        2\Re S_4.
    \end{equation*}
    To estimate any of the above terms, we must understand the contribution from the inner product. Given \(P:=\vI\times\vomega\in \cP\cap\P_k\) and \(P':=\vI'\times\vomega'\in\cP\cap\P_{k'}\) with \(k'\leq k\), we have the following identity:
    \begin{equation*}
        \ang{F\br{P},F\br{P'}}
        =
        \ang{
        \pi_{\vomega^{\prime \br{+,1}}}\br{
            \Psi_{P'}\ast\Psi_P^\ast\ast\pi_{\vomega^{\br{+,1}}}
            \br{\Psi_P\ast f}^{\br{+,1}}_P
        },
        \br{\Psi_{P'}\ast f}^{\br{+,1}}_{P'}
        }.
    \end{equation*}
    Observe that
    \begin{equation*}
        \pi_{\vomega^{\prime \br{+,1}}}\br{
            \Psi_{P'}\ast\Psi_P^\ast\ast\pi_{\vomega^{\br{+,1}}}
            g
        }
        =
        \Psi_{P',P}
        \ast g,\quad
        \Psi_{P',P}:=
        \widecheck{\phi_{\vomega^{\prime \br{+,1}}}}\ast
        \Psi_{P'}\ast\Psi_P^\ast\ast \widecheck{\phi_{\vomega^{\br{+,1}}}}.
    \end{equation*}
     On the one hand, since by \eqref{eq_def_Hdim_freq_asign_freq} \(\supp\widehat{\Psi_{P',P}}\subset \Omega_P\cap\Omega_{P'}\) we deduce
    \begin{equation*}
        \ang{F\br{P},F\br{P'}}\neq 0\implies \Omega_P\cap\Omega_{P'}\neq \varnothing.
    \end{equation*}
    On the other hand, one can easily check that \(\Psi_{P',P}\ast\chi^N_\vI\underset{N}{\lesssim}\chi^N_\vI\) and applying \textbf{Proposition \ref{prop_Hdim_wp_in_prod_psi_VP}} one has
    \begin{align}\label{eq_lem_Hdim_bessel_in_prod}
        \abs{\ang{F\br{P},F\br{P'}}}= &
        \abs{
            \ang{
                \Psi_{P',P}\ast
                \br{\Psi_P\ast f}^{\br{+,1}}_P
            ,
            \br{\Psi_{P'}\ast f}^{\br{+,1}}_{P'}
            }
        }\nonumber\\
        \underset{N}{\lesssim} &
        \nrm{\br{\Psi_P\ast f}^{\br{+,1}}}_P\nrm{\br{\Psi_{P'}\ast f}^{\br{+,1}}}_{P'}
        \min\br{1,\sqrt{\frac{\lambda^{D+1}\abs{\vI'}}{\abs{\vI}}}}
        \nrm{\chi^N_{\vI}}_{L^\infty\br{\vI'}}.
    \end{align}
    We first focus on \(S_3\). By \eqref{eq_lem_Hdim_bessel_in_prod}, we deduce:
    \begin{equation*}
        \abs{S_3}\underset{N}{\lesssim}
        \sum_i
        \sum_{\substack{
            P:=\vI\times\vomega\in\cT_i,\\
            P':=\vI'\times\vomega'\in\cT_i\\
            \abs{\vomega}=\abs{\vomega'}
        }}
        \nrm{\br{\Psi_P\ast f}^{\br{+,1}}}_P\nrm{\br{\Psi_{P'}\ast f}^{\br{+,1}}}_{P'}
        \min\br{\nrm{\chi^N_{\vI}}_{L^\infty\br{\vI'}},\nrm{\chi^N_{\vI'}}_{L^\infty\br{\vI}}}.
    \end{equation*}
    Cauchy-Schwarz inequality and assumption \eqref{eq_lem_Hdim_bessel_tree_ass} allow us to dominate the above by
    \begin{align*}
        &\sum_i
        \sum_{\substack{
            P:=\vI\times\vomega\in\cT_i
        }}
        \nrm{\br{\Psi_P\ast f}^{\br{+,1}}}_P^2
        \sum_{\substack{
            P':=\vI'\times\vomega'\in\cT_i\\
            \abs{\vomega}=\abs{\vomega'}
        }}
        \nrm{\chi^N_{\vI}}_{L^\infty\br{\vI'}}\\
        \underset{N}{\lesssim} &
        \sum_i
        \sum_{\substack{
            P:=\vI\times\vomega\in\cT_i
        }}
        \nrm{\br{\Psi_P\ast f}^{\br{+,1}}}_P^2
        =
        \sum_i \nrm{\nrm{\br{\Psi_P\ast f}^{\br{+,1}}}_P}_{\ell^2\br{P\in\cT_i}}^2
        \eqsim \sum_i\abs{\vI_{\cT_i}}\,,
    \end{align*}
thus ending the treatment of \(S_3\).

   We now shift our focus to the remaining terms \(S_2,S_4\). Suggested by the \(\lambda^\epsilon\)-strong disjointness assumption, we shall treat \(S_2\) and \(S_4\) in parallel by setting the notation 
    \begin{equation*}
       \sigma_2\br{m,n,\vI\times\vomega,\vI'\times\vomega'}:=\br{\br{\vI\times\vomega,\vI'\times\vomega'}\in \cT_m\times\cT_n} \land \br{\abs{\vomega}<\abs{\vomega'}}
    \end{equation*}
    and 
    \begin{equation*}
      \sigma_4\br{m,n,\vI\times\vomega,\vI'\times\vomega'}:= \br{\br{\vI\times\vomega,\vI'\times\vomega'}\in \cT_m\times\cT_n} \land\br{m<n} \land \br{\abs{\vomega}=\abs{\vomega'}}.
    \end{equation*}
    Thus, we slightly notational abuse, we write from now on
    \begin{equation*}
        S_\ast=
        \sum_{\substack{m,n,P,P':\\
        \sigma_\ast\br{m,n,P,P'}}}
        \ang{F\br{P},F\br{P'}},\quad\ast \in\BR{2,4}.
    \end{equation*}
    With these, using \eqref{eq_lem_Hdim_bessel_in_prod} and letting\footnote{Throughout this section we assume wlog that  \(\sigma'_\ast\br{m,n,P,P'}\not=\emptyset\).} \(\sigma'_\ast\br{m,n,P,P'}:=\sigma_\ast\br{m,n,P,P'}\land \br{\Omega_P\cap\Omega_{P'}\neq \varnothing }\not\), we have
    \begin{equation*}
        \abs{S_\ast}
        \underset{N}{\lesssim}
        \lambda^{\frac{D+1}{2}}
        \sum_{\substack{m,n,P,P':\\
        \sigma'_\ast\br{m,n,P,P'}}}
        \nrm{\br{\Psi_P\ast f}^{\br{+,1}}}_P\nrm{\br{\Psi_{P'}\ast f}^{\br{+,1}}}_{P'}
        \frac{\abs{\vI'}^{1/2}}{\abs{\vI}^{1/2}}
        \nrm{\chi^N_{\vI}}_{L^\infty\br{\vI'}},\quad \ast\in\BR{2,4},
    \end{equation*}
  
    Applying now assumption \eqref{eq_lem_Hdim_bessel_tile_ass}, we dominate the above by
    \begin{equation*}
     \lambda^{\frac{D+1}{2}}
        \sum_{\substack{m,n,P,P':\\
        \sigma'_\ast\br{m,n,P,P'}}}
        \abs{\vI'}\cdot
        \nrm{\chi^N_{\vI}}_{L^\infty\br{\vI'}}.
    \end{equation*}
    We claim that is enough to prove that for all \(m\) and \(\ast\in\BR{2,4}\) the following estimate holds:
    \begin{equation}\label{eq_lem_Hdim_bessel_reduced2}
        \sum_{\substack{
            n,P,P':\\
            \sigma'_\ast\br{m,n,P,P'}
        }}
        \abs{\vI'}\cdot
        \nrm{\chi^N_{\vI}}_{L^\infty\br{\vI'}}
        \underset{N}{\lesssim}\lambda^{\epsilon\br{1-N}\br{D+1}}\abs{\vI_{\cT_m}}.
    \end{equation}
    Indeed, assuming \eqref{eq_lem_Hdim_bessel_reduced2}, we simply take \(N\geq 1+\frac{1}{2\epsilon}\) in order to deduce that \(\abs{S_2}+\abs{S_4}\) in bounded by the right-hand side of \eqref{eq_lem_Hdim_bessel_reduced} thus completing the proof of \textbf{Lemma \ref{lem_Hdim_bessel}}.
     
     We thus focus on \eqref{eq_lem_Hdim_bessel_reduced2} and first notice that
    \begin{equation*}
        \sigma'_\ast\br{m,n,P,P'}\implies \abs{\vI}\geq \abs{\vI'}\implies
        \chi_\vI\br{\vz}\eqsim \nrm{\chi_\vI}_{L^\infty\br{\vI'}},\,\forall \vz\in\vI'.
    \end{equation*}
    This allows us to approximate the \(L^\infty\) quantities with \(L^1\) quantities:
    \begin{equation}\label{eq_lem_Hdim_bessel_reduced3}
        \sum_{\substack{
            n,P,P':\\
            \sigma'_\ast\br{m,n,P,P'}
        }}
        \abs{\vI'}\cdot
        \nrm{\chi^N_{\vI}}_{L^\infty\br{\vI'}}
        \underset{N}{\eqsim}
        \sum_{\substack{
            n,P,P':\\
            \sigma'_\ast\br{m,n,P,P'}
        }}
        \nrm{\chi^N_{\vI}}_{L^1\br{\vI'}}.
    \end{equation}
    To further reformulate the above sum, we introduce for \(\br{m,n,P}\) the following collection of intervals:
    \begin{equation*}
        \I^{\br{\ast}}_{m,n,P}:=
        \BR{\vI'\in\P\::\:
            \exists\, \vomega'\in\W,\, \sigma'_\ast\br{m,n,P,\vI'\times\vomega'}
        },\quad\ast\in \BR{2,4}.
    \end{equation*}
    With the above notion, we may rewrite \eqref{eq_lem_Hdim_bessel_reduced3}:
    \begin{equation}\label{eq_lem_Hdim_bessel_reduced4}
        \sum_{\substack{
            n,P,P':\\
            \sigma'_\ast\br{m,n,P,P'}
        }}
        \abs{\vI'}\cdot
        \nrm{\chi^N_{\vI}}_{L^\infty\br{\vI'}}
        \underset{N}{\eqsim}
        \sum_P
        \sum_n
        \sum_{\vI'\in\I^{\br{\ast}}_{m,n,P}}
        \nrm{\chi^N_{\vI}}_{L^1\br{\vI'}}.
    \end{equation}
    We claim the following two statements about \(\I^{\br{\ast}}_{m,n,P}\):
    \begin{itemize}
        \item Given \(\vI_i\in\I^{\br{\ast}}_{m,n_i,P}\), \(\vI_1\cap\vI_2\neq \varnothing\) implies \(\vI_1=\vI_2\) and \(n_1=n_2\).
        \item \(\lambda^\epsilon\vI_{\cT_m}\cap\bigsqcup \I^{\br{\ast}}_{m,n,P}=\varnothing\). \footnote{Throughout this paper, in order to simplify the exposition, we allow some notational abuse when the significance of the write-up is obvious.}
    \end{itemize}
    To demonstrate the first statement, let \(P=\vI\times\omega\) and \(\vI_i\in\I^{\br{\ast}}_{m,n_i,P}\) with \(\vI_1\cap\vI_2\neq  \varnothing\). By construction, there are \(\vomega_1,\vomega_2\in \W\) such that
    \begin{equation*}
        \sigma'_\ast\br{m,n_i,\vI\times\vomega,\vI_i\times\vomega_i},\quad i=1,2.
    \end{equation*}
    We now have two cases. Starting with the \(\ast=2\) case, we expand the definition of \(\sigma_2'\br{\cdot,\cdot,\cdot,\cdot}\) and utilize \textbf{Observation \ref{obs_Hdim_subfreq_inclus}} to derive the following implications:
    \begin{equation*}
        \sigma'_2\br{m,n_i,\vI\times\vomega,\vI_i\times\vomega_i}\implies
        \br{
        \abs{\vomega}<\abs{\vomega_i}
        \land
            \Omega_{\vI\times\vomega}\cap\Omega_{\vI_i\times\vomega_i}
            \neq \varnothing
        }
        \implies \Omega_{\vI\times\vomega}\subset\Omega_{\vI_i\times\vomega_i}.
    \end{equation*}
    As a direct consequence, \(\Omega_{\vI\times\vomega}\subset\Omega_{\vI_1\times\vomega_1}\cap \Omega_{\vI_2\times\vomega_2} \neq \varnothing\). Since \(\vI_i\times\vomega_i\in\cT_{n_i}\) and \(\vI_1\cap\vI_2\neq\varnothing\), the \(\lambda^\epsilon\)-strong disjointness between \(\cT_{n_1}\) and \(\cT_{n_2}\) forces that \(\abs{\vomega_1}=\abs{\vomega_2}\) and \(n_1=n_2\). As a result, we must have \(\abs{\vI_1}=\abs{\vI_2}\) and thus, \(\vI_1=\vI_2\). For the \(\ast=4\) case, we recall that
    \begin{equation*}
        \sigma'_4\br{m,n_i,\vI\times\vomega,\vI_i\times\vomega_i}\implies
        \br{
        \abs{\vomega}=\abs{\vomega_i}
        \land
            \Omega_{\vI\times\vomega}\cap\Omega_{\vI_i\times\vomega_i}
            \neq \varnothing
        }.
    \end{equation*}
    As a direct consequence, we have \(
    \abs{\vI}=\abs{\vI_1}=\abs{\vI_2}\) and thus, \(\vI_1=\vI_2\). To show \(n_1=n_2\), we note that the above statement also implies
    \begin{equation*}
        \implies \varnothing\neq  \Omega_{\vI\times\vomega}\cap\Omega_{\vI_i\times\vomega_i} \subset \vomega^{\br{+,1}}\cap \vomega^{\br{+,1}}_i.
    \end{equation*}
    The fact that \(\abs{\vomega}=\abs{\vomega_i}\) forces that \(\vomega=\vomega_1=\vomega_2\). This implies that \(\vI_1\times\vomega_1=\vI_2\times\vomega_2\in\cT_{n_1}\cap\cT_{n_2}\). Since \(\Omega_{\vI_1\times\vomega_1}= \Omega_{\vI_2\times\vomega_2}\) and in particular \(\Omega_{\vI_1\times\vomega_1}\cap \Omega_{\vI_2\times\vomega_2}\neq\varnothing\), the \(\lambda^\epsilon\)-strong disjointness between \(\cT_{n_1}\) and \(\cT_{n_2}\) forces that \(n_1=n_2\). This completes the proof of the first claim. As for the second statement, it follows directly from the disjointness proven in the first claim and the \(\lambda^\epsilon\)-strong disjointness assumption.

    As a direct consequence of the two claims, we have:
    \begin{equation*}
        \bigsqcup_n \bigsqcup\I^{\br{\ast}}_{m,n,P}\subset \br{\lambda^\epsilon\vI_{\cT_m}}^c.
    \end{equation*}
    This allows us to further dominate \eqref{eq_lem_Hdim_bessel_reduced4} by
    \begin{equation}\label{eq_lem_Hdim_bessel_reduced5}
        \sum_P
        \nrm{\chi^N_{\vI}}_{L^1\br{\bigsqcup_n \bigsqcup\I^{\br{\ast}}_{m,n,P}}}
        \leq
        \sum_{\vI:\vI\subset \vI_{\cT_m}}
        \nrm{\chi^N_{\vI}}_{L^1\br{\br{\lambda^\epsilon\vI_{\cT_m}}^c}}.
    \end{equation}
    Lastly, we observe that for all \(\vI\subset \vI_{\cT_m}\) with \(\vI\in\I_k\), the following approximation holds:
    \begin{equation*}
        \ang{
            \vS^{-k}\br{\vx-\vx'}
        }_\otimes^{-1}
        \eqsim
        \chi_\vI\br{\vx'},\quad
        \forall \vx\in\vI,\,\forall \vx'\notin \lambda^\epsilon\vI_{\cT_m}.
    \end{equation*}
    We can thus approximate the RHS of \eqref{eq_lem_Hdim_bessel_reduced5} by
    \begin{equation*}
        \underset{N}{\eqsim}
        \sum_k
        \sum_{\substack{
            \vI\in\I_k\\
            \vI\subset \vI_{\cT_m}
        }}
        \frac{1}{\abs{\vI}}
        \int_\vI
            \int_{\br{\lambda^\epsilon\vI_{\cT_m}}^c}
                \ang{
                    \vS^{-k}\br{\vx-\vx'}
                }_\otimes^{-N}
            d\vx'
        d\vx
    \end{equation*}
   Assuming \(\vI_{\cT_m}\in\I_{k_0}\), a direct calculation allows us to further rewrite the above expression as
    \begin{equation*}
        =
        \sum_{k\leq k_0}
        \det\br{\vS}^{-k}
        \int_{\vI_{\cT_m}}
            \int_{\br{\lambda^\epsilon\vI_{\cT_m}}^c}
                \ang{
                    \vS^{-k}\br{\vx-\vx'}
                }_\otimes^{-N}
            d\vx'
        d\vx
    \end{equation*}
which, after a change of variable, becomes
    \begin{align}\label{eq_lem_Hdim_bessel_reduced6}
        =
        \abs{\vI_{\cT_m}}
        \sum_{k=0}^\infty
        \bigg(& 2^k
        \int_{\frac{1}{2}}^{\frac{1}{2}}
            \int_{\br{-\lambda^\epsilon/2,\lambda^\epsilon/2}^c}
                \ang{
                    2^k\br{x-x'}
                }_\otimes^{-N}
            dx'
        dx
        \bigg)^D\nonumber\\
        \cdot &
        4^k
        \int_{\frac{1}{2}}^{\frac{1}{2}}
            \int_{\br{-\lambda^\epsilon/2,\lambda^\epsilon/2}^c}
                \ang{
                    4^k\br{y-y'}
                }_\otimes^{-N}
            dy'
        dy.
    \end{align}
   However, we note that for \(l\geq 1\) the following holds:
    \begin{equation*}
        l
        \int_{\frac{1}{2}}^{\frac{1}{2}}
            \int_{\br{-\lambda^\epsilon/2,\lambda^\epsilon/2}^c}
                \ang{
                    l\br{s-s'}
                }_\otimes^{-N}
            ds'
        ds
        \underset{N}{\eqsim}
        l^{1-N}\cdot \lambda^{\epsilon\br{1-N}}.
    \end{equation*}
    This allows us to further approximate \eqref{eq_lem_Hdim_bessel_reduced6} and conclude the proof of \eqref{eq_lem_Hdim_bessel_reduced2}:
    \begin{equation*}
        \sum_{\substack{
            n,P,P':\\
            \sigma'_\ast\br{m,n,P,P'}
        }}
        \abs{\vI'}\cdot
        \nrm{\chi^N_{\vI}}_{L^\infty\br{\vI'}}
        \underset{N}{\lesssim}
        \eqref{eq_lem_Hdim_bessel_reduced6}
        \underset{N}{\eqsim}
        \abs{\vI_{\cT_m}}
        \sum_{k=0}^\infty
        2^{3k\br{1-N}}\cdot \lambda^{\epsilon\br{1-N}\br{D+1}}
        \leq 
        \abs{\vI_{\cT_m}}
        \lambda^{\epsilon\br{1-N}\br{D+1}}.
    \end{equation*}
\end{proof}

\section{Single tile estimate: Proof of Proposition \ref{thm_Hdim_sing_tile_ests} subject to Lemma \ref{lem_sub_lev_2_2}}\label{singltilesestimgen}

In this section, we perform a high-resolution analysis via the LGC-method in order to reduce the proof of  \textbf{Proposition \ref{thm_Hdim_sing_tile_ests}} to the level set analysis encapsulated within \textbf{ Lemma \ref{lem_sub_lev_2_2}}.


\subsection{Initial reductions: The natural \texorpdfstring{\(L^p\)}{} improving bound}
    Due to the symmetry structures \(\Lambda_{\vI\times\vomega}^{\br{j}}\) satisfies, it suffices to take \(k=0\) and
    \begin{equation}\label{eq_thm_Hdim_sing_tile_tile_normal}
        \vI\times\vomega=
        \bR{0,1}^{D+1}\times
        \br{
            \bR{-\lambda/2,\lambda/2}\times\bR{0,\lambda}^D
        }\in\P_0\br{\lambda}.
    \end{equation}
    We claim that the \(\vu\) physical shift can be ignored during the analysis of a single tile: indeed, writing the formula for \(\Lambda^{\br{j}}_{\vI\times\vomega}\br{f,g}\) explicitly
    \begin{equation*}
        \Lambda_{\vI\times\vomega}^{\br{j}}\br{f,g}
        =
            \int
                \int
                    f^{\br{+,j}}_{\vI\times\vomega}\br{\vx-\vtau}
                    e\br{\va\br{\vx}^\top\vtau}
                d\mu_0\br{\vtau-\vu}
                \cdot 
                \overline{\1_{A^{\br{-,j}}_\vomega}\br{\vx}g\br{\vx}}
            d\vx\,,
    \end{equation*}
    after a few changes of variables, we further have
    \begin{equation*}
        =
            \int
                \int
                    f^{\br{+,j}}_{\vI\times\vomega}\br{\vx-\vtau}
                    e\br{\va\br{\vx+\vu}^\top\vtau}
                d\mu_0\br{\vtau}
                \overline{\1_{\vomega^{\br{-,j}}}\br{\va\br{\vx+\vu}}e\br{-\va\br{\vx+\vu}^\top\vu}g\br{\vx+\vu}}
            d\vx
        .
    \end{equation*}
    Recalling the following definition:
    \begin{equation*}
        \cC\cR^{\br{\va}}_{\mu_0} f\br{\vx}:=
        \int 
            f\br{\vx-\vtau}
            e\br{\va\br{\vx}^\top\vtau}
        d\mu_0\br{\vtau},
    \end{equation*}
    we realize that \eqref{eq_thm_Hdim_sing_tile_est_c} in \textbf{Proposition \ref{thm_Hdim_sing_tile_ests}} is equivalent\footnote{A similar argument works when dealing with \eqref{eq_thm_Hdim_sing_tile_est_half}.} with
    \begin{equation}\label{eq_thm_Hdim_sing_tile_est_wo_shift}
        \abs{\ang{
            \cC\cR^{\br{\va}}_{\mu_0}
                f^{\br{+,j}}_{\vI\times\vomega}
            ,
            \1_{A^{\br{-,j}}_\vomega} g
        }}\lesssim 
        \lambda^{-\sigma}
        \mu^c_\vI\br{A^{\br{-,j}}_\vomega\cap\supp g}
        \nrm{f^{\br{+,j}}}_{\vI\times\vomega}
        \nrm{\1_{A^{\br{-,j}}_\vomega}g}_{L^2\br{d\mu_\vI}}
    \end{equation}
    when we replace \(\va\) with
     \(\va_\vu:=\Tr_{-\vu}\va\) and \(g\) with \(g_\vu:=
    e\br{-\va^\top \vu}\Tr_{-\vu}g\).
     
     In other words, for the purpose of proving \textbf{Proposition \ref{thm_Hdim_sing_tile_ests}}, we may assume \(\vu=\vnull\) and thus \(\cC\cR^{\br{\va}}_{\mu_0}=\cC_0\).
    
    With these being said, we follow the treatment presented in \textsc{Section \ref{subsec_planar_sing_tile_est_bd_Lambda_c}}. Recalling \textbf{Definition \ref{def_sing_tile_est}}, we notice that in order to prove \eqref{eq_thm_Hdim_sing_tile_est_wo_shift}, and thus \eqref{eq_thm_Hdim_sing_tile_est_c}, it suffices to show the two estimates 
    \begin{equation}\label{estimnm}
        \nrm{\Lambda}_{\frac{D}{2D+2}}\lesssim 1,\quad
        \nrm{\Lambda}_0\lesssim \lambda^{-\sigma}
    \end{equation}
  and then take their corresponding geometric mean---see \eqref{eq_thm_Hdim_sing_tile_Lp_improv} and \eqref{eq_thm_Hdim_sing_tile_main_decay} below.    
    
    Now for the first inequality in \eqref{estimnm} we notice the following: recalling \eqref{eq_planar_sing_tile_triv_for_Lp_improv}, once the \(L^4\br{dx_0dx_1}\)--\(L^{4/3}\br{dx_0dx_1}\) H\"{o}lder inequality is replaced by the \(L^{2D+2}\br{d\vx}\)--\(L^{\frac{2D+2}{2D+1}}\br{d\vx}\) analog, the standard \(L^p\)-improving estimate part of \textbf{Theorem \ref{thm_Lp_improv}} yields the desired estimate
    \begin{equation}\label{eq_thm_Hdim_sing_tile_Lp_improv}
        \abs{\ang{
            \cC\cR^{\br{\va}}_{\mu_0}
                f^{\br{+,j}}_{\vI\times\vomega}
            ,
            \1_{A^{\br{-,j}}_\vomega} g
        }}
        \lesssim 
        \mu^{\frac{D}{2D+2}}_\vI\br{A^{\br{-,j}}_\vomega\cap\supp g}
        \nrm{f^{\br{+,j}}}_{\vI\times\vomega}
        \nrm{\1_{A^{\br{-,j}}_\vomega}g}_{L^2\br{d\mu_\vI}}.
    \end{equation}
    Thus, in order to complete the proof of \eqref{estimnm}, it remains to obtain the decay estimate
    \begin{equation}\label{eq_thm_Hdim_sing_tile_main_decay}
        \abs{\ang{
            \cC\cR^{\br{\va}}_{\mu_0}
                f^{\br{+,j}}_{\vI\times\vomega}
            ,
            \1_{A^{\br{-,j}}_\vomega} g
        }}\lesssim \lambda^{-\sigma}
        \nrm{f^{\br{+,j}}}_{\vI\times\vomega}
        \nrm{\1_{A^{\br{-,j}}_\vomega}g}_{L^2\br{d\mu_\vI}}
    \end{equation}
  
    \subsection{The LGC-method: proof of \texorpdfstring{\eqref{eq_thm_Hdim_sing_tile_est_c}}{} via the decay estimate \texorpdfstring{\eqref{eq_thm_Hdim_sing_tile_main_decay}}{}}
    To proceed, we continue the reduction procedure presented in \textsc{Sections \ref{subsec_planar_sing_tile_est_bd_Lambda_c} and \ref{subsec_tf_rep_lin}}.
    Recall \eqref{eq_Hdim_wp_CL_decomp}:
    \begin{equation*}
        f^{\br{+,j}}_{\vI\times\vomega}:=\chi^\Phi_\vI\pi^{\br{+,j}}_\vomega f=
        \int
            \ang{\pi^{\br{+,j}}_\vomega f,\varphi^{\br{L}}_{0,\vx,\vxi}}
            \varphi^{\br{L}}_{0,\vx,\vxi}
            \chi^{\br{L}}_{\vI}\br{\vx}
        d\vx d\vxi.
    \end{equation*}
    Again, we make the substitution \(\vxi=\sqrt{\lambda}\vzeta\), set
    \begin{equation}\label{wvpa}
        \varphi_{\vx,\vzeta}:=
        \Tr_\vx \Dil^1_{1/\sqrt{\lambda}} \Mod_{\vzeta}\big(\widecheck{\phi}\big)^{\otimes D+1}=
        \nrm{\phi}_{L^2}^{D+1}
        \lambda^{\frac{D+1}{4}}\varphi^{\br{L}}_{0,\vx,\sqrt{\lambda}\vzeta}
        ,
    \end{equation}
    and thus, reformulate
    \begin{equation*}
        f^{\br{+,j}}_{\vI\times\vomega}
        =
        \nrm{\phi}_{L^2}^{-2D-2}
        \int 
        \ang{
            \pi^{\br{+,j}}_\vomega f, 
            \varphi_{\vx,\vzeta}
        }
            \varphi_{\vx,\vzeta}
            \chi^{\br{L}}_\vI\br{\vx}
        d\vx d\vzeta
        .
    \end{equation*}
     In analogy with \eqref{eq_planar_def_dmu_I_norm}, we simply take \(d\mu\br{\vx}:=\ang{\vx}^{-N}_{\otimes}d\vx\) while in analogy with \eqref{eq_planar_def_Fxyuv} we observe
    \begin{equation*}
        \nrm{f^{\br{+,j}}}_{\vI\times\vomega}
        \eqsim
        \nrm{\ang{
            \pi^{\br{+,j}}_\vomega f, 
            \varphi_{\vx,\vzeta}
        }}_{L^2\br{\chi^{\br{L}}_\vI\br{\vx}d\vx d\vzeta}}
        =
        \nrm{
            \ang{
                \pi^{\br{+,j}}_\vomega f, 
                \varphi_{\vx,\vzeta}
            }
            \sqrt{
                \chi^{\br{L}}_\vI\br{\vx}
            }
        }_{L^2\br{d\vx d\vzeta}}.
    \end{equation*}
    This suggests that we take
    \begin{equation*}
        F\br{\vx,\vzeta}:=
        \ang{
            \pi^{\br{+,j}}_\vomega f, 
            \varphi_{\vx,\vzeta}
        }
        \sqrt{
            \chi^{\br{L}}_\vI\br{\vx}
        }.
    \end{equation*}
    Finally, notice that \(0\leq \chi^{\br{L}}_\vI\br{\vx} 
        \underset{N}{\lesssim} \ang{\vx}_\otimes^{-N}\).
    Proving \eqref{eq_thm_Hdim_sing_tile_main_decay} can thus be reduced to showing:
    \begin{equation}\label{eq_thm_Hdim_sing_tile_main_wp_decay}
        \nrm{
            F\br{\vx,\vzeta}
            \ang{\cC\cR^{\br{\va}}_{\mu_0}\varphi_{\vx,\vzeta},g}
        }_{L^1\br{
        \substack{
            d\mu\br{\vx}d\vzeta,\\
            \abs{\vzeta}\lesssim \sqrt{\lambda}\ll \zeta_j
        }
        }}
        \underset{N}{\lesssim} 
        \lambda^{-\sigma}
        \nrm{F}_{L^2}\nrm{g}_{L^2\br{d\mu}}
    \end{equation}
    under the assumption that\(\dia\br{\va\br{\supp g}}\lesssim \lambda\). Furthermore, by utilizing the following relation:
    \begin{equation*}
     \sum_{\valpha\in\Z^{D+1}}\ang{\valpha}_\otimes^{-N} \chi^N_{\valpha\boxplus \bR{0,1}^{D+1}}\br{\vx}
        \underset{N}{\lesssim} \ang{\vx}_\otimes^{-N}
        \underset{N}{\lesssim} \sum_{\valpha\in\Z^{D+1}} \ang{\valpha}_\otimes^{-N}
        \1_{\valpha\boxplus \bR{0,1}^{D+1}}\br{\vx}
    \end{equation*}
    and the translation symmetry of the formulation, it suffices to show \eqref{eq_thm_Hdim_sing_tile_main_wp_decay} with \(d\mu\br{\vx}\) on the left-hand side replaced by \(\1_{\bR{0,1}^{D+1}}\br{\vx}d\vx\). To prove the modified \eqref{eq_thm_Hdim_sing_tile_main_wp_decay}, we adopt the reduction presented in \textsc{Section \ref{subsec_tf_rep_lin}}. Omitting the standard details, the statement transforms into a high-dimensional version of \eqref{eq_the_tile_mod}:
    \begin{equation}\label{eq_Hdim_the_tile_main_mod}
    \Lambda\br{F,G}:=
        \nrm{
            \frac{
                F\br{\vx-\vX\br{\vt},\vzeta}
                G\br{\vx,\vt,\vzeta,\vw}
            }{
                \ang{
                    \br{\vzeta-\vw}^\top
                    \nabla\vX\br{\vt}
                }_\otimes^A
            }
        }_{L^1\br{
            \substack{
            d\vx d\vt d\vzeta d\cH^D_\V\br{\vw},\\
                \abs{\vx}\lesssim \abs{\vt}\eqsim 1,\\
                \abs{\vzeta},\abs{\vw}\lesssim \sqrt{\lambda}\lesssim \zeta_j-w_j
            }
        }}
    \lesssim 
        \lambda^{-\sigma}
        \nrm{F}_{L^2}\nrm{g}_{L^2\br{d\mu}},
    \end{equation}
    where 
    we recall \(\V:=\BR{0}\times\R^D\), write
    \(\vt:=\br{t_0,\dots,t_{D-1}}\in \R^D\), \(\vzeta:=\br{\zeta_0,\dots,\zeta_D}\in\R^{D+1}\), \(\vw:=\br{w_0,\dots,w_D}=\br{0,w_1,\dots,w_D}\in\V\), and 
    \begin{equation*}
        G\br{\vx,\vt,\vzeta,\vw}=
        \abs{
        \ang{
            e\br{\va^\top\vX\br{\vt}}
            g_\vw,
            \psi_{\vx,\vzeta}
        }},\quad
        g_\vw:=
        \phi^{\otimes D+1}\br{\frac{\va}{\sqrt{\lambda}}-\vw}
        g
    \end{equation*}
    with the wave-packet defined explicitly as
    \begin{equation*}
        \psi_{\vx,\vzeta}:=\Tr_\vx \Dil^1_{1/\sqrt{\lambda}} \Mod_\vzeta\psi
    \end{equation*}
    for some \(\psi\in\mathcal{S}\br{\R^{D+1}}\) satisfying the following pairs of assumptions:
    \begin{equation*}
        \abs{\psi\br{\vx}}\underset{N}{\lesssim} \ang{\vx}_\otimes^{-N},\quad
        0\leq \widehat{\psi}\br{\vxi}\lesssim \1_{\bR{-1/2,1/2}^{D+1}}\br{\vxi}.
    \end{equation*}
    To initiate the strategy presented in \textsc{Section \ref{subsec_3_endpoints}}, we first lay out some preliminary settings. We make the following harmless support assumptions:
    \begin{equation*}
        F\br{\vx,\vzeta}\neq 0\implies \abs{\vzeta}\lesssim \sqrt{\lambda};\quad
        g_\vw\br{\vx}\neq 0\implies \abs{\vw}\lesssim \sqrt{\lambda}
    \end{equation*}
    and introduce the auxiliary measure:
    \begin{equation*}
        \mu_\vx\br{S}:=
        \int
            \1_S\br{\vx+\frac{\vz}{\sqrt{\lambda}}}
        d\mu\br{\vz}
    \end{equation*}
    analogous to \eqref{eq_planar_aux_meas}. Additionally, we have the high-dimensional analog of \eqref{cor_G_bds}:
    \begin{equation}\label{eq_Hdim_Gbds}
        \nrm{G\br{\vx,\vt,\vzeta,\vw}}_{L^{p'}\br{d\vzeta}}
        \lesssim\nrm{g_\vw}_{L^p\br{d\mu_\vx}},\quad p\in\mr{1,2}.
    \end{equation}
    To complete the proof of \eqref{eq_thm_Hdim_sing_tile_main_decay}, it remains to derive the high-dimensional analogs of the three estimates \eqref{eq_S_U}, \eqref{eq_S_S}, and \eqref{eq_U_S}:
\begin{itemize}    
    
\item the \textsf{sparse-uniform} regime given by    
    
    \begin{equation}\label{eq_Hdim_S_U}
        \Lambda\br{F,G}\lesssim 
        \log\br{e+\lambda}
        \nrm{
            \nrm{
                F\br{\vx,\vzeta}
            }_{L^1\br{d\vzeta}}
        }_{L^2\br{d\vx}}
        \nrm{
            \nrm{\nrm{g_\vw}_{L^1\br{d\mu_\vx}}}_{L^\infty\br{d\cH^D_\V\br{\vw}}}
        }_{L^2\br{\abs{\vx}\lesssim 1}}
     \end{equation} 
       
\item the \textsf{sparse-sparse} regime given by   

   \begin{equation}\label{eq_Hdim_S_S}      
   \Lambda\br{F,G}\lesssim \lambda^{-\sigma}
            \nrm{
                \nrm{
                    F\br{\vx,\vzeta}
                }_{L^1\br{d\vzeta}}
            }_{L^2\br{d\vx}}
            \nrm{
                \nrm{\nrm{g_\vw}_{L^1\br{d\mu_\vx}}}_{L^1\br{d\cH^D_\V\br{\vw}}}
            }_{L^2\br{\abs{\vx}\lesssim 1}}
  \end{equation}          
    
 \item the \textsf{uniform-uniform} regime given by    
    \begin{equation}\label{eq_Hdim_U_S}    
        \Lambda\br{F,G} \lesssim  \log^{\frac{1}{2}}\br{e+\lambda}
        \nrm{
            \nrm{
                F\br{\vx,\vzeta}
            }_{L^4\br{d\vzeta}}
        }_{L^2\br{d\vx}}
        \nrm{
            \nrm{
                \mu_\vx^{\frac{1}{4}}\br{\supp g_\vw}
                \nrm{g_\vw}_{L^2\br{d\mu_\vx}}
            }_{L^1\br{d\cH^D_\V\br{\vw}}}
        }_{L^2\br{\abs{\vx}\lesssim 1}}.
    \end{equation}
\end{itemize}    
    
    The corresponding interpolation argument will be identical to the one presented in \textsc{Section \ref{subsec_interpol}}.

\subsection{Treatment of the sparse-uniform regime: proof of \texorpdfstring{\eqref{eq_Hdim_S_U}}{}}\label{subsec_Hdim_S_U}
    We start with the \(L^1\)--\(L^\infty\) H\"{o}lder inequality on \(d\vzeta d\cH^D_\V\br{\vw}\).
    \begin{equation*}
        \Lambda\br{F,G}\lesssim
        \nrm{
            \nrm{
                \frac{
                    F\br{\vx-\vX\br{\vt},\vzeta}
                }{
                    \ang{
                        \br{\vzeta-\vw}^\top
                        \nabla\vX\br{\vt}
                    }_\otimes^A
                }
            }_{L^1\br{d\vzeta d\cH^D_\V\br{\vw}}}\hspace{-5em}\cdot
            \nrm{
                G\br{\vx,\vt,\vzeta,\vw}
            }_{L^\infty\br{d\vzeta d\cH^D_\V\br{\vw}}}
        }_{L^1\br{\abs{\vx}\lesssim\abs{\vt}\eqsim 1}}.
    \end{equation*}
    By \eqref{eq_Hdim_Gbds}, we dominate the above by
    \begin{align*}
       \Bigg\Vert
            \bigg\Vert
            F\br{\vx-\vX\br{\vt},\vzeta}
            \cdot &
                \bigg\Vert
                    \ang{
                        \br{\vzeta-\vw}^\top
                        \nabla\vX\br{\vt}
                    }_\otimes^{-A}
                \bigg\Vert_{L^1\br{d\cH^D_\V\br{\vw},\abs{\vw}\lesssim \sqrt{\lambda}}}
            \bigg\Vert_{L^1\br{d\vzeta}}\\
        \cdot &
            \nrm{
                \nrm{g_\vw}_{L^1\br{d\mu_\vx}}
            }_{L^\infty\br{d\cH^D_\V\br{\vw}}}
        \Bigg\Vert_{L^1\br{\abs{\vx} \lesssim \abs{\vt} \eqsim 1}}.
    \end{align*}
    In contrast to the planar setting, the lack of a lower bound on the size of \(t_0\) gives rise to the weaker estimate:
    \begin{equation*}
        \nrm{
            \ang{
                \br{\vzeta-\vw}^\top
                \nabla\vX\br{\vt}
            }_\otimes^{-A}
        }_{L^1\br{d\cH^D_\V\br{\vw},\abs{\vw}\lesssim \sqrt{\lambda}}}
        \lesssim
        \nrm{
            \ang{
                \zeta_0+2t_0\br{\zeta_D-w_D}
            }^{-A}
        }_{L^1\br{\abs{w_D}\lesssim\sqrt{\lambda}}}
        \lesssim \sqrt{\lambda}\wedge \frac{1}{\abs{t_0}}.
    \end{equation*}
    We thus have
    \begin{equation*}
        \Lambda\br{F,G}
        \lesssim
        \nrm{
            \nrm{F\br{\vx-\vX\br{\vt},\vzeta}}_{L^1\br{d\vzeta}}
            \br{\sqrt{\lambda}\wedge \frac{1}{\abs{t_0}}}
            \nrm{
                \nrm{g_\vw}_{L^1\br{d\mu_\vx}}
            }_{L^\infty\br{d\cH^D_\V\br{\vw}}}
        }_{L^1\br{\abs{\vx}\lesssim\abs{\vt}\eqsim 1}}.
    \end{equation*}
    By Cauchy-Schwarz inequality on \(d\vx\), we dominate the above by
    \begin{equation*}
      \nrm{
            \br{\sqrt{\lambda}\wedge \frac{1}{t}}
        }_{L^1\br{0<t\lesssim 1}}
        \nrm{
            \nrm{F\br{\vx,\vzeta}}_{L^1\br{d\vzeta}}
        }_{L^2\br{d\vx}}
        \nrm{
            \nrm{
                \nrm{g_\vw}_{L^1\br{d\mu_\vx}}
            }_{L^\infty\br{d\cH^D_\V\br{\vw}}}
        }_{L^2\br{\abs{\vx}\lesssim 1}}.
    \end{equation*}
    Lastly, direct calculation on the first factor yields \eqref{eq_Hdim_S_U}.

    \subsection{Treatment of the sparse-sparse regime:  proof of \texorpdfstring{\eqref{eq_Hdim_S_S}}{}}\label{subsec_Hdim_S_S}
    One may follow the steps presented in \textsc{Section \ref{subsec_S_S}}, replace the \(L^4\)--\(L^{\frac{4}{3}}\) H\"{o}lder inequality by the \(L^{\frac{2D+2}{D}}\)--\(L^{\frac{2D+2}{D+2}}\) analog in \eqref{eq_SS_cor} and \eqref{eq_SS_func_imp}, and replace the \(L^2\)--\(L^4\) H\"{o}lder inequality by the \(L^2\)--\(L^{2D+2}\) analog in \eqref{eq_SS_last} to derive the following:
    \begin{align}\label{eq_Hdim_meas_factor}
        \Lambda\br{F,G}\lesssim 
            &\nrm{
                \ang{
                    \br{\vzeta\br{\vx}-\vw\br{\vx+\vX\br{\vt}}}^\top
                    \nabla\vX\br{\vt}
                }_\otimes^{-A}
            }_{L^{\frac{2D+2}{D}}\br{\abs{\vx}\lesssim \abs{\vt}\eqsim 1}}\\
            \cdot &
            \nrm{
                \nrm{
                    F\br{\vx,\vzeta}
                }_{L^1\br{d\vzeta}}
            }_{L^2\br{d\vx}}
            \nrm{
                \nrm{\nrm{g_\vw}_{L^1\br{d\mu_\vx}}}_{L^1\br{d\cH^D_\V\br{\vw}}}
            }_{L^2\br{\abs{\vx}\lesssim 1}},\nonumber
    \end{align}
    where \(\vzeta\br{\cdot},\vw\br{\cdot}\) are vector-valued measurable functions satisfying
    \begin{equation}\label{eq_Hdim_meas_condi}
        w_0\br{\cdot}\equiv 0,\quad
        \nrm{\vzeta\br{\cdot}}_{L^\infty},
        \nrm{\vw\br{\cdot}}_{L^\infty}\lesssim 
        \sqrt{\lambda} \lesssim 
        \inf_{\abs{\vx_0},\abs{\vx_1}\lesssim 1}\zeta_j\br{\vx_0}-w_j\br{\vx_1}.
    \end{equation}
    It remains to demonstrate a decay estimate for the factor \eqref{eq_Hdim_meas_factor}. We intend to reduce the estimate to the planar setting and apply \textbf{Lemma \ref{lem_sub_lev_2_2}}. We start by expanding the expression:
    \begin{align}\label{eq_Hdim_meas_expand}
        K\br{\vx,\vt}:= &
        \ang{
            \br{\vzeta\br{\vx}-\vw\br{\vx+\vX\br{\vt}}}^\top
            \nabla\vX\br{\vt}
        }_\otimes^{-A}\nonumber\\
        = &
        \ang{\zeta_0\br{\vx}+2t_0\br{\zeta_D\br{\vx}-w_D\br{\vx+\vX\br{\vt}}}}^{-A}\nonumber\\
        \cdot 
        \prod_{k=1}^{D-1} &
        \ang{\zeta_k\br{\vx}-w_k\br{\vx+\vX\br{\vt}}+2t_k\br{\zeta_D\br{\vx}-w_D\br{\vx+\vX\br{\vt}}}}^{-A}.
    \end{align}
    Let \(\epsilon\in\br{0,1/4}\) to be decided later. We consider two regions:
    \begin{equation*}
        S_\epsilon:=\BR{\vx\in\R^{D+1}\::\: \abs{\zeta_0\br{\vx}}\lesssim \lambda^{\frac{1}{2}-2\epsilon},\;\abs{\vx}\lesssim 1},\quad B_\epsilon:=\BR{\vx\in \R^{D+1}\setminus S_\epsilon\::\: \abs{\vx}\lesssim 1}
    \end{equation*}
    and split the factor \eqref{eq_Hdim_meas_factor} in the following manner:
    \begin{align*}
        \nrm{K\br{\vx,\vt}}_{L^{\frac{2D+2}{D}}\br{\abs{\vx}\lesssim\abs{\vt}\eqsim 1}}^{\frac{2D+2}{D}}
        \lesssim &
        \lambda^{-\epsilon}
        +
        \nrm{K\br{\vx,\vt}}_{L^1\br{d\vx d\vt, \vx\in S_\epsilon, \lambda^{-\epsilon}\lesssim \abs{t_0} \lesssim \abs{\vt}\eqsim 1}}
        +
        \nrm{K\br{\vx,\vt}}_{L^1\br{d\vx d\vt, \vx\in B_\epsilon, \abs{\vt}\lesssim 1}}\\
        =: &
        \lambda^{-\epsilon}+ K_1+K_2.
    \end{align*}
    For the treatment of \(K_1\), we look at the \(0\)-th and the \(j\)-th terms in \eqref{eq_Hdim_meas_expand}. We either have
    \begin{equation*}
        \abs{\zeta_D\br{\vx}-w_D\br{\vx+\vX\br{\vt}}}\eqsim \sqrt{\lambda}
    \end{equation*}
    and thus, the \(0\)-th factor in \eqref{eq_Hdim_meas_expand} produces decay:
    \begin{equation*}
        \ang{\zeta_0\br{\vx}+2t_0\br{\zeta_D\br{\vx}-w_D\br{\vx+\vX\br{\vt}}}}^{-A}
        \lesssim \ang{\lambda^{\frac{1}{2}-2\epsilon}-\lambda^{-\epsilon}\cdot \sqrt{\lambda}}^{-A} \eqsim\lambda^{-A\br{\frac{1}{2}-\epsilon}},
    \end{equation*}
    or the alternative:
    \begin{equation*}
        \abs{\zeta_D\br{\vx}-w_D\br{\vx+\vX\br{\vt}}}\ll \sqrt{\lambda}
    \end{equation*}
    and by assumption \eqref{eq_Hdim_meas_condi}, the \(j\)-th factor produces decay:
    \begin{equation*}
        \ang{\zeta_j\br{\vx}-w_j\br{\vx+\vX\br{\vt}}+2t_j\br{\zeta_D\br{\vx}-w_D\br{\vx+\vX\br{\vt}}}}^{-A}
        \lesssim \ang{\sqrt{\lambda}}^{-A} \eqsim\lambda^{-A/2}.
    \end{equation*}
    In combination, we deduce \(K_1\lesssim \lambda^{-A/4}\).
    As for the treatment of \(K_2\), we let:
    \begin{equation*}
        \zeta\br{\vx}:=\left\{
        \begin{aligned}
            &\zeta_0\br{\vx}, &\vx\in B_\epsilon,\\
            &\lambda^{\frac{1}{2}-2\epsilon}, & \vx \in S_\epsilon.
        \end{aligned}
        \right.
    \end{equation*}
    By design, we have the trivial estimate:
    \begin{equation*}
        K_2\leq \nrm{
            \ang{\zeta\br{\vx}+2t_0\br{\zeta_D\br{\vx}-w_D\br{\vx+\vX\br{\vt}}}}^{-A}
        }_{L^1\br{\abs{\vx},\abs{\vt}\lesssim 1}}.
    \end{equation*}
    Interchanging the order of the \(L^1\) norm equate the above with
    \begin{equation*}
        = 
        \nrm{
            \nrm{
                \ang{\zeta\br{\vx}+2t_0\br{\zeta_D\br{\vx}-w_D\br{\vx+\vX\br{\vt}}}}^{-A}
            }_{L^1\br{dx_0 dx_D dt_0,\abs{x_0},\abs{x_D},\abs{t_0}\lesssim 1}}
        }_{L^1\br{dx_1\cdots dx_{D_1} dt_1\cdots dt_{D-1}, \abs{x_k},\abs{t_j}\lesssim 1}}
    \end{equation*}
    Using the fact that \(\abs{\zeta\br{\cdot}}\gtrsim \lambda^{\frac{1}{2}-2\epsilon} \), we may apply the sub-level set result \textbf{Lemma \ref{lem_sub_lev_2_2}} in the planar setting to the inner expression and dominate the above whole expression with
    \begin{equation*}
        \lesssim
        \nrm{
            \br{\lambda^{\frac{1}{2}-2\epsilon}}^{-\delta}
        }_{L^1\br{dx_1\cdots dx_{D_1} dt_1\cdots dt_{D-1}, \abs{x_k},\abs{t_j}\lesssim 1}}
        \eqsim \lambda^{-\delta\br{\frac{1}{2}-2\epsilon}}.
    \end{equation*}
    We thus deduce \(K_2\lesssim \lambda^{-\delta\br{\frac{1}{2}-2\epsilon}}\). In conclusion, we have
    \begin{equation*}
        \nrm{
            \ang{
                \br{\vzeta\br{\vx}-\vw\br{\vx+\vX\br{\vt}}}^\top
                \nabla\vX\br{\vt}
            }_\otimes^{-A}
        }_{L^{\frac{2D+2}{D}}\br{\abs{\vx}\lesssim \abs{\vt}\eqsim 1}}
        \lesssim \br{\lambda^{-\epsilon}+K_1+K_2}^{\frac{D}{2D+2}}\lesssim \lambda^{-\sigma}
    \end{equation*}
    for some constant \(\sigma>0\) that only depends on \(D\).

    \subsection{Treatment of the uniform-uniform regime: proof of \texorpdfstring{\eqref{eq_Hdim_U_S}}{}}\label{subsec_Hdim_U_U}
    We claim that \eqref{eq_Hdim_U_S} follows from suitable adaptations of the argument in \textsc{Section \ref{subsec_U_S}}. By an identical computation, we first establish the direct analog of \eqref{eq_UU_doub}:
    \begin{equation}\label{eq_Hdim_UU_doub}
        \nrm{
            \frac{G\br{\vx,\vt,\vzeta,\vw}
            }{
                \ang{
                    \br{\vzeta-\vw}^\top
                    \nabla\vX\br{\vt}
                }_\otimes^{-A/2}
            }
        }_{L^2\br{
            \phi\br{\lambda^{-\frac{1}{2}}\zeta_D}
            d\vzeta
        }}^2\hspace{-2em}
        \lesssim \quad
        \nrm{
            \frac{
                \sqrt{\lambda}
                \prod_{k=0,1}
                    g_\vw\br{\vx-\lambda^{-\frac{1}{2}}\vx_k}
            }{
                \ang{
                    \sqrt{\lambda}
                    \br{-2\vt^\top,1}\cdot \br{\vx_0-\vx_1}
                }^N
            }
        }_{L^1\br{
            d\mu\br{\vx_0} d\mu\br{\vx_1}
        }}.
    \end{equation}
    From here, we obtain the analog of \eqref{eq_UU_doub_sqrt}
    \begin{equation*}
        \nrm{\sqrt{\eqref{eq_Hdim_UU_doub}}}_{L^2\br{\abs{\vt}\eqsim 1}}
        \leq
        \mu_{\vx}^{\frac{1}{4}}\br{g_\vw}
        \nrm{g_\vw}_{L^2\br{d\mu_\vx}}
        \nrm{
            \int_{\abs{\vt}\lesssim 1}
                \frac{
                    \sqrt{\lambda} d\vt
                }{
                    \ang{
                        \sqrt{\lambda}
                        \br{-2\vt^\top,1}\cdot
                        \vz
                    }^N
                }
        }^{\frac{1}{2}}_{L^2\br{
            \ang{\widehat{\vz}_D}_\otimes^{-N}d\vz
        }},
    \end{equation*}
    At first glance, it seems that a generalization of \textbf{Lemma \ref{lem_mu_tempering}} is needed to estimate the above expression. Yet, it is possible to estimate the above expression by the \(D=1\) counterpart. Indeed, we start with some interchanges of the \(L^1\) and \(L^2\) norms:
    \begin{align*}
        &
        \nrm{
            \int_{\abs{\vt}\lesssim 1}
                \frac{
                    \sqrt{\lambda} d\vt
                }{
                    \ang{
                        \sqrt{\lambda}
                        \br{-2\vt^\top,1}\cdot
                        \vz
                    }^N
                }
        }^{\frac{1}{2}}_{L^2\br{
            \ang{\widehat{\vz}_D}_\otimes^{-N}d\vz
        }}\\
        \leq &
        \nrm{
            \nrm{
                \nrm{
                \int_{\abs{t_0}\lesssim 1}
                    \frac{
                        \sqrt{\lambda} dt_0
                    }{
                        \ang{
                            \sqrt{\lambda}
                            \br{-2\vt^\top,1}\cdot
                            \vz
                        }^N
                    }
                }_{L^2\br{\frac{d z_0 d z_D}{\ang{z_0}^N}}}
            }_{L^1\br{d\widehat{t}_0,\abs{\widehat{t}_0}\lesssim 1}}
        }^{\frac{1}{2}}_{L^2\br{
            \ang{\widehat{z}_{0,D}}_\otimes^{-N}d\widehat{z}_{0,D}
        }}.
    \end{align*}
    We then perform the change of variables \(z:=z_D-2\sum_{k=1}^{D-1}t_k z_k\) to equate the above with its \(D=1\) counterpart
    \begin{align*}
        = &
        \nrm{
            \nrm{
                \nrm{
                \int_{\abs{t_0}\lesssim 1}
                    \frac{
                        \sqrt{\lambda} dt_0
                    }{
                        \ang{
                            \sqrt{\lambda}
                            \br{z-2t_0z_0}
                        }^N
                    }
                }_{L^2\br{\frac{d z_0 d z}{\ang{z_0}^N}}}
            }_{L^1\br{d\widehat{t}_0,\abs{\widehat{t}_0}\lesssim 1}}
        }^{\frac{1}{2}}_{L^2\br{
            \ang{\widehat{z}_{0,D}}_\otimes^{-N}d\widehat{z}_{0,D}
        }}\\
        \eqsim &
        \nrm{
            \int_{\abs{t_0}\lesssim 1}
            \frac{
                \sqrt{\lambda} dt_0
            }{
                \ang{
                    \sqrt{\lambda}
                    \br{z-2t_0z_0}
                }^N
            }
        }^{\frac{1}{2}}_{L^2\br{\frac{d z_0 d z}{\ang{z_0}^N}}}\lesssim \log^{\frac{1}{4}}\br{e+\lambda}.
    \end{align*}
    As a brief summary, we deduce:
    \begin{equation*}
        \Lambda\br{F,G}
        \lesssim \log^{\frac{1}{4}}\br{e+\lambda}
        \nrm{
            \nrm{
                \frac{F\br{\vx-\vX\br{\vt},\vzeta}}{
                    \ang{
                        \br{\vzeta-\vw}
                        \nabla\vX\br{\vt}
                    }_\otimes^{A/2}
                }
            }_{L^2\br{\substack{d\vt d\vzeta\\ \abs{\vt}\eqsim 1}}}\hspace{-4em}
            \cdot 
            \mu_{\vx}^{\frac{1}{4}}\br{\supp g_\vw}\nrm{g_\vw}_{L^2\br{d\mu_{\vx}}}
        }_{L^1\br{\substack{d\vx d\cH^D_\V\br{\vw}\\ \abs{\vx}\lesssim 1}}}.
    \end{equation*}
    To complete the rest of the adaptation of \textsc{Section \ref{subsec_U_S}}, it remains to estimate the analog of \eqref{eq_UU_2nd_cor_simp}:
    \begin{equation*}
        \nrm{
            \nrm{
                \ang{
                    \br{\vzeta-\vw\br{\vt}}
                    \nabla\vX\br{\vt}
                }_\otimes^{-A/2}
            }_{L^2\br{\abs{\vt}\lesssim 1}}
        }_{L^4\br{\phi\br{\lambda^{-\frac{1}{2}}\zeta_D}d\vzeta}},
    \end{equation*}
    where \(\vw\br{\cdot}\) is a measurable function valued in \(\V\).
    We conclude the proof by estimating the above expression by its \(D=1\) counterpart. Indeed, we interchange the ordering of the norms and dominate the above by
    \begin{align*}
        &
        \nrm{
            \nrm{
                \nrm{
                    \nrm{
                        \ang{
                            \br{\vzeta-\vw\br{\vt}}
                            \nabla\vX\br{\vt}
                        }_\otimes^{-A/2}
                    }_{L^4\br{d\widehat{\vzeta}_{0,D}}}
                }_{L^2\br{\abs{t_0}\lesssim 1}}
            }_{L^4\br{\phi\br{\lambda^{-\frac{1}{2}}\zeta_D}d\zeta_0 d\zeta_D}}
        }_{L^2\br{\abs{\widehat{\vt}_0}\lesssim 1}}\\
        \eqsim &
        \nrm{
            \nrm{
                \nrm{
                    \ang{
                        \zeta_0-2t_0\br{\zeta_D-w_D\br{t_0,\widehat{\vt}_0}}
                    }^{-A/2}
                }_{L^2\br{\abs{t_0}\lesssim 1}}
            }_{L^4\br{\phi\br{\lambda^{-\frac{1}{2}}\zeta_D}d\zeta_0 d\zeta_D}}
        }_{L^\infty\br{\abs{\widehat{\vt}_0}\lesssim 1}}\\
        \eqsim &
        \nrm{
            \eqref{eq_UU_2nd_cor_simp}
        }_{L^\infty\br{\abs{\widehat{\vt}_0}\lesssim 1}}
        \lesssim \log^{\frac{1}{4}}\br{e+\lambda}.
    \end{align*}
    This completes the proof of \eqref{eq_Hdim_U_S}.

\subsection{An upgrade of the natural \texorpdfstring{\(L^p\)}{} improving range: proof of \texorpdfstring{\eqref{eq_thm_Hdim_sing_tile_est_half}}{}}
In this situation, we cannot rely on the estimate \eqref{eq_thm_Hdim_sing_tile_Lp_improv} derived from the standard \(L^p\) improving estimate \textbf{Theorem \ref{thm_Lp_improv}}. In other words, we must utilize the frequency assumption \eqref{eq_thm_Hdim_sing_tile_freq_cond} to go beyond \(c \geq \frac{D}{2D+2}\) and extract a \(\lambda^{-\frac{1}{2}}\) decay. Recall that we may assume \eqref{eq_thm_Hdim_sing_tile_tile_normal}, and it suffices to show \eqref{eq_thm_Hdim_sing_tile_est_wo_shift} for \(c=\sigma=\frac{1}{2}\) under the frequency assumption \eqref{eq_thm_Hdim_sing_tile_freq_cond}. To proceed, we take \(\ast=C\) in the expression \eqref{eq_Hdim_wp_CL_decomp}:
\begin{equation*}
    f^{\br{+,1}}_{\vI\times\vomega}
    =
    \int_{\vomega^{\br{+,1}}} 
        \int_{\R^{D+1}}
            \ang{
                \pi_{\vomega^{\br{+,1}}} f,
                \varphi^{\br{C}}_{0,\vx,\vxi}
            }
                \varphi^{\br{C}}_{0,\vx,\vxi}
                \chi^{\br{C}}_\vI\br{\vx}
        d\vx
    d\vxi.
\end{equation*}
By \eqref{eq_thm_Hdim_sing_tile_freq_cond}, the frequency support properties of \(\varphi^{\br{C}}_{0,\vx,\vxi}\), and \eqref{eq_ass_freq_supp}, there is \(\Omega\subset\supp \br{\phi\ast\phi}_{\vomega^{\br{+,1}}}\) such that
\begin{equation*}
    f^{\br{+,1}}_{\vI\times\vomega}
    =
    \int_{\Omega} 
        \int_{\R^{D+1}}
            \ang{
                \pi_{\vomega^{\br{+,1}}} f,
                \varphi^{\br{C}}_{0,\vx,\vxi}
            }
                \varphi^{\br{C}}_{0,\vx,\vxi}
                \chi^{\br{C}}_\vI\br{\vx}
        d\vx
    d\vxi,\quad
    \dia\br{\pi_{\R\ve_0\oplus\R\ve_1}\Omega}\lesssim 1.
\end{equation*}
We then insert the above identity into the left-hand side of \eqref{eq_thm_Hdim_sing_tile_est_wo_shift} and derive the following:
\begin{equation*}
    \nrm{
            \ang{
                \pi_{\vomega^{\br{+,1}}} f,
                \varphi^{\br{C}}_{0,\vx,\vxi}
            }
            \1_\Omega\br{\vxi}
            \chi^{\br{C}}_\vI\br{\vx}
            \ang{
                C_{\mu_0}^{\br{\va}}\varphi^{\br{C}}_{0,\vx,\vxi},
                \1_{A^{\br{-,1}}_\vomega}g
            }
    }_{L^1\br{ d\vx d\vxi}}.  
\end{equation*}
Applying \(L^2\)--\(L^2\)--\(L^\infty\) H\"{o}lder's inequality, we dominate the above with
\begin{equation*}
    \leq
    \nrm{f^{\br{+,1}}}_{\vI\times\vomega}
    \abs{\Omega}^{1/2}
    \nrm{\br{\chi^{\br{C}}_\vI\br{\vx}}^{1/4}}_{L^2\br{d\vx}}
    \nrm{
        \ang{
            C_{\mu_0}^{\br{\va}}\varphi^{\br{C}}_{0,\vx,\vxi},
            \1_{A^{\br{-,1}}_\vomega} g
        }
        \1_\Omega\br{\vxi}
        \br{\chi^{\br{C}}_\vI\br{\vx}}^{1/4}
    }_{L^\infty\br{d\vx d\vxi}}.
\end{equation*}
Since \(\abs{\Omega}\lesssim\lambda^{D-1}\) and \(0 \leq \chi^{\br{C}}_\vI\br{\vx}\underset{N}{\lesssim} \ang{\vx}_\otimes^{-N}\), it suffices to prove the following uniform bound:
\begin{equation}\label{eq_thm_Hdim_sing_tile_est_half_g}
    \abs{
        \ang{
            C_{\mu_0}^{\br{\va}}\varphi^{\br{C}}_{0,\vx,\vxi},
            \1_{A^{\br{-,1}}_\vomega} g
        }
    }
    \1_\Omega\br{\vxi}
    \ang{\vx}_\otimes^{-N}
    \underset{N}{\lesssim}
    \lambda^{-\frac{D}{2}}
    \nrm{\1_{A^{\br{-,1}}_\vomega} g}_{L^1\br{d\mu}}
\end{equation}
for all \(\vx,\vxi\) and apply Cauchy-Schwarz inequality on \(\1_{A^{\br{-,1}}_\vomega} g= \1_{A^{\br{-,1}}_\vomega\cap \supp g} \cdot \1_{A^{\br{-,1}}_\vomega} g\) to conclude \eqref{eq_thm_Hdim_sing_tile_est_wo_shift} with \(c=\sigma=\frac{1}{2}\) and thus, \eqref{eq_thm_Hdim_sing_tile_est_half}. To derive \eqref{eq_thm_Hdim_sing_tile_est_half_g}, we first study the symbol formulation:
\begin{align*}
    C_{\mu_0}^{\br{\va}}\varphi^{\br{C}}_{0,\vx,\vxi}\br{\vz}= &
    \nrm{\phi}_{L^2}^{-D-1}
    \int
        \phi^{\otimes D+1}\br{\vzeta-\vxi}
        \widehat{\mu_0}\br{\vzeta-\va\br{\vz}}
        e\br{\br{\vz-\vx}^\top\vzeta}
    d\vzeta\\
    = &
    \nrm{\phi}_{L^2}^{-D-1}
    e\br{\br{\vz-\vx}^\top\vxi}
    \int
        \phi^{\otimes D+1}\br{\vzeta-\vxi}
        \widehat{\mu_0}\br{\vzeta-\va\br{\vz}}
        e\br{\br{\vz-\vx}^\top\br{\vzeta-\vxi}}
    d\vzeta\\
    = &
    \nrm{\phi}_{L^2}^{-D-1}
    e\br{\br{\vz-\vx}^\top\vxi}
    \int
        \phi^{\otimes D+1}\br{\vzeta}
        \widehat{\mu_0}\br{\vzeta+\vxi-\va\br{\vz}}
        e\br{\br{\vz-\vx}^\top\vzeta}
    d\vzeta
\end{align*}
Explicitly, we have:
\begin{align}\label{eq_thm_Hdim_sing_tile_est_half_be4_Taylor}
    &
    \int
        \phi^{\otimes D+1}\br{\vzeta}
        \widehat{\mu_0}\br{\vzeta+\vxi-\va\br{\vz}}
        e\br{\br{\vz-\vx}^\top\vzeta}
    d\vzeta\nonumber\\
    = &
    \int
        \phi^{\otimes D+1}\br{\vzeta}
        \int
            \overline{e\br{\br{\vxi-\va\br{\vz}}^\top \vX\br{\vt}}}
            \cdot
            \overline{e\br{\vzeta^\top \vX\br{\vt}}}
            \rho\br{\abs{\vt}}
            K\br{\vt}
        d\vt
        \cdot
        e\br{\br{\vz-\vx}^\top\vzeta}
    d\vzeta.
\end{align}
By performing the Taylor series expansion on the low-oscillatory term \(e\br{\vzeta^\top \vX\br{\vt}}\) similar to what has been done in \textsc{Section \ref{subsec_tf_rep_lin}}, we obtain the following estimate:
\begin{align}\label{eq_thm_Hdim_sing_tile_est_half_be4_sta}
    \abs{\eqref{eq_thm_Hdim_sing_tile_est_half_be4_Taylor}}
    \lesssim &
    \abs{
    \int
        \Phi\br{\vzeta}
        \int
            \overline{e\br{\br{\vxi-\va\br{\vz}}^\top \vX\br{\vt}}}
            \Psi\br{\vt}
            K\br{\vt}
        d\vt
        \cdot
        e\br{\br{\vz-\vx}^\top\vzeta}
    d\vzeta
    }\nonumber\\
    = &
    \abs{
    \int
        \overline{e\br{\br{\vxi-\va\br{\vz}}^\top \vX\br{\vt}}}
        \Psi\br{\vt}
        K\br{\vt}
    d\vt
    }
    \cdot
    \abs{
    \widecheck{\Phi}\br{\vz-\vx}
    }
\end{align}
for some choice of \(\Phi\in C^\infty_c\br{\R^{D+1}}\) and \(\Psi\in C^\infty_c\br{\R^D}\) satisfying
\begin{equation*}
    \abs{\partial^\valpha \Phi} \underset{\valpha}{\lesssim} \1_{\bR{-1/2,1/2}^{D+1}},\quad
    \abs{\partial^\vbeta \Psi} \underset{\vbeta}{\lesssim} \1_{\bR{-1/2,1/2}^D}.
\end{equation*}
Apply the above property and the (non-)stationary phase principle, we obtain:
\begin{equation*}
    \eqref{eq_thm_Hdim_sing_tile_est_half_be4_sta}
    \underset{N}{\lesssim}
    \abs{\vxi-\va\br{\vz}}^{-\frac{D}{2}}
    \ang{\vz-\vx}_\otimes^{-N}.
\end{equation*}
This proves the following pointwise estimate:
\begin{equation*}
    \abs{C_{\mu_0}^{\br{\va}}\varphi^{\br{C}}_{0,\vx,\vxi}\br{\vz}}
    \underset{N}{\lesssim}
    \abs{\vxi-\va\br{\vz}}^{-\frac{D}{2}}
    \ang{\vz-\vx}_\otimes^{-N}.
\end{equation*}
In combination, we have:
\begin{equation*}
    \abs{\ang{
        C_{\mu_0}^{\br{\va}}\varphi^{\br{C}}_{0,\vx,\vxi},
        \1_{A^{\br{-,1}}_\vomega} g
    }}
    \1_\Omega\br{\vxi}\ang{\vx}_\otimes^{-N}
    \underset{N}{\lesssim}
    \nrm{
        \abs{\vxi-\va\br{\vz}}^{-\frac{D}{2}}
        \1_\Omega\br{\vxi}
        \ang{\vz-\vx}_\otimes^{-N}
        \ang{\vx}_\otimes^{-N}
        \1_{A^{\br{-,1}}_\vomega}\br{\vz}g\br{\vz}
    }_{L^1\br{d\vz}}.
\end{equation*}
By \textbf{Lemma \ref{lem_jap_mul_diff}} and the fact that \(\abs{\vxi-\va\br{\vz}}\geq \dist\br{\Omega,\vomega^{\br{-,1}}}\geq \dist\br{\supp \br{\phi\ast\phi}_{\vomega^{\br{+,1}}},\vomega^{\br{-,1}}}\gtrsim \lambda\), we complete the proof of \eqref{eq_thm_Hdim_sing_tile_est_half_g} and thus, \eqref{eq_thm_Hdim_sing_tile_est_half} by dominating the above by
\begin{equation*}
    \lesssim
    \nrm{
        \lambda^{-\frac{D}{2}}
        \ang{\vz}_\otimes^{-N}
        \1_{A^{\br{-,1}}_\vomega}\br{\vz}
        g\br{\vz}
    }_{L^1\br{d\vz}}
    \eqsim
    \lambda^{-\frac{D}{2}}
    \nrm{\1_{A^{\br{-,1}}_\vomega}g}_{L^1\br{d\mu}}.
\end{equation*}


\section{Sub-level set analysis: Proof of Lemma \ref{lem_sub_lev_2_2}}\label{sec_sub_lev_est}

In this section, we provide the proof of \textbf{Lemma \ref{lem_sub_lev_2_2}} which follows almost verbatim the proof of \textbf{Lemma 4.8} in \cite{HsuL24}.
However, for the reader's convenience, we provide here the main ingredients.

We start by denoting
\begin{equation*}
    K\br{x,y,t}:=\ang{u\br{x,y}+t\br{v\br{x,y}+w\br{x+t,y+t^2}}}^{-N}
\end{equation*}
and set\footnote{The implicit constant in the definition of \(\cI\) may change from line to line throughout the argument due to the change of variables that we will employ in several key steps.}
\(
    \cI:=\set{
        x\in\R
    }{
        \abs{x}\lesssim 1
    }
\).
It suffices to show
\begin{equation*}
    \nrm{K\br{x,y,t}}_{L^1\br{dxdydt,\cI^3}}\lesssim \lambda^{-\delta}.
\end{equation*}
The proof is subdivided into four subsections, as follows:

\subsection{Decoupling of smooth variables from rough variables}\label{subsec_decoup_smooth_from_rough}
We first observe that the $L^1$ norm is stable under the parabolic shift\footnote{When the context is clear, we drop the domain \(\cI^3\) dependency in the \(L^p\) norm.}
\begin{equation}\label{eq_K_less_K_shift_123}
    \nrm{K}_{L^1}\leq
    \nrm{K\br{x-t,y-t^2,t}}_{L^1\br{dxdydt}}.
\end{equation}
Then, we square the above quantity and apply Cauchy-Schwarz 
\begin{align*}
   \eqref{eq_K_less_K_shift_123}^2
    \lesssim
    \nrm{
        \nrm{
            K\br{x-t,y-t^2,t}
        }_{L^1\br{dt}}
    }_{L^2\br{dxdy}}^2
    =
    \nrm{
        \prod_{j=0,1}K\br{x-t_j,y-t_j^2,t_j}
    }_{L^1\br{dxdydt_0dt_1}}
\end{align*}
which, after undoing the shift, gives
\begin{equation*}
    \nrm{K}_{L^1}^2\lesssim
    \nrm{
        K\br{x,y,t_+}K\br{x+t_+-t_-,y+t_+^2-t_-^2,t_-}
    }_{L^1\br{dxdydt_+dt_-}}.
\end{equation*}
Once here, we notice that
\begin{equation*}
    K\br{x+t_+-t_-,y+t^2_+-t^2_-,t_-}
    =\ang{
        w\br{x+t_+,y+t_+^2}t_-+
        \cdots
    }^{-N}.
\end{equation*}
shares the term \(w\br{x+t_+,y+t_+^2}\) with \(K\br{x,y,t_+}\).

Denoting by
\begin{equation*}
    X:=x+t_+-t_-,\quad
    Y:=y+t_+^2-t_-^2.
\end{equation*}
we apply \textbf{Lemma \ref{lem_jap_mul_diff}} in order to deduce
\begin{equation}\label{eq_KK_les_t2utv}
    K\br{x,y,t_+}K\br{X,Y,t_-}
    \lesssim
    \ang{
        \det
        \begin{pmatrix}
            t_+ & u\br{x,y}+t_+v\br{x,y}\\
            t_- & u\br{X,Y}+t_-v\br{X,Y}
        \end{pmatrix}
    }^{-N}
    =\ang{
        \bm{t}\vec{1}\wedge
        \br{\vec{u}+\bm{t}\vec{v}}
    }^{-N},
\end{equation}
where we introduce the alternating tensor product \(\wedge\) and the matrix/vector notations:
\begin{equation*}
    \bm{t}:=
    \begin{pmatrix}
        t_+ & 0\\
        0 & t_-
    \end{pmatrix}
    ,\quad
    \vec{1}:=
    \begin{pmatrix}
        1\\
        1
    \end{pmatrix}
    ,\quad
    \vec{u}:=
    \begin{pmatrix}
        u\br{x,y}\\
        u\br{X,Y}
    \end{pmatrix}
    ,\quad
    \vec{v}:=
    \begin{pmatrix}
        v\br{x,y}\\
        v\br{X,Y}
    \end{pmatrix}
\end{equation*}
to express the determinant in a more succinct form. 

Bounding \eqref{eq_KK_les_t2utv} trivially by \(1\) when \(\abs{t_+-t_-}\lesssim \lambda^{-2\delta}\) gives the following estimate:
\begin{equation*}
    \nrm{K}^2_{L^1}
    \lesssim \lambda^{-2\delta}+
    \nrm{
        \ang{
            \bm{t}\vec{1}\wedge
            \br{\vec{u}+\bm{t}\vec{v}}
        }^{-N}
    }_{L^1\br{dxdydt_+dt_-,\lambda^{-2\delta} \lesssim \abs{t_+-t_-} \lesssim 1 }}.
\end{equation*}
Thus, to prove \textbf{Lemma \ref{lem_sub_lev_2_2}}, it suffices to show
\begin{equation}\label{eq_K_2_bdd}
    \nrm{
        \ang{
            \bm{t}\vec{1}\wedge
            \br{\vec{u}+\bm{t}\vec{v}}
        }^{-N}
    }_{L^1\br{dxdydt_+dt_-,\lambda^{-2\delta} \lesssim \abs{t_+-t_-} \lesssim 1 }}\lesssim \lambda^{-2\delta}.
\end{equation}
We now notice that
\begin{equation}\nonumber
    t_\pm=\frac{Y-y}{2\br{X-x}}\pm\frac{X-x}{2},\quad dt_+dt_-=\frac{dXdY}{2\br{X-x}}\,,
\end{equation}
and thus, by setting
\(\bm{j}:=
    \begin{pmatrix}
        1 & 0\\
        0 & -1
    \end{pmatrix}\), 
we can write\footnote{For the simplicity of the exposition, we allow the following notational abuse: in what follows, if $\bm{A}$ is a matrix and $a\in\R$ then $\bm{A}+a:=\bm{A}+a \bm{I}$ where $\bm{I}$ stands for the identity matrix.} \(\bm{t}=\frac{Y-y}{2\br{X-x}}+\bm{j}\frac{X-x}{2}\)
and obtain
\begin{equation*}
    \ang{
        \bm{t}\vec{1}\wedge
        \br{\vec{u}+\bm{t}\vec{v}}
    }^{-N}
    =\ang{
        \br{
            \frac{Y-y}{2\br{X-x}}+
            \bm{j}\frac{X-x}{2}
        }\vec{1}
        \wedge
        \br{
            \vec{u}
            +
            \br{
                \frac{Y-y}{2\br{X-x}}+
                \bm{j}\frac{X-x}{2}
            }
            \vec{v}
        }
    }^{-N}.
\end{equation*}
Let \(K_2\br{x,y,X,Y}\) be the right-hand-side expression of the above equality; it then remains to show
\begin{equation*}
    \nrm{K_2\br{x,y,X,Y}}_{
        L^1\br{
            \frac{dxdX}{\abs{X-x}}dydY,
            \lambda^{-2\delta} \lesssim \abs{X-x} \lesssim 1
        }
    }
    \lesssim \lambda^{-2\delta}.
\end{equation*}
Our next intermediate goal is to eliminate the \(u\br{x,y}\) and \(v\br{x,y}\) terms in the \(K_2\br{x,y,X,Y}\) expression.

For this, we proceed as follows: firstly, we apply H\"older inequality
\begin{align*}
    &\nrm{K_2\br{x,y,X,Y}}_{
        L^1\br{
            \frac{dxdX}{X-x}dydY,
            \lambda^{-2\delta} \lesssim \abs{X-x} \lesssim 1
        }
    }^3\\
    \lesssim &
    \nrm{
        \frac{
            \1_{\cI^3}\br{x,y,X}
        }{
            X-x
        }
    }^3_{
        L^{\frac{3}{2}}
        \br{
            dxdXdy,
            \lambda^{-2\delta}\lesssim \abs{X-x} \lesssim 1
        }
    }
    \cdot
    \nrm{
        \nrm{
            K_2\br{
                x,y,X,Y
            }
        }_{L^1\br{dY}}
    }^3_{L^3
        \br{
            dxdXdy
        }
    }\\
   \lesssim & \lambda^{2\delta}\cdot \nrm{
        \prod_{j=1}^3
            K_2\br{x,y,X,Y_j}
    }_{L^1\br{
        dxdXdydY_1dY_2dY_3
    }}.
\end{align*}
Next, we apply \textbf{Lemma \ref{lem_jap_mul_diff}} on the inner expression \(\prod_{j=1}^3K_2\br{x,y,X,Y_j}\); for this, we shall first describe the choice of the \(O\br{1}\) bounded coefficients we choose.
Observe that
\begin{equation*}
    K_2\br{x,y,X,Y_j}=\ang{
        \cdots
        -
        \br{
            \frac{Y_j-y}{2\br{X-x}}
            -\frac{X-x}{2}
        }
        \br{
            u\br{x,y}
            +
            \br{
                \frac{Y_j-y}{2\br{X-x}}
                +\frac{X-x}{2}
            }
            v\br{x,y}
        }
    }^{-N}.
\end{equation*}
There exists a natural choice of coefficients $P_j$ given by
\begin{align*}
    \begin{pmatrix}
        P_1\\
        P_2\\
        P_3
    \end{pmatrix}
    := &
    \br{
        \frac{\bm{Y}-y}{2\br{X-x}}
        -\frac{X-x}{2}
    }
    \vec{1}
    \times
    \br{
        \frac{\bm{Y}-y}{2\br{X-x}}
        -\frac{X-x}{2}
    }
    \br{
        \frac{\bm{Y}-y}{2\br{X-x}}
        +\frac{X-x}{2}
    }
    \vec{1}\\
    = &
    \br{
        \frac{\bm{Y}-y}{2\br{X-x}}
        -\frac{X-x}{2}
    }
    \vec{1}
    \times
    \br{
        \frac{\bm{Y}-y}{2\br{X-x}}
        -\frac{X-x}{2}
    }^2
    \vec{1},
\end{align*}
where \(\times\) stands for the cross product in \(\R^3\) and \(\bm{Y}:=
    \begin{pmatrix}
        Y_1 & 0 & 0\\
        0 & Y_2 & 0\\
        0 & 0 & Y_3
    \end{pmatrix}\),
    \(\vec{1}:=
    \begin{pmatrix}
        1\\
        1\\
        1
    \end{pmatrix}\).
By construction,
\begin{equation}\nonumber
    \sum_{j=1}^3
        P_j\cdot
        \br{
            \frac{Y_j-y}{2\br{X-x}}
            -\frac{X-x}{2}
        }
        \br{
            u\br{x,y}
            +
            \br{
                \frac{Y_j-y}{2\br{X-x}}
                +\frac{X-x}{2}
            }
            v\br{x,y}
        }
    =0.
\end{equation}
Applying now \textbf{Lemma \ref{lem_jap_mul_diff}} for the choice above, we obtain
\begin{align*}
    &
    \prod_{j=1}^3
        K_2\br{x,y,X,Y_j}\\
    \lesssim &
    \ang{
        \sum_{j=1}^3
            P_j \cdot
            \br{
                \frac{Y_j-y}{2\br{X-x}}
                +\frac{X-x}{2}
            }
            \br{
                u\br{X,Y_j}
                +
                \br{
                    \frac{Y_j-y}{2\br{X-x}}
                    -\frac{X-x}{2}
                }
                v\br{X,Y_j}
            }
    }^{-N}\nonumber\\
    = &
    \Bigg<
        \br{
            \frac{\bm{Y}-y}{2\br{X-x}}+\frac{X-x}{2}
        }
        \br{
            \vec{U}
            +
            \br{
                \frac{\bm{Y}-y}{2\br{X-x}}-\frac{X-x}{2}
            }
            \vec{V}
        }\\
        &\phantom{
            \br{
                \frac{\bm{Y}-y}{2\br{X-x}}+\frac{X-x}{2}
            }
        }
        \wedge
        \br{
            \frac{\bm{Y}-y}{2\br{X-x}}-\frac{X-x}{2}
        }
        \vec{1}
        \wedge
        \br{
            \frac{\bm{Y}-y}{2\br{X-x}}-\frac{X-x}{2}
        }^2
        \vec{1}
    \Bigg>^{-N}\\
    =&: K_6\br{x,y,X,\bm{Y}},
\end{align*}
where here we define the two vectors as
\begin{equation*}
    \vec{U}=\vec{U}\br{X,\bm{Y}}:=
    \begin{pmatrix}
        u\br{X,Y_1}\\
        u\br{X,Y_2}\\
        u\br{X,Y_3}
    \end{pmatrix}
    ,\quad
    \vec{V}=\vec{V}\br{X,\bm{Y}}:=
    \begin{pmatrix}
        v\br{X,Y_1}\\
        v\br{X,Y_2}\\
        v\br{X,Y_3}
    \end{pmatrix}\,.
\end{equation*}

We thus reduced \eqref{eq_K_2_bdd}
 to showing the following:
\begin{equation*}
    \nrm{K_6\br{x,y,X,\bm{Y}}}_{L^1\br{
        dxdXdydY_1dY_2dY_3
    }}
    \lesssim \lambda^{-8\delta}.
\end{equation*}
In order to simplify our expression, we perform a change of variables with the corresponding Jacobian \(\eqsim 1\):
\begin{align*}
   & K_6\br{X-x,2y-x^2,X,2Y}
    =
    \ang{
        \br{
            \frac{\bm{Y}-y}{x}+x
        }
        \br{
            \vec{\mathcal{U}}
            +
            \frac{\bm{Y}-y}{x}
            \vec{\mathcal{V}}
        }
        \wedge
        \br{
            \frac{\bm{Y}-y}{x}
        }
        \vec{1}
        \wedge
        \br{
            \frac{\bm{Y}-y}{x}
        }^2\vec{1}
    }^{-N}\\
    \lesssim &
    \ang{
        \br{
            \bm{Y}-y+x^2
        }
        \br{
            x
            \vec{\mathcal{U}}
            +
            \br{\bm{Y}-y}
            \vec{\mathcal{V}}
        }
        \wedge
        \br{\bm{Y}-y}
        \vec{1}
        \wedge
        \br{\bm{Y}^2-y^2}
        \vec{1}
    }^{-N}=:\mathcal{K}_6\br{x,y,X,\bm{Y}},
\end{align*}
where we modify the two vectors
\(\vec{\mathcal{U}}=\vec{U}\br{X,2\bm{Y}}\) and \(\vec{\mathcal{V}}=\vec{V}\br{X,2\bm{Y}}\)
accordingly. 

This way we have created a polynomial with rough coefficients having as \emph{smooth} variables \(x,y\):
\begin{equation}\label{polyn}
    \br{
        \bm{Y}-y+x^2
    }
    \br{
        x
        \vec{\mathcal{U}}
        +
        \br{\bm{Y}-y}
        \vec{\mathcal{V}}
    }
    \wedge
    \br{\bm{Y}-y}
    \vec{1}
    \wedge
    \br{\bm{Y}^2-y^2}
    \vec{1}
    =\sum_{0\leq j,k\leq 5}
    C_{j,k}\br{X,\bm{Y}}x^jy^k
\end{equation}

Thus, our goal becomes now the following estimate:
\begin{equation}\label{eq_det_smooth_-8bd}
    \nrm{
        \ang{
            \sum_{0\leq j,k\leq 5}
            C_{j,k}\br{X,\bm{Y}}x^jy^k
        }^{-N}
    }_{L^1\br{dxdXdydY_1dY_2dY_3,\cI^6}}
    \lesssim \lambda^{-8\delta}.
\end{equation}

\subsection{A Van der Corput type estimate}
To prove \eqref{eq_det_smooth_-8bd}, we rely on the following lemma whose simple proof is presented in Section 4.3.3 of \cite{HsuL24} and thus will be omitted here:

\begin{lemma}\label{lem_van_der_cor_type}
     Let \(n,d\in\N\). Then there exists a constant\footnote{The exponent  $\delta\br{d}$ may be taken $\frac{1}{2d+2}$.} \(\delta\br{d}\in\bR{0,1}\) such that for any real-valued multivariate polynomial \(P\) on \(\R^n\) given by
    \begin{equation*}
        P\br{x}:=\sum_{\alpha:
        \abs{\alpha}\leq d}c_\alpha x^\alpha,
    \end{equation*}
    we have the following Van der Corput type estimates:
    \begin{equation*}
        \nrm{
            \ang{P\br{x}}^{-1}
        }_{L^1\br{dx,\cI^n}}
        \underset{n,d}{\lesssim} \ang{\nrm{c_\alpha}_{\ell^1\br{\alpha}}}^{-\delta\br{d}}.
    \end{equation*}
\end{lemma}

\subsection{Control over the coefficients of the polynomial}
By the multilinearity of the alternating tensor product, we may express coefficients \(C_{j,k}\br{X,\bm{Y}}\) in \eqref{polyn} as the result of linear maps applied on the vectors \(\begin{pmatrix}
    \vec{\mathcal{U}}\\
    \vec{\mathcal{V}}
\end{pmatrix}\). In particular, there is a \(3\times 3\) matrix \(\bm{M}_{\bm{Y}}\) consisting of \(O\br{1}\) bounded polynomial entries \(M_{i,j}\br{Y_1,Y_2,Y_3}\in\R\Br{Y_1,Y_2,Y_3}\) such that
\begin{equation*}
    \vec{\mathcal{C}}\br{X,\bm{Y}}:=
    \begin{pmatrix}
       C_{3,0}\br{X,\bm{Y}}\\
       C_{3,1}\br{X,\bm{Y}}\\
       C_{1,0}\br{X,\bm{Y}}
    \end{pmatrix}
    =
    \bm{M}_{\bm{Y}}
    \vec{\mathcal{U}}
    .
\end{equation*}
Now, in order to utilize \textbf{Lemma \ref{lem_van_der_cor_type}},
it suffices to derive lower bound on \(\nrm{\vec{\mathcal{C}}\br{X,\bm{Y}}}\) since
\begin{equation}\nonumber
    \nrm{\vec{\mathcal{C}}\br{X,\bm{Y}}}\lesssim
    \nrm{C_{j,k}\br{X,\bm{Y}}}_{\ell^1\br{j,k}}.
\end{equation}
This is achieved by showing that \(\bm{M}_{\bm{Y}}\) is non-singular ``for most" of the values of \(\br{Y_1,Y_2,Y_3}\in \cI^3\).
To understand the structure of the matrix \(\bm{M}_{\bm{Y}}\), we may successively take \(\vec{\mathcal{U}}=\vec{1},\;\bm{Y}\vec{1},\;\bm{Y}^2\vec{1}\) in \eqref{polyn} to derive
\begin{equation*}
    \bm{M}_{\bm{Y}}
    \cdot
    \begin{pmatrix}
        \vec{1} & \bm{Y}\vec{1} & \bm{Y}^2\vec{1}
    \end{pmatrix}
    =
    \begin{pmatrix}
        \bm{Y}^{0,1,2} & 0 & 0\\
        0 & \bm{Y}^{0,1,2} & 0\\
        0 & 0 & \bm{Y}^{1,2,3}
    \end{pmatrix},
\end{equation*}
where we introduce the notation
\begin{equation*}
    \bm{Y}^{a,b,c}:=\bm{Y}^a\vec{1}\wedge \bm{Y}^b\vec{1}\wedge\bm{Y}^c\vec{1}=
    \det\begin{pmatrix}
        Y_1^a & Y_1^b & Y_1^c\\
        Y_2^a & Y_2^b & Y_2^c\\
        Y_3^a & Y_3^b & Y_3^c
    \end{pmatrix}
    .
\end{equation*}
Taking the determinant on both sides
\begin{equation*}
    \det \bm{M}_{\bm{Y}}\cdot \bm{Y}^{0,1,2}=\bm{Y}^{1,2,3}\br{\bm{Y}^{0,1,2}}^2,
\end{equation*}
we deduce that
\begin{equation*}
    \det \bm{M}_{\bm{Y}}=\det\bm{Y}\cdot \br{\bm{Y}^{0,1,2}}^2=
    Y_1Y_2Y_3 \br{Y_1-Y_2}^2\br{Y_1-Y_3}^2\br{Y_2-Y_3}^2\,.
\end{equation*}
With this done,  we can now estimate the product of all the singular values of \(\bm{M}_{\bm{Y}}\) as follows:
\begin{equation*}
    \prod_{k=1}^3\mu_k\br{\bm{Y}}=
    \abs{
        \det\bm{M}_{\bm{Y}}
    }
    =
    \abs{Y_1 Y_2 Y_3}
    \br{Y_1-Y_2}^2\br{Y_1-Y_3}^2\br{Y_2-Y_3}^2.
\end{equation*}
Once at this point, we remark that we trivially have the following bounds on all singular values:
\begin{equation*}
    \max_k\nrm{\mu_k\br{\bm{Y}}}_{L^\infty\br{\cI^3}}
    \lesssim
    \max_{i,j}\nrm{M_{i,j}\br{Y_1,Y_2,Y_3}}_{L^\infty\br{\cI^3}}
    \lesssim 1\,.
\end{equation*}
Thus, if we further assume a lower bound on their product
\begin{equation*}
    \prod_{k=1}^3\mu_k\br{\bm{Y}}=
    \abs{
        \det\bm{M}_{\bm{Y}}
    }\gtrsim\lambda^{-\epsilon}\,,
\end{equation*}
we deduce that
\begin{equation*}
    \lambda^{-\epsilon}
    \lesssim
        \min_k\mu_k\br{\bm{Y}}
    \lesssim
        \max_k\mu_k\br{\bm{Y}}
    \lesssim 1
\end{equation*}
which immediately implies
\begin{equation*}
    \nrm{\bm{M}_{\bm{Y}}^{-1}}\lesssim \lambda^\epsilon.
\end{equation*}
With all these, we now have
\begin{equation*}
    \lambda
    \lesssim
    \nrm{
        \vec{\mathcal{U}}
    }
    \leq
     \nrm{
        \bm{M}_{\bm{Y}}^{-1}
    }
    \cdot
    \nrm{\vec{\mathcal{C}}\br{X,\bm{Y}}}\lesssim
    \lambda^\epsilon
    \nrm{C_{j,k}\br{X,\bm{Y}}}_{\ell^1\br{j,k}}\,,
    \end{equation*}
Thus, in summary, we proved the following:
\begin{equation}\label{eq_detM_R_dich}
    \lambda^{-\epsilon}\lesssim
    \abs{
        \det\bm{M}_{\bm{Y}}
    }
    \implies
    \lambda^{1-\epsilon}\lesssim
    \nrm{C_{j,k}\br{X,\bm{Y}}}_{\ell^1\br{j,k}}
    .
\end{equation}

\subsection{Concluding our lemma: the proof of \texorpdfstring{\eqref{eq_det_smooth_-8bd}}{}}
Take now \(\epsilon>0\) a constant to be chosen later, and consider the set
\begin{equation*}
    \mathcal{S}_3:=\set{\br{Y_1,Y_2,Y_3}\in\cI^3}{\lambda^\epsilon
        \det\bm{M}_{\bm{Y}}
    \in \cI}.
\end{equation*}
Focusing on the left-hand-side of \eqref{eq_det_smooth_-8bd}, we first apply the decomposition
\begin{align}
    &
    \nrm{
        \ang{
            \sum_{0\leq j,k\leq 5} 
            C_{j,k}\br{X,\bm{Y}}x^jy^k
        }^{-N}
    }_{L^1\br{dxdXdydY_1dY_2dY_3,\cI^6}}\nonumber\\
    \leq &
    \nrm{
        \ang{
            \sum_{0\leq j,k\leq 5} 
            C_{j,k}\br{X,\bm{Y}}x^jy^k
        }^{-N}
    }_{L^1\br{dxdXdy,\cI^3}\otimes L^1\br{dY_1dY_2dY_3,\mathcal{S}_3}}\label{eq_R_est_Sing}\\
    + &
    \nrm{
        \ang{
            \sum_{0\leq j,k\leq 5} 
            C_{j,k}\br{X,\bm{Y}}x^jy^k
        }^{-N}
    }_{L^1\br{dxdXdy,\cI^3}\otimes L^1\br{dY_1dY_2dY_3,\cI^3\setminus\mathcal{S}_3}}\label{eq_R_est_Poly}.
\end{align}%

Once at this point, we use the trivial bound for \eqref{eq_R_est_Sing}:
\begin{equation*}
    \eqref{eq_R_est_Sing}\leq \abs{\mathcal{S}_3}  \leq
    \nrm{
        \ang{
            \lambda^\epsilon
            \det\bm{M}_{\bm{Y}}
        }^{-1}
    }_{L^1\br{dY_1dY_2dY_3,\cI^3}}
    \lesssim \lambda^{ -\epsilon \delta\br{9} }\,,
\end{equation*}
while for \eqref{eq_R_est_Poly}, we use the fact that
\begin{equation*}
    \br{Y_1,Y_2,Y_3}\in \cI^3\setminus \mathcal{S}_3\implies \lambda^{-\epsilon}\lesssim \abs{\det\bm{M}_{\bm{Y}}},
\end{equation*}
together with \eqref{eq_detM_R_dich} and Lemma \ref{lem_van_der_cor_type} in order to deduce
\begin{equation*}
    \eqref{eq_R_est_Poly}\leq
    \nrm{
        \nrm{
            \ang{
                \sum_{0\leq j,k\leq 5} C_{j,k}
                    \br{X,\bm{Y}}x^jy^k
            }^{-1}
        }_{L^1\br{dxdy,\cI^2}}
    }_{L^\infty\br{dXdY_1dY_2dY_3,\cI\times \br{\cI^3\setminus \mathcal{S}_3}}}\hspace{-5em}
    \lesssim
    \lambda^{\br{\epsilon-1} \delta\br{5}}.
\end{equation*}
Finally, employing a standard variational argument, we optimize \(\epsilon=\frac{\delta\br{5}}{\delta\br{5}+\delta\br{9}}\) which gives
\begin{equation*}
    \nrm{
        \ang{
            \sum_{0\leq j,k\leq 5} C_{j,k}\br{X,\bm{Y}}x^jy^k
        }^{-N}
    }_{L^1\br{dxdydXdY_1dY_2dY_3,\cI^6}}
    \lesssim \lambda^{-\frac{\delta\br{5}\delta\br{9}}{\delta\br{5}+\delta\br{9}}}
    =:\lambda^{-8\delta}\,,
\end{equation*}
thus providing us with the desired estimate
\begin{equation*}
    \nrm{K}_{L^1}\lesssim\lambda^{-\frac{\delta\br{5}\delta\br{9}}{8\delta\br{5}+8\delta\br{9}}}.
\end{equation*}
This concludes the proof of \textbf{Lemma \ref{lem_sub_lev_2_2}} with \(\delta=\frac{\delta\br{5}\delta\br{9}}{8\delta\br{5}+8\delta\br{9}}\).

\section{Final Remarks}\label{finrem}

In this final section of our paper, we focus on providing more context to \textbf{Remark \ref{MotV}} and \textbf{Remark \ref{ext}}. With respect to the latter,  we will provide a comparative view between our current paper and the various relevant works \cite{c1,sj2,Bcarlrad,BGH24,CGGHIW24,PY19}.

\subsection{On the choice of \texorpdfstring{\(\V\)}{} in \texorpdfstring{\eqref{degone}}{}: optimality  within the linear resonance regime}\label{optV}

In this first subsection, we provide a rigorous argumentation to the informal claim\footnote{For a related, more succinct discussion, please see \textsc{Section 1.5} in \cite{Bcarlrad}.} made in \textbf{Remark \ref{MotV}}. Indeed, in what follows, we will prove that under the natural symmetry requirements \eqref{eq_Mod_2_V_symtize} and \eqref{eq_Dil_to_V_scaling} imposed on a non-empty subset\footnote{Notice here that we are not requiring any additional structure on \(\V\).} \(\V\subset \R^{D+1}\), restriction \eqref{degone} embodies the optimal condition---\emph{i.e.}, the largest possible $\V$---which allows linear but no higher order modulation invariance for $ CR^\ast_\V$ defined by \eqref{def_C_ast_V_K}.
       
With these being said, for \(\V\subset \R^{D+1}\) a non-empty subset, we introduce the following classes of symmetries: 
\begin{itemize}
    \item \textbf{Modulation symmetry:} Guided by \eqref{Modl2}, it is natural to require
    \begin{equation}\label{eq_Mod_2_V_symtize}
      \V-\va=\V,\quad \forall \va\in\V.
    \end{equation}
    \item \textbf{Kernel rescaling symmetry:} Departing from the observation that \eqref{eq_kernel_cond} is closed under the rescaling \(K\mapsto \Dil^1_s K\) for all \(s>0\) and based on the parabolic scaling property of \eqref{eq_paraboloid}
    \begin{equation*}
        \vX\br{s\cdot\vt}=
        \begin{pmatrix}
            s\vI_D & 0\\
            0 & s^2
        \end{pmatrix}
        \cdot
        \vX\br{\vt},
    \end{equation*}
   we are led to the following: if \(K\) is homogeneous of degree \(-D\),
    the natural dilation symmetry\footnote{This is the analogous to the dilation symmetry the classical Carleson operator satisfies.}
    \begin{equation*}
        CR^\ast_\V \Dil^p_{\vX_s}f=\Dil^p_{\vX_s} CR^\ast_\V f,\quad
        \vX_s:=\begin{pmatrix}
            s\vI_D & 0\\
            0 & s^2
        \end{pmatrix}
    \end{equation*}
    enforces the condition
    \begin{equation}\label{eq_Dil_to_V_scaling}
        \begin{pmatrix}
            s\vI_D & 0\\
            0 & s^2
        \end{pmatrix}
        \V=\V,\quad
        \forall s>0.
    \end{equation}
\end{itemize}
Given that the subset \(\V\subset \R^{D+1}\) satisfies \eqref{eq_Mod_2_V_symtize} and \eqref{eq_Dil_to_V_scaling}, we now demonstrate that \(\V\) can be decomposed into a direct sum of two vector subspaces of \(\R^{D+1}\),  as follows:
\begin{equation}\label{eq_sym_t_V_struct}
    \V=\widehat{V}_D\oplus V_D
    ,\quad 
    \widehat{V}_D\leq \R^D\times\BR{0},\quad
    V_D\leq \BR{0}^D\times\R.
\end{equation}
Indeed, we begin by observing that
\begin{equation*}
    \eqref{eq_Mod_2_V_symtize}
    \implies \V-\V=\V \implies \V=-\V.
\end{equation*}
In particular, we have that \(\V\) is a subgroup of \(\ang{\R^{D+1},+}\). Explicitly, we have
\begin{equation}\label{eq_V_abelian}
    0\in\V,\quad \pm \V \pm \V \pm \V \pm \cdots \pm \V=\V.
\end{equation}
Next, combining \eqref{eq_V_abelian} with the scaling property \eqref{eq_Dil_to_V_scaling}, we deduce that
\begin{equation}\label{eq_V_proj_vert}
    s
    \begin{pmatrix}
        \bm{O}_D & 0\\
        0 & 1
    \end{pmatrix}
    \va
    =
    \begin{pmatrix}
        \br{s+\frac{1}{2}}\vI_D & 0\\
        0 & \br{s+\frac{1}{2}}^2
    \end{pmatrix}
    \va -
    \begin{pmatrix}
        s\vI_D & 0\\
        0 & s^2
    \end{pmatrix}
    \va
    -
    \begin{pmatrix}
        \frac{1}{2}\vI_D & 0\\
        0 & \frac{1}{4}
    \end{pmatrix}
    \va
    \in\V,\quad
    \forall \va\in\V,\,\forall s>0.
\end{equation}
From \eqref{eq_V_proj_vert}, we further deduce
\begin{equation}\label{eq_V_proj_hori}
    s
    \begin{pmatrix}
        \vI_D & 0\\
        0 & 0
    \end{pmatrix}
    \va
    =
    \begin{pmatrix}
        s\vI_D & 0\\
        0 & s^2
    \end{pmatrix}
    \va
    -
    s^2
    \begin{pmatrix}
        \bm{O}_D & 0\\
        0 & 1
    \end{pmatrix}
    \va
    \in
    \V,\quad
    \forall \va\in\V,\,\forall s>0.
\end{equation}
The two properties \eqref{eq_V_proj_vert} and \eqref{eq_V_proj_hori} suggest that we set
\begin{equation*}
    \widehat{V}_D:=\begin{pmatrix}
        \vI_D & 0\\
        0 & 0
    \end{pmatrix}\V \subset \V,\quad
    V_D:=
    \begin{pmatrix}
        \bm{O}_D & 0\\
        0 & 1
    \end{pmatrix}
    \V\subset \V.
\end{equation*}
Routine verification shows that \eqref{eq_V_proj_vert} and \eqref{eq_V_proj_hori} respectively imply \(\widehat{V}_D\) and \(V_D\) are cones in \(\V\).
Moreover, since both \(\widehat{V}_D\) and \(V_D\) are subgroups of \(\ang{\V,+}\), we thus deduce that
\begin{equation*}
    \V\supset
    \widehat{V}_D + V_D
    =
    \widehat{V}_D\oplus V_D
    ,\quad 
    \widehat{V}_D\leq \R^D\times\BR{0},\quad
    V_D\leq \BR{0}^D\times\R.
\end{equation*}
The reverse inclusion holds trivially. This concludes \eqref{eq_sym_t_V_struct}.

Now, under the  assumption \eqref{eq_sym_t_V_struct}, we know from the discussion in \textsc{Section \ref{GoalCR}} on \eqref{q1} that whenever \(\widehat{V}_D=\R^D\) the operator \(CR^\ast_\V\) exhibits dependent linear-quadratic resonances analogous to \eqref{qi12}:
\begin{equation*}
    CR^\ast_\V M^0_{2,b} M^1_{2,b}\cdots M^{D-1}_{2,b} M^D_{1,b}= CR^\ast_\V,\quad \forall b\in\R.
\end{equation*}
From here, we conclude that if it is to rule out quadratic (higher) order resonances, requirement \eqref{degone} in our main result \textbf{Theorem \ref{thm_main}} becomes, indeed, optimal.

\subsection{The classical Carleson operator \texorpdfstring{\cite{c1}}{} and its higher-dimensional generalization \texorpdfstring{\cite{sj2}}{}}
We first claim that our \textbf{Theorem \ref{thm_main1D}} implies the \(L^p\)-boundedness of the classical Carleson operator
\begin{equation*}
    Cf\br{x}:=\sup_a
    \abs{
    \pv
        \int
            f\br{x-t}
            e\br{at}
        \frac{dt}{t}
    }.
\end{equation*}
In fact, we can deduce a slightly stronger result, that is, the boundedness of the following
\begin{equation*}
    C_\ast f\br{x}:=
    \sup_{\substack{
        0<r<R<\infty\\
        a\in\R
    }}
    \abs{
        \int_{r<\abs{t}\leq R}
            f\br{x-t}
            e\br{at}
        \frac{dt}{t}
    }.
\end{equation*}
Indeed, we first notice that via some changes of variables, we identify the above object with
\begin{equation*}
    =\frac{1}{2}\cdot
    \sup_{\substack{
        0<r<R<\infty\\
        a\in\R
    }}\abs{
        \int_{r<\abs{t}\leq R}
            f\br{x-t\abs{t}}
            e\br{a t\abs{t}}
        \frac{dt}{t}
    }
    .
\end{equation*}
However, the latter relates to the following Carleson-Radon transform
\begin{equation*}
    CR_\ast f\br{x_0,x_1}:=
    \sup_{\substack{
        0<r<R<\infty\\
        a\in\R
    }}\abs{
        \int_{r<\abs{t}\leq R}
            f\br{x_0-t,x_1-t\abs{t}}
            e\br{at\abs{t}}
        \frac{dt}{t}
    },
\end{equation*}
which is exactly \eqref{par} with $t^2$ replaced by $t|t|$. One easily sees that the proof of \textbf{Theorem  \ref{thm_main1D}} applies with trivial modifications in order to deduce that
\begin{equation}\label{eq_CR_tabst_bd}
    \nrm{CR_\ast f}_{L^p}\underset{p}{\lesssim}\nrm{f}_{L^p}.
\end{equation}
From here on, the boundedness of \(Cf\), follows via some standard reasonings involving suitable projection and limiting arguments. Indeed, one first
introduces the truncated operators
\begin{equation*}
    C_N f\br{x}=
    \sup_{\substack{
        0<r<R<N\\
        a\in\R
    }}\abs{
        \int_{r<\abs{t}\leq R}
            f\br{x-t\abs{t}}
            e\br{a t\abs{t}}
        \frac{dt}{t}
    },
\end{equation*}
\begin{equation*}
    CR_N f\br{x_0,x_1}:=
    \sup_{\substack{
        0<r<R<N\\
        a\in\R
    }}\abs{
        \int_{r<\abs{t}\leq R}
            f\br{x_0-t,x_1-t\abs{t}}
            e\br{at\abs{t}}
        \frac{dt}{t}
    }.
\end{equation*}
Next, given \(f\in L^p\br{\R}\), we observe the following identity
\begin{equation*}
    \1_{\br{-N,N}}\br{x_0}
    CR_N \br{\1_{\br{-2N,2N}}\otimes f}\br{\vx}
    =
    \1_{\br{-N,N}}\br{x_0} C_N f\br{x_1}
\end{equation*}
By direct computation and \eqref{eq_CR_tabst_bd}, we deduce that
\begin{equation*}
    N^{1/p}\nrm{C_N f}_{L^p\br{\R}}\underset{p}{\eqsim}\nrm{CR_N \br{\1_{\br{-2N,2N}}\otimes f}}_{L^p\br{\br{-N,N}\times \R}}\leq \nrm{CR_N \br{\1_{\br{-2N,2N}}\otimes f}}_{L^p\br{\R^2}}
\end{equation*}
\begin{equation*}
    \leq \nrm{CR_\ast \br{\1_{\br{-2N,2N}}\otimes f}}_{L^p\br{\R^2}}
    \underset{p}{\lesssim}
    \nrm{
        \1_{\br{-2N,2N}}\otimes f
    }_{L^p\br{\R^2}}
    \underset{p}{\eqsim} N^{1/p}\nrm{f}_{L^p\br{\R}}.
\end{equation*}
This implies \(L^p\) boundedness of \(C_N\) with the implicit constant independent of the choice of \(N\). Finally, we recover the classical Carleson's theorem by taking \(N\to \infty\) and applying the monotone convergence theorem.

As for the high-dimensional setting \cite{sj2}, we note that a naive application of the above projection trick on \eqref{def_C_ast_V_K} would not provide the desired result since tensor products of singular kernels satisfying \eqref{eq_kernel_cond} do not themselves satisfy \eqref{eq_kernel_cond}.
Due to this technicality, we instead perform a similar projection trick on \textbf{Theorem \ref{thm_log_shift_est_D_dim}} adapted to a suitable compact frequency support to deduce the full \(L^p\) boundedness range of the key model operators associated with the high-dimensional generalization of the Carleson operator. To elaborate, let \(j\in\BR{1,\dots,D-1}\), \(\underline{l}\in\Z\), \(\vb:\R^{D-1}\to\R^{D-1}\), and \(\overline{l}:\R^{D-1}\to\Z\cap\br{\underline{l},\infty}\). For any choice of coefficients \(\br{\varepsilon_k}_{k\in\Z}\in\ell^\infty\br{\C}\) and vector \(\vv\in\R^{D-1}\), consider the following model operator:
\begin{equation*}
        C^{\br{\vb,\varepsilon,\underline{l},\overline{l}}}_{j,\vv}f\br{\vz}
        := 
        \sum_{k=\underline{l}}^{\overline{l}\br{\vx}}
            \varepsilon_k
        \int_{\R^{D-1}}
        \widehat{f}\br{\vxi}
        \br{
            \Mod_{-\vv}
            \widetilde{\chi}_j
        }
        \br{
            2^k \br{\vzeta -\vb\br{\vz}}
        }
        e\br{\vzeta^\top\vz}
        d\zeta,
    \end{equation*}
where \(\widetilde{\chi}\) is given by modifying \eqref{eq_Hdim_chi_j_def}:
\begin{equation*}
    \widetilde{\chi}_j\br{\vzeta}:=
    \tpsi\br{\zeta_j}
    \varphi^{\otimes D-2}\br{\widehat{\vzeta}_j}
    .
\end{equation*}
Similar to the discussion for the classical Carleson theorem, we may assume \(\overline{l}\leq N\) for some large number \(N\) and leverage the key identity:
\begin{equation*}
    \cC^{\br{\br{0,\vb,0},\varepsilon,\underline{l},\overline{l}}}_{\vI,j,\br{0,\vv,0}}
    \br{\br{\Dil^1_{2^N}\widehat{\phi}} \otimes f \otimes\br{\Dil^1_{2^{2N}}\widehat{\phi}}}\br{x_0,\vz,x_D}
    =
    \br{\Dil^1_{2^N}\widehat{\phi}}\br{x_0}
    C^{\br{\vb,\varepsilon,\underline{l},\overline{l}}}_{j,\vv}f\br{\vz}
    \br{\Dil^1_{2^{2N}}\widehat{\phi}}\br{x_D}
\end{equation*}
to infer from \textbf{Theorem \ref{thm_log_shift_est_D_dim}} the following chain of estimates:
\begin{equation*}
    2^{{\color{red}{-}}3N/p'}\nrm{C^{\br{\vb,\varepsilon,\underline{l},\overline{l}}}_{j,\vv}f}_{L^p\br{\R^{D-1}}}
    \underset{p}{\eqsim}
    \nrm{
        \cC^{\br{\br{0,\vb,0},\varepsilon,\underline{l},\overline{l}}}_{\vI,j,\br{0,\vv,0}}
    \br{\br{\Dil^1_{2^N}\widehat{\phi}} \otimes f \otimes\br{\Dil^1_{2^{2N}}\widehat{\phi}}}
    }_{L^p\br{\R^{D+1}}}
\end{equation*}
\begin{equation*}
    \underset{p}{\lesssim}
    \log^2\br{e+\abs{\br{0,\vv,0}}}\nrm{\varepsilon}_{\ell^\infty}
    \nrm{
        \br{\Dil^1_{2^N}\widehat{\phi}} \otimes f \otimes\br{\Dil^1_{2^{2N}}\widehat{\phi}}
    }_{L^p\br{\R^{D+1}}}
\end{equation*}
\begin{equation*}
    \underset{p}{\eqsim}
    2^{{\color{red}{-}}3N/p'}
    \log^2\br{e+\abs{\vv}}\nrm{\varepsilon}_{\ell^\infty}
    \nrm{
        f
    }_{L^p\br{\R^{D-1}}}.
\end{equation*}
We thus conclude via the limiting argument mentioned previously that
\begin{equation*}
    \nrm{C^{\br{\vb,\varepsilon,\underline{l},\overline{l}}}_{j,\vv}f}_{L^p\br{\R^{D-1}}}
    \underset{p}{\lesssim}
    \log^2\br{e+\abs{\vv}}\nrm{\varepsilon}_{\ell^\infty}
    \nrm{
        f
    }_{L^p\br{\R^{D-1}}}.
\end{equation*}
We omit the standard argument deducing from the last relation above the result in \cite{sj2}.

\subsection{The degree one resonant higher-dimensional Carleson-Radon transform \texorpdfstring{\cite{Bcarlrad}}{}: discussion on the planarity, subspace \texorpdfstring{\(\V\)}{}, and \texorpdfstring{\(L^p\)}{} range conditions}\label{LB}
In \cite[\textbf{Theorem 1}]{Bcarlrad}, the \(L^p\) boundedness of \(C^\ast_\V\) is provided under the following three restrictions:
\begin{itemize}
    \item \textbf{\(L^p\) range}:
    \begin{equation}\label{eq_Lp_rang_Lars}
        \frac{D^2+4D+2}{\br{D+1}^2}<p<2\br{D+1},
    \end{equation}
    \item \textbf{Non-planarity:}
    \begin{equation}\label{eq_non_planar}
        D\geq 2,
    \end{equation}
    \item \textbf{Subspace requirement:}
    \begin{equation}\label{eq_subspace_Lars}
        \V\lneq \R^D\times\BR{0}\quad\text{or}\quad
        \V=\BR{0}^D\times\R.
    \end{equation}
\end{itemize}
The above restrictions have been briefly explained in the original paper \cite[\textsc{Section 1.5}]{Bcarlrad}. Below, we provide further context on the above via a comparative discussion between the original approach and ours.
\smallskip

We start our discussion with a more philosophical comment: a common feature of \cite{Bcarlrad} and our current work is the adaptation of the time-frequency analysis framework to frequency-\(\lambda\)-dilated systems of tiles in order to capture/exploit the high-frequency information.  However, though cosmetic at the core, the two works involve different choices of the time-frequency analysis frameworks. To be specific, \cite{Bcarlrad} is based on the approach introduced by Fefferman in \cite{f}, while our present work is an elaboration of the setting proposed by Lacey and Thiele in \cite{lt3}. The two frameworks are believed to be mutually transferable on Carleson-type problems with \emph{only linear resonances}. However,  the subtle differences between these two frameworks manifest in plain sight when passing to the higher order modulation invariance setting: indeed, all the known progress made in  Carleson-type problems where higher degree resonances are involved utilize Fefferman's framework as a stepping stone for the implementation of the relational time-frequency analysis---see \cite{LVQuadCarl,LVPolynCarl}, and, more recently, \cite{LVHidimPolyCarl,ZK,B24,BDJST}; consequently, if one aims to approach the degree two (or higher) resonant Carleson-Radon type problems---see the discussion in \textsc{Section \ref{GoalCR}} and in particular \textbf{Problems} $\bf{1}$ and $\bf{2}$ therein---one expects to involve, yet again, Fefferman's framework.

\subsubsection{Organization of the high-frequency information: impact on the \texorpdfstring{\(L^p\)}{} range}
In this section we discuss a first major difference addressing the $L^p$ range restriction---see \eqref{eq_Lp_rang_Lars}. At the origin of this lies the organization of the high-frequency information. Indeed, in \cite{Bcarlrad}, 
contributions of the high-frequency information at various \(\lambda\)-levels are controlled altogether in the key estimates; we however, control each \(\lambda\)-level high-frequency contribution individually with a \(\lambda^{-\delta}\) polynomial decay factor. This distinction reflects into the general \(L^p\) bounds for the high-frequency components. 

For context, without restricting to a single \(\lambda\)-level high-frequency component, the operator associated to a single tree is morally a modulated maximally truncated singular Radon transform. Due to the lack of a ``true convolution" behavior in a Radon-type operator, one cannot infer the moral constancy at small enough physical scales for a singular Radon transform. This precludes the standard argument in \cite[\textbf{Lemma 3}]{f} to derive the control on a single tree with a density factor. To bypass this issue, \cite[\textsc{Section 6.4}]{Bcarlrad} utilizes implicitly the \(L^p\) improving \textbf{Theorem \ref{thm_Lp_improv}} in the form of sparse bounds on the associated singular Radon transforms. As a result, via an extrapolation argument in the spirit of \cite{MR3148061}, the final \(L^p\) range is restricted to the regime given by \eqref{eq_Lp_rang_Lars}.
See \cite[(80), \textsc{Section 8}]{Bcarlrad} for more details.

In our case, the analysis on a single \(\lambda\)-level high-frequency component allows us to exploit locally full  \(L^p\) bounds. To be precise, we reinterpret the associated operator as a superposition of shifted Carleson operators. Each such shifted Carleson operator does possess a ``true convolution" behavior within the defining formula. This allows us to obtain the full \(L^p\) range control up to a tame \(\log^C\br{\lambda}\) growth for a single \(\lambda\)-level high-frequency component. Coupled with the weak-type \(\br{2,2}\) bound with a \(\lambda^{-\delta}\) decay factor provided by the \(L^2\) theory, we recover via interpolation the full range \(L^p\) bounds with a \(\lambda^{-\delta_p}\) decay factor which suffices for addressing summability issue in order to add together all the \(\lambda\)-level high-frequency contributions.

\subsubsection{Treatment of the stationary phase contribution, planarity, and the subspace condition}\label{}
The second distinction arises from the treatment of the stationary phase contribution of the symbol
\begin{equation}\label{sph}
    \widehat{\mu_0}\br{\vxi}:=
    \int
        \overline{e\br{\vxi\cdot \vX\br{\vt}}}
        \rho\br{\abs{\vt}}
        K\br{\vt}
    d\vt
    \approx
    e\br{-\frac{\sum^{D-1}_{j=0}\xi^2_j}{4\xi_D}}
    \cdot
    \abs{\vxi}^{-\frac{D}{2}}.
\end{equation}
%
%
%
%
%
%
More precisely, referring to the two limitations \eqref{eq_non_planar} and \eqref{eq_subspace_Lars}, we have:
\begin{itemize}
\item the analysis in \cite{Bcarlrad} involves only the \emph{size} of the amplitude in \eqref{sph};

\item in our present analysis, we use both the \emph{size} and the \emph{phase} in the stationary phase approximation evoked by \eqref{sph}.

\end{itemize}

In terms of layout, requirements \eqref{eq_non_planar} and \eqref{eq_subspace_Lars} in \cite[\textbf{Theorem 1}]{Bcarlrad} arise from the treatment of the high-frequency part of the anti-chain in \cite[\textsc{Section 5.3}]{Bcarlrad}. This part of the argument corresponds to our single tile decay estimates \textbf{Theorem \ref{thm_Hdim_sing_tile_ests}} and our control \eqref{eq_lem_Hdim_gen_tree_large_scale_split} on the \(\br{+,1}\)-tree in the proof of the general tree estimate \textbf{Lemma \ref{lem_Hdim_gen_tree}}.

To provide more context, the presence of the two assumptions \eqref{eq_non_planar} and \eqref{eq_subspace_Lars} in \cite{Bcarlrad} are directly related to the proof of \cite[\textbf{Lemma 6}]{Bcarlrad}, which,  making only use of the amplitude information in \eqref{sph}, involves the following relation:
\begin{equation}\label{eq_sum_symb_est}
    \sum_{\substack{
        \vzeta\in\Z^{D+1}\cap\V;\\
        \abs{\vxi+\zeta}\eqsim \lambda
    }}
    \abs{
        \widehat{\mu_0}\br{\vxi +\vzeta}
    }^2
    \lesssim
    \lambda^{-\delta},\qquad \forall\: \lambda\gtrsim 1\,.
\end{equation}
Once at this point, we notice that this strategy cannot address the particular case \(\V=\BR{0}\times\R^D\) covered by \textbf{Theorem \ref{thm_main}} since the amplitude information only gives
\begin{equation}\label{eq_sum_symb_est_strong}
    \sum_{\substack{
        \vzeta\in\Z^{D+1}\cap\V:\\
        \abs{\xi_j+\zeta_j}\eqsim \abs{\xi_D+\zeta_D}\eqsim \lambda,\\
        \forall j<D 
    }}
    \abs{
        \widehat{\mu_0}\br{\vxi +\vzeta}
    }^2
    \lesssim
    \lambda^{\dim \V-D}
\end{equation}
and sees no decay in \(\lambda\) when \(\dim \V=D\), which is always the case in the planar setting. However, \eqref{eq_sum_symb_est_strong} does suggest that the approach in \cite{Bcarlrad} could in fact cover all the cases verifying \(\dim \V\leq D-1\).


In comparison, our LGC-based approach utilizes both the amplitude information and the oscillatory nature of the symbol \(\widehat{\mu_0}\br{\vxi}\) to bypass the summability issue in \eqref{eq_sum_symb_est_strong}.
More specifically, an argument using only the amplitude information provides the more classical\footnote{In \cite{lt3}, the power of the mass/density term is one, which corresponds to \(c=\frac{1}{2}\) in our formulation.} single tile estimate \eqref{eq_thm_Hdim_sing_tile_est_half} under the specific frequency condition \eqref{eq_thm_Hdim_sing_tile_freq_cond}. Without this condition, the power decay in \eqref{eq_thm_Hdim_sing_tile_est_c} requires the careful analysis that truly makes use of the oscillatory nature of the symbol \(\widehat{\mu_0}\br{\vxi}\). Once we obtain the two forms of the single tile estimates \eqref{eq_thm_Hdim_sing_tile_est_c} and \eqref{eq_thm_Hdim_sing_tile_est_half} in \textbf{Theorem \ref{thm_Hdim_sing_tile_ests}}, the control \eqref{eq_lem_Hdim_gen_tree_large_scale_split} on a \(\br{+,1}\)-tree is achieved via a classical treatment\footnote{See the treatment of the \(1\)-tree in \cite[\textsc{Section 6}]{lt3}.} using \eqref{eq_thm_Hdim_sing_tile_est_half} for the central parts and an orthogonality argument using \eqref{eq_thm_Hdim_sing_tile_est_c} for the outer parts.




\subsection{The non-resonant one-dimensional Carleson-Radon transform  via Kakeya/Nikodym non-compression phenomenon \texorpdfstring{\cite{BGH24}}{}}\label{Dec}

In this section, we discuss some relationships between the methods employed in \cite{BGH24} and the general Carleson-Radon theme. Our discussion adapts to the two main features of the central theme discussed in this paper:
\begin{itemize}
\item the non-zero curvature feature encoded in the Radon behavior, and

\item the zero-curvature feature encoded in the linear or/and higher order resonances part of the (genuine) Carleson behavior. 
\end{itemize}

Now, as also evoked by the authors of \cite{BGH24}, the methods therein can only deal with the former item and are not suitable for the latter. 
With these being said, diving deeper into our analysis, we reveal the following:

The key feature of the formal time-frequency framework introduced in our current paper in \textsc{Section \ref{sec_formal_tf_ana}} is that the linear resonance effect and the curvature phenomenon are morally \emph{decoupled}. The former enters only within the low-resolution analysis, namely in \textbf{Lemma \ref{lem_-_tree}}, \textbf{Lemma \ref{lem_gen_tree_c_half}}, \textbf{Lemma \ref{lem_u_shift_mass_sel}}, and \textbf{Lemma \ref{lem_cE_2_bd}}, while the latter manifests only within the high-resolution analysis in which the single-tile decay estimate \textbf{Proposition \ref{singtileestimdec}} exploits the smoothing effect generated by the curvature. This latter aspect is precisely the phenomenon addressed in both \cite[\textbf{Theorem 1.4}]{HsuL24} and \cite[\textbf{Proposition 3.4} and \textbf{3.5}]{BGH24}.  Thus it comes as no surprise, that the method developed in \cite{BGH24}---via local smoothing estimates beyond universal bounds for variable-coefficient Schr\"{o}dinger operators\footnote{See also the companion paper \cite{CGGHIW24}.}---may, in principle, provide an alternate treatment of \textbf{Proposition \ref{singtileestimdec}}.

Indeed, to proceed as in \cite{BGH24}, one begins with the frequency localization \(\abs{\xi_0},\abs{\xi_1-a}\eqsim \lambda\gg 1\) and applies a stationary phase approximation for the symbol:
\begin{equation*}
    \widehat{\mu_0}\br{\xi_0,\xi_1-a}:=
    \int
        \overline{e\br{
            \xi_0 t+\br{\xi_1-a} t^2
        }}
        \rho\br{t}
    dt
    \approx
    e\br{-\frac{\xi^2_0}{4\br{\xi_1-a}}}/\sqrt{\lambda}.
\end{equation*}
With the following normalization
\begin{equation*}
    \br{a;\vxi} \mapsto \br{\lambda\br{w+3}; 2\lambda \xi_0, \lambda \xi_1},
\end{equation*}
one may further perform some simple computation and show that the phase function
\begin{equation*}
    \phi\br{w;\vxi}:=\xi^2_0/\br{\xi_1-w+3}
\end{equation*}
satisfies the Nikodym non-compression hypothesis \cite[\textbf{Definition 1.3}]{BGH24}
\begin{equation*}
    \partial^2_\vxi \phi\br{w;\vxi}\equiv \vnull,\quad
    \abs{
        \det
        \begin{pmatrix}
            \partial_w \partial^2_{\xi_0\xi_0}\phi\br{w;\vxi} & \partial_w \partial^2_{\xi_0\xi_1}\phi\br{w;\vxi} & \partial_w \partial^2_{\xi_1\xi_1}\phi\br{w;\vxi}\\
            \partial^2_w \partial^2_{\xi_0\xi_0}\phi\br{w;\vxi} & \partial^2_w \partial^2_{\xi_0\xi_1}\phi\br{w;\vxi} & \partial^2_w \partial^2_{\xi_1\xi_1}\phi\br{w;\vxi}\\
            \partial^3_w \partial^2_{\xi_0\xi_0}\phi\br{w;\vxi} & \partial^3_w \partial^2_{\xi_0\xi_1}\phi\br{w;\vxi} & \partial^3_w \partial^2_{\xi_1\xi_1}\phi\br{w;\vxi}
        \end{pmatrix}
    }\eqsim \frac{\abs{\xi_0}^3}{\abs{\xi_1-w+3}^{12}}\eqsim 1
\end{equation*}
for all \(\br{w,\xi,\eta}\in \br{-1,1}^3\). Consequently, one may proceed as in \cite{BGH24} to first derive a \(\lambda^{-\delta}\) decay \(L^p\) estimate associated to a single scale operator for some \(2<p<4\) and then further interpolate it with trivial estimates to recover the \(L^2\) decay estimate and lastly localizes the physical domain of the input function via standard tricks to derive \textbf{Proposition \ref{singtileestimdec}}.

Finally, returning to the second item in the introductory paragraph concerning the key Carleson behavior, we would like to explain why none of the methods in \cite{HsuL24} and \cite{BGH24} are able to address the \emph{degree two resonance regime}---see \emph{e.g.} the case represented by \(CR_{II,1}\) defined in \eqref{q1}.

Indeed, focusing first on the approach in \cite{HsuL24}, we notice that the statement analogous to \textbf{Proposition \ref{singtileestimdec}} is false for the model form associated with \(CR_{II,1}\). The exact point of failure for our current argument based on \cite{HsuL24} is the failure of the sub-level set estimate analogous to \textbf{Lemma \ref{lem_sub_lev_2_2}}.

Turning our attention towards the local smoothing argument in \cite{BGH24}, one failure point lies in the lack of Nikodym non-compression phenomenon \cite{CGGHIW24} for the phase function:
\begin{equation}\label{eq_phi_for_CR_II_Nico}
    \phi\br{w;\vxi}:=\br{\xi_0-w+3}^2/\xi_1
\end{equation}
associated with the symbol
\begin{equation*}
    \widehat{\mu_0}\br{\xi_0-a,\xi_1}:=
    \int
        \overline{e\br{
            \br{\xi_0-a} t+\xi_1 t^2
        }}
        \rho\br{t}
    dt
    \approx
    e\br{-\frac{\br{\xi_0-a}^2}{4\xi_1}}/\sqrt{\lambda}
\end{equation*}
in the region \(\abs{\xi_0-a},\abs{\xi_1}\eqsim \lambda\gg 1\) corresponding to the stationary phase contribution. More precisely, direct computation shows that the phase function \eqref{eq_phi_for_CR_II_Nico} does not satisfy the Nikodym non-compression hypothesis \cite[\textbf{Definition 1.3}]{BGH24}:
\begin{equation*}
    \partial^2_\vxi \phi\br{w;\vxi}\equiv \vnull,\quad
    \det
    \begin{pmatrix}
        \partial_w \partial^2_{\xi_0\xi_0}\phi\br{w;\vxi} & \partial_w \partial^2_{\xi_0\xi_1}\phi\br{w;\vxi} & \partial_w \partial^2_{\xi_1\xi_1}\phi\br{w;\vxi}\\
        \partial^2_w \partial^2_{\xi_0\xi_0}\phi\br{w;\vxi} & \partial^2_w \partial^2_{\xi_0\xi_1}\phi\br{w;\vxi} & \partial^2_w \partial^2_{\xi_1\xi_1}\phi\br{w;\vxi}\\
        \partial^3_w \partial^2_{\xi_0\xi_0}\phi\br{w;\vxi} & \partial^3_w \partial^2_{\xi_0\xi_1}\phi\br{w;\vxi} & \partial^3_w \partial^2_{\xi_1\xi_1}\phi\br{w;\vxi}
    \end{pmatrix}
    \equiv 0.
\end{equation*}

\subsection{The non-resonant higher-dimensional Carleson-Radon transform \texorpdfstring{\cite{PY19,AMPY24}}{}}\label{PY}
The work\footnote{Here for expository reasons we will only focus on the contribution in \cite{PY19}, though most of the comments below apply to \cite{AMPY24}) as well.} \cite{PY19}, initiating the study of Carleson-Radon transform in \eqref{CRdD}, provides \(L^p\) bounds for
\begin{equation}\label{eq_CR_cP}
    CR_{\cP}f\br{\vx,y}:=\sup_{P\in\cP}
    \abs{
        \int_{\R^D}
            f\br{\vx-\vt,y-\abs{\vt}^2}
            e\br{P\br{\vt}}
            K\br{\vt}
        d\vt
    },\quad
    \br{\vx,y}\in \R^D\times\R
\end{equation}
under the following three restrictions:
\begin{itemize}
    \item \textbf{Non-planarity:} \(D\geq 2\);
    \item \textbf{Non-resonance:} \(\cP \cap \operatorname{span}\BR{t_0,\dots,t_{D-1},\abs{\vt}^2} = \BR{0}\);
    \item \textbf{Degree-uni-homogeneity:}
    \(\cP\) takes the form\footnote{Recall the definition of $\cQ_{d,D}$ in \eqref{CRdD}.}
    \begin{equation}\label{eq_Poly_PY_cond}
        \cP:=\operatorname{span}\BR{P_j}_{j=2,\dots,D}
        \leq \cQ_{d,D},\quad P_j\br{\lambda \vt}=\lambda^j P_j\br{\vt}.
    \end{equation}
\end{itemize}
Our works \cite{HsuL24} and the present one are conceived in order to
address the first two restrictions, respectively. 

Regarding the third one, it is fair to say the following: while we did not identify any \emph{conceptual} obstacle in the applicability of our methods, a  naive, minimal time-investment adaptation of these methods faces several  \emph{technical} difficulties which prevent us from fully recovering the original result in \cite{PY19}. 

However, it is worth noticing that even involving just minimal adaptations, we are able to address some higher-dimensional non-resonant cases for which \(\cP\) is not required to satisfy the \textbf{degree-uni-homogeneity} property \eqref{eq_Poly_PY_cond}. In fact, we suspect that with only small changes our approach provides \(L^p\) bounds for \(CR_\cP\) when \(\cP\) verifies the \textbf{non-resonance} property and the following condition:
\begin{equation*}
    \cP\leq \cQ_{2,D},\quad \dim \cP\leq D.
\end{equation*}

Having established these considerations, for the sake of concreteness, we provide below an outline for approaching the explicit example
\begin{equation*}
    \cP:=\operatorname{span}\BR{
        t^2_1, \dots, t^2_{D-1}
    }.
\end{equation*}
Aligned with the works in \cite{PY19}, \cite{HsuL24} and \cite{BGH24}, the key to the \(L^p\) boundedness result in the non-resonant setting lies in the analysis of
\begin{equation*}
    \cC\cR^{\br{\va}} f\br{\vx,y}:= 
    \int_{\R^D}
        f\br{\vx-\vt,y -\abs{\vt}^2}
        e\br{\va\br{\vx,y}^\top \vQ\br{\vt}}
        \rho\br{\abs{\vt}}
    d\vt,
\end{equation*}
where we write the column vectors
\begin{equation*}
    \va\br{\vx,y}:=\br{a_1\br{\vx,y},\dots,a_{D-1}\br{\vx,y}}^{\top},\quad
    \vQ\br{\vt}:=\br{t^2_1,\dots,t^2_{D-1}}^{\top}.
\end{equation*}
We claim that the approach in our present paper can be adapted in order to provide the following smoothing inequalities/decay estimates:
\begin{equation}\label{eq_smooth_ineq_cross}
    \abs{
    \ang{
        \cC\cR^{\br{\va}} f,
        g
    }
    }
    \lesssim \lambda^{-\delta}
    \nrm{f}_{L^2}\nrm{g}_{L^2}
\end{equation}
whenever the following spatial/frequency localization is met
\begin{equation*}
    \br{\widehat{f}\br{\vxi,\eta},g\br{\vx,y}}\neq \br{0,0}
    \implies 
    \abs{\br{\vxi,\eta;\va\br{\vx,y}}}\eqsim \lambda. 
\end{equation*}
Indeed, we can achieve this by modifying the argument presented in \textsc{Section \ref{sec_est_of_sing_tile}} and \textsc{Section \ref{singltilesestimgen}}.
To see this, we start from the associated time-frequency representation
\begin{equation*}
    \Lambda\br{F,G}=
    \nrm{
        \frac{F\br{\vx-\vt,y-\abs{\vt}^2,\vu,v}
        G\br{\vx,y,\vt,\vu,v,\vw}}{
        \ang{
            \nabla \vQ\br{\vt}^\top\vw -\vu-2\vt v
        }^A_{\otimes}
        }
    }_{L^1\br{\substack{
        \abs{\vx},\abs{y}\lesssim 1\eqsim \abs{\vt},\\
        \abs{\br{\vu,v,\vw}}\eqsim \sqrt{\lambda}
    }}}
\end{equation*}
and proceed with deriving the three estimates analogous to \eqref{eq_Hdim_S_U}, \eqref{eq_Hdim_S_S}, and \eqref{eq_Hdim_U_S}, respectively: 
\begin{itemize}
    \item the \textsf{sparse-uniform} regime: one proceeds as in \textsc{Section \ref{subsec_Hdim_S_U}}; in the current setting, the log factor in \eqref{eq_Hdim_S_U} is replaced by a similar quantity
    \begin{equation*}
        \nrm{\nrm{\ang{
            \nabla\vQ\br{\vt}^\top\vw-\vu-2\vt v
        }^{-A}_\otimes}_{L^1\br{
            \substack{
                d\vw d\vt\\
                \abs{\vw}\lesssim \sqrt{\lambda},\abs{\vt}\lesssim 1
            }
        }}}_{L^\infty\br{\vu,v}}
        \lesssim \log^{D-1}\br{e+\lambda}
    \end{equation*}
    \item the \textsf{sparse-sparse} regime: one proceeds as in \textsc{Section \ref{subsec_Hdim_S_S}}; in the current setting, the power decay factor in \eqref{eq_Hdim_S_S} arises also from a sub-level set estimate
    \begin{equation*}
    \nrm{
        \ang{\vu\br{\vx,y}+2\vt v\br{\vx,y}-\nabla\vQ\br{\vt}^\top\vw\br{\vx+\vt,y+\abs{\vt}^2}}^{-A}_\otimes
    }_{L^1\br{\abs{\vx},\abs{y},\abs{\vt}\eqsim 1}}\lesssim \lambda^{-\sigma}
    \end{equation*}
    where we assume either \(\abs{\br{\vu,v}}\eqsim \lambda\) or \(\abs{\vw}\eqsim\lambda\). Via reductions similar to those presented in \textsc{Section \ref{subsec_Hdim_S_S}}, the above estimate can be inferred from \textbf{Lemma \ref{lem_sub_lev_2_2}}.
    variables
    \item the \textsf{uniform-uniform} regime: one proceeds as in \textsc{Section \ref{subsec_Hdim_U_U}}; we note that \(\nabla\vQ\) effectively disappears at the step analogous to \eqref{eq_Hdim_UU_doub}.
\end{itemize}
Putting all of the above together, a standard interpolation argument yields
\begin{equation*}
    \Lambda\br{F,G}\lesssim \lambda^{-\delta}\nrm{F}_{L^2}\nrm{g}_{L^2},
\end{equation*}
which proves the desired decay estimates \eqref{eq_smooth_ineq_cross}. We leave further details to the interested reader.

\end{document}